\RequirePackage{fix-cm}
\documentclass[smallcondensed]{svjour3}     % onecolumn (ditto)
\smartqed  % flush right qed marks, e.g. at end of proof
\usepackage{blindtext}

\usepackage{graphicx}
\usepackage{lineno,hyperref}
\modulolinenumbers[5]
\usepackage{mathtools}
\usepackage{amsfonts}
\usepackage{algorithm}
\usepackage{algorithmic}
\usepackage{amssymb}
\usepackage{amsmath}
\usepackage{subfig}
\usepackage{xcolor}
\usepackage[nocompress]{cite}
\usepackage{appendix}
\DeclareMathOperator*{\argmin}{argmin}
\DeclareMathOperator*{\argmax}{argmax}

\newcommand\hadis[1]{\textcolor{black}{#1}}
\newcommand\new[1]{\textcolor{black}{#1}}

\usepackage{booktabs}
\usepackage{lscape} 
\usepackage{graphicx,rotating,booktabs}

\usepackage{graphicx}
\usepackage{subfig}
\begin{document}
%\tableofcontents

\title{High-dimensional Black-box Optimization Under Uncertainty
%TK-MARS: An Efficient Approach for Optimizing Unknown Complex Functions with Limited Labeled Data
%\thanks{Grants or other notes
%about the article that should go on the front page should be
%placed here. General acknowledgments should be placed at the end of the article.}
}
%\subtitle{Do you have a subtitle?\\ If so, write it here}

%\titlerunning{Short form of title}        % if too long for running head

\author{Hadis Anahideh         \and
        Jay Rosenberger \and
        Victoria Chen%etc.
}

%\authorrunning{Short form of author list} % if too long for running head

\institute{H. Anahideh \at
              University of Illinois - Chicago \\
%              Tel.: +123-45-678910\\
%              Fax: +123-45-678910\\
              \email{hadis@uic.edu}           %  \\
%             \emph{Present address:} of F. Author  %  if needed
           \and
           J. Rosenberger \at
              University of Texas at Arlington\\
              \email{jrosenbe@uta.edu}
          \and 
            V. Chen \at
              University of Texas at Arlington\\
              \email{vchen@uta.edu}
}

\date{Received: date / Accepted: date}
% The correct dates will be entered by the editor

\maketitle
\begin{abstract}

%Limited data remain a challenge for optimizing expensive black-box systems. 

%Learning from limited data and finding a set of variables that optimizes an expected output arises in many real-world problems.
%practically in  design problems. %everywhere from molecular structure design for drug discovery to a deep Neural Network tuning.
% In such situations, the underlying function is complex yet unknown, a large number of variables are involved though not all of them are important and the outcome is under uncertainty.
%, and the interactions between the variables may be significant.
%On the other hand, it is usually expensive to collect more data and the outcome is under uncertainty. Unfortunately, despite such situations being prevalent in real-world problems, existing research has not addressed these issues jointly.
% In this research, we propose a new surrogate optimization approach to address these challenges. 
\new{Optimizing expensive black-box systems with limited data is an extremely challenging problem. As a resolution, we present a new surrogate optimization approach by addressing two gaps in prior research---unimportant input variables and inefficient treatment of uncertainty associated with the black-box output. 
We first design a new flexible non-interpolating parsimonious surrogate model using a partitioning-based multivariate adaptive regression splines approach, Tree Knot MARS (TK-MARS). The proposed model is specifically designed for optimization by capturing the structure of the function, bending at near-optimal locations, and is capable of screening unimportant input variables.
Furthermore, we develop a novel replication approach called \emph{Smart-Replication}, 
to overcome the uncertainty associated with the black-box output. The Smart-Replication approach identifies promising input points to replicate and avoids unnecessary evaluations of other data points. Smart-Replication is agnostic to the choice of a surrogate and can adapt itself to an unknown noise level. 
Finally to demonstrate the effectiveness of our proposed approaches we consider different complex global optimization test functions from the surrogate optimization literature.
The results indicate that TK-MARS outperforms original MARS within a surrogate optimization algorithm and successfully detects important variables. The results also show that although non-interpolating surrogates can mitigate uncertainty, replication is still beneficial for optimizing highly complex black-box functions. 
The robustness and the quality of the final optimum solution found through Smart-Replication are competitive with that using no replications in environments with low levels of noise and using a fixed number of replications in highly noisy environments.}
% \PACS{PACS code1 \and PACS code2 \and more}
% \subclass{MSC code1 \and MSC code2 \and more}
\end{abstract}

\section{Introduction}\label{sec:intro}
% Optimization involves finding a best solution among a set of feasible solutions.
Challenging optimization problems appear in many areas of science, technology, and industry. Examples include optimizing designs of wind farms~\cite{wilson2018evolutionary}, autonomous vehicle control systems~\cite{su2018autonomous}, green buildings~\cite{nguyen2014review}, vehicle safety systems~\cite{hamza2005vehicle}, molecular structures in pharmaceuticals~\cite{pyzer2018bayesian}, and material structures~\cite{griffiths2017constrained}. 
%finite element design~\cite{jones1998efficient},
Most of these problems include a complex system of inputs and outputs without well-known information about the underlying system behavior. 
Optimizing such problems is sometimes referred to as \emph{Black-Box Optimization (BBO)}. In this situation, the underlying function is complex yet unknown, a large number of input variables are involved, and there are substantial interactions between the input variables. 
%though not all of them are important and the outcome is under uncertainty. 
On the other hand, %it is usually expensive to collect more data.%Simulators can be used to study the black-box systems. Indeed, a large number of variables are involved, and 
%This results in a large combinatorial search space.
the evaluation process includes costly experiments, % due to the complexity of the systems, 
which can be either based upon computer simulators, such as finite element simulation tools~\cite{farkas2010fuzzy}, or actual experiments, such as crash simulators~\cite{gu2001comparison}. Achieving a near-optimal solution of a high-dimensional expensive black-box function within a limited number of function evaluations is the primary goal of this and other BBO research. %Therefore, finding the best sampled mean solution (BSMS) in a fewer number of function evaluations is a matter of concern.

\hadis{The BBO problem formulation we address in this research has a single objective, (\ref{eq:objfunc}), and a box constraint, (\ref{eq:boxconst}), and is given by:}
%similar to a conventional optimization problem with an objective function, Equation~(\ref{eq:objfunc}), and constraints defining the feasible region. 
% In this research, we only consider the box-constraints, Equation~(\ref{eq:boxconst}). 
% Therefore, the formulation given in Equations (\ref{eq:objfunc}) and (\ref{eq:boxconst}) is the BBO problem  we address in this research:
\begin{align}
\vspace{-3mm}
\label{eq:objfunc} &\min f(x) \\
\nonumber &\mbox{ s.t. }\\
%\label{eq:const} &g_i(x) \leq 0, \forall i=1\ldots m\\
\label{eq:boxconst} &\;a \leq x \leq b, \forall x \in \mathbb{R}^d.
\end{align}
Here $f(x)$ is the \emph{black-box function}, and the goal is to obtain a global optimal solution of $f(x)$ in the \emph{feasible input space} $D$, where $D =\{ x \in \mathbb{R}^d: a_j \leq x_{j} \leq b_j, \forall j=1,\ldots,d\}$. $a, b \in \mathbb{R}^d$ represent the lower bound and the upper bound of $x$, which is a $d$-dimensional vector of \emph{input variables}.

Since computing $f(x)$ is often expensive, % to be embedded directly within a global optimization method, 
a practical BBO approach is to employ a surrogate model that approximates $f(x)$ but is less expensive to evaluate. A \emph{surrogate model} is a mathematical approximation of the relationship between the input and the output variables. %Several types of statistical models including 
\hadis{Gaussian Processes (GP), also referred to as Kriging~\cite{krige1951statistical}}, 
Radial Basis Function (RBF)~\cite{powell1985radial}, Regression Trees~\cite{breiman1984classification}, Multivariate Adaptive Regression Splines (MARS)~\cite{friedman1991multivariate}, Artificial Neural Networks~\cite{papadrakakis1998structural}, and Support Vector Regression~\cite{clarke2005analysis} are examples of surrogate models based upon statistical modeling.
Optimizing the cheap to evaluate surrogate models for BBO is one of the existing well-known derivative-free techniques called \emph{surrogate optimization}. Surrogate optimization requires careful selection of the candidate points for evaluation to simultaneously improve the approximation performance and find a potential optimum, within a few evaluations. An exploration and exploitation trade-off is required to find the most promising data points for black-box function evaluation. In \S\ref{sect:literature}, we will elaborate on different input selection strategies for surrogate optimization.

\new{Most of the existing surrogate optimization methods assume that there is no uncertainty in the black-box system, especially in higher dimensional space, and the set of important input variables is known \emph{a priori}~\cite{jones1998efficient,regis2005constrained,muller2014influence}.} Both of these assumptions are often unrealistic in real-world problems~\cite{picheny2013quantile}. \hadis{One
such example, is the hyper-parameter tuning of deep learning for image processing \cite{feurer2019hyperparameter}, which is a stochastic complex system. }
\hadis{In one other example, Pilla et. al.~\cite{pilla2007robust} shows that, in a fleet assignment case study, only 42 out of 1264 variables remain significant.}
This research introduces a new surrogate optimization paradigm to address these two primary concerns.
% 1) How can we handle unimportant input variables (decision variables of the black-box function) in surrogate optimization?
% 2) How can we handle the uncertainties associated with the black-box outcome in surrogate optimization?
One surrogate model that has the potential to overcome both of these concerns is MARS \cite{friedman1991multivariate}. 
MARS is a parsimonious and non-interpolating model, implicitly capable of screening unimportant input variables. \hadis{Moreover, as a non-interpolating model, MARS is able to mitigate noise effects that are present in the objective values.}
%and handling the uncertainties that is present in the objective values. 
MARS has rarely been used in the surrogate optimization context. \hadis{M{\"u}ller et al. \cite{muller2014influence} studied the performance of MARS as a single surrogate and showed that MARS in its original form is inferior compared with other surrogates, like Radial Basis Functions and Gaussian Processes.}
However, in this research, we develop a partitioning-based MARS, called \emph{Tree-Knot Multivariate Adaptive Regression Splines (TK-MARS)}, to tailor MARS for the purpose of surrogate optimization \hadis{and improve its performance}. In addition, we develop a smart replication strategy to mitigate the uncertainty of the black-box output. \hadis{The robustness and the quality of the final optimum solution found through our proposed smart replication approach is competitive with using no replications in environments with low levels of noise and using a fixed number of replications in highly noisy environments. There are two notable advantages of the smart replication approach. One, it is agnostic to the choice of a surrogate model, and two, it adapts to the noise level that is present in the data. Hence, if the noise level is low it does not replicate, but if the noise level is high it replicates extensively.}

The remainder of this paper is organized as follows. We introduce a general surrogate optimization algorithm and related literature in \S\ref{sec:background} to highlight the gaps in prior research.
In \S\ref{sec:contribution}, we describe TK-MARS and the proposed replication approach for black-box functions under uncertainty. %This section explains the notion of unimportant variables in the context of surrogate optimization.
Finally, in \S\ref{sec:results}, we provide results of computational experiments for the performance of the proposed approach.
\vspace{-3mm}

\section{Background}
\label{sec:background}
\hadis{In this section, we first provide some technical background on surrogate optimization. We then further describe its core components that we focus on in this research---choosing an appropriate surrogate model (interpolating and non-interpolating models), as well as an exploration-exploitation Pareto sampling technique.
%method that is used in this research to select input points for the function evaluation. 
We also provide a brief background on Classification and Regression Trees (CART)~\cite{breiman1984classification}, which are used as a partitioning technique in our proposed modified MARS approach \S~\ref{sec:tkmars}.}
\vspace{-5mm}

% \subsection{Surrogate Optimization}
% \vspace{-3mm}
%Surrogate optimization is a derivative-free optimization technique for the expensive black-box functions. It is a sequential adaptive approach to discover the solution space, and find a global optimum. 
%Figure~\ref{fig:srg} shows the specific steps in a generic surrogate optimization approach. First of all, a small set of input points are initialized using the design of experiments methods (DOE) to be evaluated by the black-box function. The solution space is well represented by the most popular DOE space-filling sequences, including Latin Hypercube Design (LHD), Sobol~\cite{sobol1967distribution, sobol1976uniformly}, and Orthogonal Arrays (OA)~\cite{hedayat2012orthogonal}. Next is the metamodeling step, which fits a response surface model to the initial data set; in Figure~\ref{fig:srg} we use a cubic splines interpolation method. The following step is to apply an optimizer, which determines new candidate points by solving an auxiliary optimization problem.
% \begin{figure*}[!tb]
% \centering
% \includegraphics[width=0.75\textwidth]{figures/srg}
% \caption{Surrogate Optimization Framework Illustration}
% \label{fig:srg}
% \end{figure*}

\subsection{General Surrogate Optimization Framework}
\vspace{-3mm}

\hadis{Surrogate optimization is a derivative-free optimization technique, which applies cheap to evaluate surrogate models to optimize expensive black-box functions. \new{Surrogate optimization is referred to as \emph{one-shot optimization} when it optimizes over a fixed well-designed pool of sample points \cite{kennedy2001bayesian}. In other settings, surrogate optimization uses sequential model-based adaptive sampling \cite{jones1998efficient}}. Algorithm~\ref{alg:SurOpt} depicts a generic sequential surrogate optimization algorithm, which is the focus of this research.} In this setting, we assume the black-box function is deterministic. % noise-free.
In Step~\ref{step1}, a set of $N$ input points is generated using a space-filling technique, such as Latin Hypercube Design (LHD), Sobol sequence~\cite{sobol1967distribution, sobol1976uniformly}, or Orthogonal Array (OA)~\cite{hedayat2012orthogonal}. Each input point $x^i$ is evaluated to obtain a function value $f(x^i)$ in Step~\ref{step2}. 
Next, the algorithm fits a surrogate model to predict the output of the system using the already evaluated input points, Step~\ref{step4}. %In each iteration the fitted surrogate model is updated and a new input point is selected to be evaluated. 
\new{The surrogate fitting can be performed in multiple ways, using an interpolating or a non-interpolating surrogate model.}
Step~\ref{step5} determines, or samples, promising new input points using exploration, exploitation, or both. The exploration-exploitation trade-off involves identifying points that are in unexplored regions of the input space and points with minimum predicted output values. The sampled input points are evaluated in the black-box system in Step~\ref{step6}. A solution with the current best objective value, referred to as a \emph{best known solution (BKS)} is determined in Step \ref{step8}. Surrogate optimization repeats these steps until some termination criteria are met. Finally, it returns a BKS. \hadis{In this research, we first propose a modified version of MARS to be used in Step 4 of Algorithm 1 and then will investigate the impact of different non-/interpolating surrogate models in \S~\ref{sec:results}. We also propose a replication strategy for the cases where the uncertainty is present and will evaluate the performance of different surrogates within the proposed framework.}
%the available budget is exhausted. 
% Ideally, surrogate optimization attempts to find a global optimum of the surrogate model $\hat{f}$ and evaluate it with the black-box function to determine the actual objective value $f$. It adds the newly evaluated point to the initial data set and iteratively minimizes the $\hat{f}-f$ gap.
\begin{algorithm}[!tb]
\caption{{\bf Surrogate Optimization}}
\begin{algorithmic}[1]
\label{alg:SurOpt}
% \STATE Generate $N$ input points $x^i \in D$ with a DOE method. 
%\STATE $\{x^1, \ldots, x^N\} =$ a set of $N$ input points in $D$, selected with a DOE method \label{step1}
\STATE $I=\{x^1, \ldots, x^N\}$, a set of $N$ input points in $D$, selected with a Design of Experiment (DOE) method \label{step1}
% \STATE $I=\{\}$
% \FOR{$i=1$ to $N$}
%     \STATE evaluate $x^i$ and obtain $f(x^i)$
%     \STATE $I = I \cup \{(x^i,f(x^i))\}$
% \ENDFOR
%\STATE Evaluate each $x^i$ to obtain $f(x^i)$, and let $I=\{(x^i,f(x^i)) ~|~\forall i=1,\ldots,N\}$. \label{step2}
\STATE $\mathcal{F}= \{f(x^i)|x^i\in I\}$ \label{step2}
%\STATE $t=1$
\WHILE {Termination criteria not satisfied}
%; $N^t~\leq N_{\mathit{Max}}$
 \STATE \emph{Surrogate}: Fit a surrogate model $\hat{f}$ on $(I,\mathcal{F})$\label{step4} %$f(x^i)$ in the response variable. \label{step4} %$|I|$ evaluated points
 \STATE \emph{Sampling}: Determine new candidate points, $P$\label{step5}
 \STATE \emph{Evaluation}: $\mathcal{F}_P=\{f(x^i)|x^i\in P\}$ %Evaluate selected candidate points, $P$, $f(x^k), \forall x^k \in P$
 \label{step6}
  \STATE $I= I\cup P$; $\mathcal{F}=\mathcal{F}\cup \mathcal{F}_P$
%  \STATE Update the collection of evaluated input points, $I=I\cup P$
 %\STATE $t=t+1$
 \STATE \emph{Best Known Solution (BKS)}: $x^o\in\argmin_{x\in I} f(x)$\label{step8}
\ENDWHILE
\STATE Return $x^*=x^o$
\end{algorithmic}
\end{algorithm}
% The stopping criteria can be either a maximum number of expensive evaluations of the black-box function or the expected improvement of the best sampled mean solution (BSMS). 
%BSMS is the best known solution after $i$ iterations.
\vspace{-5mm}
\subsection{Interpolating vs. Non-interpolating Surrogate Models}\label{sect:intvnonint}
\begin{figure}[tb!]
\vspace{-5mm}
\centering
\subfloat[]{\includegraphics[width=0.35\textwidth]{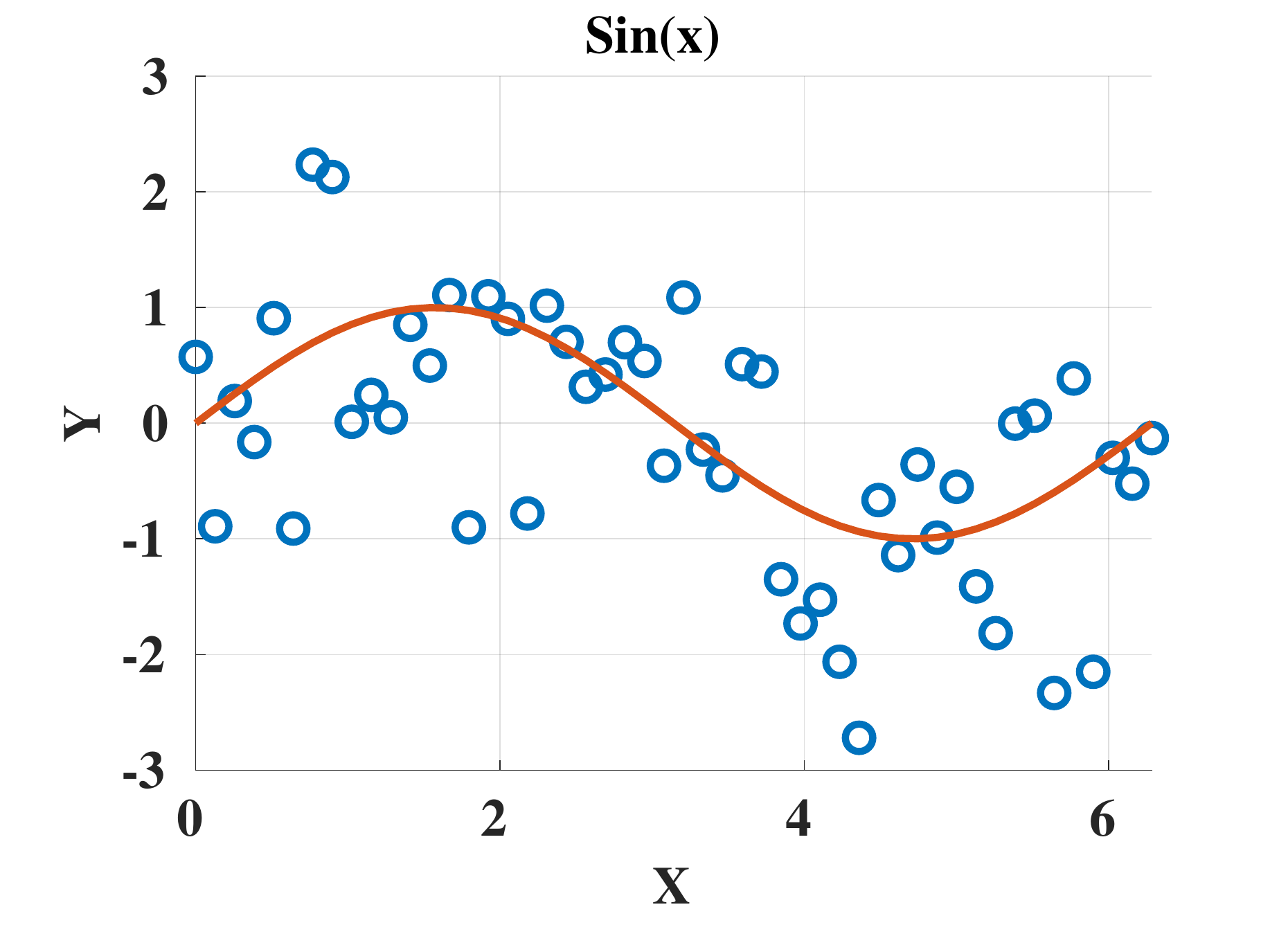}}
\subfloat[]{\includegraphics[width=0.35\textwidth]{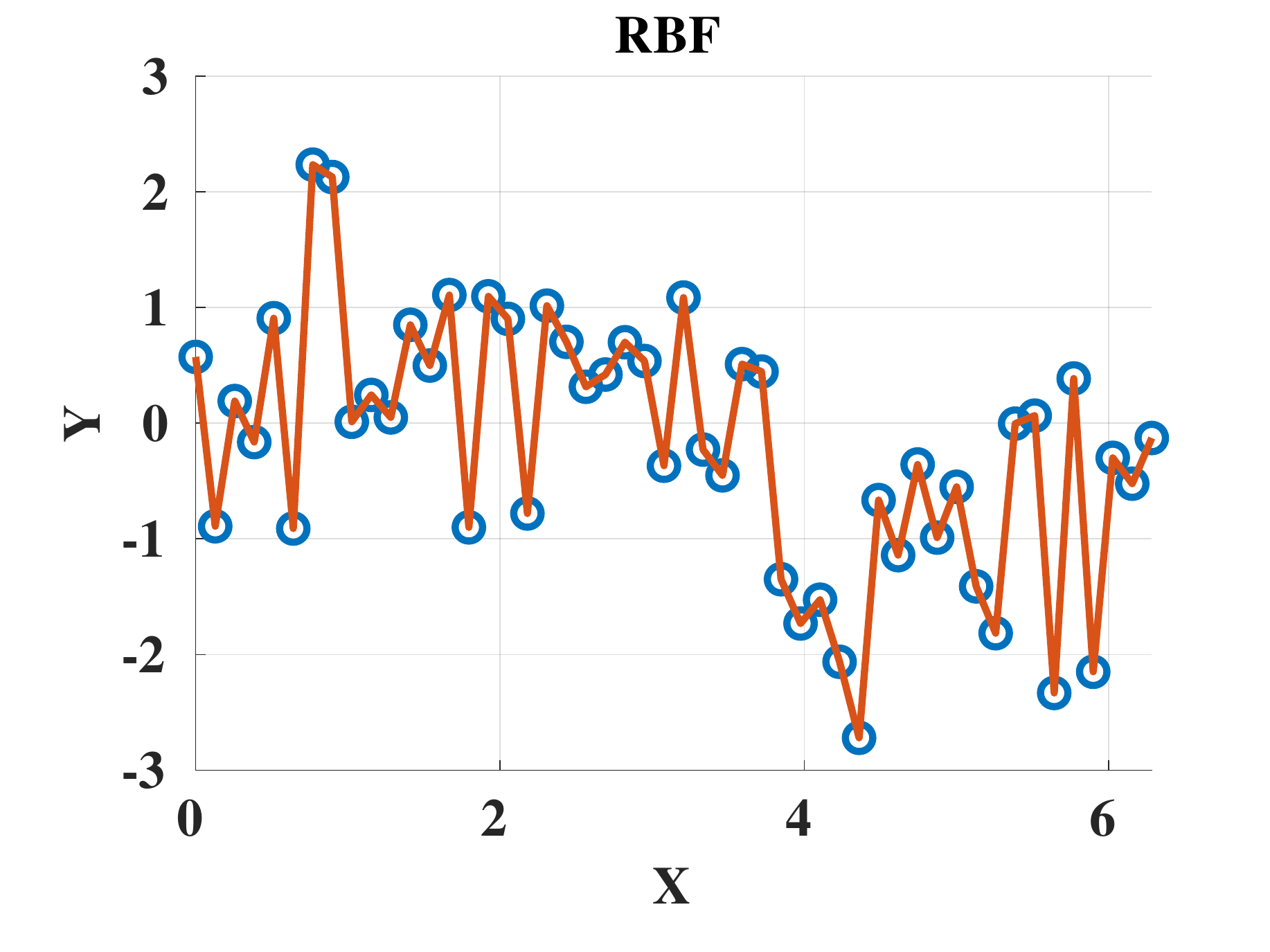}}
\subfloat[]{\includegraphics[width=0.35\textwidth]{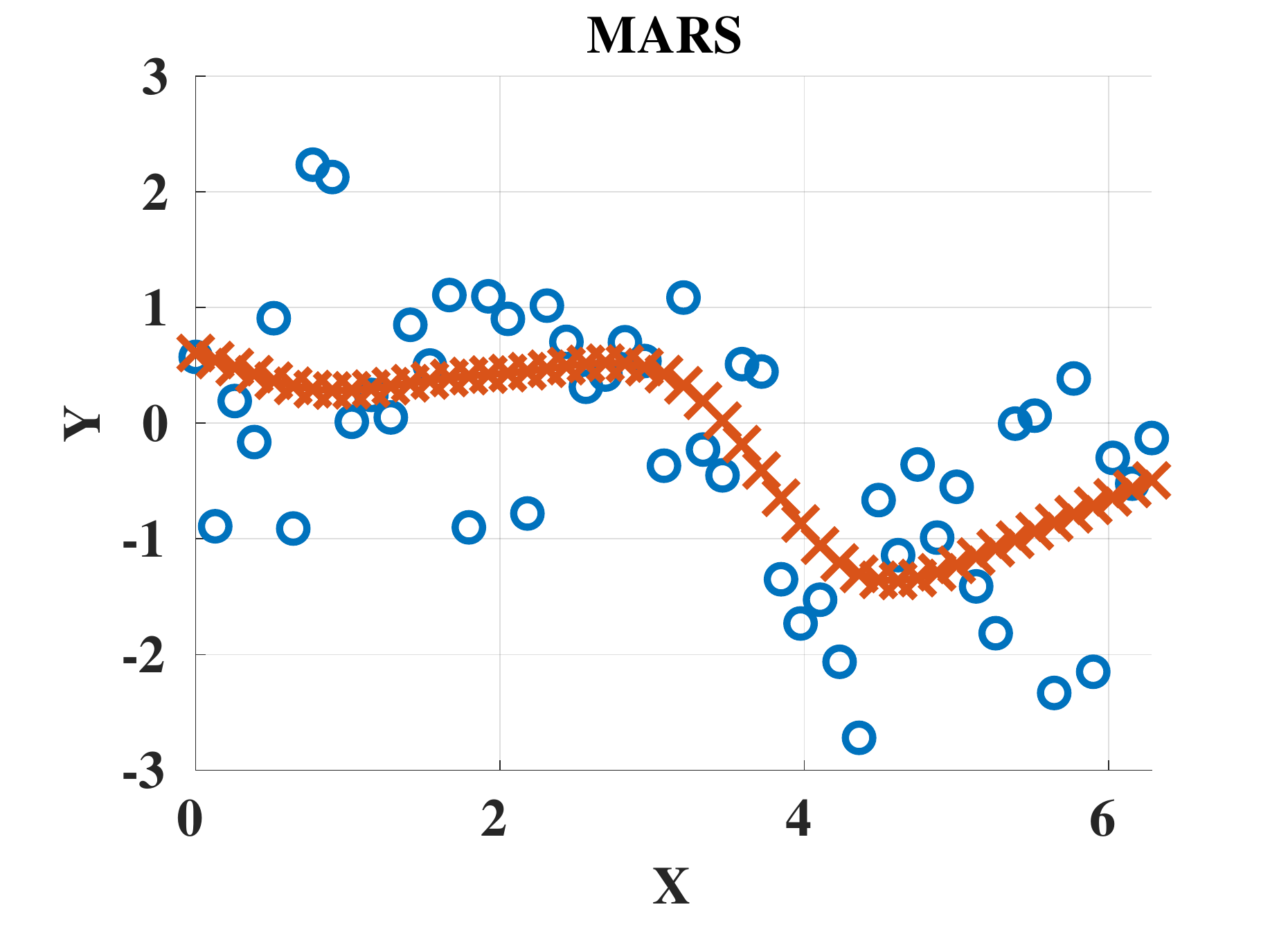}}
\caption{Interpolating versus non-interpolating models for noisy $\sin$ function}
\label{fig:sin}
\hspace{\fill}
\vspace{-7mm}
\end{figure}

\vspace{-3mm}
Choosing an appropriate surrogate model, $\hat{f}$ (Step~\ref{step4}), is extremely dependent on %the performance of distinct methods under different 
the characteristics of the black-box function. Surrogate models are classified as interpolating (e.g., \hadis{often} \hadis{GP}~\cite{krige1951statistical}, RBF~\cite{powell1985radial}) or non-interpolating (e.g., polynomial regression models~\cite{neter1996applied}, MARS~\cite{friedman1991multivariate}).
Interpolating models are the most common surrogates used in the surrogate optimization literature. They predict deterministic function values accurately by traversing the output of the input points. By contrast, non-interpolating models smoothly approximate the output under uncertainty \cite{hwang2018fast}. % Non-interpolating surrogate models are mainly based on the least squared error of some predetermined functional form. 
Unlike interpolating surrogate models, non-interpolating surrogate models do not necessarily traverse all the input points to capture the exact behavior of the function.
%In terms of bias-variance trade-off, non-interpolating surrogate models may have higher bias and lower variance than interpolating surrogates.
Therefore, in highly uncertain, or noisy, systems, non-interpolating models are preferred to avoid oscillations caused by interpolation~\cite{caballero2008algorithm}.
%Stochastic black-box systems are rarely addressed in surrogate optimization literature as we discussed before.
\hadis{For example, we show an interpolating RBF model in Figure~\ref{fig:sin}(b) and a non-interpolating MARS model in Figure~\ref{fig:sin}(c); both of these models approximate a simple $\sin (x)$ function under uncertainty.}
Observe that to capture a given input with two or more different outputs, the interpolation method needs to traverse the points with an extremely large slope. 
Consequently, interpolation-based surrogate models result in highly fluctuating approximations for the output of a black-box function with uncertainty. \hadis{In \S\ref{sec:results}, we demonstrate the performance of our proposed TK-MARS versus interpolating RBF, non-interpolating RBF, and non-interpolating GP in the surrogate optimization context.}

\vspace{-5mm}
\subsection{Radial Basis Function}\label{sec:rbf}
\vspace{-3mm}
Radial Basis Function (RBF) is one of the most common interpolating surrogate models~\cite{powell1985radial}. \new{Assuming $N$ distinct already evaluated input points, $x^1,x^2,\ldots,x^N \in \mathbb{R}^d$, the RBF interpolant is of the form, $\hat{f}(x)=\sum_{i=1}^{N} \lambda^i B(||x-x^i||)+p(x)$, $\forall x \in \mathbb{R}^d$, where the coefficients $\lambda^i$, $\forall i=1,\ldots,N$, are real numbers, $B(x)$ is a basis function, and $p(x)$ is a low degree polynomial, which is added to avoid singularity~\cite{wright2003radial}. Although there are multiple forms of the basis function $B(x)$, we use Multiquadric (MQ) RBF models\footnote{MQ had the highest performance among different RBF basis function types in our preliminary analysis.} with the basis function of the form $B(r)=\sqrt{r^2+\omega^2}$, where $r=||.||$ refers to the Euclidean $L_2$ norm, and $\omega$ is a constant called the \emph{shape parameter}. A larger shape parameter corresponds to a flatter basis function.} \hadis{In this research, the shape parameter is selected using cross-validation}. %It is worth mentioning that MQ guarantees non-singularity~\cite{wright2003radial}.
 %In this work we used one of the most common basis functions are Multiquadric (MQ), $B(r)=\sqrt{r^2+\omega^2}$, Gaussian (G) $B(r)=e^{\frac{-r^2}{\omega^2}}$, Cubic (C) $ B(r)=r^3$, and Thin Plate Splines (TPS) $B(r)=r^3$.

%The shape parameter $\omega$, and the number of points affect the accuracy and stability of RBF, with 
%For multiquadric and Gaussian basis functions, a unique interpolant is guaranteed when solving the RBF equation. However, for cubic and thin plate spline, the matrix of basis values might be singular. Rocha~\cite{rocha2009selection} adds a low degree polynomial, $p(x)$, to RBF interpolations to avoid singularity. 

\vspace{-5mm}
\subsection{Non-interpolating Radial Basis Function}\label{sec:nonRBF_LR}
\vspace{-3mm}
\new{Jakobsson et al.~\cite{jakobsson2010method} proposed a modified version of RBF as a non-interpolating surrogate to handle noisy objective functions, $\tilde{f}(x)=f(x)+\varepsilon$, and avoid oscillating interpolation. The proposed approximation allows the model to deviate from the data points. 
Specifically, the non-interpolating RBF 
is determined by minimizing a weighted objective that considers both the standard interpolating RBF objective, or surrogate norm, and the sum of squared errors, or deviation. 
Let $\hat{f}$ be the interpolating RBF function with coefficients $\lambda^i$. The squared semi norm is defined as $|S|^2=\sum_{i=1}^N \lambda^i \hat{f}(x^i)$. The new optimization problem finds the coefficients $\lambda^i$ that minimize $\eta|S|^2+(1-\eta) ||e||^2$ subject to $e_i=\hat{f}(x^i)-\tilde{f}(x^i), \forall i=1,\dots,N$, where $||e||$ denotes the $L_2$ norm of $(e_1,\dots,e_N)$.}
The trade-off between the surrogate norm and deviation is controlled by $\eta$, a 
hyperparameter that can be calibrated via cross-validation. Note that the proposed model cannot handle replicated data points and requires taking an average response if there are any.
\vspace{-5mm}

\subsection{Gaussian Processes for Noisy Observations}\label{sec:nonRBF_LR}
\vspace{-3mm}

Gaussian processes, also known as Kriging, is one of the most common interpolating surrogate models that has been used in the BBO literature. \new{A Gaussian process is a probabilistic model based on Gaussian distributions that are defined over a function $f$ with a mean function $\mu$ at each $x^i$ and a covariance function $\Sigma$ at each pair of points $x^i$ and $x^j$ \cite{rasmussen2003gaussian}.
%GP models are nonparametric and could be used for nonlinear problems.
\new{Let $X_N = [x^1, …, x^N]$ represent the matrix of evaluated input data points, let $\tilde{f}$ be a noisy function, and let $y_N = [\tilde{f}(x^1), …, \tilde{f}(x^N)]$ be the sampled output. The posterior distribution of $f$ (the true function being approximated) at any point $x$ is given as $f\sim\mathcal{GP}(\mu,\sigma^2)$, where $\mu(x)=\Sigma(X_N,x)^T[\Sigma(X_N,X_N)+\sigma^2_f \mathbb{I}_N]^{-1} y_N$, and $\sigma^2(x)=\Sigma(x,x)-\Sigma(X_N,x)^T[\Sigma(X_N,X_N)+\sigma^2_f\mathbb{I}_N]^{-1}\Sigma(X_N,x)$.
$\Sigma(x,x)$ represents a kernel function over the input space, $\Sigma(X_N,X_N)$ is an $N\times N$ matrix of kernel values calculated on the input data, $X_N$, $\Sigma(X_N,x)$ is a $N\times 1$ vector of kernel values, and $\sigma_f$ is the standard deviation of $\varepsilon$.}
There are different possibilities for the kernel function, such as exponential, squared exponential, and Matern class of covariance functions~\cite{rasmussen2003gaussian}. The most popular kernel is the squared exponential kernel $\Sigma(x,x')=\sigma_f^2\exp(-0.5||x-x'||^2/\sigma^2_m)$, where $\sigma_m$ represents the length scale for predictors. The kernel hyperparameters can also be calibrated using cross-validation.}

\vspace{-5mm}
\subsection{Multivariate Adaptive Regression Splines}\label{sec:mars_LR}
\vspace{-3mm}
Multivariate Adaptive Regression Splines (MARS) was introduced by Friedman~\cite{friedman1991multivariate}. MARS is a non-parametric non-interpolating surrogate model. The structure of the MARS model is based on basis functions. The MARS algorithm utilizes these basis functions to construct a piecewise continuous function in the following
%The structure of a MARS model is based on basis functions. 
% A MARS model is %algorithm utilizes splines to construct 
% a piecewise continuous function of model the dependent variables with the following 
%The basis functions are either univariate truncated linear functions or the product of different truncated linear functions for interaction terms. The truncated linear functions are of the form $b^{+}(x;k) = [x-k]_{+} $ or $b^{-}(x;k) = [k-x]_{+}$, where $\left[ q \right]_{+}= \max\{0,q\}$, $x$ is a single independent variable, and $k$ is the corresponding univariate knot, where the approximation bends to capture curvatures. %A knot is the location of an intersection between two splines and therefore is an important concept associated with MARS. 
%The eligible knot locations are selected from the set of training points. 
%Each dimension value of the training data set is an eligible knot for the corresponding independent variable $\{x^1,\ldots ,x^N\}$. 
%The MARS approximation has the 
form:
%\vspace{-3mm}
\begin{equation}
\label{eq:MARS_approx_form}
\hat{f}(x,\beta) = \beta_{0}+ \sum^{M_{\mathit{max}}}_{m=1} \beta_{m}B_{m}(x).\\
%\vspace{-3mm}
\end{equation}
where $x$ is a $d$-dimensional vector of input variables, $\beta_{0}$ is an intercept coefficient, $M_{\mathit{max}}$ is the maximum number of linearly independent basis functions, $\beta_{m}$ is the coefficient for the $m$th basis function, and $B_{m}(x)$ is a basis function that is either univariate or multivariate with interaction terms. 
These interaction terms have the following form:
%\vspace{-3mm}
\begin{equation}
\label{eq:MARS_int_basis_fun}
B_{m}(x) = \prod^{L_{m}}_{l=1} [s_{ml}(x_{j(m,l)}-t_{ml})]_{+}.
%\vspace{-3mm}
\end{equation}
where $L_{m}$ is the number of interaction terms in the $m$th basis function, $x_{j(m,l)}$ is the $j$th input variable of the $l$th truncated linear function in the $m$th basis function, and $t_{ml}$ is a knot location where the MARS function bends. The constant $s_{ml}$ is the direction of the truncated linear basis function and is either +1 or -1. 

\hadis{Within the context of modeling complex systems, MARS has two significant benefits.}
%in the sense of modeling complex systems.} 
%has two major advantages over interpolating models, such as RBF and Kriging, which are prevalent in the surrogate optimization literature~\cite{jones1998efficient, sacks1989design, simpson1998comparison, huang2006global, regis2007stochastic, regis2005constrained, muller2015ch, muller2017gosac, datta2016surrogate, krityakierne2016sop, krityakierne2014global}. 
First, it can mitigate uncertainty as discussed in \S\ref{sect:intvnonint}. Second, MARS is intended to be parsimonious and is able to screen unimportant input variables, which frequently occur in real-world complex systems~\cite{craig2005automotive}. 
The MARS algorithm includes 
a forward selection and a backward elimination procedures. In forward selection, MARS adds the basis functions in pairs for different dimensions, which gives the maximum reduction in the sum of squared error. The process of adding continues until it reaches the maximum number of basis functions, $M_{max}$. By removing the least effective term at each step, the backward elimination process avoids overfitting. MARS has an embedded dimension reduction technique, which is very useful for BBO where there is no previous understanding of the input variables and their impact on the output. As a consequence, the final MARS model includes only important input variables. %The advantages of MARS over other statistical models, such as linear regression models lie in its ability to handle curvature in high-dimensional space and produce easier-to-interpret models. 

Nonetheless, the MARS algorithm was developed for approximation purposes but not for optimization. 
%To do this, we must analyze the knot locations where the MARS model bends and develop a new eligible knot optimization for the latter purpose.
%To do this, the knot locations, where the MARS model bends, must be evaluated, and a new eligible knot optimization built.
Specifically, in the traditional version of MARS described in Friedman~\cite{friedman1991multivariate}, each input point in the training set can be an eligible knot location. As the size of the training set increases, the number of eligible knot locations increases. %Adding more input points increases the size of eligible knot locations. 
MARS tends to interpolate the input points and loses its flexibility as the number of eligible knots increases.
Interpolating a set of limited input points early leads to a false assessment of function behavior and multiple local optima for highly noisy functions, as described in \S\ref{sect:intvnonint} and in Koc and Iyigun \cite{koc2014restructuring}. 
%The primary reason we consider MARS in surrogate optimization is the fact that MARS does not interpolate but finds the general structure of the function. So we have to avoid including all the eligible knot locations in the model.
%For highly noisy functions, interpolating leads to high local variance in function estimation~\cite{koc2014restructuring}.
This may cause difficulties for surrogate optimization. In addition, when a large number of eligible knots are selected, multicollinearity can occur between basis functions with knots that are close to each other. 
%Choosing a subset of eligible knot locations for basis function construction and slowly increase the subset as the number of training points increases, the complexity and flexibility of MARS increases.

Eligible knot selection techniques are developed to mitigate local variance and multicollinearity issues~\cite{koc2014restructuring}. One of the most common techniques selects evenly spaced knot locations within the range of the input points~\cite{chen1999applying, martinez2013variants, cao2010penalized, huang2004functional, song2010oracally}. %Though this method may not capture the true patterns in the data. %, and a small number of knots may result in underfitting. 
Friedman~\cite{friedman1991multivariate} proposes a minimum span (MinSpan) approach to minimize the local variability. In the MinSpan approach, for each independent variable, a local search around its current knot location is designed to reduce the number of eligible knot locations. 
Miyata and Shen~\cite{miyata2005free} presents a simulated annealing approach to choose eligible knot locations. 
Koc and Iyigun~\cite{koc2014restructuring} develops a mapping strategy by transforming the original input points into a network of nodes through a nonlinear mapping. %Each node has a specific topological position in the lattice and is represented by a weight vector. 
The nodes in the mapped network act as %basically 
a reference for choosing the new eligible knot locations.
%The size and structure of the grid are preset in advance. The nodes are located with equal space on the grid. 
%The mapping hyperparameters may have an impact on the accuracy and time efficiency of the model.
%The grid size and a threshold value are self-organized mapping parameters, which have an impact on the accuracy and time efficiency of the model. Decreasing the grid size and increasing the threshold value lead to fewer candidate knot locations and decrease CPU time and model accuracy. Hence, the mapping parameters need to be optimized beforehand.
\new{Traditional knot selection approaches are mainly designed to fit the entire approximation over the whole space, which is different from what is best for surrogate optimization. 
Moreover, Müller and Shoemaker \cite{muller2014influence} conducted a comprehensive study on the influence of ensemble surrogate models and sampling strategies, and they showed that using MARS as a surrogate for optimization was inferior to using both RBF and GP. These facts motivate modifying MARS to be more effective for surrogate optimization.
To the best of our knowledge, there exists no other method for MARS in surrogate optimization. In this work, we propose a partitioning-based approach using classification and regression trees to tailor MARS for surrogate optimization.}

\vspace{-5mm}
\subsection{Classification and Regression Trees}\label{sec:cart}
\vspace{-3mm}
Classification and Regression Trees (CART) use a non-parametric decision tree as a nonlinear predictive modeling method~\cite{breiman1984classification}. 
The CART algorithm recursively partitions the input variable space into two smaller sub-regions. The goal is to find partitions that minimize the sum of squared error, SSE, of the output in the resulting sub-regions. Greedy recursive binary splitting approach is performed to identify an optimal splitting variable and an optimal cutpoint that leads to the greatest possible reduction in SSE. \new{The splitting process continues until a stopping
criterion is reached; for instance, CART may continue until each node has fewer than some
minimum number of observations, or the depth of the tree reaches a maximum. The final partitions after the CART algorithm stops are known as leaves or \emph{terminal nodes}}. As a result of the binary recursive splitting approach, the input points should have similar function values within one terminal node, and there should be a significant change in function values between different partitions.
% resulting CART predictive model is a decision tree, in which the the aforementioned binary splits define a tree logic. A predicted response is determined by following the tree logic to determine a final sub-region, which is referred to as a terminal node. 
The predicted output is the average output of the input points within each terminal node. %Because of the tree logic, the terminal nodes partition the input variable space. Moreover, because the same average response is used for prediction within each terminal node, the resulting CART predictive model is a step function.
In \S\ref{sec:tkmars}, we elaborate on our proposed eligible knot location approach for TK-MARS using CART. %The resulting TK-MARS approach is specifically designed for surrogate optimization.
%so that the interactions are manageable. %It is a classification method when the response variable is categorical, and in the case of having a continuous numerical response variable, CART predicts the response value in each terminal node. 
%For regression predictive modeling problems, CART chooses binary splits by minimizing the sum of the squared error, Equation~\ref{eq:LSE}, across all input points that fall within each partition. %The partitions discover the structure of the unknown function as we progress over time.
%\vspace{-3mm}
%\begin{align}\label{eq:LSE}
%LSE=\frac{1}{N}\sum_{v\in V}\sum_{k=1}^{K_v}(y^k_v-\bar{y}_v)^2,
%\end{align}
%where $y$ is the response value for each observations and $\bar{y}_v$ is the average response value of the observations in terminal $v$.
\vspace{-5mm}
\subsection{Exploration-Exploitation Pareto Approach}\label{sec:eepa}
\vspace{-3mm}

As we discussed in Algorithm~\ref{alg:SurOpt}, sampling input points in Step \ref{step5} is critically important to find promising points to evaluate in the black-box function and avoid unnecessary expensive evaluations. \hadis{In this research, we apply an Exploration-Exploitation Pareto Approach (EEPA)~\cite{dickson2014exploration, bischl2014moi, krityakierne2014global, dong2015kind, krityakierne2016sop} for sampling in Step~\ref{step5}.} \hadis{Let $P\subset R$, where $P$ is a set of new candidate points, and $R$ is a fixed pool of randomly generated candidate points.} 
To determine a set of new input points from $R$ to be evaluated, % $P$ in , 
%EEPA first generates a fixed random pool of sample points $R$. 
EEPA creates a Pareto frontier on the predicted function value using the surrogate model, as one dimension, and the distance of the candidate points from the already evaluated input points, as the second dimension. The first dimension exploits near-optimal areas, and the second dimension explores undiscovered regions of the input space.
In particular, for each input point $x\in R$, the exploration metric is $\delta(x)=\min_{\tilde{x}\in I}\|x-\tilde{x}\|$, and the exploitation metric is $\hat{f}(x)$. The first metric should be maximized, while the second should be minimized. The non-dominated Pareto set is given by $F = \{x \in R~|~\nexists~\tilde{x} \in R, \hat{f}(\tilde{x}) \leq \hat{f}(x), \delta(\tilde{x})\geq \delta(x)\}$. 
There might be candidate points on the Pareto set that are close to each other. \hadis{To find the final set of candidate points from the Pareto set, 
%a \emph{maximin} exploration technique is applied, so $P\subseteq F$. 
EEPA first selects the point with the minimum predicted function value, then it applies a \emph{maximin} function on the distance of the Pareto set points from the evaluated data points, to maximize the exploration further.} Dickson~\cite{dickson2014exploration} shows that EEPA outperforms pure exploration and exploitation methods. For more details, we refer to the pseudocode of EEPA in Appendix~\ref{app:EEPA}. 
\vspace{-5mm}
\subsection{Related Literature}\label{sect:literature}
\vspace{-3mm}
%Figure~\ref{fig:lr1} demonstrates in chronological order 
Various surrogate models have been used in the BBO literature. A few of them
employ non-interpolating models, such as MARS~\cite{costas2014multi, crino2007global, dickson2014exploration} and polynomial regression~\cite{moore2000q2}. %The rest focus on interpolating surrogates where all the input variables are considered to be significant to the output. 
Despite the aforementioned shortcomings of interpolation methods, such as Kriging~\cite{sacks1989design, jones1998efficient, simpson1998comparison, huang2006global, picheny2013quantile, dong2015kind} and RBF~\cite{liu2021trust,regis2016multi,regis2011stochastic, regis2014constrained, regis2014evolutionary, krityakierne2016sop, datta2016surrogate, muller2016miso}, their use within BBO is far more prevalent. %\new{HADIS WILL ADD: There are few research used non-interpolating RBF and GP}

%Kriging~\cite{sacks1989design, jones1998efficient, simpson1998comparison, huang2006global, picheny2013quantile, dong2015kind} is one the most prevalent traditional surrogate models used in the literature for black-box optimization. RBF, an interpolating technique, is the most common surrogate model that has been given attention lately~\cite{regis2011stochastic, regis2014constrained, regis2014evolutionary, krityakierne2016sop, datta2016surrogate, jakobsson2010method, muller2014influence, muller2016miso}. 
In addition, these methods consider every variable as an important input to the output, which is a weakness in the surrogate optimization literature. \hadis{A regularized RBF has been developed in the PYSOT toolbox \cite{eriksson2019pysot} to ensure a well-conditioned system for RBF approximation. However, its ability to screen important variables in surrogate optimization has not been discussed in the literature.} By contrast, Crino and Brown~\cite{crino2007global} demonstrate that MARS is capable of screening and reducing input variables using the parsimonious nature of MARS. 

% \begin{figure}[tb]
% \centering
% \includegraphics[width=\textwidth]{figures/LR1}
% \caption{Metamodeling Literature}
% \label{fig:lr1}
% \end{figure}

%From all the research displayed in Figure~\ref{fig:lr1},
In some expensive computer simulations, there is inherent noise associated with the system. Most of the existing methods cannot handle uncertainty~\cite{rojas2020survey,davis2007adaptive, grill2015black, moore2000q2}. \new{While the majority of surrogate optimization literature ignores uncertainty, some research has taken it into account \cite{wang2017novel,  crino2007global, costas2014multi, picheny2013quantile, jakobsson2010method, huang2006global}}. %Kriging and RBF do not inherently handle uncertainty, and some modifications are required.
Huang et al.~\cite{huang2006global} develops a Kriging-based surrogate optimization approach and applies it to low-dimensional test problems that include a low-level of noise. The computational effort of fitting Kriging increases for higher-dimensional problems. The proposed method for highly fluctuating functions under higher noise levels requires further investigation.
Picheny, Wagner, and Ginsbourger ~\cite{picheny2013benchmark} adds Gaussian noise with a fixed independent variance to the output of low-dimensional optimization test problems. The results show the relative poor modeling performance of Kriging. A large part of the variability that cannot be explained by the model is due to the observational noise during optimization. %Further, analysis of variance of modeling parameters is presented to show the significance of the modeling parameters, such as initial design, different infill criteria, the noise level of test functions, and covariance kernel. 
Jakobsson et al.~\cite{jakobsson2010method} and Picheny et al.~\cite{picheny2013quantile} apply RBF and Kriging based surrogate optimization methods on low-dimensional test problems under low levels of noise. 
%MARS and regression require no revision to handle uncertainty. 
Costas et al.~\cite{costas2014multi} shows MARS is preferable in real-world applications due to slope discontinuities and uncertainties, even though they do not study the effect of noise. %We show TK-MARS is able to handle uncertainties associated with the black-box systems.

A few other studies in surrogate optimization concentrate on sampling. \hadis{The main approach to select a new candidate point to sample is to solve an optimization sub-problem that is subject to exploration and exploitation constraints ~\cite{li2021surrogate,moore2000q2, regis2005constrained, knowles2006parego, muller2017gosac}. In these researches,
a weighted score of response surface predictions and a distance metric is minimized to choose the next sample point. %However, the methods solve the optimization sub-problem is computationally costly, and the complexity increases as the input dimension increases.
M\"uller et al.~\cite{muller2016miso, muller2015ch} applied a stochastic sampling approach by perturbing the variables of a BKS with a perturbation probability distribution.}
% \begin{figure}[tb]
% \centering
% \includegraphics[width=\textwidth]{figures/LR2}
% \caption{Surrogate Optimization Literature}
% \label{fig:lr2}
% \end{figure}
The Pareto-based candidate point sampling approach has recently begun %taken into consideration 
\cite{dickson2014exploration, bischl2014moi, krityakierne2014global, dong2015kind, krityakierne2016sop}. 
\vspace{-5mm}
\subsection{Contributions}
\vspace{-3mm}
In this research, we propose a new flexible, parsimonious, and non-interpolating surrogate model called TK-MARS, which can identify important input variables. We specifically design TK-MARS for surrogate optimization of black-box functions using a partitioning technique. Besides, we propose a smart replication strategy to mitigate the uncertainty associated with the black-box output. To evaluate the performance of surrogate optimization under uncertainty, we develop a new metric, called the Maximal True Function Area Under the Curve. We demonstrate the effectiveness of the proposed TK-MARS model and the smart replication approach using various complex global optimization test functions..% through a comprehensive set of experiments.
% We demonstrate that the modified surrogate model is capable of identifying significant variables and handling uncertainty.
% Furthermore, a centroid-based dynamic pool generation approach is developed to improve the candidate selection procedure. We propose an alternative approach to identify subregions of interest, which nicely surround promising points for subsequent function evaluations. 
% We introduce uncertainty associated with the black-box functions and propose effective strategies to handle it in the surrogate optimization procedure.
%\vspace{-5mm}
\section{Technical Description}\label{sec:contribution}
In this section, we describe the aforementioned contributions to surrogate optimization.
\vspace{-5mm}
\subsection{Tree-Knot Multivariate Adaptive Regression Splines}\label{sec:tkmars}
\vspace{-3mm}
%To address the two major issues in black-box optimization, namely uncertainty and insignificant inputs, a flexible, non-interpolating and parsimonious surrogate model $\hat{f}$ (Step 4 in Algorithm~\ref{alg:SurOpt}) can be used.
As we have discussed so far, in practice, there are unimportant input variables that are unknown \emph{a priori} in black-box systems. Consequently, we would like to use MARS as a surrogate model in Step~\ref{step4-3} of Algorithm \ref{alg:SurOpt}, since it can identify important input variables and screen unimportant ones. However, the original version of MARS is not customized for optimization. In this section, we develop Tree-Knot MARS (TK-MARS), a new version of MARS that is more efficient for surrogate optimization. More precisely, TK-MARS has fewer eligible knot locations but ones that are more promising for optimization. Moreover, developing a MARS model that uses fewer eligible knot locations avoids interpolation causing oscillations for highly noisy systems, as discussed in \S\ref{sect:intvnonint}.
% within the context of surrogate optimization.
%to identify significant variables and address the uncertainties associated with a black-box function. 
% In this section, we develop a modified version of MARS that is capable of identifying significant variables and specifically designed for black-box optimization. 

%The existing eligible knot location approaches has not to do with the optimization and are mainly focused on the model precision.
The motivation behind the new approach is to provide eligible knots around potential optimum input points. %In addition, MARS is sensitive to the order of the data points. The regression tree therefore helps bring in the near-optimal data points earlier, which accelerates the efficiency of MARS for optimization.
%So a new eligible knot location approach is required for this purpose. %The proposed approach considers identifying the new eligible knot locations plus sampling the new candidate points, simultaneously. 
%Along these lines, appropriate eligible knot locations are identified, and the approximation model is constructed gradually as more data are collected. Therefore, there is no need for many basis function evaluations and costly function evaluations. 
%The proposed eligible knot location approach is customized to locate knots for the optimization purpose by capturing peaks and valleys. Locating knots, where the response changes more frequently or significantly, identifies an optimal solution faster without focusing on the overall accuracy of the function estimation. %TK-MARS does not necessarily focus on function behavior but an optimum point. 
The proposed TK-MARS uses CART as a partitioning technique to capture the function structure and identify potentially near-optimal knot locations in each partition. 

%Specifically, the least squares error is given by Equation \ref{LSE}. 
%Classification and Regression Tree, CART~\cite{breiman1984classification} (CART) is a greedy partitioning technique to divide the input space with recursive binary splitting. All input variables and all possible split points are evaluated and chosen in a greedy manner. %As the data set is updated dynamically over time, more information from the black-box function structure is available. 

\begin{figure}[!tb]
\centering
\subfloat[Sin(x)]{\includegraphics[width=0.25\textwidth]{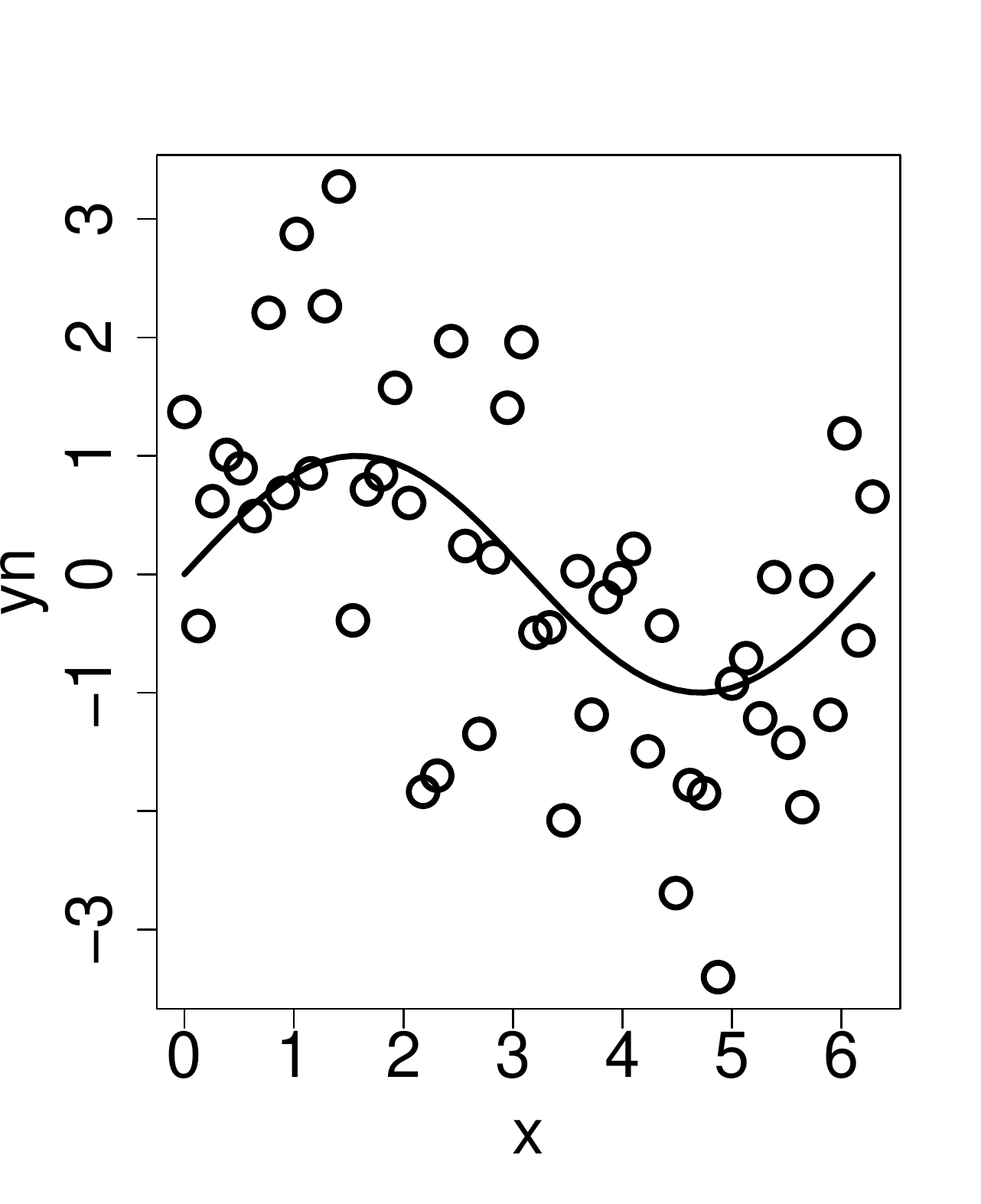}}
\subfloat[CART partitions]{\includegraphics[width=0.25\textwidth]{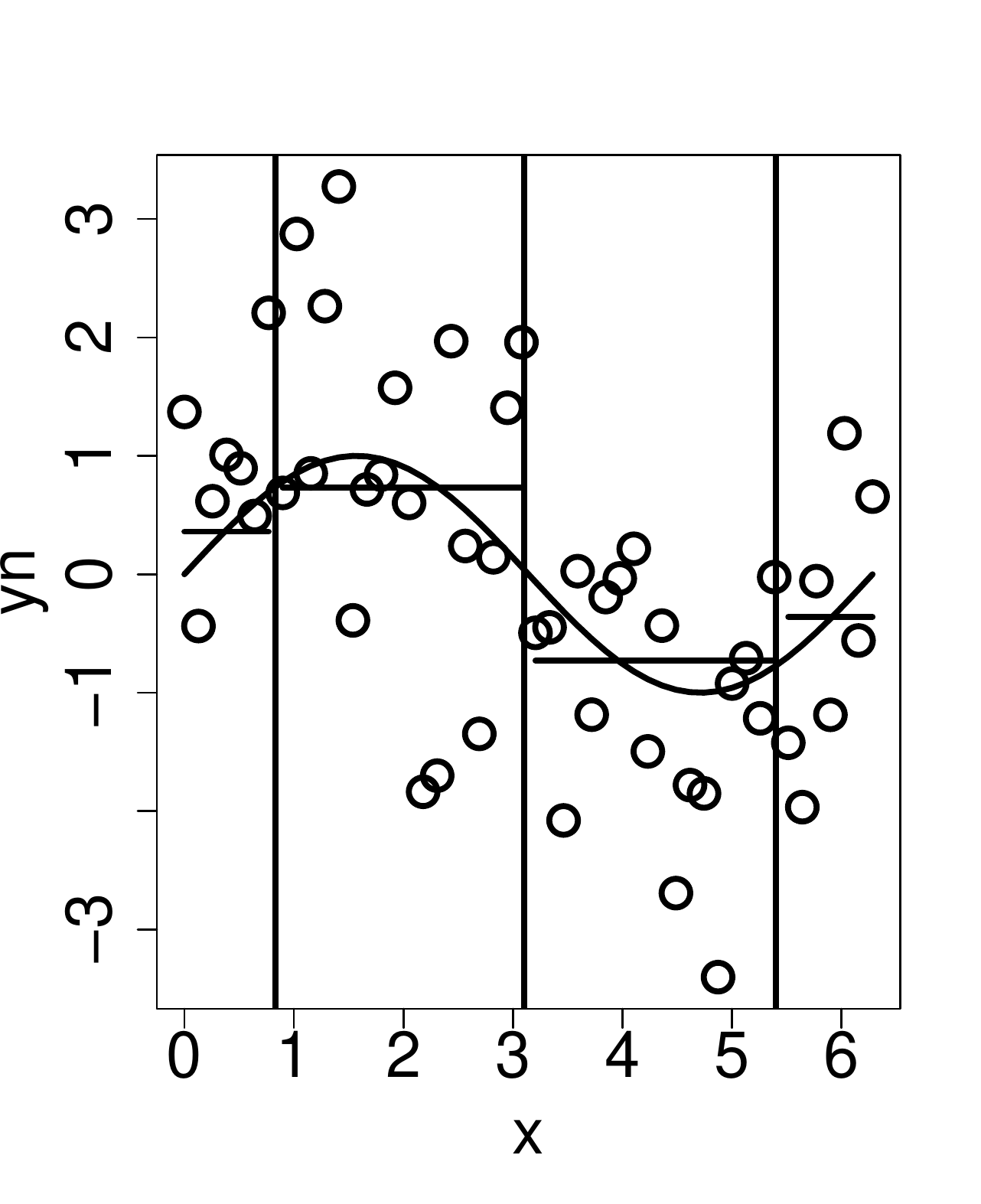}}
\subfloat[Centroids]{\includegraphics[width=0.25\textwidth]{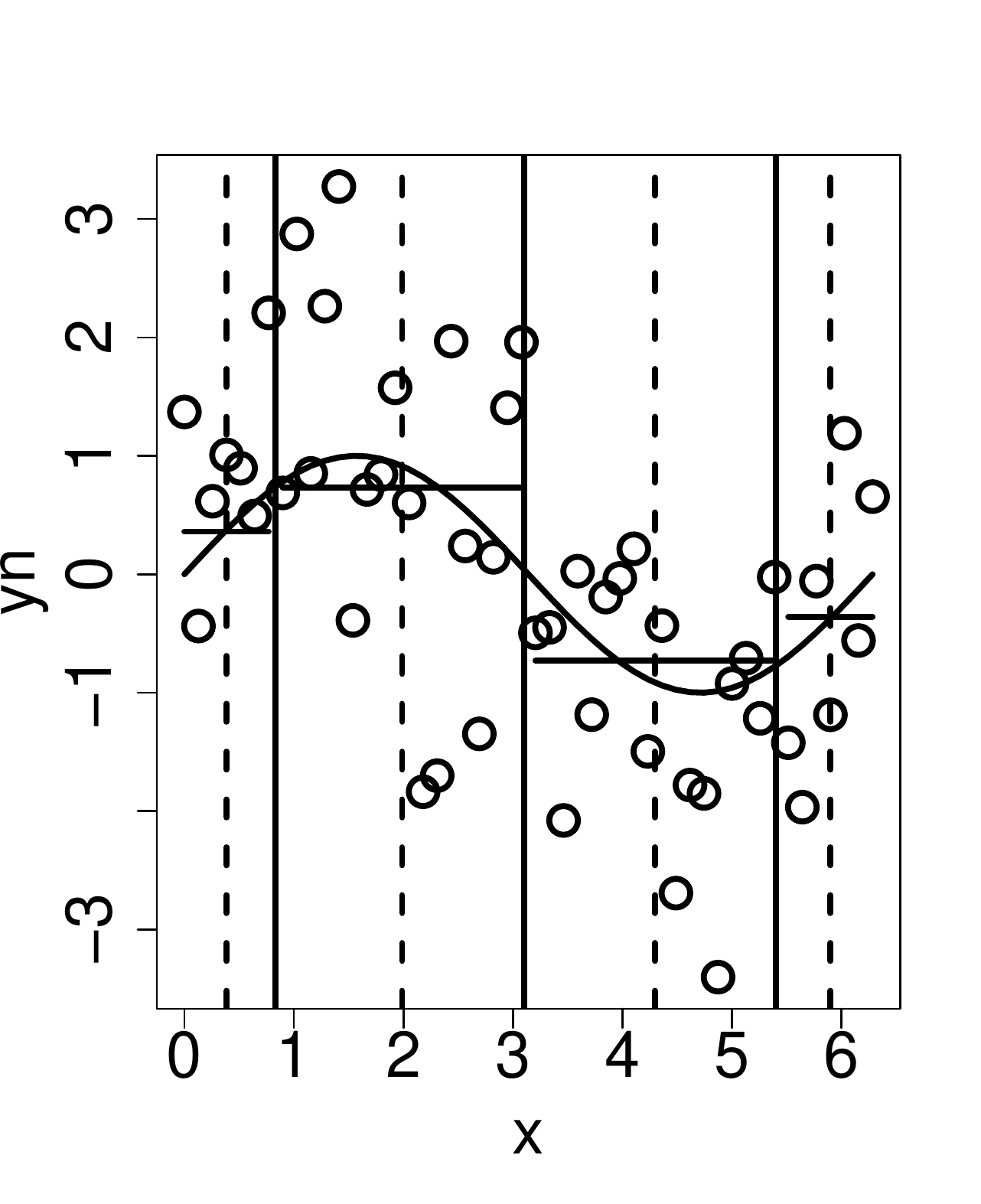}}
\subfloat[Knot locations]{\includegraphics[width=0.25\textwidth]{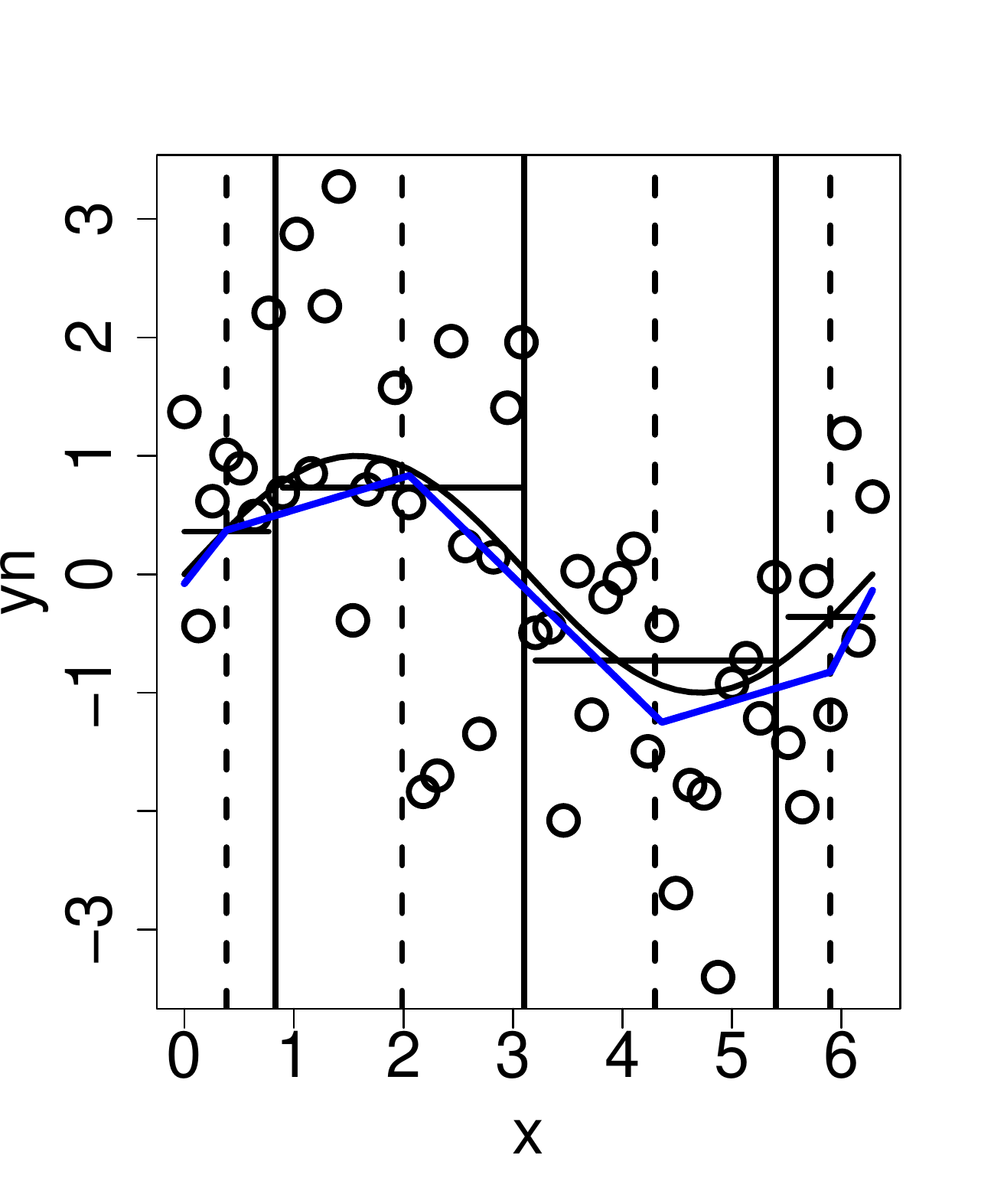}}
\caption{CART partitioning for eligible knot selection of TK-MARS}
\label{fig:cart}
\hspace{\fill}
% \vspace{-7mm}
\end{figure}

Before giving a detailed description of TK-MARS, consider the following high-level example with
%Consider %Suppose we have 
a data set generated from the function $f(x)=\sin(x)+\varepsilon$, where $\varepsilon$ is a Gaussian noise term with a mean of $0$, as in Figure~\ref{fig:cart}(a). 
Consequently, $E[f(x)] = \sin(x)$.
TK-MARS first uses CART to split the input points into four partitions, as in Figure~\ref{fig:cart}(b). The solid vertical lines in Figure~\ref{fig:cart}(b) represent the boundaries of the partitions, and the horizontal lines show the predicted average function value in each partition. \hadis{Note that CART partitions the data based on the split where there is a significant change in the underlying function $f$, i.e., points with a similar objective value (closer to a peak or closer to a valley) fall into a single partition. As a result, the peaks and valleys of $f$ occur near the centers of the two middle partitions in the example as in Figure~\ref{fig:cart}. As we discussed in \S~\ref{sec:cart}, the final partitions obtained based on the CART binary recursive splitting algorithm using a stopping criterion are terminal nodes.} Next, TK-MARS identifies the centroid of each terminal node and selects the closest input points to the centroids as eligible knot locations. The selected knot locations are represented by the dashed vertical lines as in Figure~\ref{fig:cart}(c). 
%We want to point out that the centroids are close to peaks and valleys. 
Since these centroids tend to be close to the peaks and valleys where optima lie, selecting points near them as eligible knot locations facilitates optimization.
Given the new set of eligible knot locations, TK-MARS approximates $f$, as $\hat{f}$, which is represented by the dotted function in Figure~\ref{fig:cart}(d). \hadis{The selected knot locations are still represented by the dashed vertical lines in Figure~\ref{fig:cart}(d).} The surrogate model $\hat{f}$ bends at near-optimal locations. 
%Note that the resulting TK-MARS model $\hat{f}$ bends at input points near the centroids, which include points near the peaks and valleys of $f$.
%the function as in  Figure~\ref{fig:cart}(c). 

%TK-MARS has a new efficient, eligible knot selection approach that is proposed for the purpose of surrogate optimization in this study. 
%Regarding the forward-backward procedure, MARS tends to have more knots where we have more information (dense regions) by minimizing the LOF.  As a result, it captures the function structure defining more knots in dense regions.
%Partitioning the data set in a way that captures the structure of the underlying function also finds promising eligible knot locations for MARS.

%The advantages of using a regression tree for partitioning data are:  (1) CART finds unbiased splits, (2) CART identifies the significant variables, and (3) CART identifies the interactions. 
%Partitioning the data set using CART, we are able to prioritize the eligible knot locations. %The data set is limited in the first iteration where there is not enough knowledge about eligible knot locations. 
%Updating the initial data set by adding more candidate points over time, the number of eligible knot locations increases.
%TK-MARS starts with a small set of eligible knot locations and slowly grows as the size of the data sets grows. This yields a more stable MARS model. 
%The training data set is updated gradually. As a result, the dense regions, where we have more information about the function, are identified by partitioning, and the structure of the underlying black-box function can be recognized. 

Specifically, let $I$ be the set of input points, let $x_j$ be the input variable in dimension $j$, $j=1,\ldots,d$, and let $V$ be the set of terminal nodes from CART. For each terminal node $v\in V$, let $ K_v\subset \{x^i|i = 1,\ldots, N\}$ be the set of observations in terminal node $v$. \new{Depending on the number of terminal nodes $|V|$, the $N$ number of points in $I$ is partitioned into $|V|$ subsets.} 
%We define \textbf{points}$(v)$ as the function that returns the observations in terminal node $v$}. 
Let $c^v_j$ be the $j$th dimension of the centroid in terminal node $v$. Let $x^k_j$ be the $j$th dimension of the $k$th observation in terminal $v$. Let $t^v_j$ be the index of the nearest input point to the centroid in terminal node $v$ for dimension $j$ and $t$ be the set of eligible knot locations, $\mathcal{T}=\{t_j^v, \forall v \in V, j=1,\ldots,d\}$. Algorithm~\ref{alg:TKMARS} presents the TK-MARS framework. % The representatives of each partition are considered as a reference point for eligible knot locations. 
As shown in Algorithm~\ref{alg:TKMARS}, CART partitions the data set based upon the function structure. 
%It splits more where the function structure changes to minimize the least square error. CART does not split where there are no significant changes in the response. Hence, there are more partitions in highly structured regions. As a consequence, TK-MARS defines more knots in highly structured regions. 
%In the beginning, the density is equal across terminal nodes since we start with a limited set of initial points that are uniformly distributed within the domain. The density changes as we add more points in the promising regions, over time.
%The peaks and valleys of the black-box function are often in the middle of the space defined by the tree logic at the terminal nodes rather than the edges since CART splits where the response values significantly change. Hence, 
TK-MARS then calculates
the centroid of each terminal node $v\in V$ 
%as calculated 
by Equation~(\ref{eq:centroid}).
%, are appropriate locations for the eligible knots
\vspace{-3mm}
\begin{align}
c^v_j=\frac{1}{|K_v|}\sum_{k=1}^{|K_v|} x^k_j\quad\forall v\in V, j=1,\ldots,d.
\label{eq:centroid}
\vspace{-5mm}
\end{align}
% where $c^v_j$ is the centroid of terminal $v$ for dimension $j$, and  $x^k_j$ is the $j$th dimension of observations in terminal $v$.
MARS considers a set of univariate knots in each dimension, so TK-MARS determines the input points near the centroids in each dimension $j=1,\ldots, d$ of each terminal node $v\in V$. %Knots in MARS are in one-dimensional space. %, i.e., each independent variable has a set of possible values for the eligible knot locations. The centroids of terminal nodes are in multi-dimensional space.
%A transformation is therefore needed from multi-dimensional to one-dimensional space. 
Specifically, Equation (\ref{eq:transformation}) describes the index $t^v_j$ of the input variable that is closest to the centroid of each terminal node $v\in V$ in each dimension $j=1,\ldots,d$. 
%represents a transformation function. The eligible knot location for each dimension is the closest point to the terminal centroid $v$ in the same dimension.
\begin{align}
\vspace{-3mm}
 t^v_j \in \arg\min_{k=1,\ldots,|K_v|}|x^k_j-c^v_j|\quad\forall v\in V, j=1,\ldots,d.
\label{eq:transformation}
 \vspace{-10mm}
\end{align}
% where $t^v_j$ is the nearest point to the centroid in terminal node $v$ for dimension $j$. 
When there are ties for the nearest input points, $t^v_j$ is the smallest index. Finally, we fit a MARS model $\hat{f}$ using the previously evaluated data set $(I, 
\mathcal{F})$ with the potential knot locations $\mathcal{T}$ to predict the output of the black-box function. % $t=\{x^{t_j^v}_j, \forall v\in V, \forall j=1\ldots, d\}$.
Note that TK-MARS uses $|V|$ potential knots in each dimension $j=1,\ldots d$.

\begin{algorithm}[!tb]
\caption{{\bf TK-MARS}\\
\small{{\bf input}: $(I,\mathcal{F})$}}
\hadis{
\begin{algorithmic}[1]
\label{alg:TKMARS}
    \STATE Fit CART on $(I,\mathcal{F})$
    \STATE $t=\emptyset$
    \FOR{$v\in V $}
    % \STATE $K_v=$ {\bf points}$(v)$ {\small \tt // points in node $v$}
        \FOR{$j=1$ to $d$}
            \STATE $c^v_j=\frac{1}{|K_v|}~\sum_{k=1}^{|K_v|} x^k_j$
            \STATE $t^v_j\in\argmin_{k=1,\ldots,|K_v|}|x^k_j-c^v_j|$
            \STATE $\mathcal{T}=\mathcal{T}\cup t^v_j$
        \ENDFOR
    \ENDFOR
    \STATE $\hat{f}(x)=MARS((I,\mathcal{F}),\mathcal{T})$
    \STATE {\bf return} $\hat{f}(x)$
\end{algorithmic}}
\end{algorithm}

% \begin{algorithm}[!tb]
% \caption{{\bf TK-MARS}\\
% \small{{\bf input: $I$} }}
% \begin{algorithmic}[1]
% \label{alg:TKMARS}
%     \STATE Execute the CART algorithm on  the data set;\\
%      $I=\{x^i\in D~|~\forall i=1,\ldots,N\}$
%     \STATE Find the centroid for each terminal node $v\in~V$\\% the set of terminal nodes \\
%     $c^v_j=\frac{1}{|K_v|}~\sum_{k=1}^{|K_v|} x^k_j\quad\forall v\in V, j=1,\ldots,d$
%     \STATE Determine the index of the closest point to centroids in each dimension;\\
%     $t^v_j\in\argmin_{k=1,\ldots,|K_v|}|x^k_j-c^v_j|\quad\forall v\in V, j=1,\ldots,d$\\
% % * <jaymrosenberger@gmail.com> 2018-07-24T16:40:58.274Z:
% % 
% % > t^v
% % The way you have this written is that it is the index, not the knot value. Is this really what you want? It makes Step 4 look weird. hadis: Yes I need the index
% % 
% % ^.
%     %(in the case of having ties randomly select one; smallest index)
%     \STATE Fit a MARS model $\hat{f}$ using eligible knot locations $x_j^{t_j^v}, \forall v \in V, j = 1, \ldots, d$.
%     % \STATE Fit MARS using the selected tree-based knots $t^v_j\quad\forall v=1,\ldots,V, j=1,\ldots,d$.
% \end{algorithmic}
% \end{algorithm}

In addition to assisting with developing the surrogate $\hat{f}$ (Step~\ref{step4} of Algorithm \ref{alg:SurOpt}) using TK-MARS, the centroids $C = \{c^v| v\in C\}$ can assist in the sampling of new candidate points $P$ (Step~\ref{step5-3}). As discussed in \S\ref{sec:eepa} and Appendix~\ref{app:EEPA}, EEPA uses a large fixed pool, $R$, from which it determines a Pareto set, $F$, from which it samples a set of new candidate points, $P$. However, EEPA and the quality of a best known solution found using EEPA depends on the fixed pool $R$,~\cite{dickson2014exploration}.
As shown in Fig.~\ref{fig:cart}, good candidate solutions are often near centroids, consequently, in this research, we dynamically augment the pool $R$ with the centroids; specifically $R=R\cup C$.

\vspace{-5mm}
\subsection{Smart-Replication}\label{subsec:smartRep}
\vspace{-3mm}
%\textcolor{blue}{In \S~\ref{sec:contribution}, we assume that the system is noise-free.} However, this assumption is not valid in many real-world applications. 
% In some expensive computer simulations, there is inherent noise associated with the system. Most of the existing methods cannot handle uncertainty~\cite{davis2007adaptive, grill2015black, moore2000q2}. 
We relax the assumption that the black-box function is deterministic in this section. Consequently, the black-box function output in Step~\ref{step6} includes uncertainty. Specifically, $\tilde{f}(x)=f(x)+\varepsilon$, where $\tilde{f}(x)$ is the output of the black-box system, and $\varepsilon$ is a random variable that follows an unknown probability distribution with a mean of 0 and a variance of $\sigma^2$. 
This implies that
%the black-box function output is a stochastic function that differs from the true value of the function.
each time a candidate point is evaluated a different output is obtained. 
%Therefore, for a single input, we have different outputs. 
Nonetheless, the goal is still to minimize the true objective function $f(x)$, Equation~(\ref{eq:objfunc}). %as in Equation~\ref{eq:objfunc_noisy}; however, only the simulated noisy function values (actual objective value) are accessible. 
Since there is uncertainty associated with the black-box output, a single evaluation of a candidate point might not be representative of the true output value. 
Consequently, the deterministic approach, which we also refer to as \emph{No-Replication}, may not be adequate to handle the uncertainty and, therefore, mislead the optimization process.
%This approach is sufficient for the systems with a low level of uncertainty and cannot handle a higher level of uncertainties. The deterministic approach may not be robust to the uncertainties associated with the black-box system.
%An alternative approach to cancel out the uncertainty effect is replicating the function evaluation for each point. 
%It may not be reliable and precise to evaluate a point once with a stochastic black-box function. 
Replicating the same input point multiple times and taking the mean of the outputs provides the most efficient estimator of the true objective value of $f$ at the input point. 
%We propose two distinct replication strategies as the alternative approaches for stochastic black-box optimization to minimize the output uncertainty. 
%An alternative approach is to replicate samples of input points to reduce the uncertainty.
%Assume that the actual objective function of the black-box is $\tilde{f}(x)$ where $\tilde{f}(x)\sim (f(x),\sigma^2)$. 
%Suppose the number of replications is $r$. 
Let $x^i$ be the $i$th input point in the set of evaluated input points $I$. 
%Let $I(x^i)$ be a set of the points that are equal to $x^i$. 
\new{Let $\varepsilon^1_i, \varepsilon^2_i, \ldots, \varepsilon^{r^i}_i$ be a random sample of size $r^i$ of the noise term $\tilde{f}(x^i)-f(x^i)$. We assume that the sample $\varepsilon^1_i, \varepsilon^2_i, \ldots, \varepsilon^{r^i}_i$ is independent and identically distributed. Consequently, we can estimate $f(x^i)$ as the sampled mean $\bar{f}(x^i)$ using Equation (\ref{eq:fbar}).}
%
%Namely, let \[
%I(x^i)=\left\{ x^k\left|x^k = x^i, \forall x^k\in %I\right.\right\}, \forall x^i\in I.
%\]
%Let us define $r(x^i)$ as the number of replications for $x^i$, which is $r(x^i) = \left|I(x^i)\right|, \forall x^i\in I$. For simplicity we will refer to $r(x^i)$ as $r^i$.
%In other words, $r(.)$ indicates the number of replications of sample points from the set of already evaluated points, $I$.
%Let $\tilde{f}(x)$ be the output of the black-box function with uncertainty, $\tilde{f}(x^i)=f(x^i)+\epsilon,~\forall i=1,\ldots,|I|$. 
%Let $\bar{f}(x^i)$ be the sample mean of $\tilde{f}(x^i)$ after $r(x^i)$ replications.
\vspace{-3mm}
\begin{align}\label{eq:fbar}
\new{\bar{f}(x^i)=f(x^i)+\dfrac{\sum_{k=1}^{r^i}\varepsilon^k_i}{r^i},~\forall x^i\in I.}
\end{align}
%\vspace{-3mm}
$\bar{f}(x^i)$ in Equation~(\ref{eq:fbar}) follows \new{an unknown distribution, with $f(x^i)$ mean and $\frac {\sigma^2} {r^i}$ variance},
based on the central limit theorem. 
%$\epsilon(x^k)$ is the uncertainty component of $\bar{f}(x^k)$. 
%The second component in Equation~(\ref{eq:fbar}) corresponds to the mean of replications' uncertainties for an input point. 
%$(x^j=x^i)$ is an indicator function that counts the number of replications for a single input point. 
If $r^i=r, \forall x^i \in I$, we refer to this approach as \emph{Fixed-Replication}, where we evaluate input points a fixed number of times, $r$. 
For noisy black-box systems, we no longer have a best known solution because we do not have the \new{true objective value} of $f(x)$. Consequently, we refer to a solution $x$ with the best $\bar{f}(x)$ as the \emph{best sampled mean solution (BSMS)} in Step~\ref{step8} of Algorithm~\ref{alg:SurOpt}.

%In the replication strategies, we track the BKS based on $\bar{f}$. 
%Consequently, we define a \emph{best sample mean solution (BSMS)} to refer to the BKS after replication. %For simplicity, in the remainder of the paper we use $BSMS$ to also refer to best known solution in deterministic case, where $r=0$.

%Based on the level of uncertainty and the maximum number of function evaluations, one may move along between exploring more points with single evaluation and exploiting a few points by replicating more than once. 
%Although the effect of noise decreases with replication and the sampled mean converges to the true function value, when the level of uncertainty is high, it needs more function evaluations.
%As a result, the replication effect is not significant; exploration, on the other hand, is more effective.
%The trade off between replicating the same point and evaluating a new point is the key idea in our Smart-Replication approach for noisy cases.

%\textcolor{red}{Hence, when the uncertainty level is high, the replication does not reduce the variance significantly.
%As a result, exploring more single evaluated points will be a better option.}

%So far, we explained No-Replication and Fixed-Replication approaches. 
To avoid unnecessary expensive evaluations in Step~\ref{step6} of Algorithm~\ref{alg:SurOpt}, we propose a novel strategy called \emph{Smart-Replication}, which replicates not all the candidate points but only the promising points for optimization % using hypothesis testing 
based on confidence intervals around the sampled means.
%Following this strategy, in this section, we propose an approach for replication, using hypothesis testing based on confidence intervals.
%Smart-Replication replicates only the promising points close to BSMS. %The number of replications is determined by the algorithm itself, following the hypothesis rule.
%Smart-Replication 
%and automatically chooses the number of replications for each selected candidate point following the hypothesis rule.
Let $x^o$ be a BSMS, and let $x^i\in P$ be a candidate point, 
%For each selected candidate point $x^i\in P$, Smart-Replication considers the following null hypothesis and stops replicating $x^i$, if it can reject the null hypothesis.%
%\vspace{0.02in}
%\medskip\noindent
%\framebox[\columnwidth]{\parbox{0.9\columnwidth}{ %\textbf{\textsc{}}
%\\
%{\bf H$_o$}: if $\bar{f}(x^o) < \bar{f}(x^i)$, where $x^i\neq x^o,~\forall i=1,\ldots,|I|$, then $f(x^o)>f(x^i)$. \\
%{\bf H$_o$}: $f(x^o) > f(x^i)$ versus 
%{\bf H$_1$}: if $\bar{f}(x^o) < \bar{f}(x^i)$, where $x^i\neq x^o,~\forall i=1,\ldots,|I|$, then $f(x^o)<f(x^i)$. \\
%{\bf H$_1$}: $f(x^o) \leq f(x^i)$.\\
%{\bf Decision rule}: If $CI_{low}(f(x^i))~\geq ~CI_{up}(f(x^o)), \forall x^i \in I~$, Reject $H_o$.\\
%{\bf Decision rule}: If $CI_{low}(f(x^i))~\geq ~CI_{up}(f(x^o))$, Reject $H_o$.\\
%{\bf Conclusion}: Even though the objective value of BSMS, $\bar{f}(x^o)$, is smaller than the objective value of $x^i$, $\bar{f}(x^i)$, the true objective value of $x^i$, $f(x^i)$, is smaller than the true objective value of $x^o$, $f(x^o)$. Hence, we are $(1-\alpha)\%$ confident that more replications are required for $x^i$.
%}}\\
% \vspace{-2mm}
\hadis{ and let $r^o$ denote the number of replications for $x^o$}. Suppose $r^o,r^i\geq 2$. Then, we are $100(1-\alpha)\%$ confident that $f(x^i)$ is within the interval $\bar{f}(x^i)\pm t_{\alpha/2, r^i-1} \frac{s^i}{\sqrt(r^i)}$, where $s^i$ the sample variance of $f(x^i)$, which is given by the Equation (\ref{stddev}), $t_{\alpha/2}$ is the critical value of the student's t-distribution, and $\alpha$ is the significance level. \hadis{We denote the sample variance of $\tilde{f}(x^o)$ by $s^o$.}
\begin{equation}\label{stddev}
    \new{s^i = \sqrt{\frac{1}{r^i-1} \sum_{k=1}^{r^i} \left(f(x^i)+\varepsilon^k_i- \bar{f}(x^i) \right)^2 }.}
\end{equation}

\new{Let $\mathit{CI}^i_{\mathit{low}}$ and $\mathit{CI}^i_{\mathit{up}}$ be the lower and upper limits of the aforementioned confidence interval of $f(x^i)$. Consequently, we can stop replicating a candidate point $x^i$ if $\mathit{CI}^i_{\mathit{low}}\geq \mathit{CI}^o_{\mathit{up}}$, i.e., if the lower bound of the objective value for point $x^i$ is smaller than the upper bound of the best sampled mean solution, BSMS, there is no need to replicate
$x^i$.}

\begin{figure}[H]
\centering
        \includegraphics[width=0.5\textwidth]{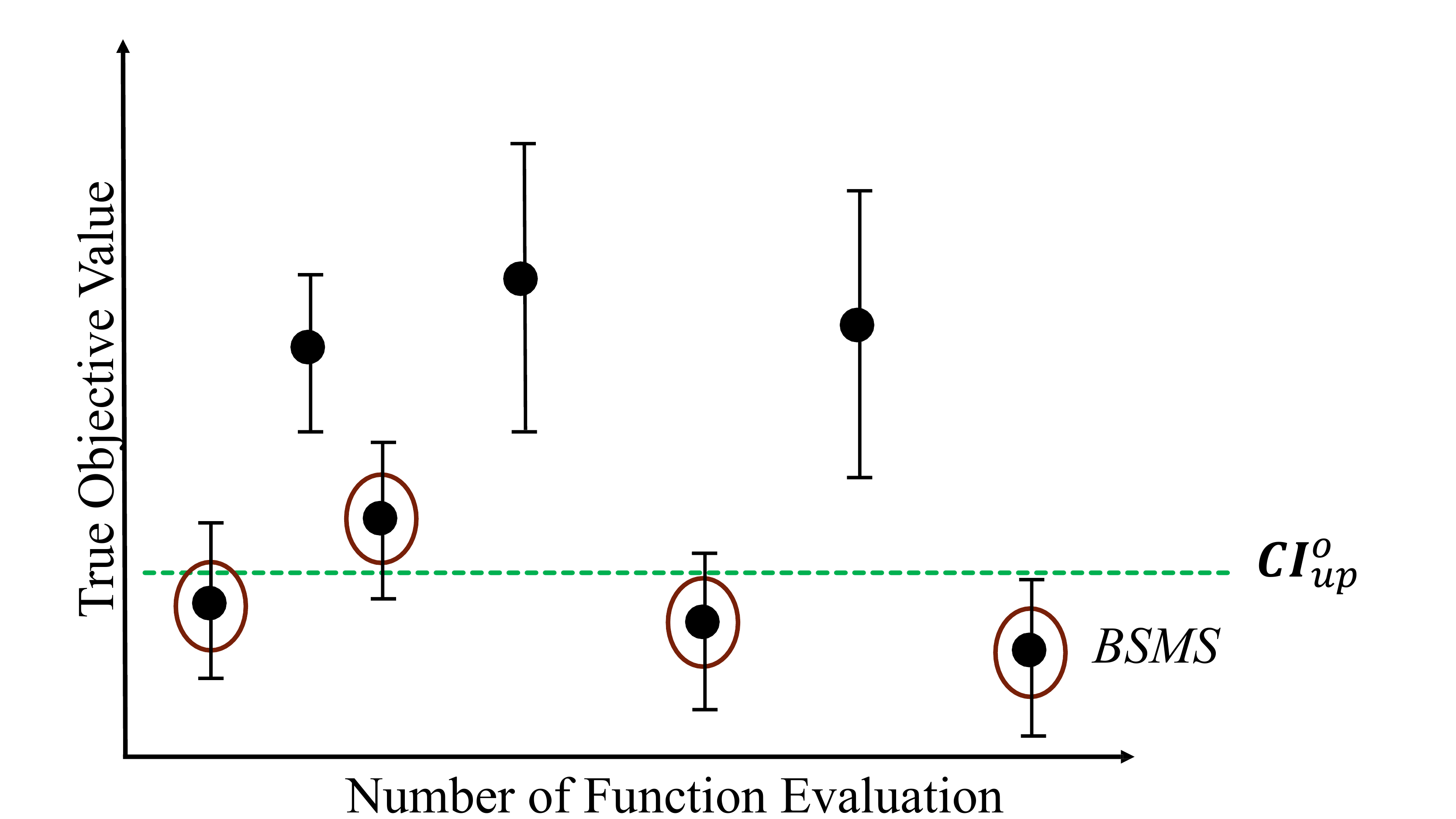}
    \caption{Confidence interval threshold mechanism of the Smart-Replication Approach}
    \label{fig:smartrep}
    \hspace{\fill}
    \vspace{-10mm}
\end{figure}

%[DESCRIBE FIGURE 3.]
Figure~\ref{fig:smartrep} presents an example of the Smart-Replication approach. The black bars show the confidence intervals for \hadis{the true objective value}, $f$, at each input point in $I$. The dashed horizontal line is at the $\mathit{CI}_{\mathit{up}}^o$ value. The promising points are the candidate points that have confidence intervals that overlap with that of the current BSMS, i.e., $\mathit{CI}^i_{\mathit{low}}< \mathit{CI}^o_{\mathit{up}}$. Consequently, the circled points are replicated further. 

Since black-box function evaluations are expensive, for a given candidate point $x^i$, Smart-Replication limits the number of replications $r^i$ to a maximum number of replications denoted as $r_\mathit{max}$. 
Specifically, Algorithm~\ref{alg:tksrg} shows the proposed surrogate optimization approach, which uses Smart-Replication (Algorithm~\ref{alg:smartrep}) for evaluation, Step~\ref{step6}.

%For simplicity we refer to $f(x^i)$ as $f^i$. For each point that is evaluated~$r^i$~times, $r^i \geq 2$, we calculate the standard deviation $s$ and the confidence interval $CI$ with a significance level of $\alpha$: $s(f^i)=\sqrt{\frac{1}{r^i-1}\sum_{k=1}^{r^i}(\tilde{f}(x^k)-\mu)^2}$, where $\mu=\frac{1}{r^i}\sum_{k=1}^{r^i}\tilde{f}(x^k)$. $\tilde{f}(x^k)$ is the noisy evaluation.

%Now, we calculate the confidence interval, $CI(f^i)=(\mu~\mp~t_{\frac{\alpha}{2}}~\frac{s(f^i)}{\sqrt[]{r^i}}), \forall i=1,\ldots,|I|$.
%To select the promising points, we quit replicating the candidate points based on the lower bound of the confidence interval $CI_{low}$, as it can be observed in Figure~\ref{fig:smartrep}. The horizantal dashed line indicates the threshold  $threshold=CI_{up}(BSMS)$. If the $CI_{low}$ of a point is less than the threshold, it is selected as a promising point to be replicated. 
%The promising points are the candidate points that are close to the current BSMS, i.e., below the threshold. The selected points are the circled points in Figure~\ref{fig:smartrep}. 
%Since the number of function evaluations is limited, and the experiments are costly, we consider a maximum number of replications for the promising points $r_\mathit{max}$. 

% \begin{algorithm}[!tb]
% \caption{\bf Smart-Replication Approach}
% \begin{algorithmic}[1]
% \label{alg:smartrep}
% \STATE Sample initial design space $I=\{x^i_j \in D ~|~ \forall i=1,\ldots ,d, j=1,\ldots ,N\}$
%     \STATE Randomly generate m uniform random points in the box region;\\
%     $R=\{u_{j}^{i}~|~a \leq u_{j}^{i} \leq b,~\forall~i=1,\ldots,d, j=1,\ldots,M\}$
% \WHILE {Termination criteria is not satisfied}
% 	\STATE $\mu^i=\frac{1}{r^i}\sum_{k=1}^{r^i}\tilde{f}^k$
% 	\STATE $s(f^i)=\sqrt{\frac{1}{r^i-1}\sum_{k=1}^{r^i}(\tilde{f}^k-\mu^i)^2}$ 
%     \STATE $CI(f^i)=(\mu~\mp~t_{\frac{\alpha}{2}}~\frac{s(f^i)}{\sqrt[]{r^i}}), \forall i=1,\ldots,M$
%     \WHILE {$CI_{low}(f(x^i))~\geq ~CI_{up}(f(x^o))~\&~r^i\leq r_\mathit{max} $}
%     \STATE Evaluate $x^i$, $r^i=r^i+1$
%     \STATE Update $std(x^i), \mu^i, CI(x^i)$, 
%     \ENDWHILE
%   	\STATE Construct CART on initial input set $I$\\
%     \STATE Find the centroids for terminal nodes\\
%     $c^v_j=\frac{1}{|K_v|}~\sum_{k=1}^{|K_v|} x^k_j,~\forall v\in V, j=1,\ldots,d$
%     \STATE Construct a surrogate model on $I$~(TK-MARS/RBF)

%     \STATE Add centroids to the $R$;$~R={R~\cup~C}$ where $C=\{c^v~|~\forall v\in V\}$
%     \STATE Optimizer (EEPA) on $R$ to determine new candidate set of points $P$
%     \STATE Update initial input set $I=I\cup P$
%     \STATE Find BSMS, $(x^o,f(x^o))$
%     \STATE t:=t+1
% \ENDWHILE
% \end{algorithmic}
% \end{algorithm}
%Algorithm~\ref{alg:tksrg} shows the proposed surrogate optimization approach. It uses Algorithm~\ref{alg:smartrep} for evaluation, Step~\ref{step6-3}.

\begin{algorithm}[!tb]
\caption{\bf Proposed Surrogate Optimization Approach}
\begin{algorithmic}[1]
\label{alg:tksrg}
    \STATE $I=\{x^1, \ldots, x^N\}$, a set of $N$ input points in $D$, selected with a DOE method \label{step1-3}
    \STATE $\mathcal{F}=\{\tilde{f}(x^i)|x^i \in I\}$ \label{step2-3}
    % \STATE Randomly generate $M$ uniform random points in $D$;\\
    % $R=\{u_{j}^{i}~|~a \leq u_{j}^{i} \leq b,~\forall~i=1,\ldots,d, j=1,\ldots,M\}$
    \WHILE {Termination criteria is not satisfied}
    	\STATE \emph{Surrogate}: $\hat{f}(x)=$TK-MARS$(I,\mathcal{F})$\label{step4-3}
        \STATE \emph{Sampling}: $P = EEPA(R, I,\hat{f}(x))$  \label{step5-3}
        % \STATE \emph{Evaluation}: $\mathcal{F}_P=\{\tilde{f}(x^i)|x^i \in P\}$\\ \label{step6-3}
        % \hspace{16mm}$I= I\cup P, \mathcal{F}=\mathcal{F}\cup\mathcal{F}_P$\\
        % \hspace{15mm} SmartRep$(I,\mathcal{F})$ \label{step6-3}
         \STATE \emph{Evaluation}: Smart-Replication$(I,\mathcal{F},P)$ \label{step6-3}
        \STATE  \emph{Best Sampled Mean Solution (BSMS)}: $x^o\in\argmin_{x\in I} \bar{f}(x)$ \label{step7-3}
    \ENDWHILE
    \STATE {\bf return} $x^*=x^o$
\end{algorithmic}
\end{algorithm}

\begin{algorithm}[!tb]
\caption{\bf Smart-Replication\\
\small{{\bf input}: $(I,\mathcal{F},P)$ }}
\begin{algorithmic}[1]
\label{alg:smartrep}
\STATE $\mathcal{F}_P=\{\tilde{f}(x^i)|x^i \in P\}$\\ \label{step6-3}
\STATE $I= I\cup P; \mathcal{F}=\mathcal{F}\cup\mathcal{F}_P$\\
\STATE \hadis{$\mathit{CI}_{\mathit{up}}^o={f}(x^o)$}%+~t_{\frac{\alpha}{2}}~\frac{s^o}{\sqrt[]{r^o}}$
\STATE \hadis{$\mathit{CI}_{\mathit{low}}^i=-\infty, \mathit{CI}_{\mathit{up}}^i=+\infty ,r^i=1, s^i=1, \forall i\in I$.}
\FOR{$i=1$ to $|I|$}
    %\WHILE {$CI_{low}(f(x^i))~\geq ~CI_{up}(f(x^o))~\&~r^i\leq r_\mathit{max} $}
     \WHILE {$\mathit{CI}_{\mathit{low}}^i~< ~\mathit{CI}_{\mathit{up}}^o~\&~r^i\leq r_\mathit{max} $}
        \STATE %Evaluate $x^i$, $\tilde{f}(x^i)$ 
        \emph{Evaluation}: $\mathcal{F} = \mathcal{F} \cup \{\tilde{f}(x^i)\}$  
        \STATE $r^i=r^i+1$
        \STATE \new{$\bar{f}(x^i)=f(x^i)+\dfrac{\sum_{k=1}^{r^i}\varepsilon^k_i}{r^i}$}
    	\STATE \new{$s^i = \sqrt{\frac{1}{r^i-1}\sum_{k=1}^{r^i}\left(f(x^i)+\varepsilon^k_i-\bar{f}(x^i) \right)^2}$}
        \STATE $\mathit{CI}_{\mathit{low}}^i=\bar{f}(x^i)~-~t_{\frac{\alpha}{2}}~\frac{s^i}{\sqrt[]{r^i}}$
         \STATE $\mathit{CI}_{\mathit{up}}^i=\bar{f}(x^i)~+~t_{\frac{\alpha}{2}}~\frac{s^i}{\sqrt[]{r^i}}$
                % \STATE Update $std(x^i), \mu^i, CI(x^i)$, 
    \ENDWHILE
    \STATE $x^o\in\argmin_{x\in I} \bar{f}(x)$
    \STATE $\mathit{CI}_{\mathit{up}}^o=\bar{f}(x^o)+~t_{\frac{\alpha}{2}}~\frac{s^o}{\sqrt[]{r^o}}$
   % \STATE $f(x^i) = \bar{f}(x^i)$
\ENDFOR
\end{algorithmic}
\end{algorithm}

Smart-Replication is similar to the No-Replication approach when the uncertainty level is low and similar to the Fixed-Replication approach when the uncertainty level is high. Specifically, 
for a high noise level, the variance is larger, so more replications are required, and most candidate points are replicated up to the maximum number of replications $r_{\mathit{max}}$.
%, which needs to be large and determined in advance.
Consequently, Smart-Replication adjusts its behaviour based on the uncertainty level of the system. %and is efficient for optimizing systems with either low or high levels of noise. %That is because it recognizes if the system is deterministic or stochastic and determines how to replicate dynamically.
\vspace{-5mm}
\subsection{Performance Metrics for Surrogate Optimization}\label{mesureofperformance}
% \vspace{-3mm}
%In this section, we propose two measures of performance for deterministic and stochastic systems.
%Minimizing the number of function evaluations is the goal of the expensive black-box optimization. 
Since black-box function evaluations are expensive, the number of such evaluations performed before finding an optimal solution is a common metric to test a surrogate optimization algorithm.
However, obtaining a global optimum cannot be guaranteed. 
A metric that can quantify the improvement in the true BSMS objective value over fewer black-box function evaluations is more appropriate in the context of BBO. \new{Area Under the Curve has been used in black-box optimization to measure how quickly an algorithm improves solution quality compared to other algorithms \cite{dewancker2016stratified}. AUC is also is commonly used in other optimization domains \cite{memarian2019optimization} and in statistics to compare models with the two metrics sensitivity and specificity \cite{hanley1982meaning, delong1988comparing}.} Specifically, consider an executed surrogate optimization algorithm, and let $\mathcal{I}$ be the total number of black-box function evaluations conducted during the execution. The \emph{Area Under the Curve (AUC)} is given by definition \ref{def:auc}.

%\subsection{Area Under the Curve: AUC}\label{sec:auc}
%First, we define the area under the curve (AUC) metric in Definition~\ref{def:auc}.
%In order to make the AUC comparable for different functions, we scale the objective value of BSMS.

\begin{definition}[Area Under the Curve -- AUC]\label{def:auc} 
Let $x^{oi}$ be a BSMS found after $i$ black-box function evaluations, and let $f(x^{oi})$ be the true objective value of the BSMS.
Let
$f^{min}$ be the \new{optimal} objective value of the optimization problem (\ref{eq:objfunc})-(\ref{eq:boxconst}), and let 
$f^{max}=\max\{f(x^{oi})| \forall i = 1\ldots \mathcal{I}\}$.
The normalized objective value of BSMS is:
$$\textcolor{blue}{{\check{f}}}(x^{oi})=\frac{f(x^{oi})-f^{min}}{f^{max}-f^{min}},~\forall i=1,\ldots,\mathcal{I}.$$
The AUC is given by
\vspace{-3mm}
\begin{align}
\label{eq:auc} & \mathit{AUC}=\sum_{i=1}^{\mathcal{I}} {\check{f}}(x^{oi}). 
\end{align}
\end{definition}
%AUC comprises both the quality of BSMS, as well as the time that the algorithm finds it.
%The total number of function evaluations is limited by financial, time or resource constraints in black-box optimization. 
AUC measures the performance of a surrogate optimization algorithm well in a \emph{deterministic environment} with no noise (i.e., $\sigma^2=0$), since the objective value of the BSMS monotonically decreases over each black-box function evaluation.
However, this monotonicity may not hold for black-box systems with uncertainty, as the true objective value of the BSMS may oscillate.
Even though oscillations in early iterations are tolerable, we would like stable behavior towards the end of the algorithm execution. %Consequently, the stability of the true objective value of the BSMS represents the robustness of the algorithm.
%A good algorithm is the one with fewer and shorter oscillations, in which after a reasonable number of function evaluations, the results are reliable. 
Consequently, we propose a metric for black-box systems with uncertainty that penalizes instability by considering the maximum true objective value of all BSMS found at the current and subsequent black-box function evaluations of the algorithm. 
%\subsection{Maximal True Function Area Under the Curve: MTFAUC}\label{sec:mtfauc}
%To consider the instability in the metric, we consider the maximum objective value of BSMS obtained among all BSMS found forward.
%Subsequently, 
Specifically, we define the \emph{Maximal True Function Area Under the Curve (MTFAUC)} in definition \ref{def:mtfauc}.

\begin{definition}[Maximal True Function Area Under the Curve -- MTFAUC]\label{def:mtfauc}
%Given the best sampled mean solution (BSMS) found by an algorithm after each function evaluation. 
%Let $f(x^{oi})$ be the true objective value of BSMS, where $x^{oi}$ is the BSMS after $i$ function evaluations.
%Let $f^{min}=\min_{i=1,\ldots |I|}(f(x^{oi}))$ and %$f^{max}=\max_{i=1,\ldots |I|}(f(x^{oi}))$.
%Using the normalized objective value of BSMS 
%$${f}(x^{oi})=\frac{f(x^{oi})-f^{min}}{f^{max}-f^{min}}~\forall i=1,\ldots,|I|,$$
%$\bar{f}(x)$ (Equation~\ref{eq:fbar}) is the sample mean of the observed objective value after $r$ replications. 
%Let $x^{oi}= \argmin_{j=1,\ldots,i} \bar{f}(x^j)$ be the BSMS after $i$ number of function evaluations and $\bar{f}(x^{oi})$ be its estimated objective value.
Let $\hat{j}(i)$ be the index of the black-box function evaluation of the maximum true objective value among all BSMS found in function evaluations after function evaluation $i-1$; that is,
\begin{equation}
\hat{j}(i) \in \argmax_{j=i,\ldots,\mathcal{I}}\check{f}(x^{oj}).
\end{equation}
The MTFAUC is given by
\vspace{-5mm}
\begin{equation}
\mathit{MTFAUC}=\sum \limits_{i=1}^{\mathcal{I}}\frac{\check{f}\left(x^{\hat{j}(i)}\right)+\check{f}\left(x^{\hat{j}(i-1)}\right)}{2}.
\end{equation}
\end{definition}
MTFAUC penalizes the instability of the true objective value of BSMS, by using $\hat{j}(i)$ to consider subsequent oscillations.
Hence, MTFAUC increases if the true objective value oscillates. 
A surrogate optimization algorithm that is less sensitive to uncertainty has a lower MTFAUC value than one that is more sensitive. MTFAUC is equal to AUC in a deterministic environment.
%\vspace{-5mm}

\begin{figure}[H]
\centering
\subfloat[]{\includegraphics[width=0.35\textwidth]{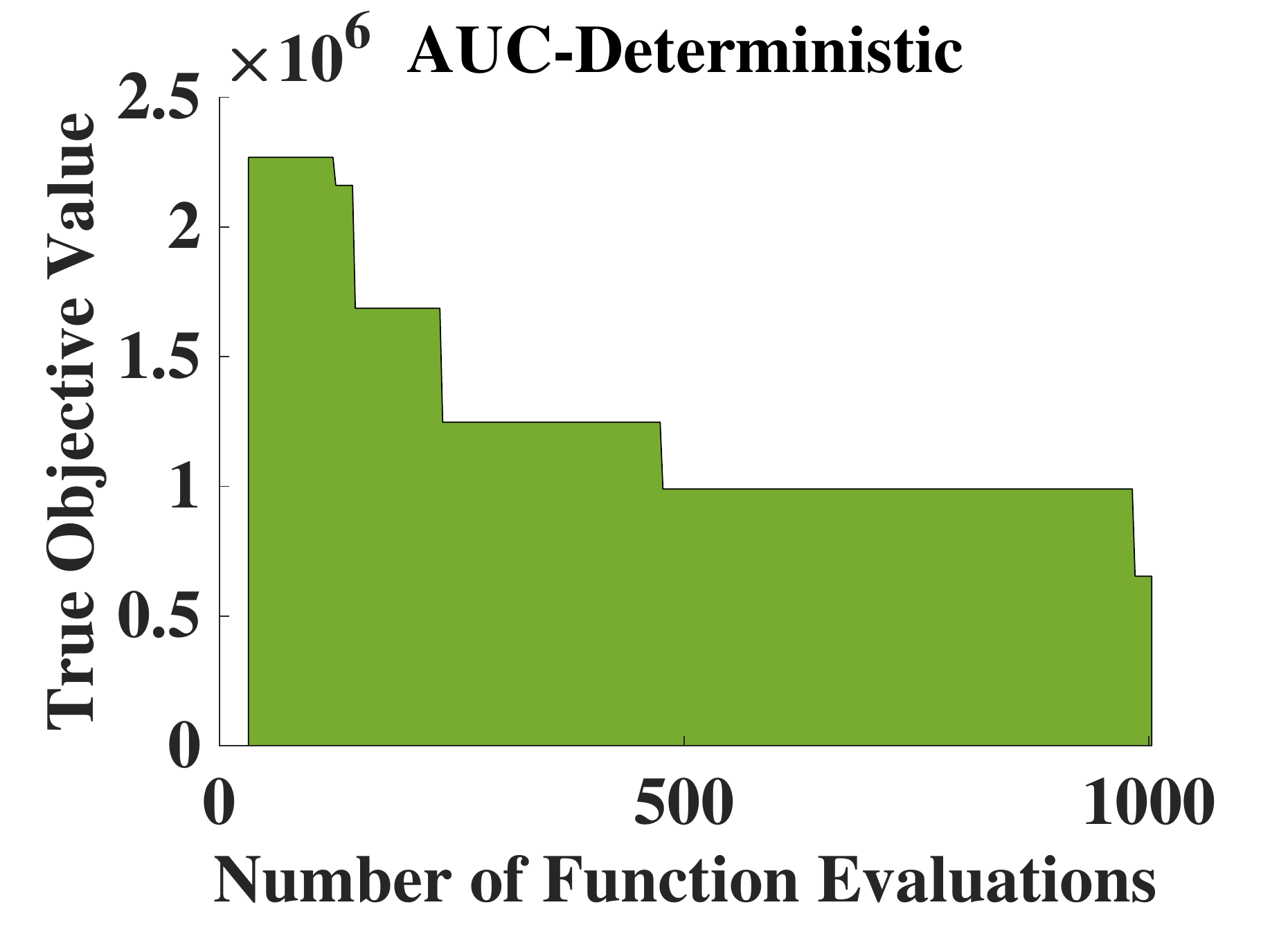}}
\subfloat[]{\includegraphics[width=0.35\textwidth]{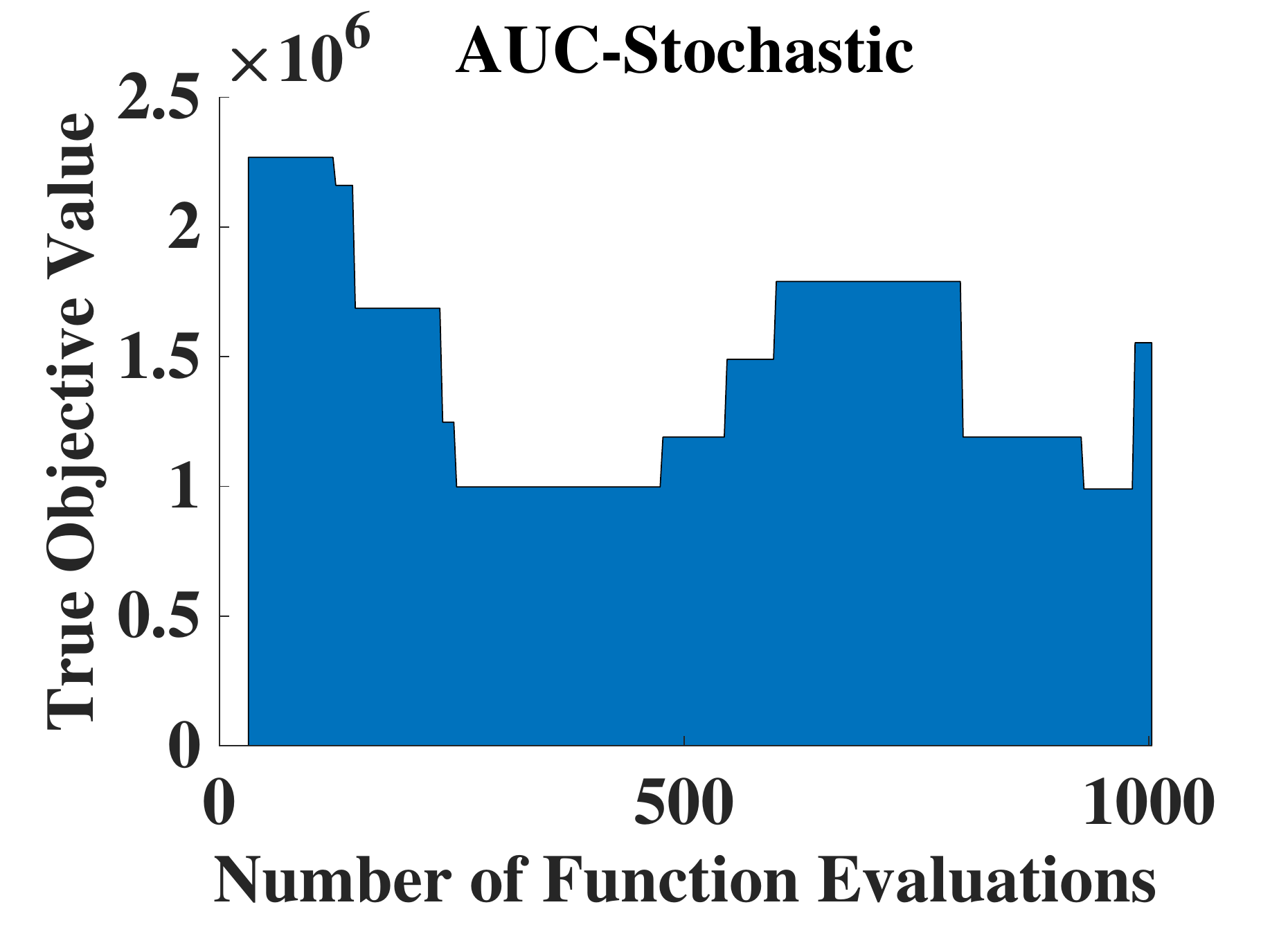}}
\subfloat[]{\includegraphics[width=0.35\textwidth]{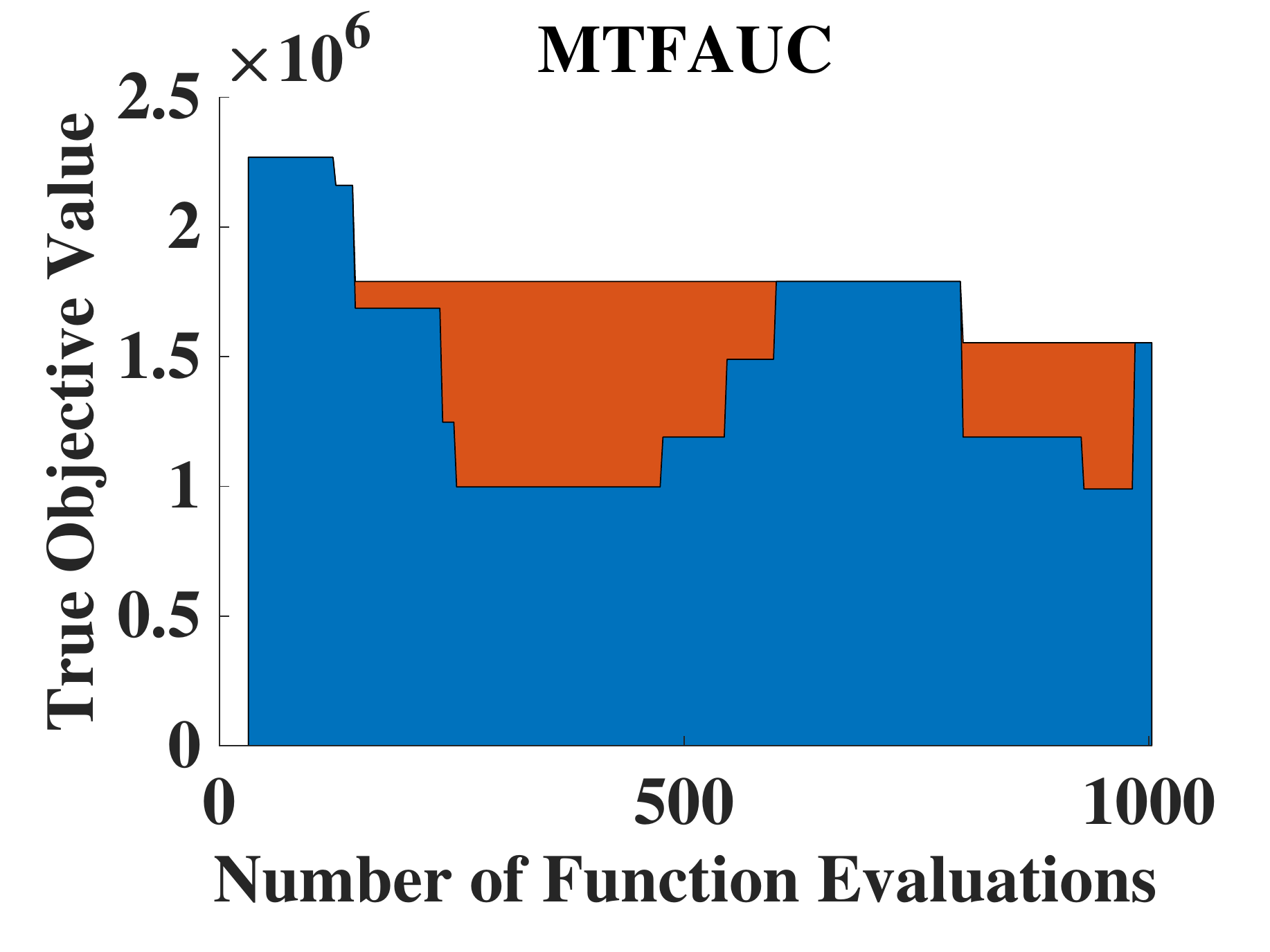}}
\caption{MTFAUC vs. AUC}
\label{fig:mtfauc-eg}
\hspace{\fill}
% \vspace{-5mm}
\end{figure}

\hadis{Figure~\ref{fig:mtfauc-eg} shows the difference between AUC and MTFAUC. \new{Figure~\ref{fig:mtfauc-eg} (a) demonstrates a deterministic test case in which the AUC curve is non-increasing}. As we can observe in Figure~\ref{fig:mtfauc-eg} (b), the AUC curve of the true function is no longer non-increasing since there are uncertainties associated with the sampled function values. It is worth mentioning that the AUC in the stochastic case in (b) is 0.35, which is smaller than that in (a) with 0.37. However, we can observe that the final BSMS reported by (b) has a higher true objective value than solutions found earlier in the surrogate optimization process. The oscillations are caused by solutions with relatively good sampled mean objective value but not a good true objective value. Our proposed MTFAUC
% When there are uncertainties associated with the black-box outcome, MTFAUC 
penalizes the unexpected oscillations caused by the noisy objective values and is a new non-increasing curve. Figure~\ref{fig:mtfauc-eg} (c) demonstrate the same test case in Figure~\ref{fig:mtfauc-eg} (b).
While the AUC true function is calculated with the blue area under the curve, MTFAUC penalizes the area painted in red and blue. 
Note that the BSMS in the final steps is critical as it is the output of the surrogate optimization framework.
%, while the initial oscillations in the objective value are less important due to the few numbers of evaluated data points. 
A robust surrogate would have fewer and shorter oscillations in the final iterations, which results in a lower overall MTFAUC value.
Note that this cannot be captured by mean and variance, as those are affected by oscillations at the beginning and end equally. 
}
%In the presence of noise, the variance is not a proper metric since as it is shown in Figure ~\ref{fig:mtfauc} the plot which is penalized has a smaller variance of the objective values of the BSMS but is not stable in the final iteration. In contrast, the plot with a decreasing pattern on the objective values of BSMS has higher variance but more stable behavior, especially at the final iterations. Similarly, mean does not capture the unexpected oscillations in the final iterations. MTFAUC focuses on stability in the final steps.}

% \begin{figure}[H]
% \centering
%         \includegraphics[width=1\textwidth]{figures/mtfauc.pdf}
%     \caption{MTFAUC vs. AUC}
%     \label{fig:mtfauc-eg}
%     \hspace{\fill}
%     \vspace{-10mm}
% \end{figure}

\new{Indeed the AUC and MTFAUC evaluation metrics are designed to evaluate the performance of different algorithms across many different types of cases. AUC and MTFAUC are designed for system testing purposes, and they are not intended to be used in a real-world engineering problem. The proposed measures are used after a full execution of an algorithm and are not involved as a convergence criterion inside the algorithm. For example, if we execute two different algorithms for the same number of iterations, then we can use MTFAUC to compare their performance. Even when the global minimum is unknown, in the deterministic case, AUC can be applied directly to evaluate the performance of algorithms with the true objective minimum substituted by the objective value of the minimum BKS of all the algorithm executions being evaluated. Under uncertainty, we would collect sampled response values for each data point over all the algorithm executions and apply MTFAUC using the mean sampled response values in place of the true function value.}
\vspace{-15mm}
\section{Experimental Results}\label{sec:results}
%\vspace{-3mm}
We evaluate the performance of our proposed approach on global optimization test functions that have proven to be challenging for BBO~\cite{simulationlib,hansen2009real}. Table~\ref{tab:testfuncchar} describes the characteristics of the selected test functions, \hadis{and $d$ denotes the number of variables.}
In this research, we first assume that the black-box function is deterministic, and later we relax this assumption in \S\ref{subsec:smartresults}. We use MTFAUC throughout the experiments for performance evaluation. \new{Recall MTFAUC is equal to AUC for deterministic cases.}

\hadis{We evaluate the performance of surrogate optimization using TK-MARS, RBF interpolant, non-interpolating RBF (nonRBF), and non-interpolating GP (non-GP) as surrogates. We fixed the budget to $B=1000$ function evaluations. The benchmark test functions presented in Table~\ref{tab:testfuncchar} are used along with the cross-validated hyperparameters presented in Table~\ref{tab:hyp} for different surrogates. \new{For nonGP implementation, we use the GP package in Matlab.\footnote{\url{https://www.mathworks.com/help/stats/gaussian-process-regression-models.html}} 
We use Bayesian 
hyperparameter optimization (which is an option in the Gaussian Processes package in Matlab) to tune the GP hyperparameters using cross-validation. We perform this hyperparameter tuning % fine-tuning of the surrogate model periodically
% The nonGP hyperparameters are tuned using the the Bayesian optimization option of the package.
once on the first iteration of our surrogate optimization algorithm, and once after the algorithm completes 500 function evaluations. 
In other iterations, we fit the nonGP model using the same set of hyperparameters as in the previous iteration. This is a common strategy in Bayesian optimization literature \cite{ju2017designing,song2018optimizing,norouzzadeh2021deep}. Moreover, in other preliminary experiments, we saw no significant difference between tuning the GP hyperparamters on every iteration and tuning them on just these two iterations}. We refer to the surrogate optimization with different surrogates as the name of surrogates in the plots to save space. For example, we mark the surrogate optimization algorithm with RBF as ``RBF''.}
%, where tuning the embedded model's hyperparameters at collecting a few datapoints are not significantly changing the performance of the fitted model ''.

% Table generated by Excel2LaTeX from sheet 'Sheet1'
\begin{table}[H]
  \centering
%   \vspace{-5mm}
  \small
  \caption{Hyperparameter settings of different surrogates}
    \begin{tabular}{rrll}
    \multicolumn{1}{l}{Model} & \multicolumn{3}{l}{Hyperparameters} \\
    \midrule
          &       &       &  \\
    \multicolumn{1}{l}{\textbf{TK-MARS}:} & \multicolumn{1}{l}{CART:} & minsplit=20 & maxdepth=30 \\
          & \multicolumn{1}{l}{MARS:} & $M_{max}=\left[\frac{2n + 3}{2+3}\right]$ & $L_m=1$ \\
    \multicolumn{1}{l}{\textbf{RBF}:} &       & $kernel=MQ$ & $\omega=2$ \\
    \multicolumn{1}{l}{\textbf{nonRB}:} &       & $kernel=MQ$ & $\omega=2$, $\eta=0.0001$ \\
    \multicolumn{1}{l}{\textbf{nonGP}:} & \multicolumn{1}{l}{Bayesopt} & kernel & $\in \{matern32, matern52, squaredexponential, exponential\}$ \\
          &       & $\sigma$ & $\in [[0.01*std(\tilde{f}(x)),10*std(\tilde{f}(x)]] $ \\
          &       & $\sigma_f$ & $\in [1,10]$ \\
    \end{tabular}%
  \label{tab:hyp}%
\end{table}%

\begin{table}[H]
% \vspace{-8mm}
  \centering
  \small
  \caption{Test Functions Definition}
    \begin{tabular}{|l|l|c|c|}
    \hline
    \textbf{Function } & \multicolumn{1}{l|}{\textbf{Formulation}} & \textbf{Range } & \multicolumn{1}{l|}{\textbf{Global Min}} \\
    \hline
    Rosenbrock & $f(x)=\sum_{i=1}^{d-1}[100(x_{i+1}-x_i^2)^2+(x_i-1)^2]$   & [-5,10] & $f(x*)=0$ \\
    &&&$x*=(1,\ldots,1) $\\
    \hline
    Rastrigin  & $f(x)=10d+\sum_{i=1}^{d}[x_i^2-10cos(2 \pi x_i)]$     & [-5.12,5.12] & $f(x*)=0$ \\
      &&&$x*=(0,\ldots,0) $ \\
      \hline
    Levy & $f(x)=\sin^2(\pi w_1)$  & [-10,10] & $f(x*)=0$ \\
    &$\qquad+\sum_{i=1}^{d-1}(w_i-1)[1+10 \sin^2(\pi w_i+1)]$&&$x*=(1,\ldots,1) $\\
    &$\qquad +(w_d-1)^2[1+\sin^2(2\pi w_d)]$&&\\
    & where,  $w_i=1+\frac{x_i-1}{4}$&&\\
    \hline
        Ackley & $f(x_0 \cdots x_n) = -20 exp(-0.2 \sqrt{\frac{1}{d} \sum_{i=1}^d x_i^2}) $  & [-32.768, 32.768] & $f(x*)=0$ \\
   &$\quad - exp(\frac{1}{d} \sum_{i=1}^d cos(2\pi x_i)) + 20 + e$&&$ x*=(1,\ldots,1) $\\
    
    \hline
        Zakharov & $f(x)=\sum_{i=1}^d x_i^2+ (\sum_{i=1}^d 0.5ix_i)^2$  & [-5,10] & $f(x*)=0$ \\
    &$\qquad+(\sum_{i=1}^d 0.5ix_i)^4$&&$x*=(0,\ldots,0) $\\
    \hline
    \end{tabular}%
  \label{tab:testfuncchar}
%   \vspace{-15mm}
\end{table}%

% \vspace{-10mm}

\subsection{Comparison of MARS and TK-MARS}

\begin{figure}[!ht]
% \vspace{-3mm}
\centering
    \begin{minipage}{\linewidth}
        \subfloat[]{\includegraphics[width=0.35\textwidth]{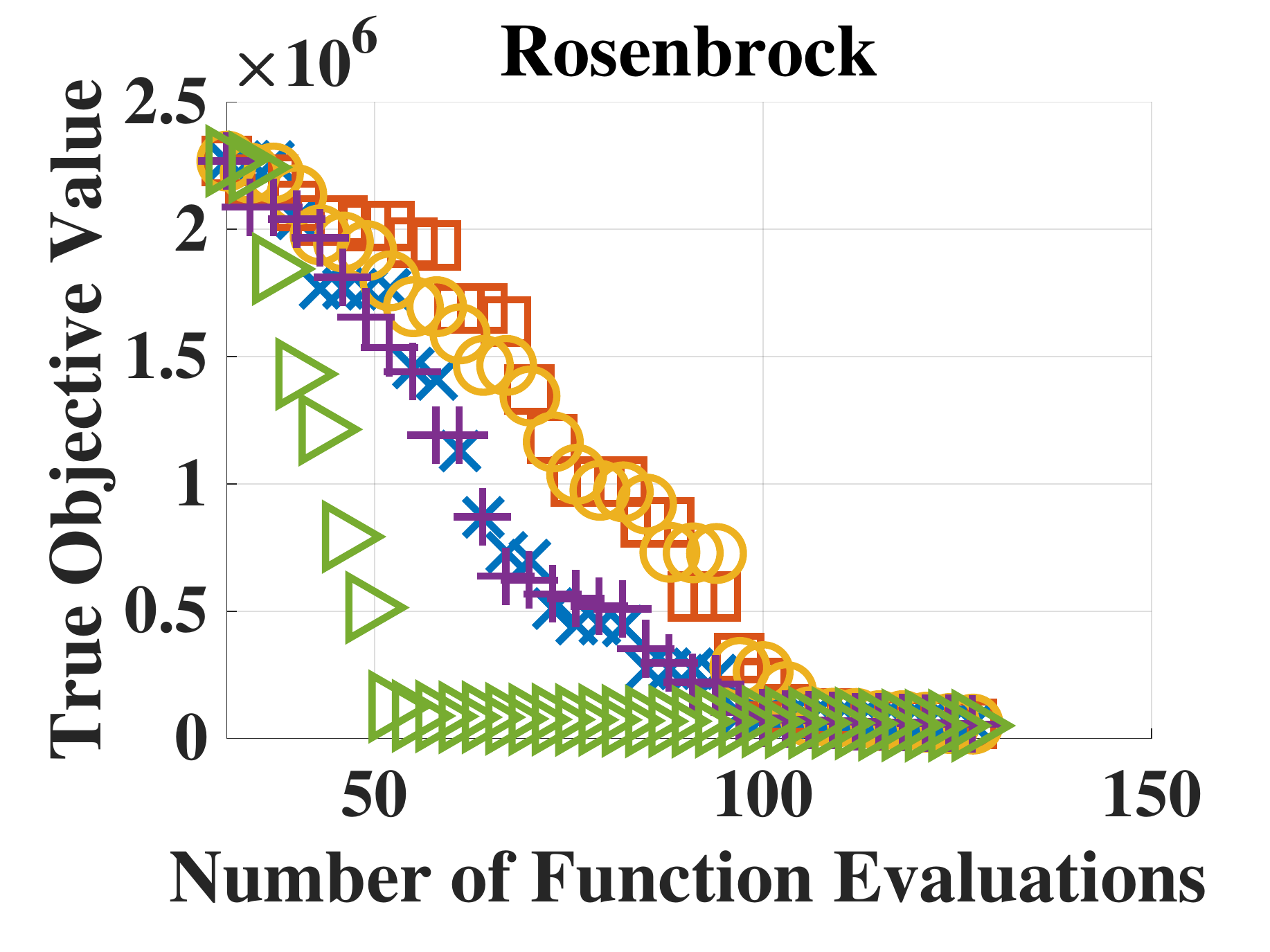}}
        \subfloat[]{\includegraphics[width=0.35\textwidth]{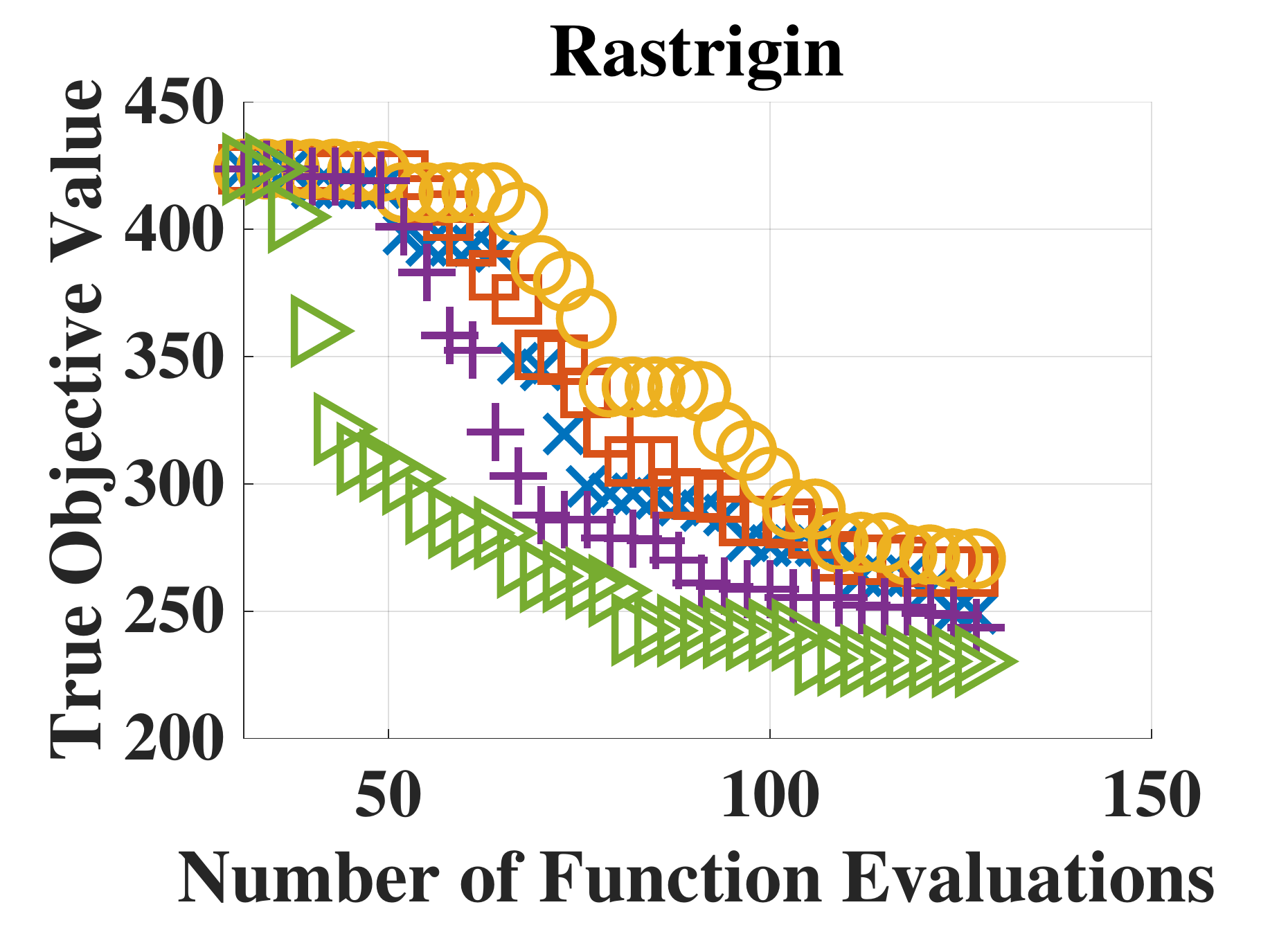}}
        \subfloat[]{\includegraphics[width=0.35\textwidth]{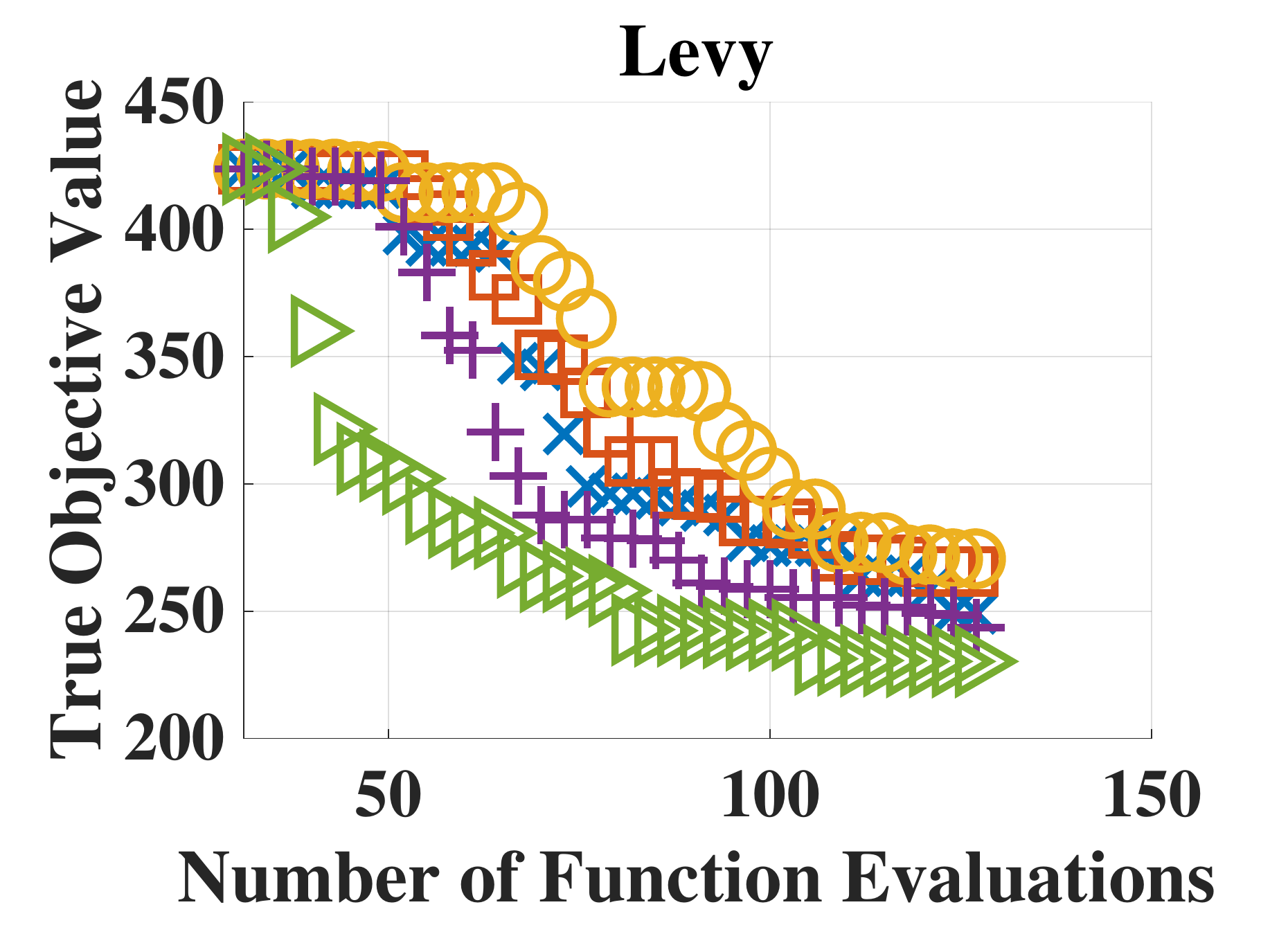}}
    \end{minipage}
    \begin{minipage}{\linewidth}
     \subfloat[]{\includegraphics[width=0.35\textwidth]{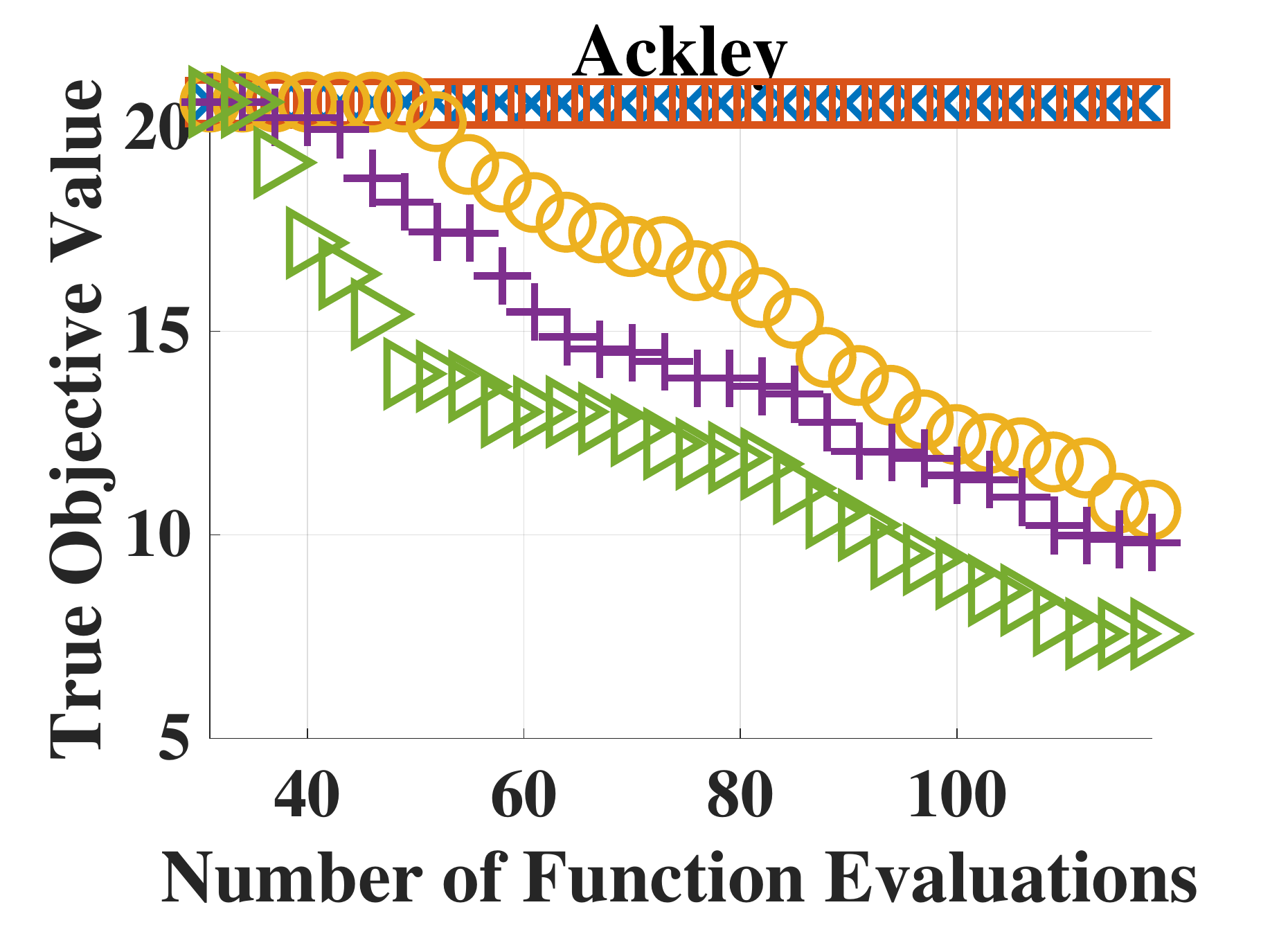}}
    \subfloat[]{\includegraphics[width=0.35\textwidth]{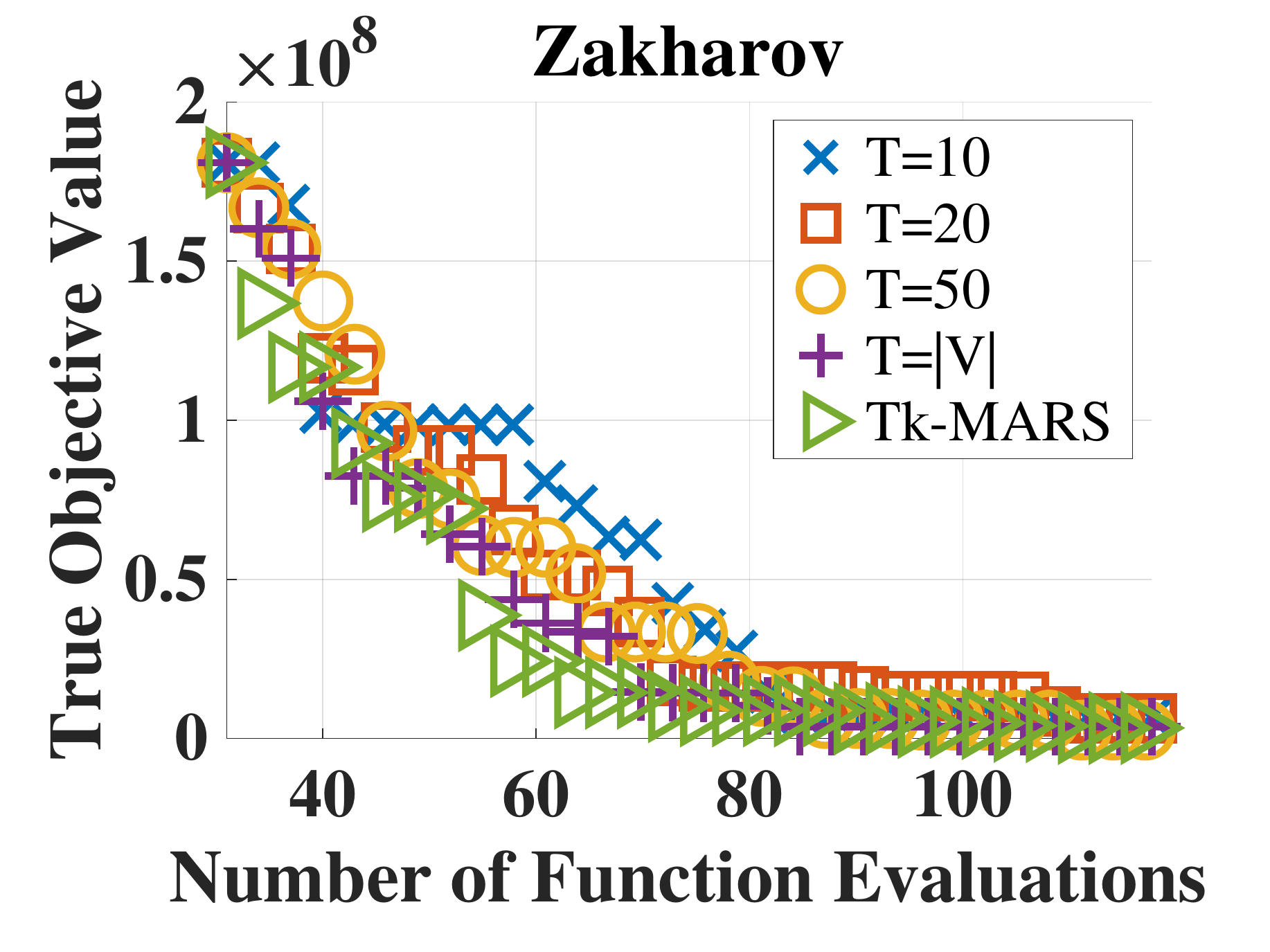}}
    \end{minipage}
    \caption{TK-MARS vs. MARS with different number
of eligible knot locations}
    \label{fig:treeknots}
    \vspace{-3mm}
\end{figure}

First, we evaluate the proposed TK-MARS versus original MARS with evenly spaced knot locations in the context of surrogate optimization. The number of eligible knot locations, $T$, has to be preset for original MARS. We examine a surrogate optimization algorithm that uses original MARS with $T=10$, $T=20$, and $T=50$. Since TK-MARS sets $T$ based on the number of terminal nodes, $|V|$, to ensure a fair comparison between the two approaches, we consider $T=|V|$ for original MARS, as well. \hadis{It should be noted that as the number of knots increases, the computation time of MARS increases, even though the accuracy may increase.}
%It is more appropriate to have MARS and TK-MARS dynamically set $T$ based on the number of terminal nodes of CART,$V$, to ensure a fair comparison between the two approaches. Yet the eligible knot location is distinct in each approach. As a consequence, we demonstrate that TK-MARS eligible knot locations promise near optimal sites.
We initialize $I$ with $N=31$ points selected with a LHD for $d=30$ independent input variables.
Figures~\ref{fig:treeknots}(a)-(e) show the average true objective value of the BSMS versus the number of black-box function evaluations for each of the test functions listed in Table~\ref{tab:testfuncchar}. The results are presented for each surrogate optimization algorithm using the aforementioned versions of MARS as well as our proposed TK-MARS over 30 different executions. %In addition, Table~\ref{tab:tkmarsauc} shows MTFAUC for them as well. 
\hadis{Note that using TK-MARS as a surrogate within a surrogate optimization algorithm outperforms that using the original MARS with different number of knots for the various test functions. The green triangles corresponds to TK-MARS, and the faster convergence leads to smaller AUC.} %except for Zakharov. The plate-shaped valley of Zakharova creates a complex functional behavior that is in general, difficult to model by surrogates. Having more Knots would help MARS to better fit the Zakharov function.} 
%The initial set of 35 points with $d=27$ independent variables designed with LHD for this experiment. As one can see, the tree-based knots approach improves the optimization process, significantly.

% \begin{table}[!ht]
% \vspace{-3mm}
% \centering
% \caption{MTFAUC values for TK-MARS vs. MARS \hadis{revise}}
% \begin{tabular}{|l|c|c|c|}
% \hline
%      & Rosenbrock & Rastrigin & Levy \\ \hline
% T=10 & 0.05      & 0.09      & 0.06 \\ \hline
% T=20 & 0.08      & 0.10      & 0.12 \\ \hline
% T=50 & 0.05      & 0.05      & 0.09 \\ \hline
% T=V  & 0.05      & 0.04      & 0.09 \\ \hline
% TK   & 0.02      & 0.04      & 0.03 \\ \hline
% \end{tabular}
% \label{tab:tkmarsauc}
% \end{table}

\hadis{Next, we compare the performance of our surrogate optimization using TK-MARS versus deterministic RBF (``RBF''), non-interpolating RBF (``nonRBF''), and non-interpolating GP for noisy observations (``nonGP'') on a new class of test functions that include unimportant input variables.}
%We believe that the existing test functions differ from many real-world applications in which certain input variables are unimportant. 
%Hence, we design a new class of test functions that mimic the real-world problems in the context of surrogate optimization~\cite{addis2007new,balasundaram2005constructing,jamil2013literature}.
%To design a new class of test functions, we consider two major factors of real-world black-box functions, (1) fraction of important variables and (2) uncertainties associated with the black-box systems.
Specifically, the new test functions only consider a fraction of the input variables, $\mathit{fiv}$ \hadis{for the function evaluation. In particular, we select the first $\mathit{fiv}$\% of the variables for the evaluation. For example, if $\mathit{fiv}=0.5$ and $d=30$ only
$\{x_1,\dots,x_{15}\}$ are considered for the function evaluation.} %Then, we test the performance of the proposed method in identifying the significant variables.

\begin{figure}[H]
\vspace{-8mm}
\centering
    \begin{minipage}{\linewidth}
        \subfloat[]{\includegraphics[width=0.35\textwidth]{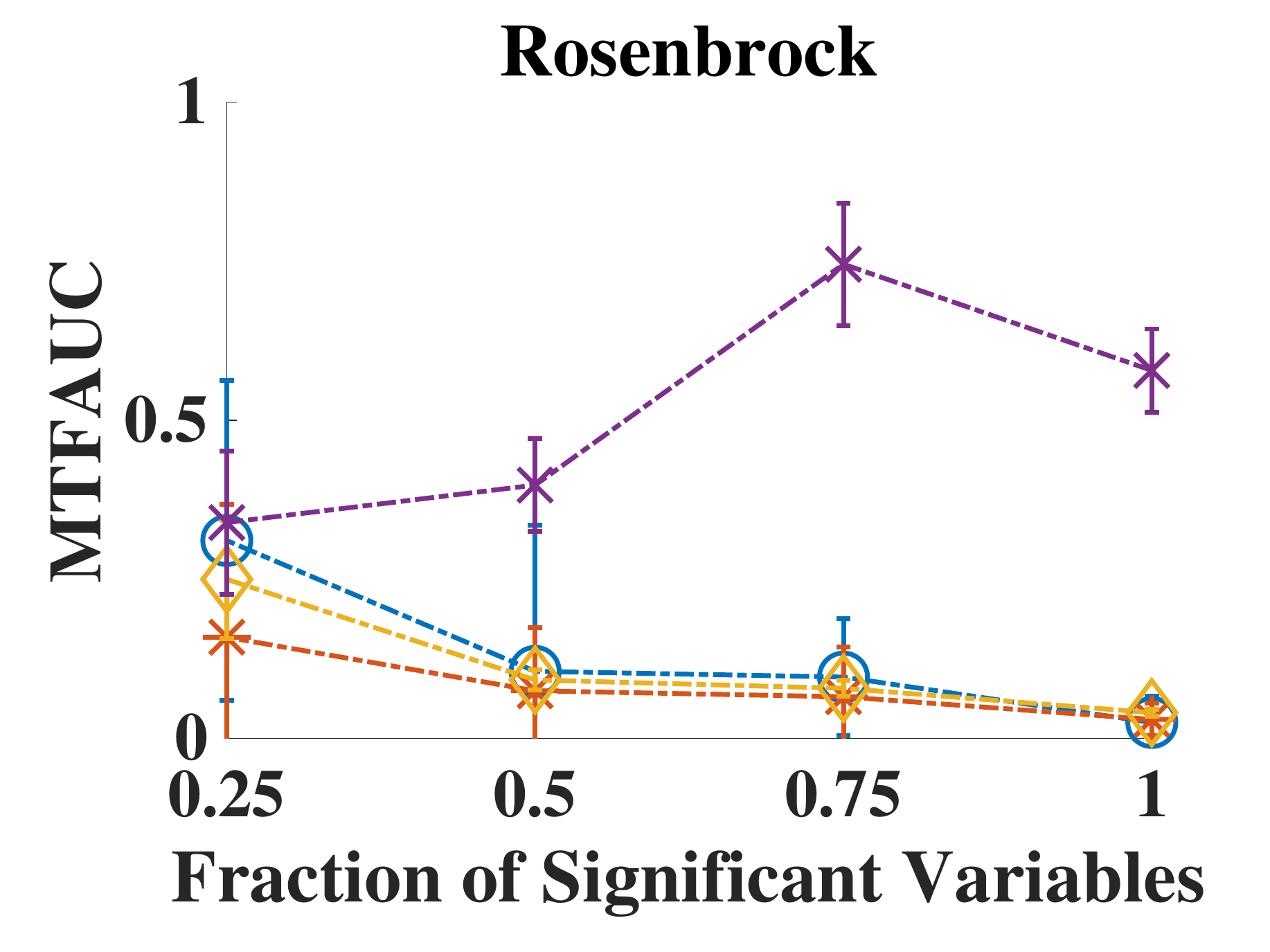}}
        \subfloat[]{\includegraphics[width=0.35\textwidth]{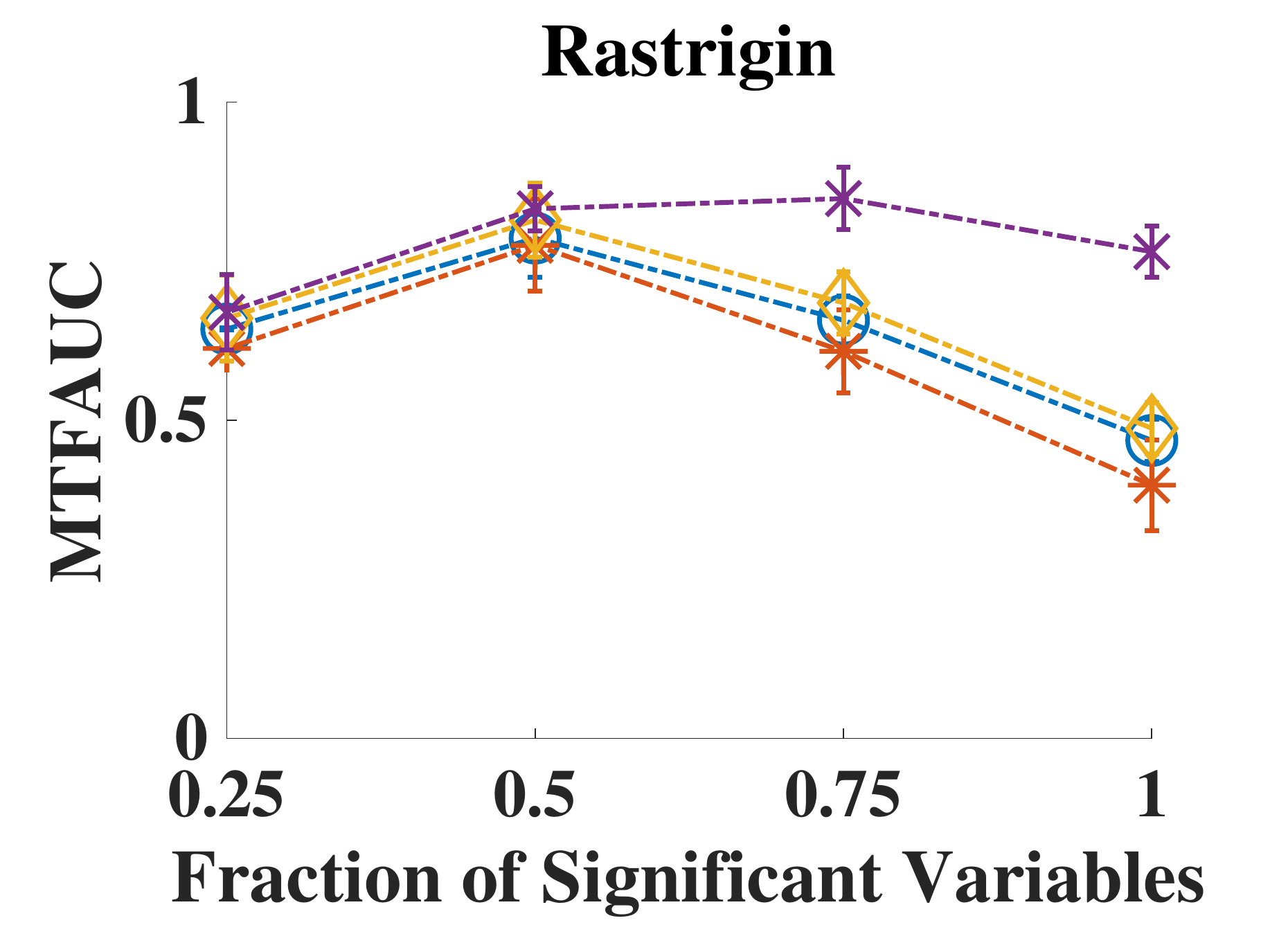}}
        \subfloat[]{\includegraphics[width=0.35\textwidth]{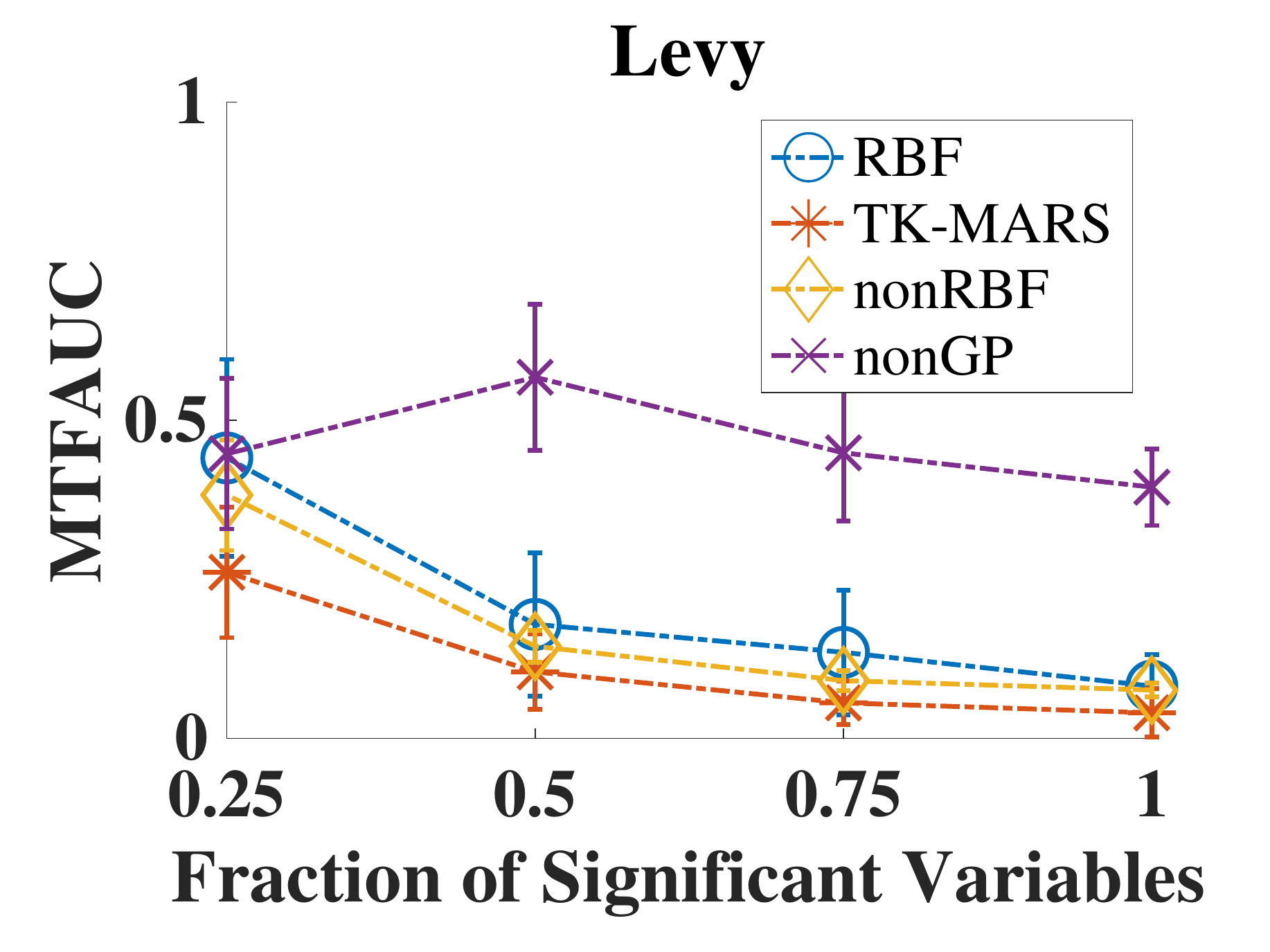}}
    \end{minipage}
    \begin{minipage}{\linewidth}
        \subfloat[]{\includegraphics[width=0.35\textwidth]{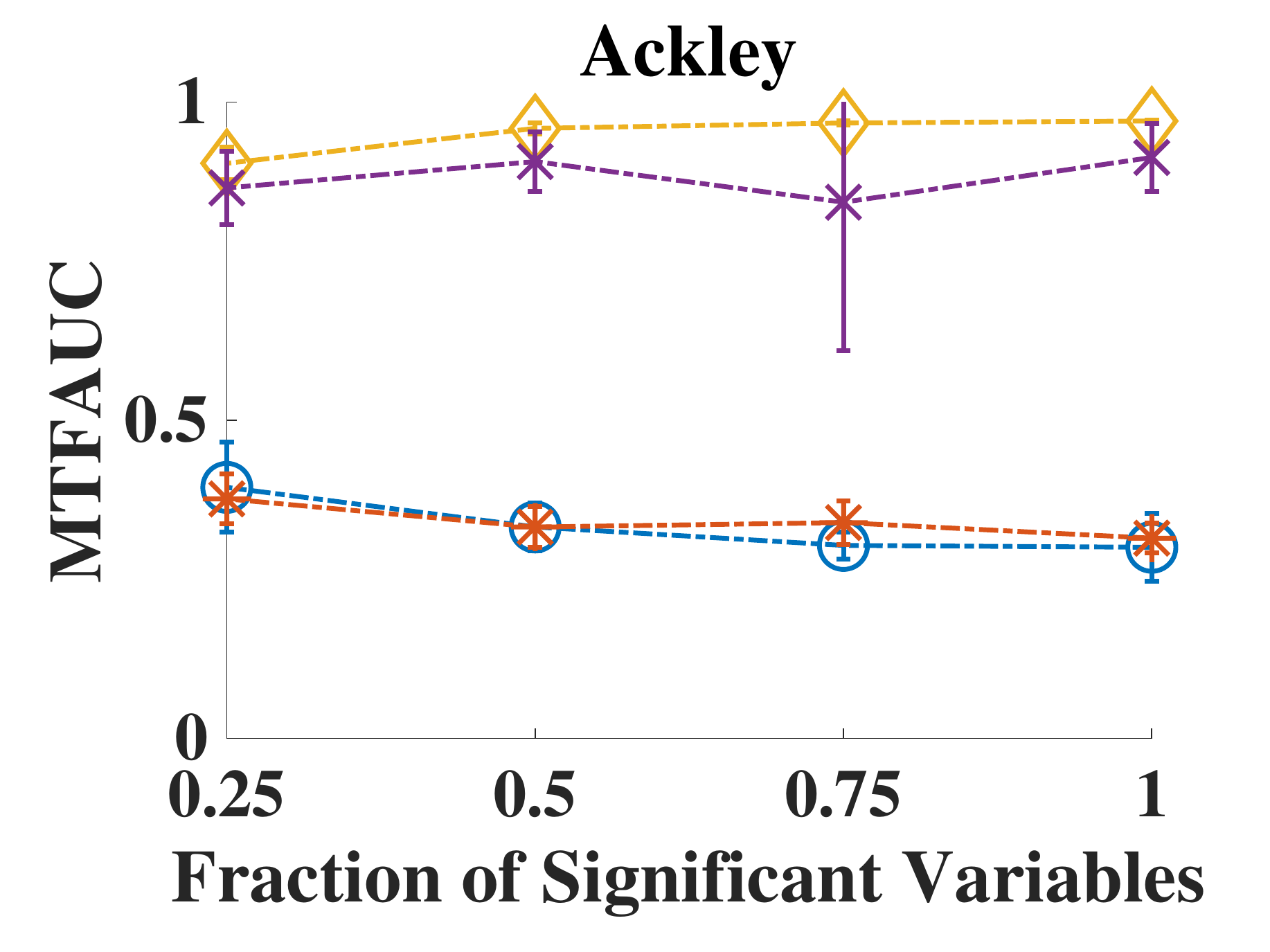}}
        \subfloat[]{\includegraphics[width=0.35\textwidth]{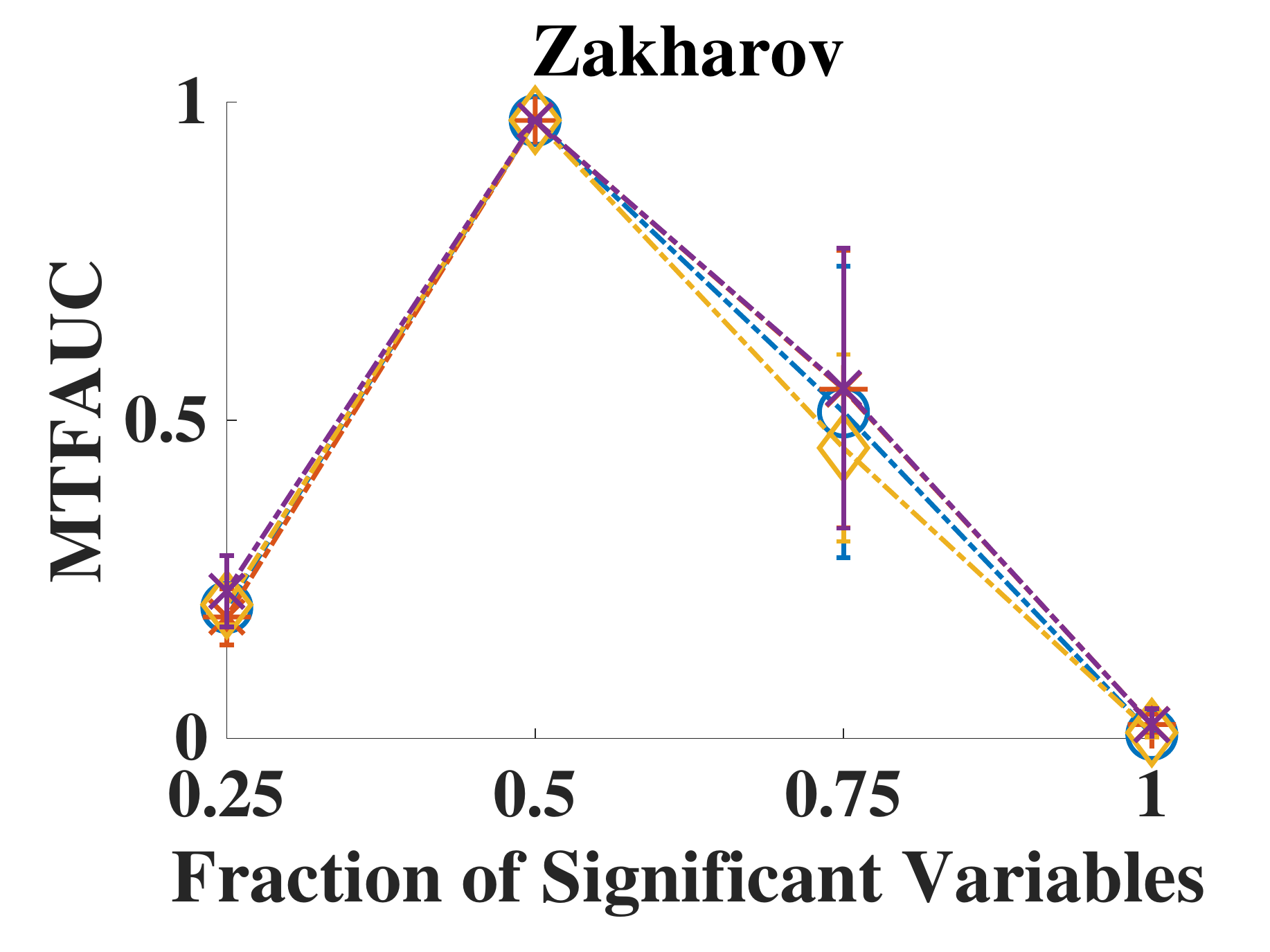}}
    \end{minipage}
    % \begin{minipage}{\linewidth}
    %     \subfloat[]{\includegraphics[width=0.35\textwidth]{figures/fracmodel_rosen_25.pdf}}
    %     \subfloat[]{\includegraphics[width=0.35\textwidth]{figures/fracmodel_rast_25.pdf}}
    %     \subfloat[]{\includegraphics[width=0.35\textwidth]{figures/fracmodel_levy_25.pdf}}
    % \end{minipage}
    \caption{TK-MARS vs. other surrogate models across different $\mathit{fiv}$}
    \label{fig:tkmars vs rbf}
    \vspace{-3mm}
\end{figure}

Figure~\ref{fig:tkmars vs rbf} shows the average MTFAUC of 30 different executions for each test function with four different levels of $\mathit{fiv}$, 0.25, 0.5, 0.75, and 1.0. \hadis{The error bar shows the standard deviation of the average MTFAUC of 30 different executions.}

\hadis{Note that MTFAUC is higher when fewer variables are important}, suggesting that the surrogates struggle to determine which input variables are important when few of them are. \hadis{
%Observe that surrogate optimization using TK-MARS performs at least as well as that using the other considered surrogate models at each level of $\mathit{fiv}$. 
The benefit of using TK-MARS is notable for the Rosenbrock, Rastrigin, and Levy functions.
%, and in general it is largest when there are the fewest number of important input variables. 
This difference is likely because MARS is parsimonious, while others use all the input variables in their functional form. We do not observe the same pattern for Zakharov function since the plate-shaped valley of Zakharov function makes the surrogate fitting extremely difficult for metamodels \cite{tsoukalas2016surrogate}. For the Ackley function, the difference is not visible in the plot because of the higher range of MTFAUC for nonRBF and nonGP, however, there is 2\% difference between the performance of \new{TK-MARS} compared with RBF when $\mathit{fiv}=0.25$.}

\begin{figure}[H]
\vspace{-5mm}
\centering
    \begin{minipage}{\linewidth}
        \subfloat[$\mathit{fiv}=1$]{\includegraphics[width=0.45\textwidth]{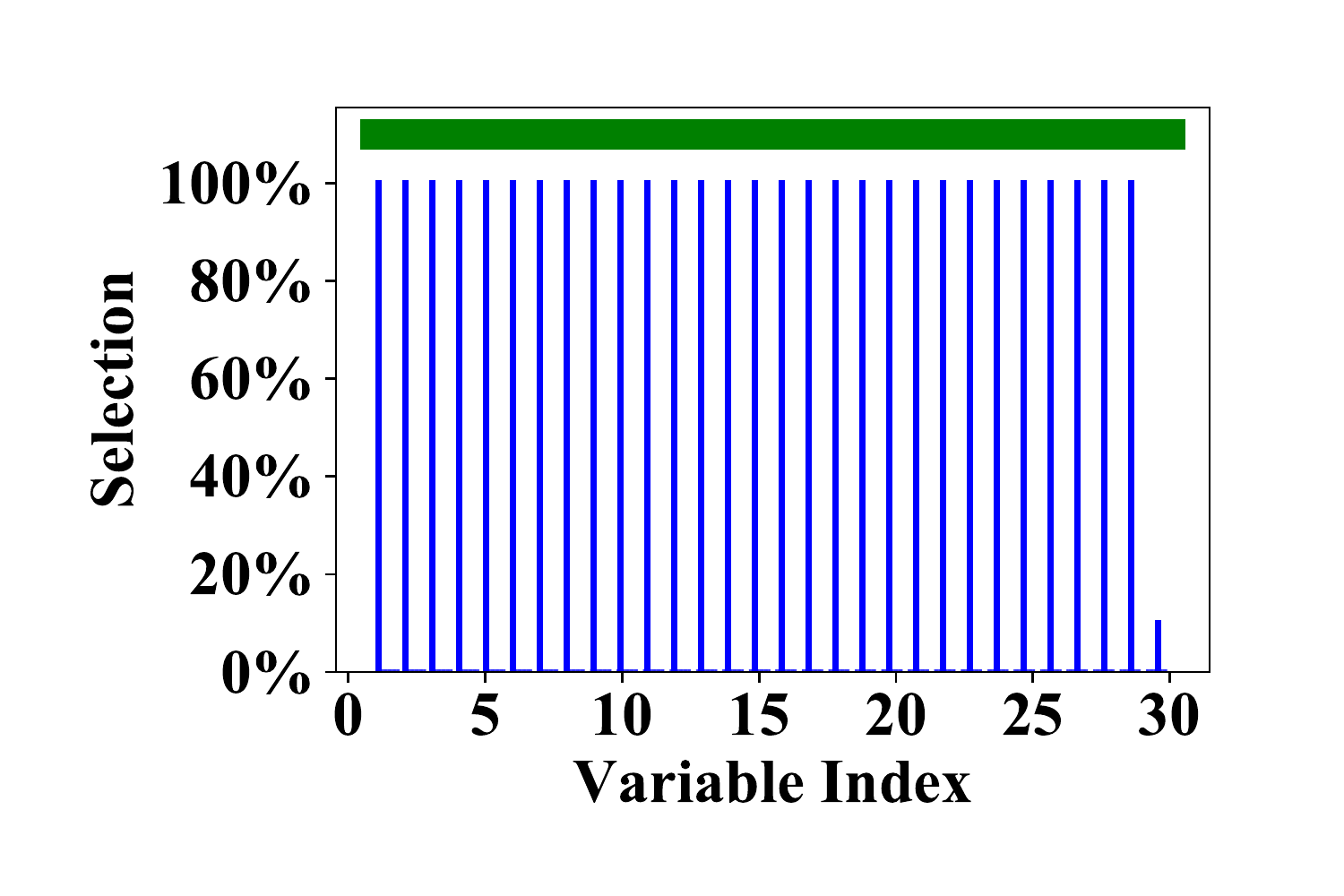}}
        \subfloat[$\mathit{fiv}=0.75$]{\includegraphics[width=0.45\textwidth]{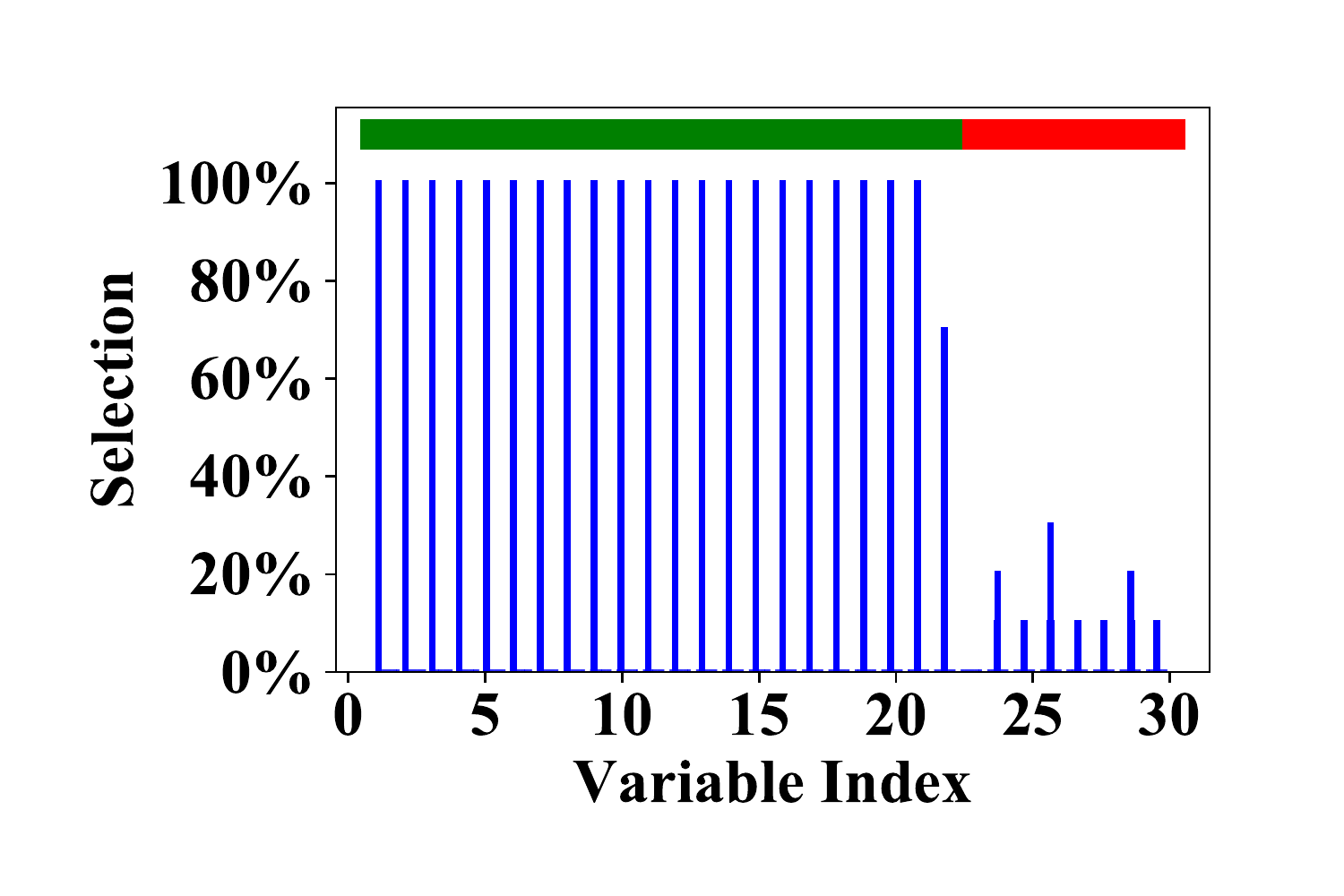}}
            \vspace{-3mm}
    \end{minipage}
    \begin{minipage}{\linewidth}
     \subfloat[$\mathit{fiv}=0.5$]{\includegraphics[width=0.45\textwidth]{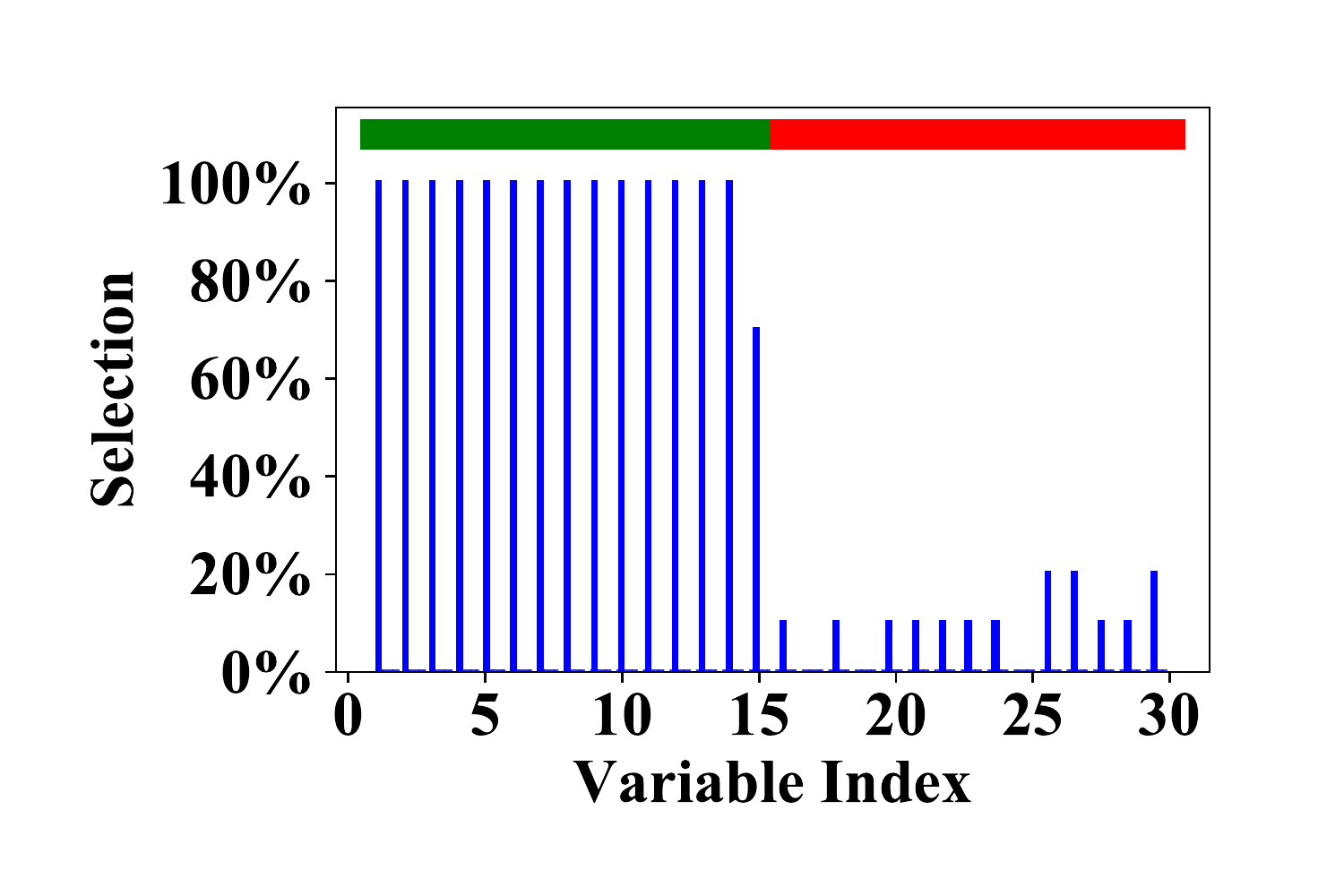}}
    \subfloat[$\mathit{fiv}=0.25$]{\includegraphics[width=0.45\textwidth]{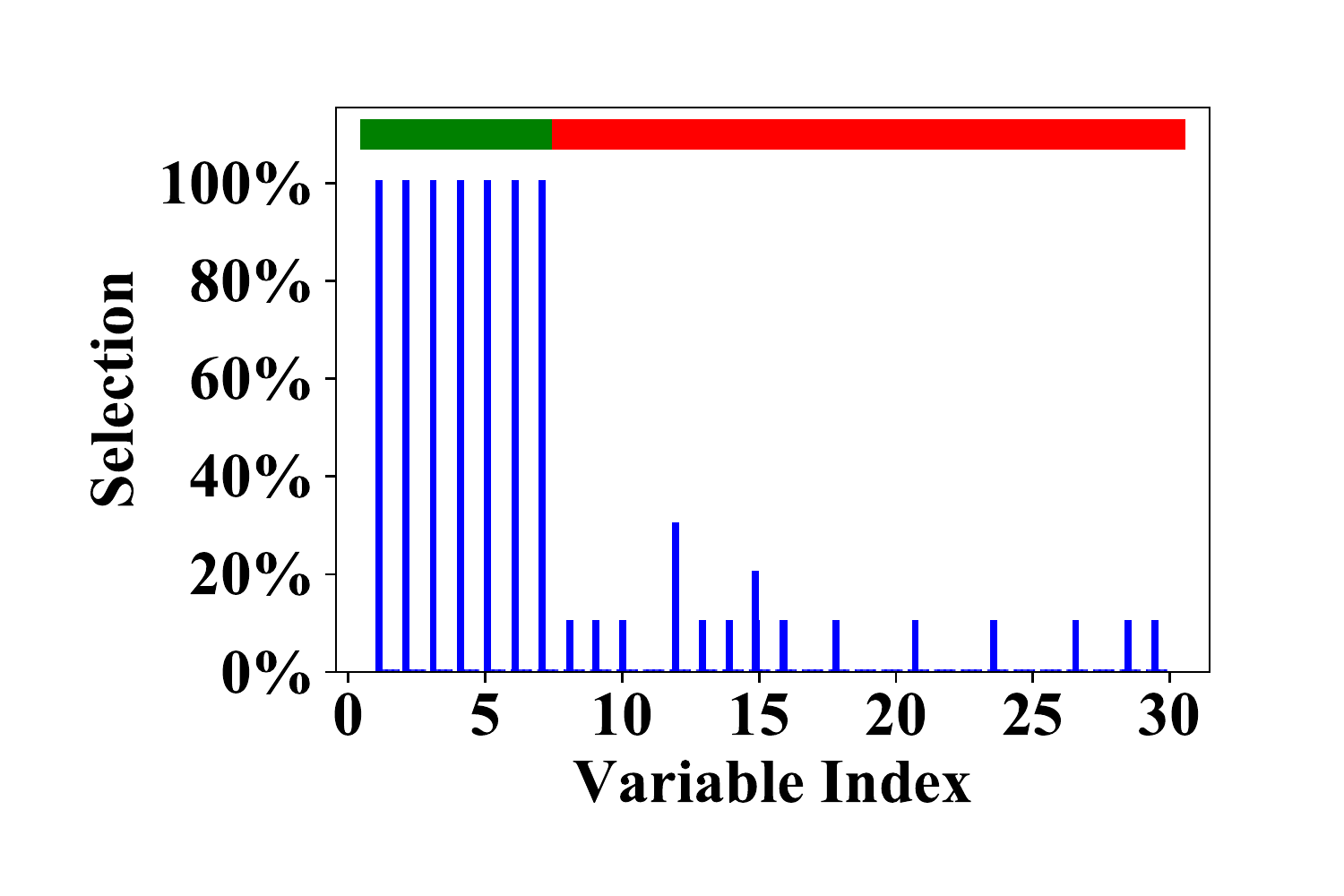}}
    \end{minipage}
    \caption{Percentage of variables selected by TK-MARS across 10 different executions}
    \label{fig:vars}
    \vspace{-5mm}
\end{figure}

\hadis{To demonstrate the performance of TK-MARS in correctly detecting the important variables, Figure~\ref{fig:vars} presents the variables selected by TK-MARS for Rosenbrock function, at the final iteration of surrogate optimization. The bars show the percentage of a variable being selected in 10 different executions (different pools). The green horizontal bar indicates the range of the important variables considered for the function evaluation, and the red horizontal bar shows the unimportant range of the variables. As we can observe, in Figure~\ref{fig:vars}(a), TK-MARS correctly selects all the important variables (except 30) in all 10 executions. In Figure~\ref{fig:vars}(b), TK-MARS in the majority of cases selects all the important variables, which are $x_1,\dots, x_{22}$, and do not prefer the unimportant variables \new{$x_{23},\dots, x_{30}$}. One can confirm the same pattern in Figure~\ref{fig:vars}(c) and Figure~\ref{fig:vars}(d). Although TK-MARS selects few unimportant variables, because the variable selection of TK-MARS also depends on the input-output of the training dataset, TK-MARS involves the significant variables for its basis function construction.}

\hadis{While the performance of surrogate optimization is usually based upon the number of expensive function evaluations, it is worth comparing the non-evaluation time of the algorithm considering different
surrogate models. In our experiments, for Rosenbrock function with $d=30$,
RBF requires an average of 0.06 seconds, original MARS and TK-MARS use 0.04 seconds, non-interpolating RBF needs 0.07 seconds, and non-interpolating GP uses 0.22 seconds. Although TK-MARS requires negligible time more than original MARS for CART to select the candidate knot locations, this step is completed in 0.006 seconds. 
%
%TK-MARS are compared to other surrogates for %Rosenbrock function with $d=30$ on the fixed dataset; RBF: 1.7163 secs,
%MARS: 1.4518 secs, and fitting a CART for TK-MARS adds 0.012 secs,
%nonRBF:     2.4117 secs, and
%nonGP: 	     7.2904 secs. 
Moreover, the total computation time of the No-Replication algorithm with using TK-MARS, MARS with $k=20$, and $k=50$ on Rosenbrock function is 1148 secs, 927 secs, and 1346 secs, respectively. Consequently, the difference in the computation times of TK-MARS and original MARS is unlikely to be a concern to a user given how well TK-MARS outperforms MARS in terms of AUC as shown in Figure~\ref{fig:treeknots}.} \new{The experiments were performed on a 2.3 GHz Core i9 machine with 16 GB DDR4 memory.}
%  The approximation quality decreases with only a fraction of variables available, and the convergence becomes slower.
%It is more difficult for the algorithm to recognize significant variables when the output is evaluated using fewer variables.
%From Table~\ref{tab:norep-anova-ff}, TK-MARS statistically outperforms RBF at the significance level of $\alpha=0.1$.
%We see that the noise effect has not been mitigated in this approach. So, next we will discuss the replication strategy effectiveness in handling noise effect.
\vspace{-5mm}
\subsection{ Optimization Under Uncertainty}\label{subsec:smartresults}
\vspace{-3mm}
% Besides, a component of uncertainty is added to the real objective value for each function evaluation to represent stochastic black-box systems.
To evaluate the performance of the proposed approach under uncertainty, \hadis{we assume that the unknown noise associated with the black-box outcome follows a Gaussian distribution and simulate the noise accordingly.} As a result, we add a Gaussian noise to the function values, i.e.,
%In this research, the experiments are provided with a Gaussian noise,
$\tilde{f}(x)=f(x)+\varepsilon$, where ~$\varepsilon~\sim~N(0, \mathit{np}*\sigma_0,)$, where $\mathit{np}\in(0, 1)$ is the noise percentage level, and $\sigma_0$ is given by
$$\sigma_0 = \max\{f(x^1),\ldots, f(x^N)\} - \min\{f(x^1),\ldots, f(x^N)\}.
$$
%range(f(x_{x\in I}))
\hadis{In a preliminary analysis, we designed a set of experiments using an Orthogonal Array,~\cite{bose1952orthogonal,plackett1946design}, on the full list of algorithm and problem parameters (Table~\ref{tab:params-OA} in the Appendix) and performed an ANOVA (Table~\ref{tab:wrep-anova-oa} in Appendix B) to focus on the significant parameters, which are given in Table \ref{tab:wrep-OA}.}
%we observed that the most appropriate settings to conduct experiments are those given in Table \ref{tab:wrep-OA}.
% We designed a set of experiments using an Orthogonal Array,~\cite{bose1952orthogonal,plackett1946design}, on the parameters listed in Table~\ref{tab:wrep-OA}. The notation of each parameter is given the parentheses. 
The Replication number is $r$ when using Fixed-Replication and $r_\mathit{max}$ for Smart-Replication. %There are several techniques to construct orthogonal arrays with multiple factor levels~\cite{bose1952orthogonal,plackett1946design}. OA designs are common for design of experiments in industries where there is a large number of factors to be studied but only a few of them affect the output. 
%OA designs save runs and are flexible for combinations of factor levels. SAS offers a complete library of Orthogonal Arrays with different numbers of runs, factors, and levels~\cite{sasOA}.
%Table~\ref{tab:norep-OA} and 
%Performing an ANOVA, Table~\ref{tab:norep-anova-oa}, on the experiments, we observe that the No-Replication approach can not handle uncertainty in the system and has an important impact on the algorithm efficiency. On the other hand, the fraction of important variables is statistically significant in almost all the analysis. 
%For simplicity and focusing on noise, we continue to study the replication strategies with fixing some of the parameters to their most significant level based upon this ANOVA; Specifically, $d=30$, $\mathit{fiv}=0.5$, $N=d+1$, DOE method to LHD, EEPA distance metric set to Euclidean, and $K'=3$. 
%We refer to Appendix?? for more details on the ANOVA results.
%In the Figures~\ref{fig:tkmars_mtfauc}-\ref{fig:rbf_bks_box_rosen},
%, \ref{fig:auc_box_rosen}, \ref{fig:bks_box_rosen}, \ref{fig:rbf_mtfauc},  \ref{fig:rbf_auc_box_rosen}, and 
\hadis{We refer to No-Replication as $\mathit{Norep}$, Fixed-Replication with $r$ replications as $\mathit{Fixedrep},r$, and Smart-Replication with $r_{\mathit{max}}$ as $\mathit{Smartrep},r_{\mathit{max}}$, throughout the experiments.}
%Table~\ref{tab:wrep-OA} represents the problem setting and algorithm parameters for an extensive set of experiments.
\begin{table*}[!tb]
% \begin{minipage}[t]{0.49\textwidth}
%     \centering
%   \caption{No-Replication parameters and levels for OA design}
%   \begin{scriptsize}
%     \begin{tabular}{|@{}l@{}@{}l@{}|}
%     \textbf{Problem Parameters} & \textbf{levels} \\ 
%     Test function & Rosenbrock, \\
%          & Rastrigin, Levy \\
%           Dimension & 10, 20, 30 \\
%           Fraction of import. vars. & 0.25, 0.50, 0.75, 1 \\
%           Noise level (\%) & 0, 10, 25 \\
%     \textbf{Algorithm Parameters} & \textbf{levels} \\ 
%      Initial pool size & $d+1$, $2(d+1)$ \\
%         DOE method & LHD, Sobol \\
%         EEPA distance & Euclidean, Cosine \\
%         EEPA \# candidates & 3, 6 \\
%         Model & RBF, TK-MARS \\
%     \end{tabular}%
%     \end{scriptsize}
%   \label{tab:norep-OA}%
% \end{minipage}
% \hfill
%\begin{minipage}[t]{0.49\textwidth}
    \centering
  \caption{Replication parameters and levels}
  \begin{scriptsize}
   \begin{tabular}{|l|l|}
   \hline
	\textbf{Problem Parameters} & \textbf{Levels} \\ 
    Test function & Rosenbrock, %\\
         %&
         Rastrigin, Levy, Ackley, Zakharov  \\
%          Dimension ($d$) & 10, 20, 30 \\
          Dimension ($d$) & 30 \\
%Fraction of important variables ($\mathit{fiv}$)& 0.25, 0.50, 0.75, 1.0 \\
Fraction of important variables ($\mathit{fiv}$)& 0.50 \\
          Noise level  ($\mathit{np}$) & 0, 0.05, 0.10, 0.25 \\
          \hline
    \textbf{Algorithm Parameters} & \textbf{Levels} \\ 
%	Initial pool size ($N$) & $d+1$, $2(d+1)$ \\
	Initial pool size ($N$) & $d+1$ \\
          DOE method in Step~\ref{step1} of Algoritm~\ref{alg:SurOpt} & LHD \\
%          EEPA distance metric ($\delta$ in Appendix~\ref{app:EEPA})& Euclidean, Cosine  \\
%          EEPA number of candidates ($K'$  in Appendix~\ref{app:EEPA}) & 3, 6 \\
          EEPA number of candidate points ($K'$  in Appendix~\ref{app:EEPA}) & 3 \\
Replication type & No-Replication, Fixed-Replication,\\ & Smart-Replication \\
          Replication number ($r$ or $r_\mathit{max}$)& 5, 10 \\
          Surrogate Model  & RBF, TK-MARS \\
          \hline
    \end{tabular}%
    \end{scriptsize}
  \label{tab:wrep-OA}%
%\end{minipage}
\end{table*}%

\begin{figure}[H]
\vspace{-7mm}
\centering
\subfloat[Linear function]{\includegraphics[width=0.35\textwidth]{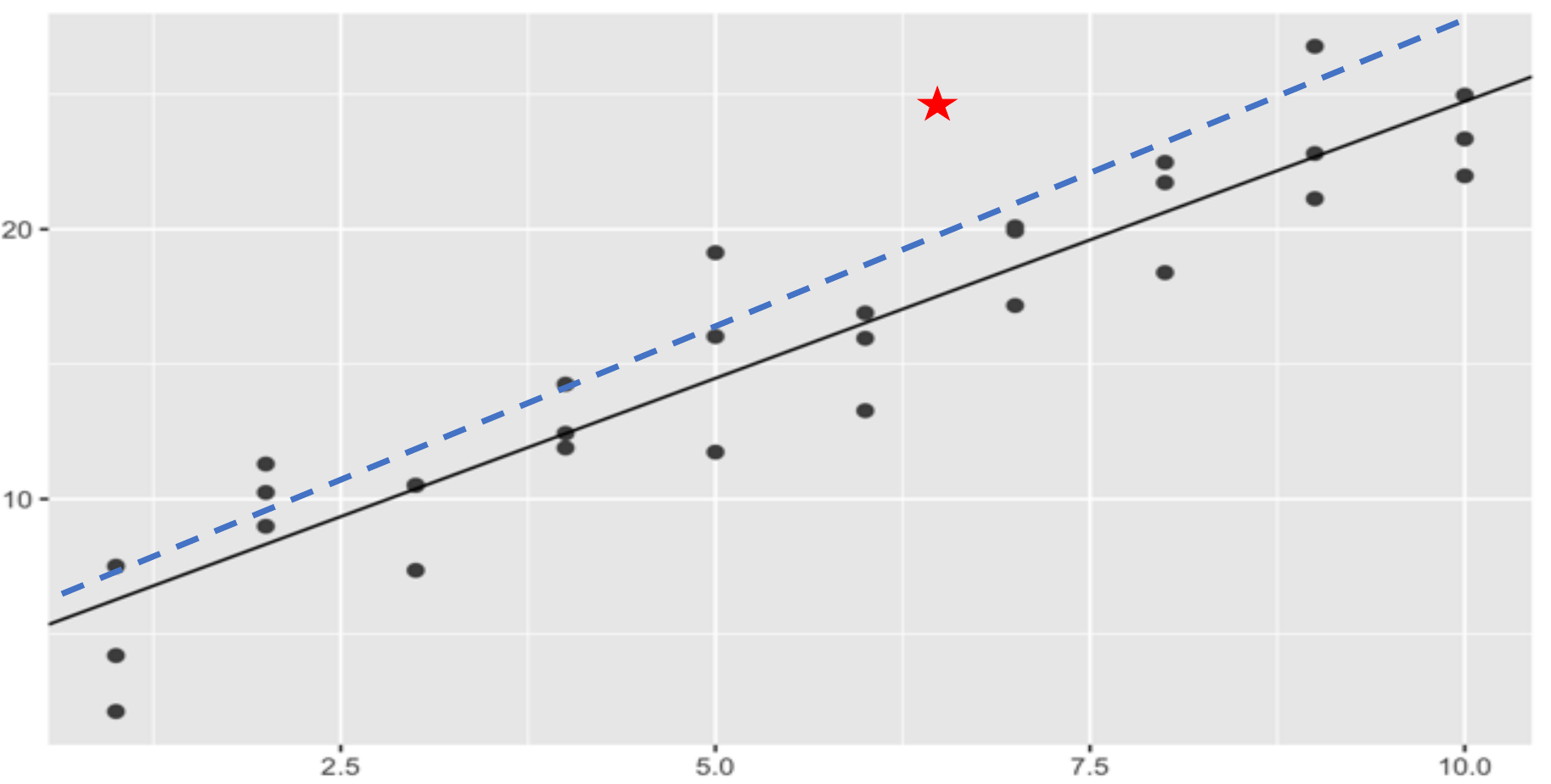}}
\subfloat[Nonlinear ]{\includegraphics[width=0.35\textwidth]{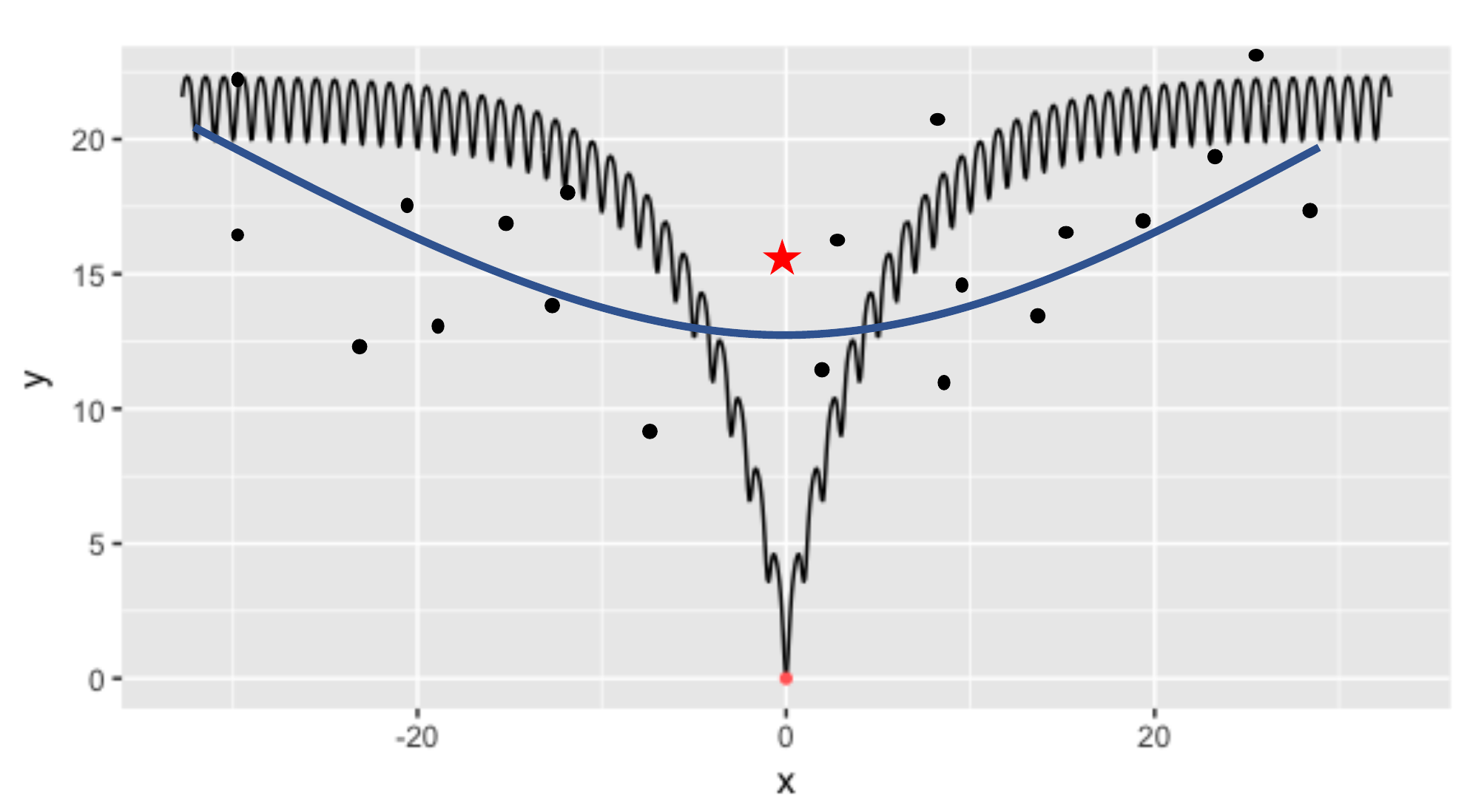}}
\subfloat[Nonlinear + replication]{\includegraphics[width=0.35\textwidth]{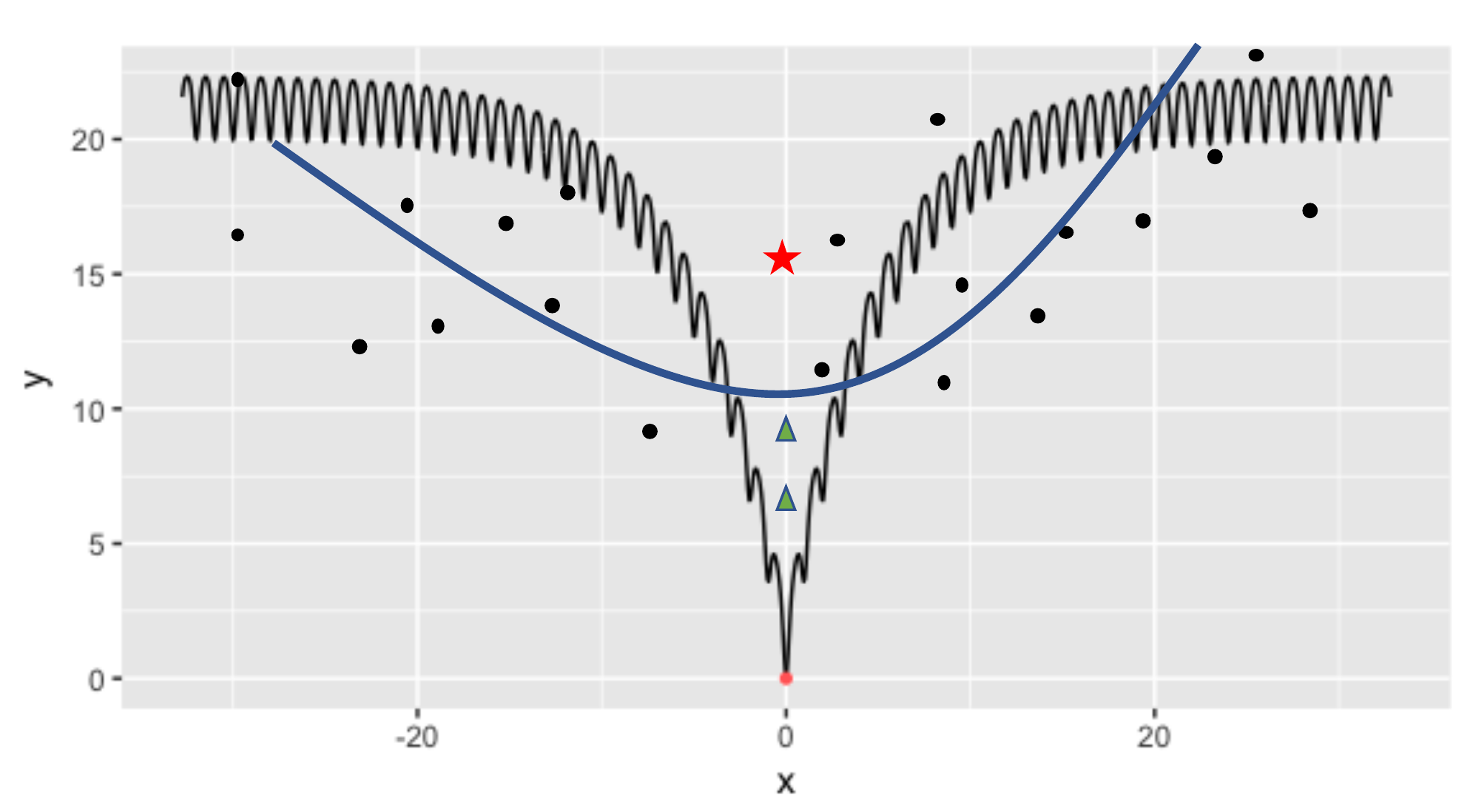}}
\caption{Toy example for non-interpolating wiht/out replication under uncertainty}
\label{fig:toy}
\hspace{\fill}
\vspace{-7mm}
\end{figure}

\begin{figure}
\begin{tabular}{@{\hspace{-8mm}}l@{\hspace{-3mm}}l@{\hspace{-3mm}}l@{\hspace{-3mm}}l@{\hspace{-3mm}}}
\subfloat[RBF]{\includegraphics[width = 1.5in]{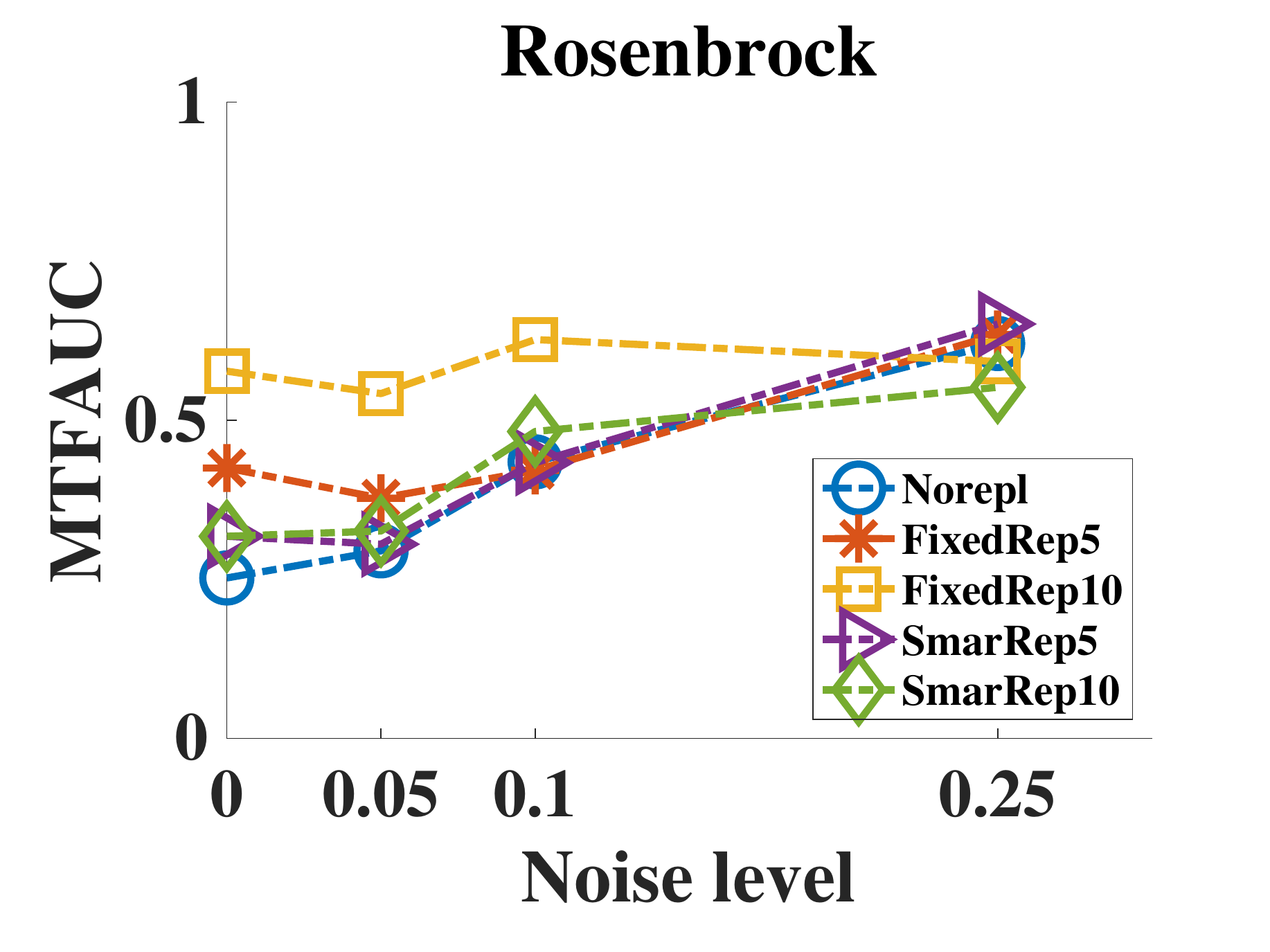}} &
\subfloat[TK-MARS]{\includegraphics[width = 1.5in]{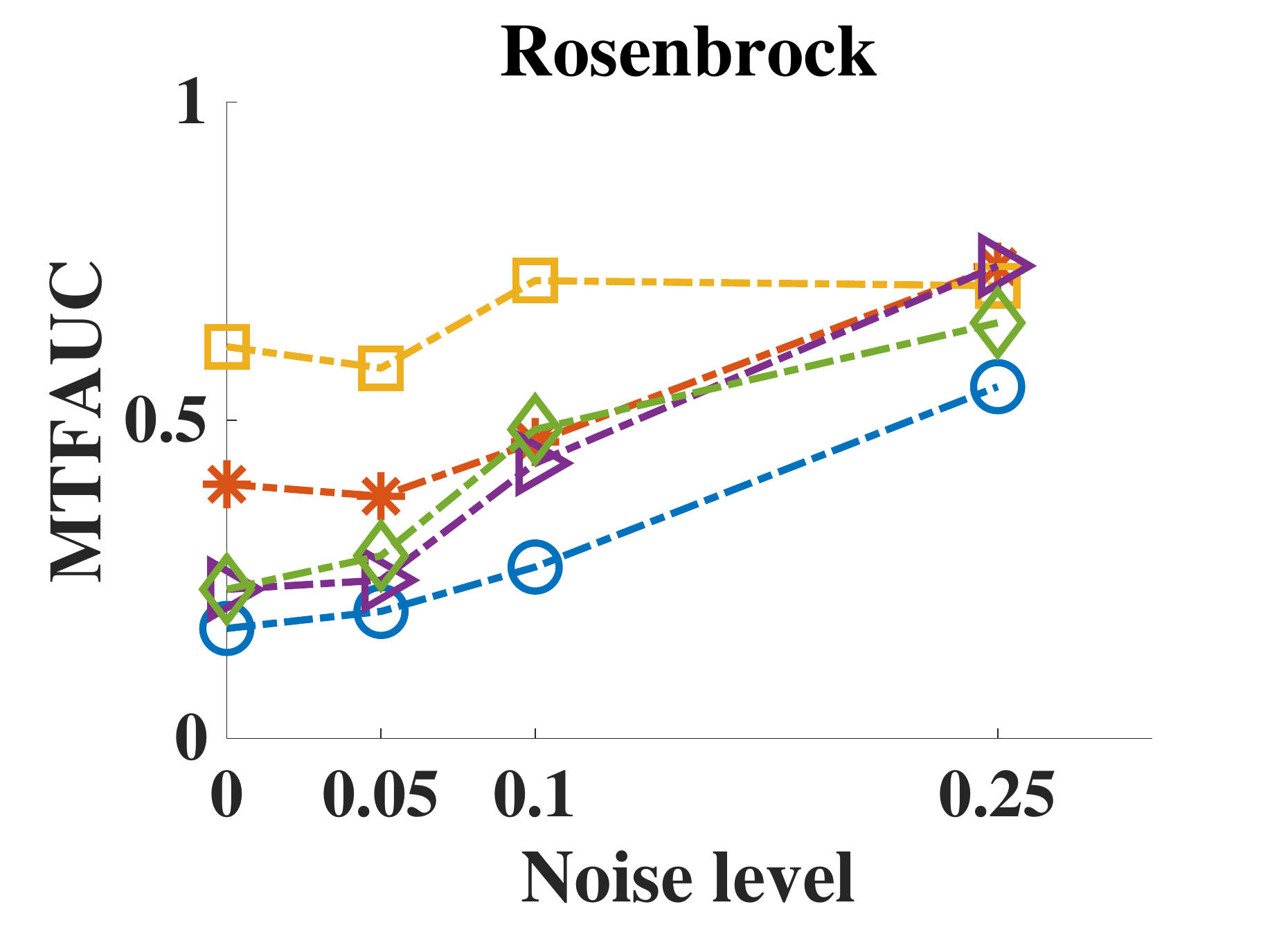}} &
\subfloat[nonRBF]{\includegraphics[width = 1.5in]{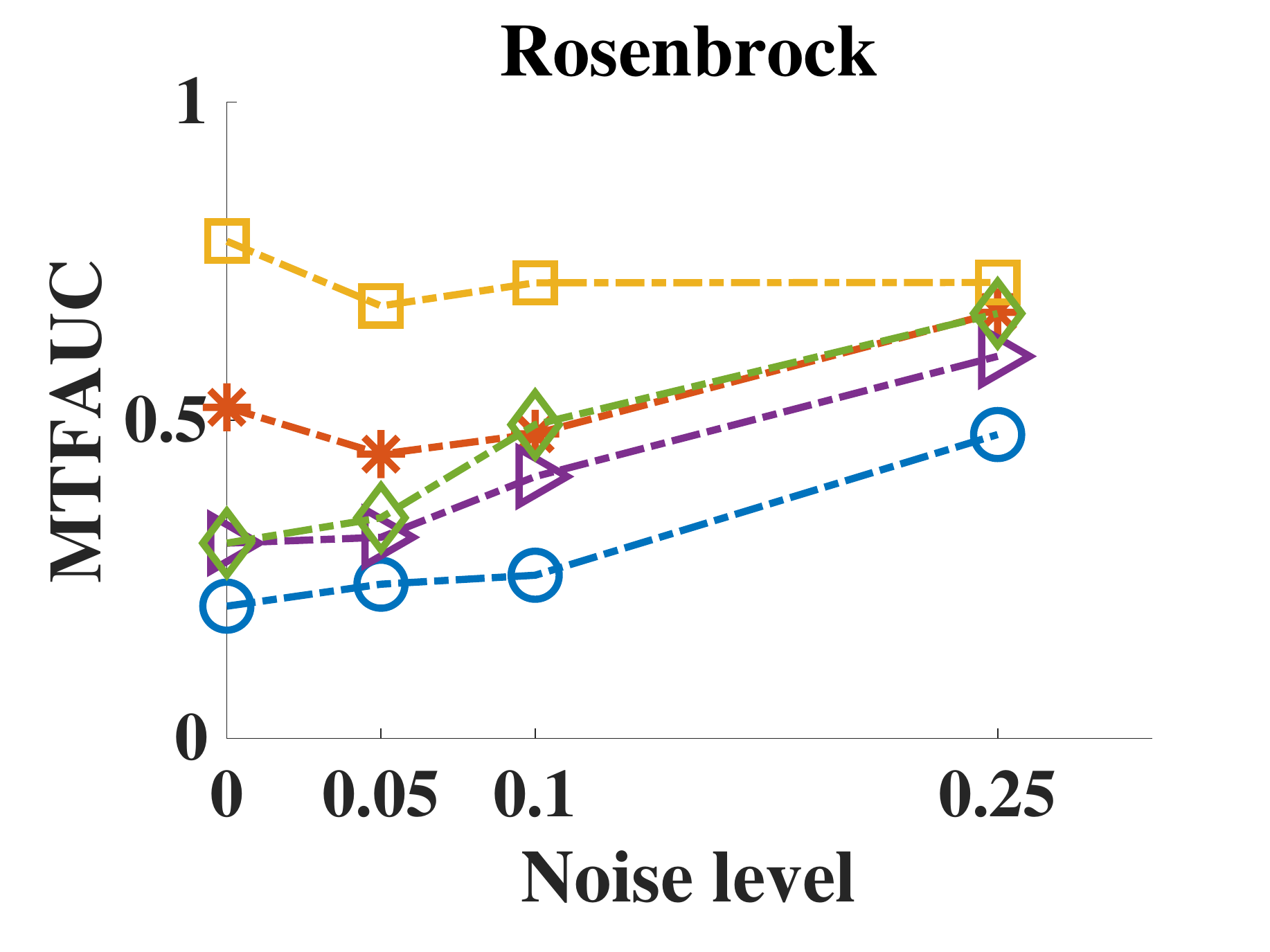}} &
\subfloat[nonGP]{\includegraphics[width = 1.5in]{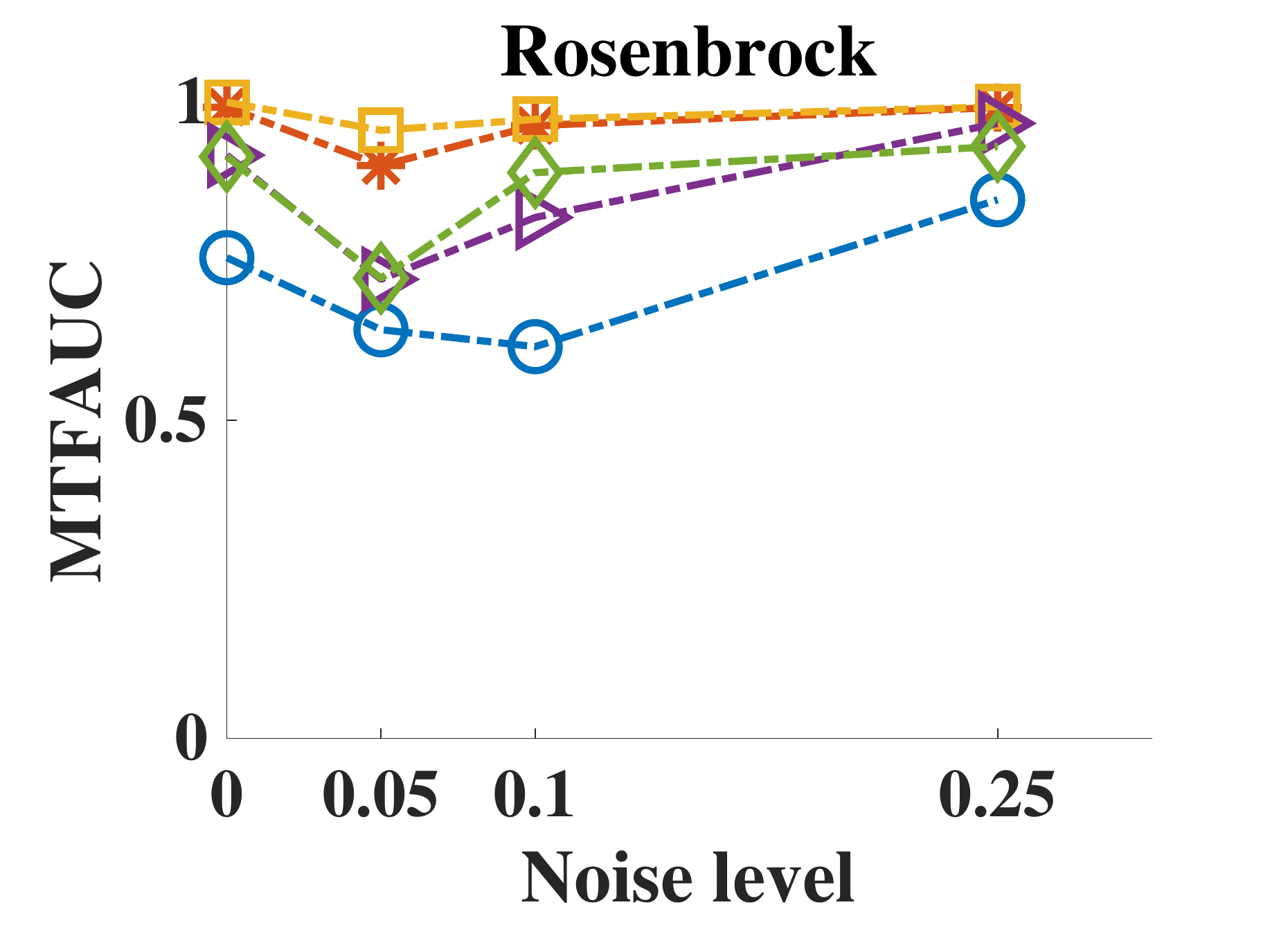}}\\
\subfloat[RBF]{\includegraphics[width = 1.5in]{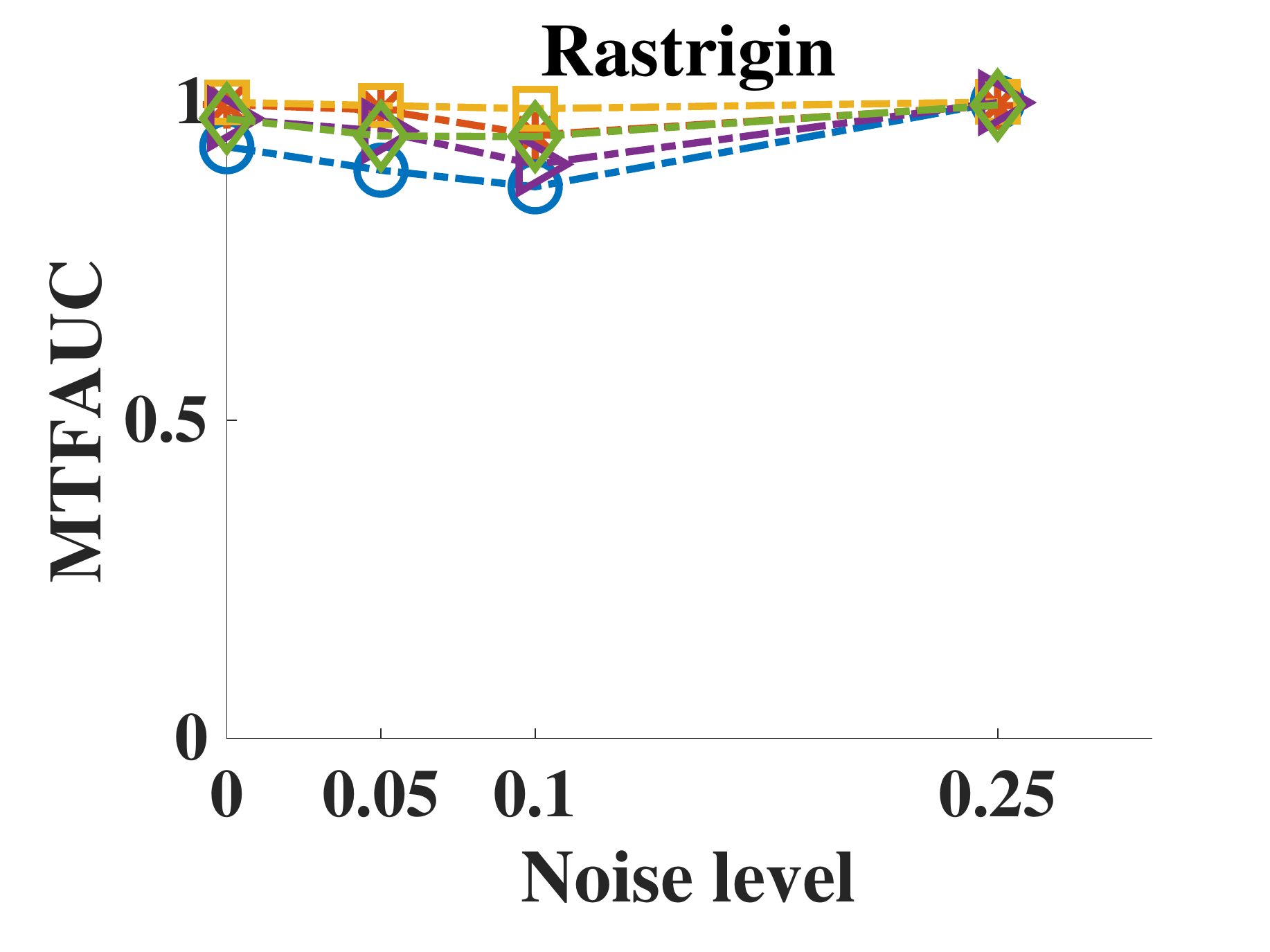}} &
\subfloat[TK-MARS]{\includegraphics[width = 1.5in]{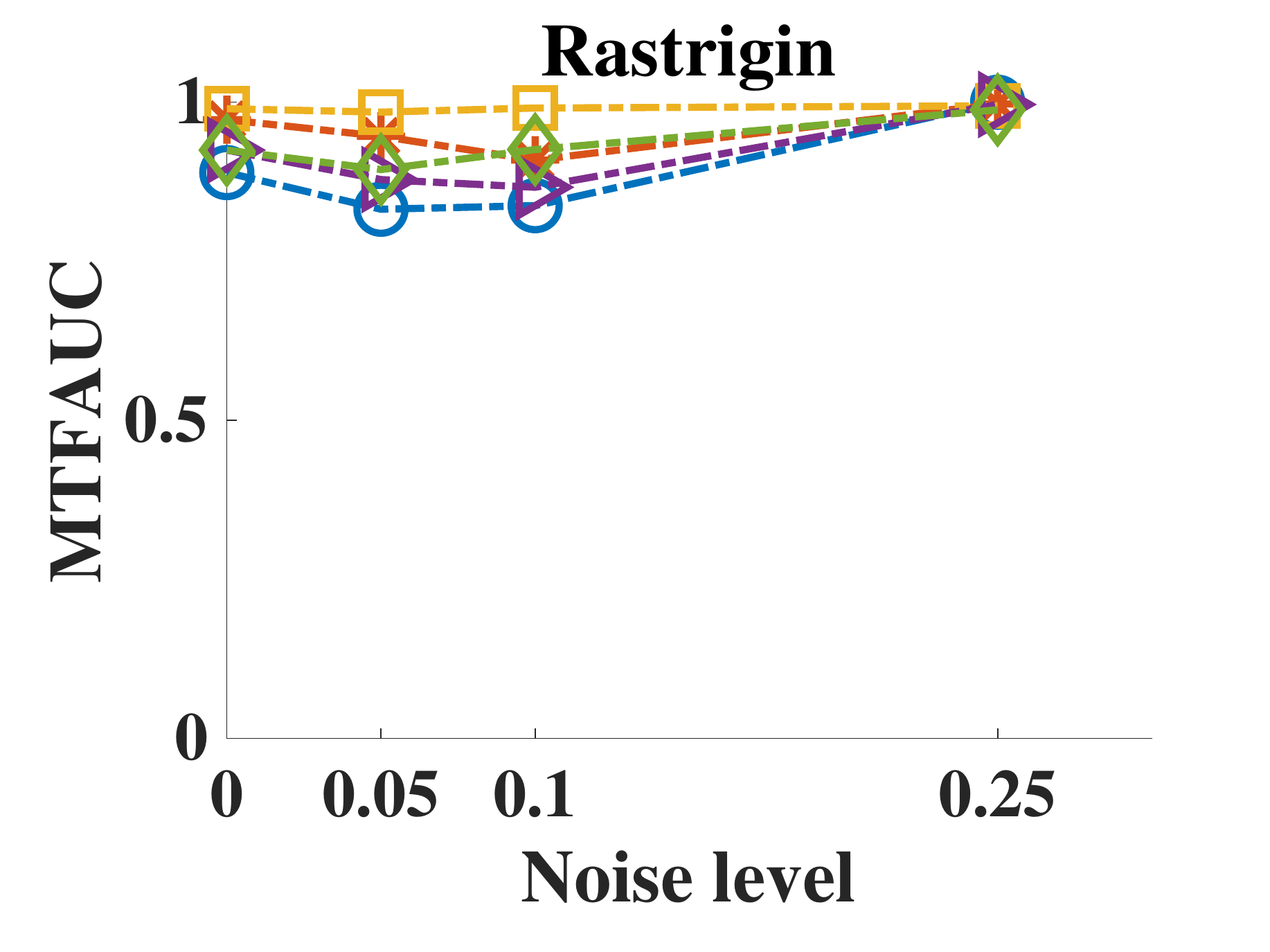}} &
\subfloat[nonRBF]{\includegraphics[width = 1.5in]{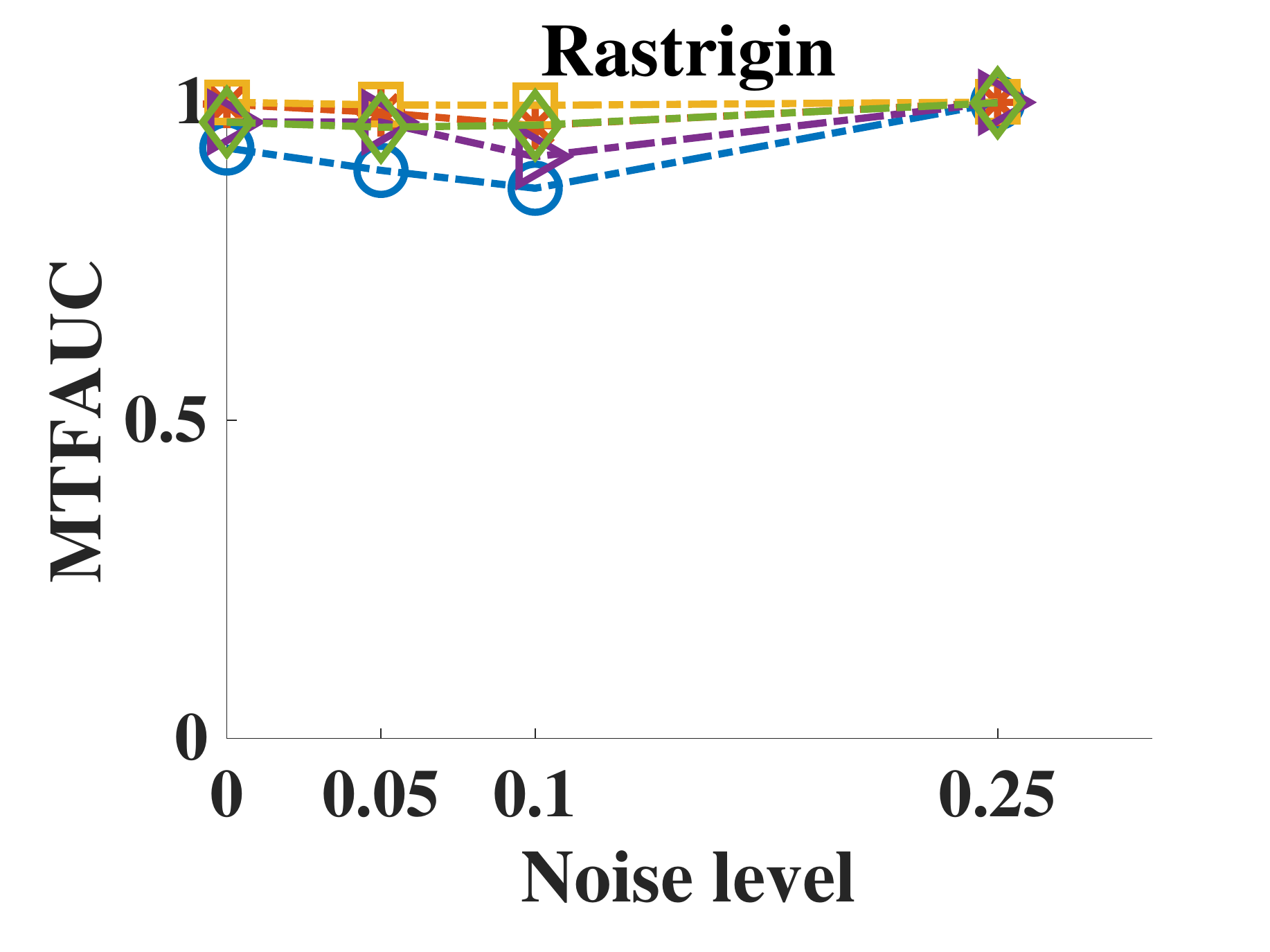}} &
\subfloat[nonGP]{\includegraphics[width = 1.5in]{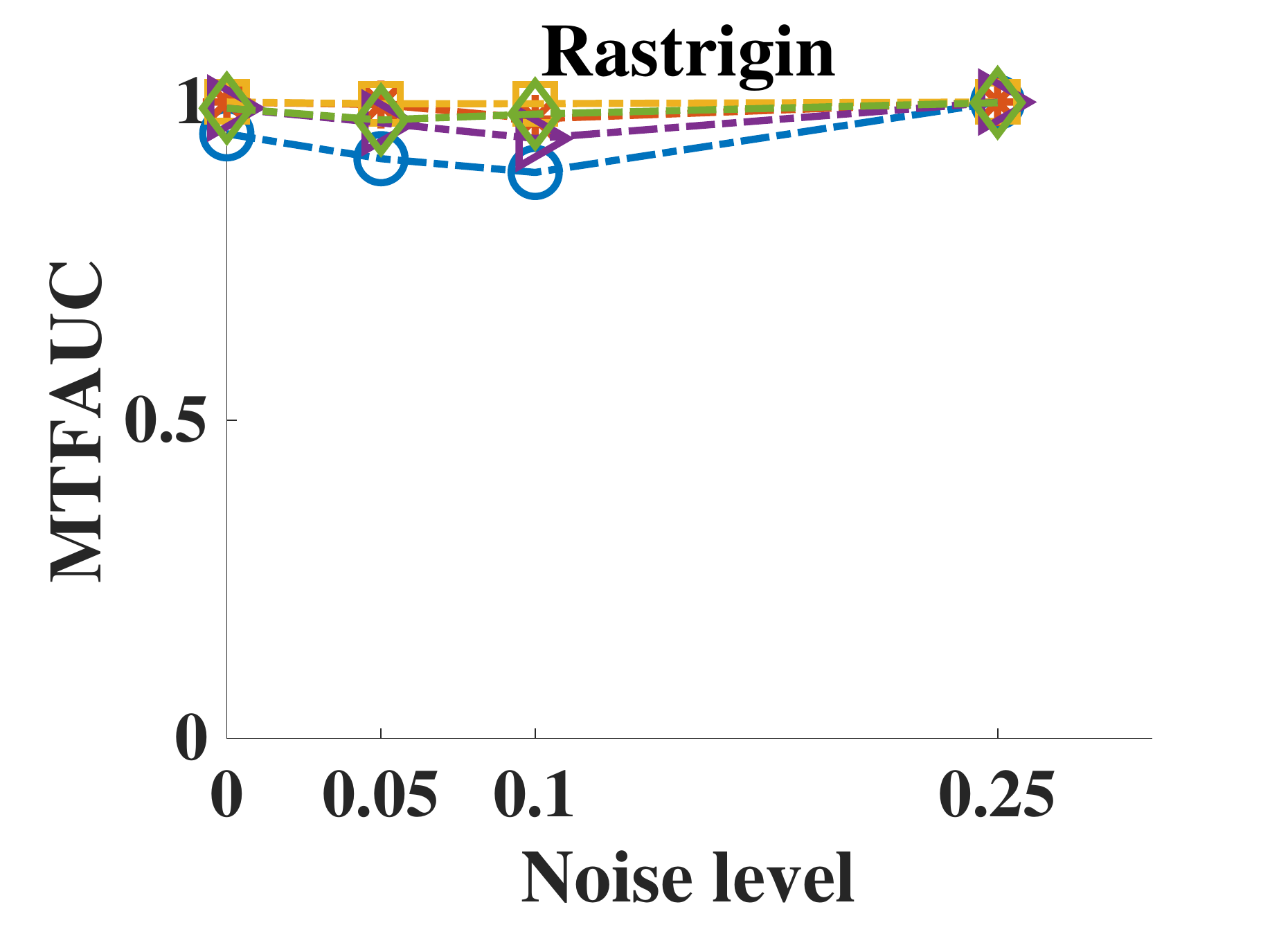}} \\
\subfloat[RBF]{\includegraphics[width = 1.5in]{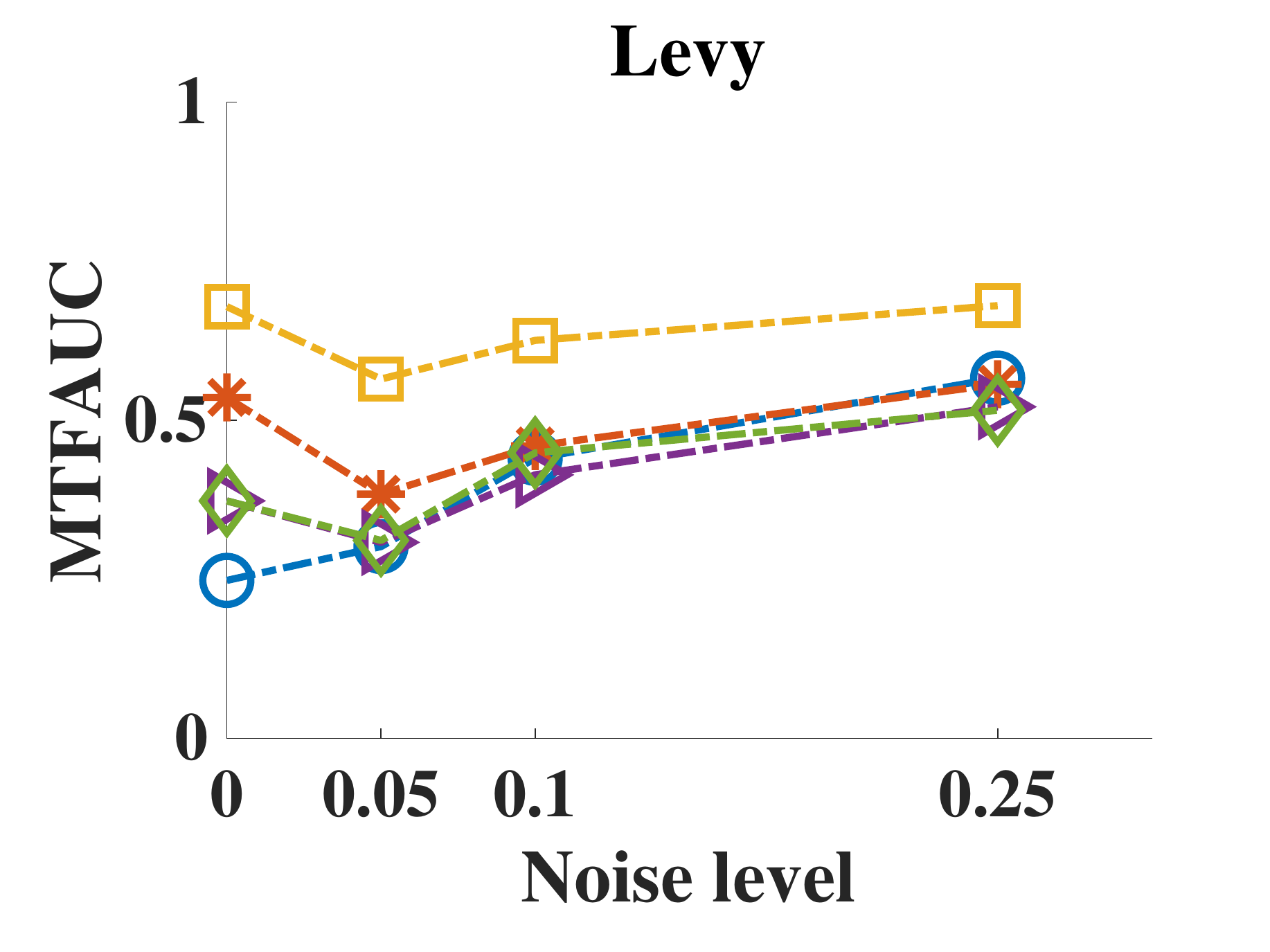}} &
\subfloat[TK-MARS]{\includegraphics[width = 1.5in]{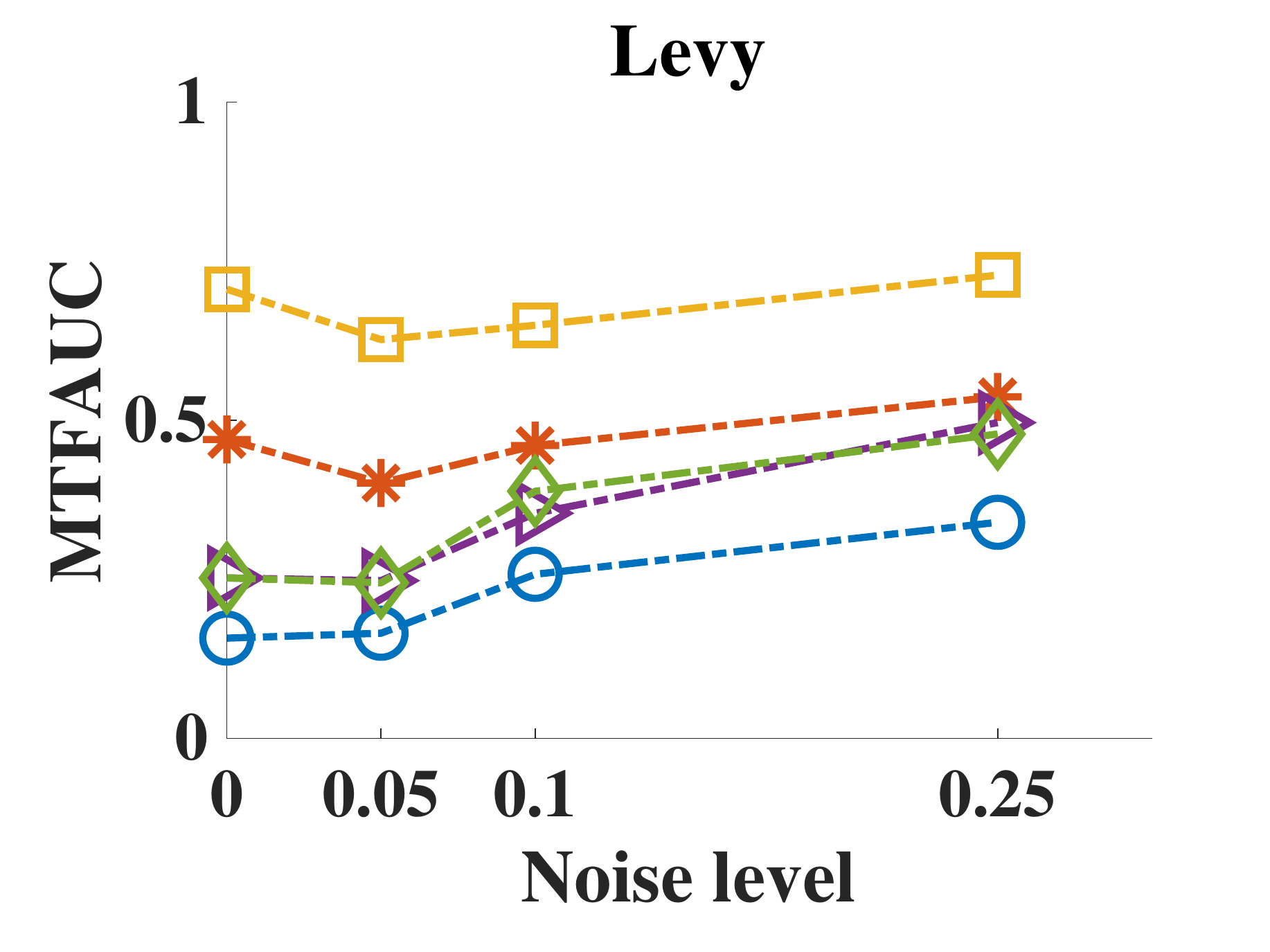}} &
\subfloat[nonRBF]{\includegraphics[width = 1.5in]{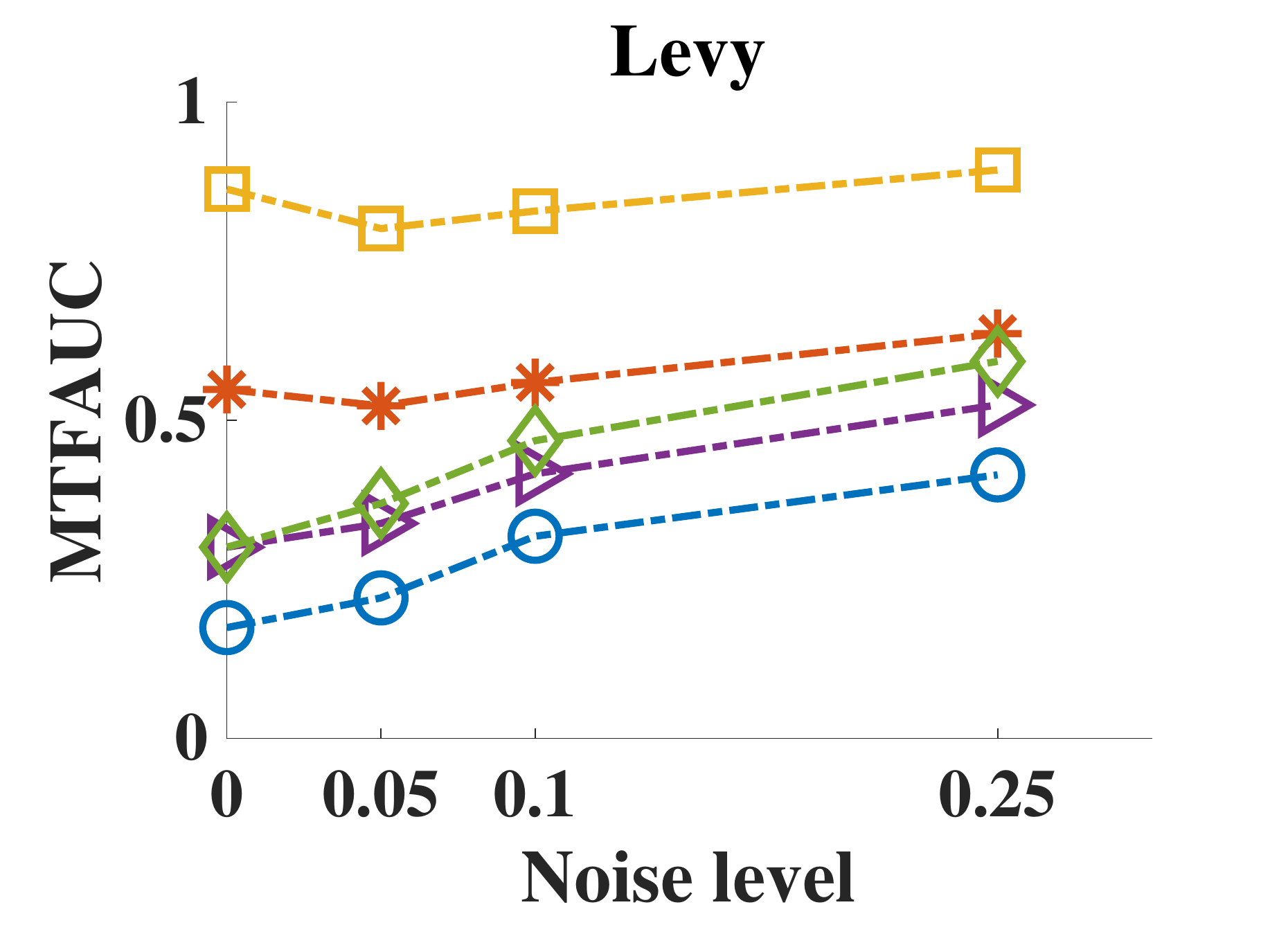}} &
\subfloat[nonGP]{\includegraphics[width = 1.5in]{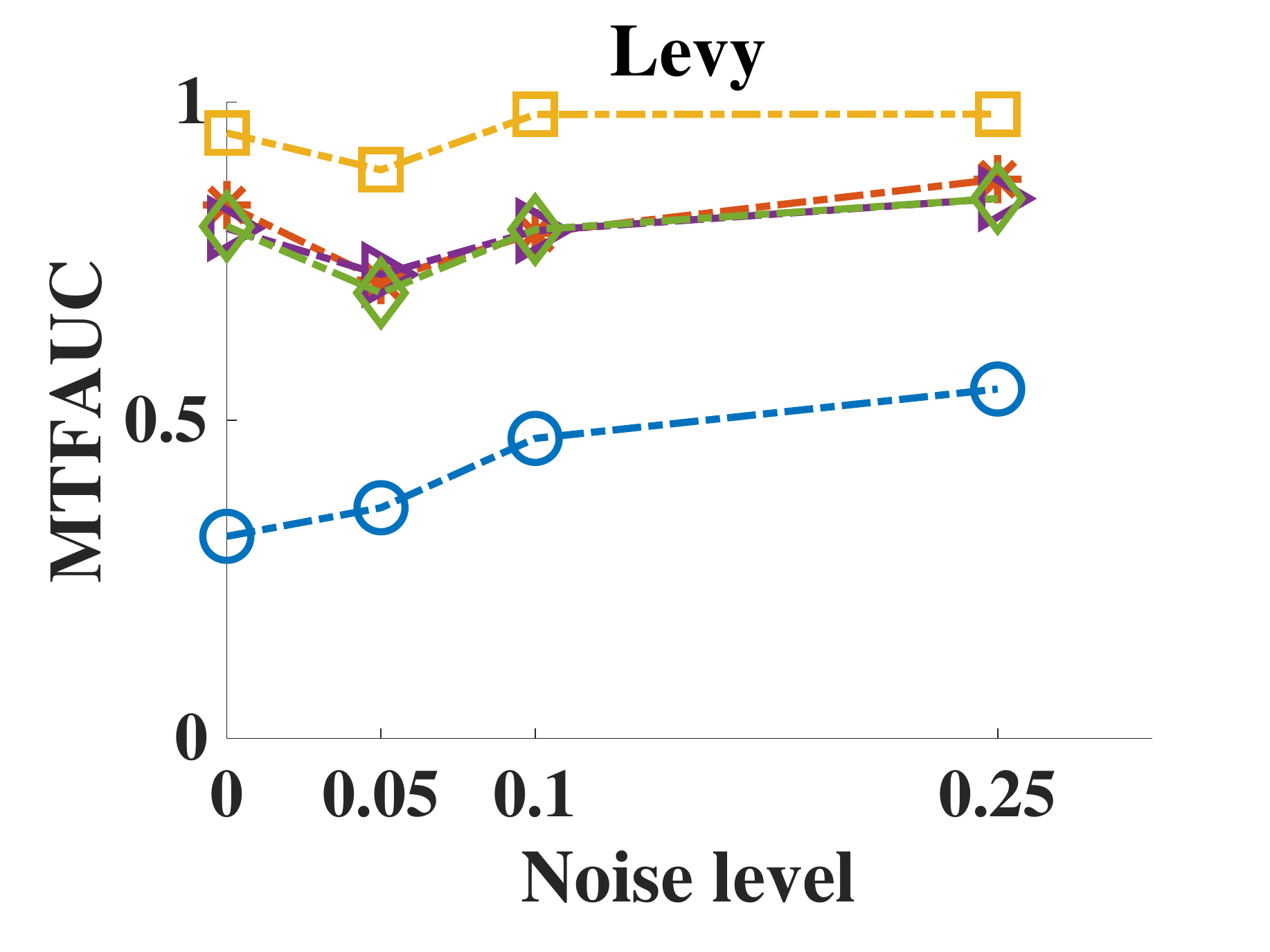}}\\
\subfloat[RBF]{\includegraphics[width = 1.5in]{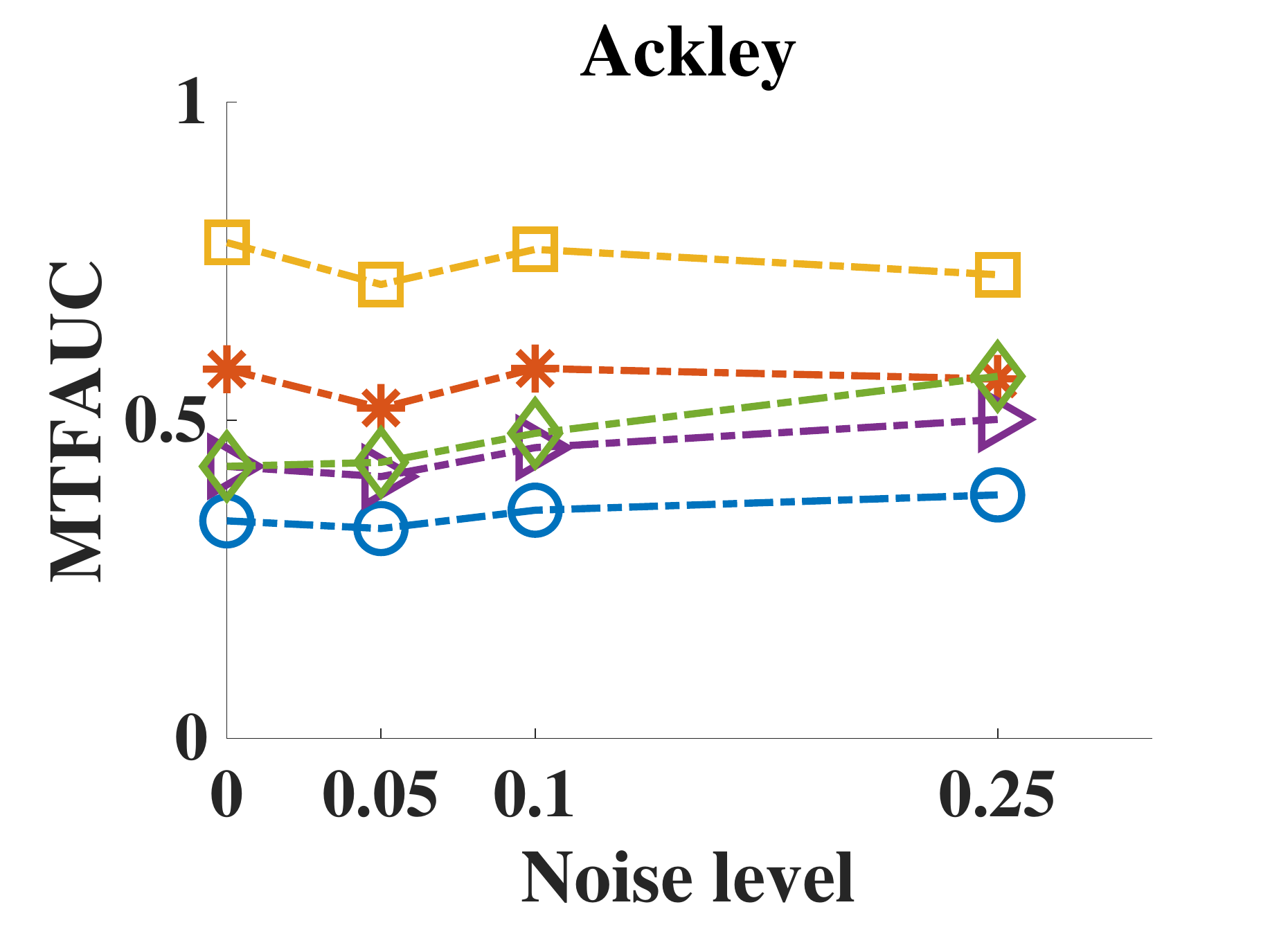}} &
\subfloat[TK-MARS]{\includegraphics[width = 1.5in]{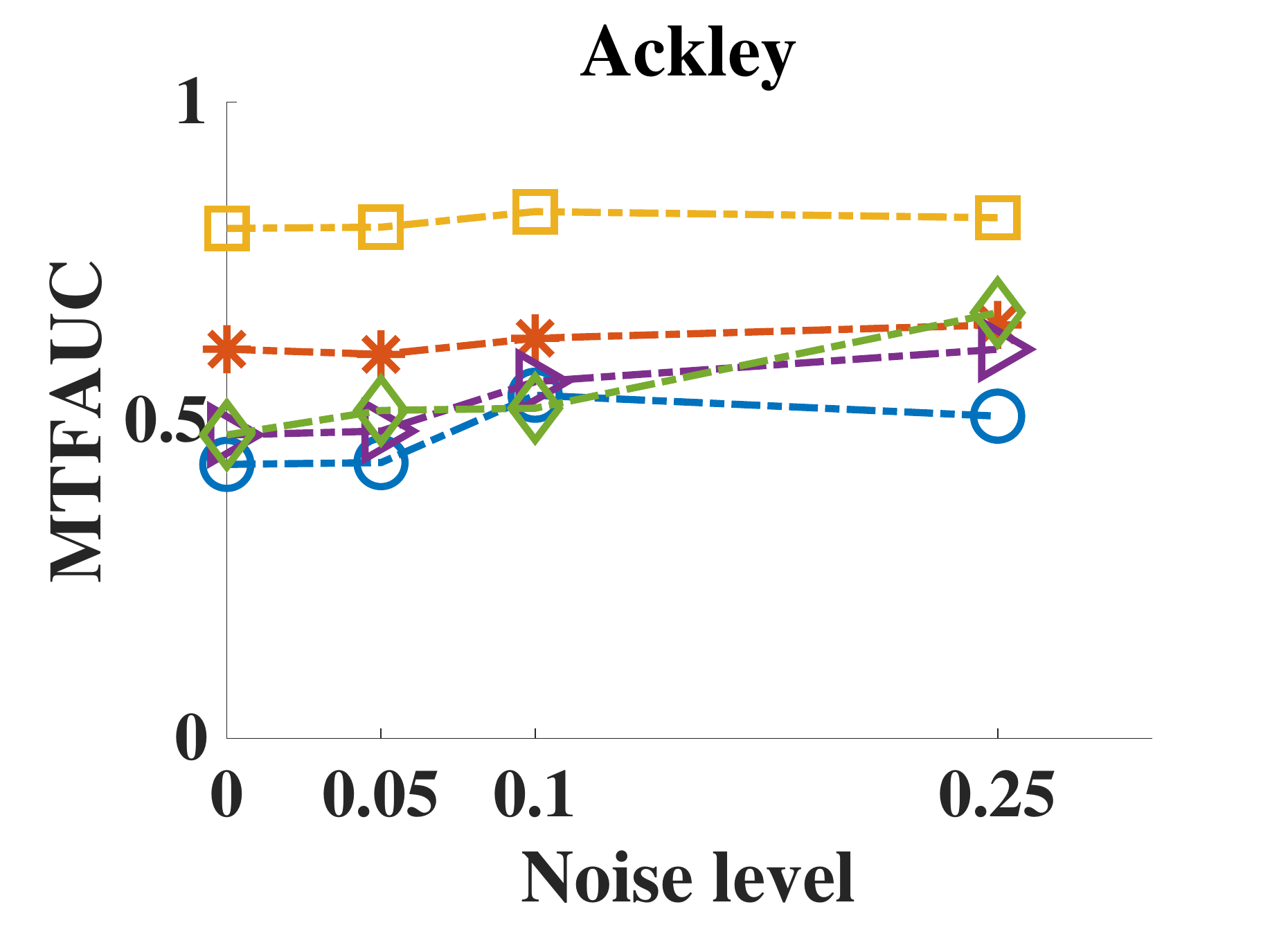}} &
\subfloat[nonRBF]{\includegraphics[width = 1.5in]{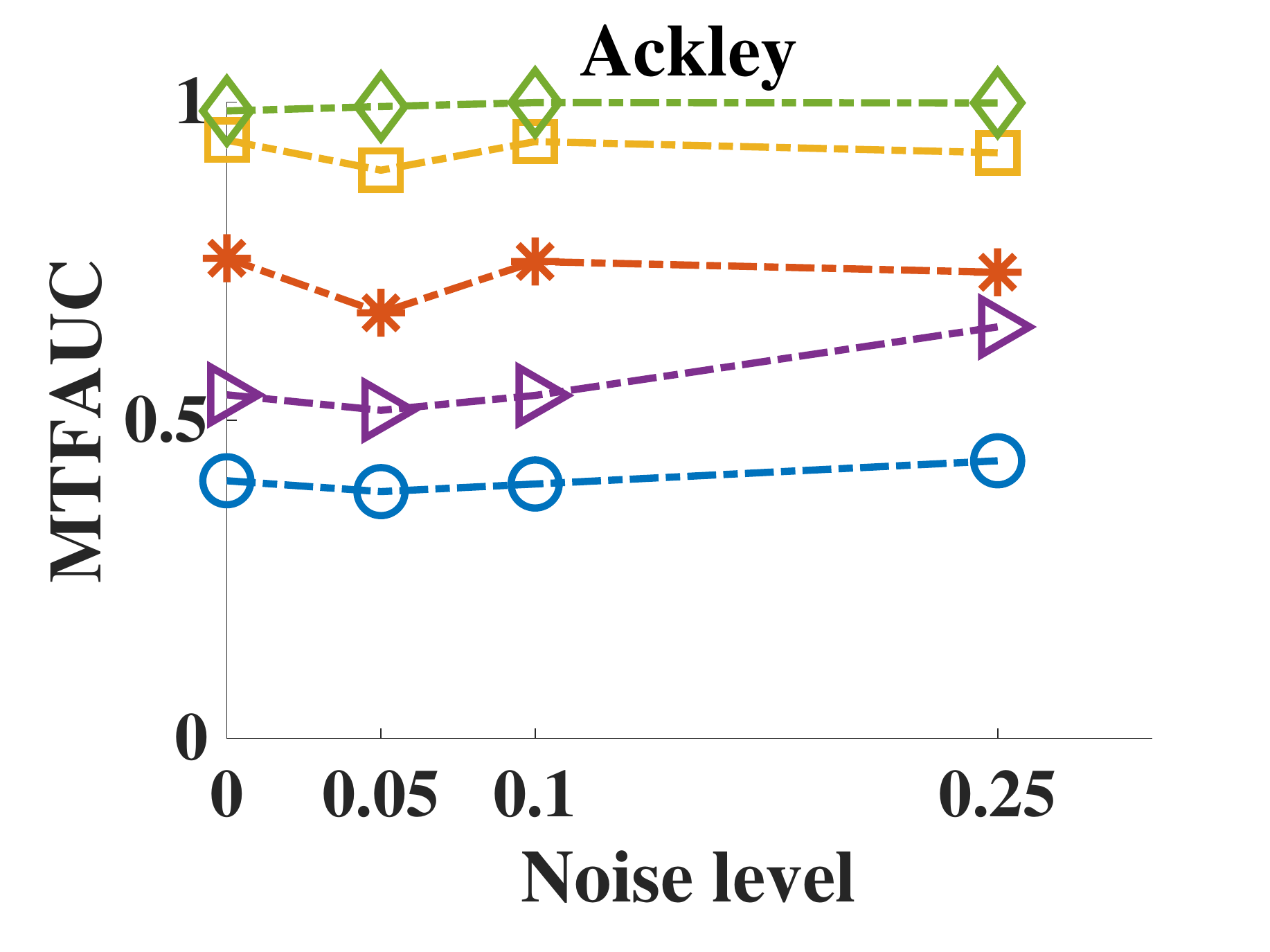}} &
\subfloat[nonGP]{\includegraphics[width = 1.5in]{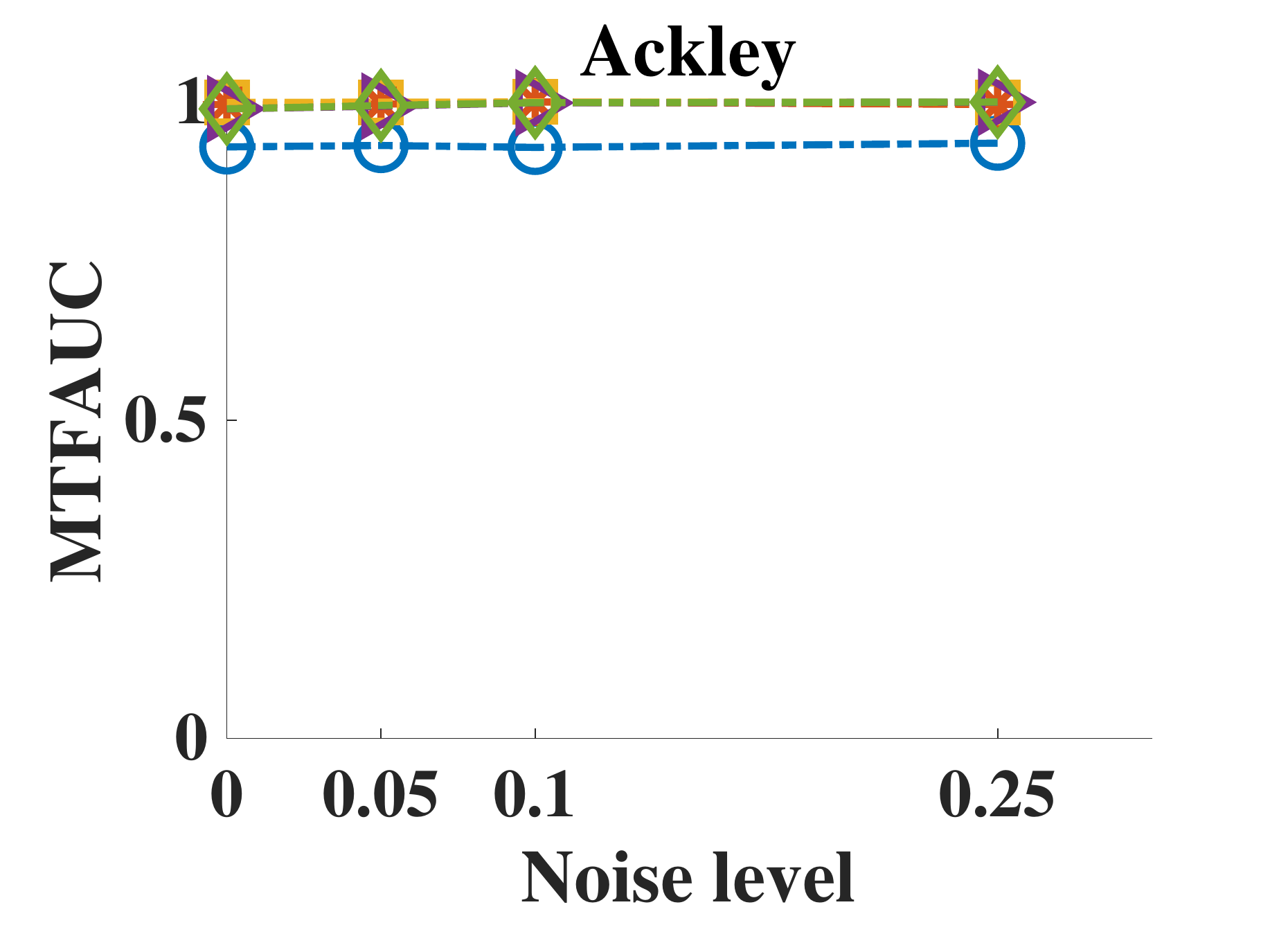}}\\
\subfloat[RBF]{\includegraphics[width = 1.5in]{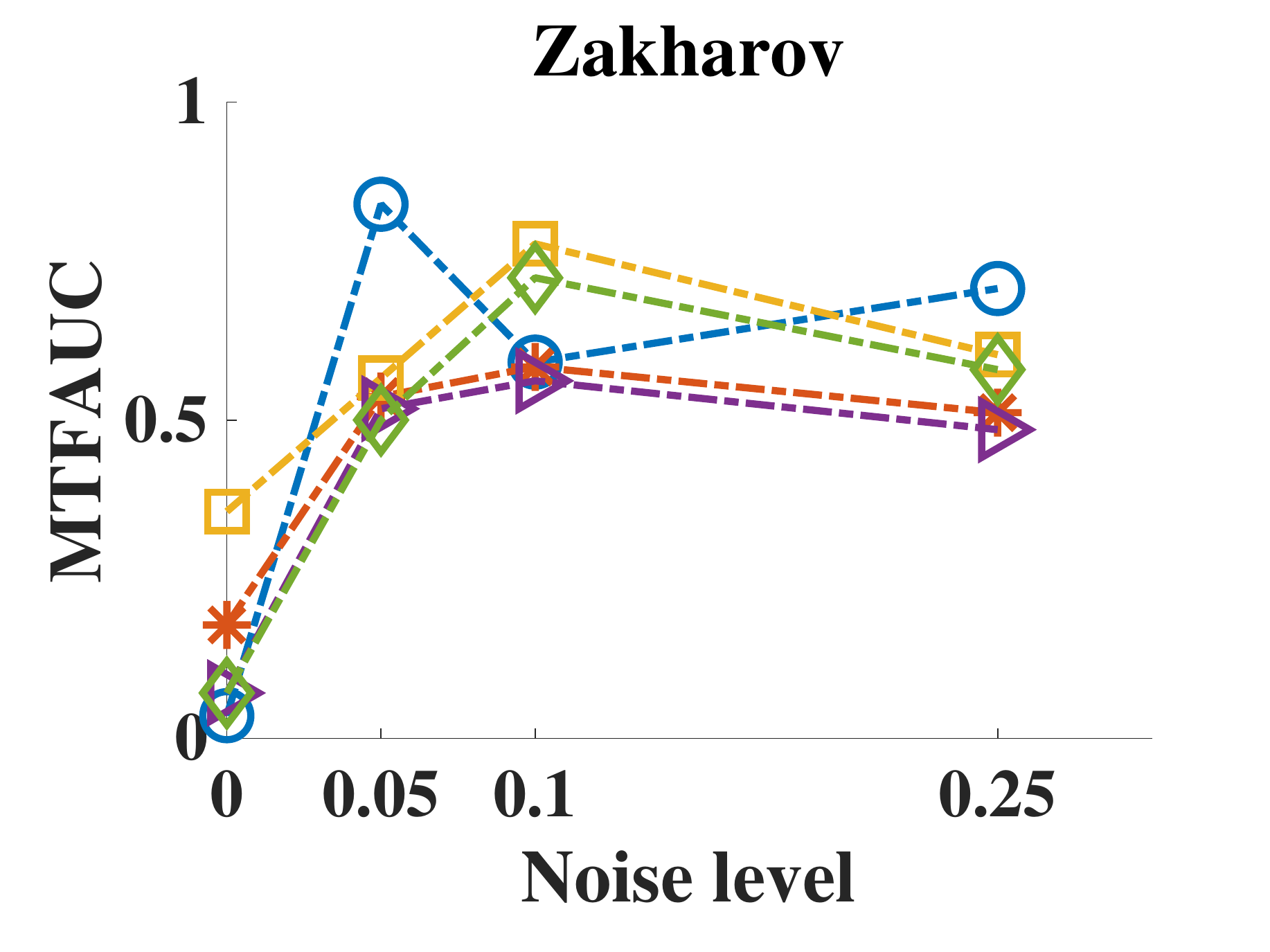}} &
\subfloat[TK-MARS]{\includegraphics[width = 1.5in]{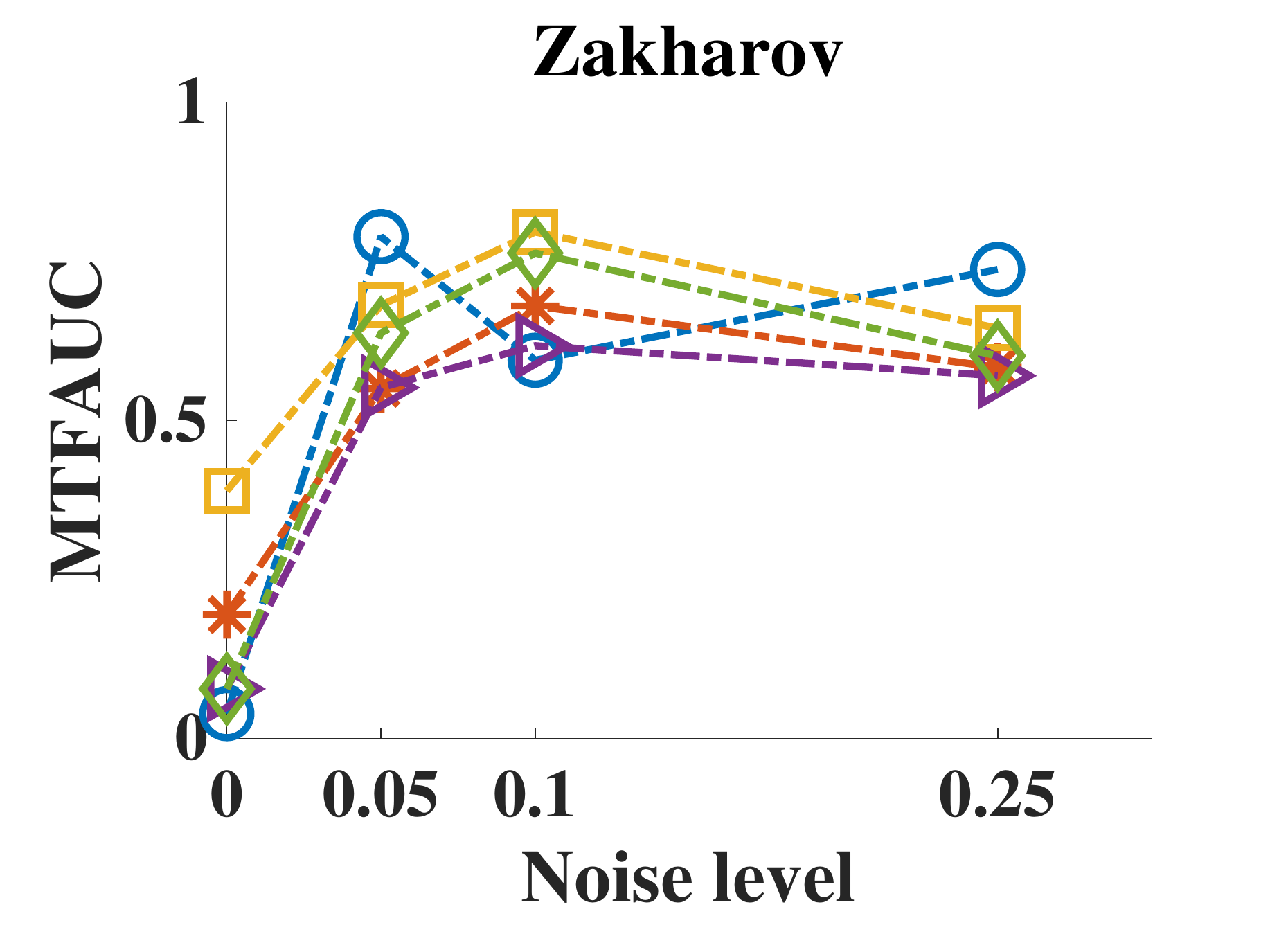}} &
\subfloat[nonRBF]{\includegraphics[width = 1.5in]{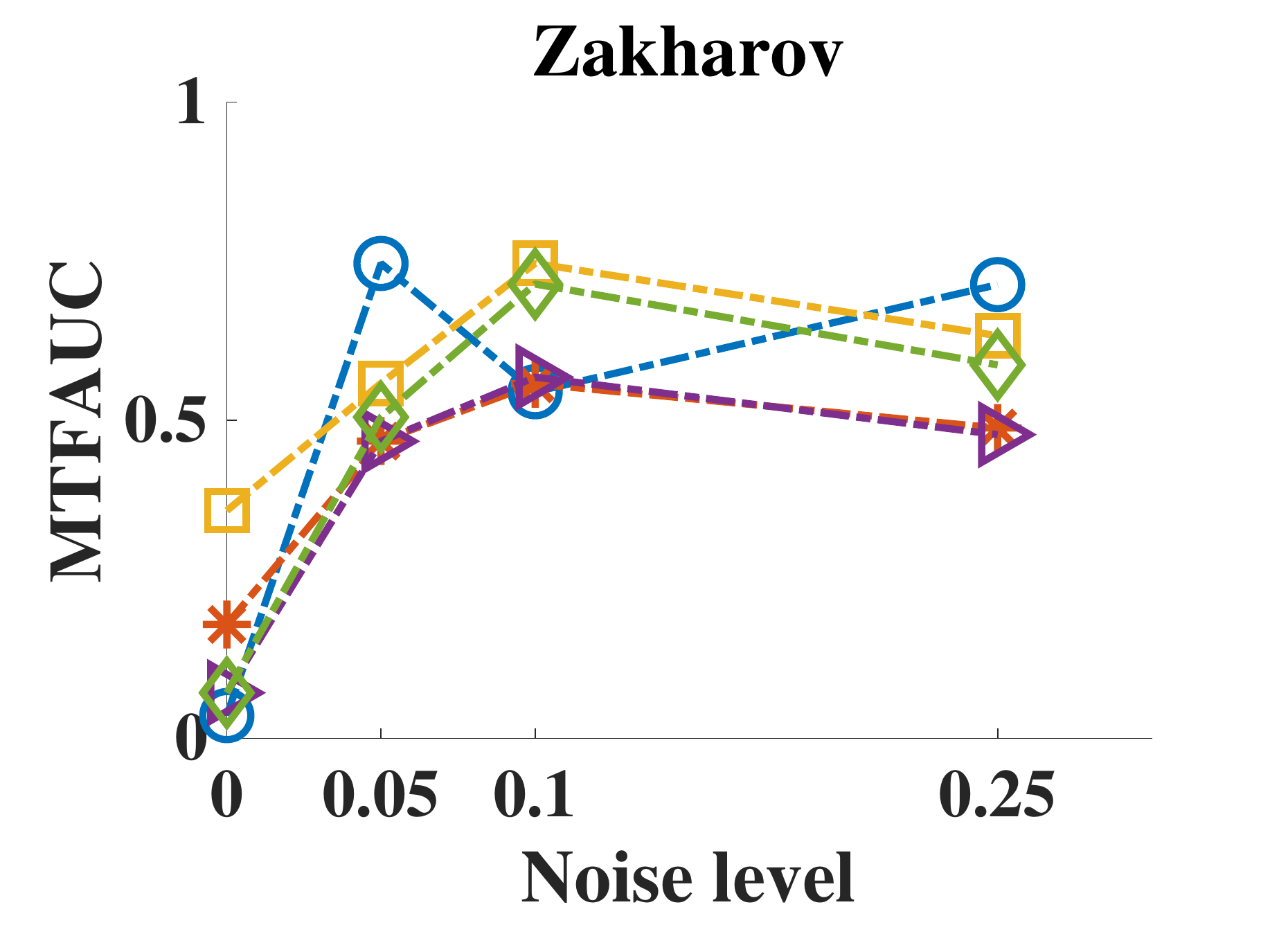}} &
\subfloat[nonGP]{\includegraphics[width = 1.5in]{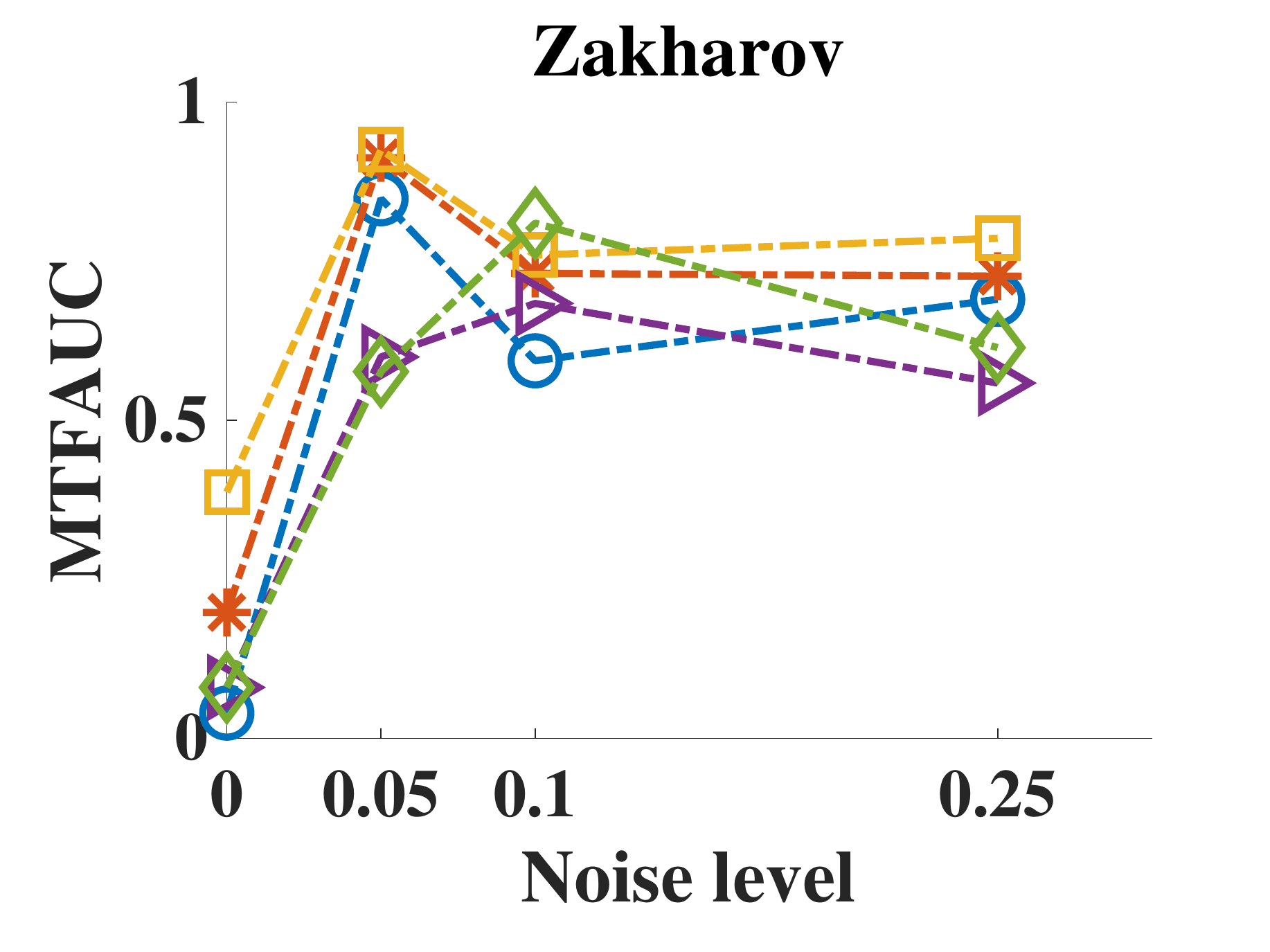}}\\

\end{tabular}
\caption{Average MTFAUC across different noise levels}
\label{fig:mtfauc}
\end{figure}

\new{Figures~\ref{fig:mtfauc}(a)-(t) show the average performance of No-Replication, Fixed- Replication, and Smart-Replication for 30 executions on the test functions in Table~\ref{tab:testfuncchar}. Each row represents a test function, and each column corresponds to a surrogate model.} Tables~\ref{tab:mtfauc1},~\ref{tab:mtfauc2}, and~\ref{tab:mtfauc3} in Appendix B, show corresponding MTFAUC values.
% Table generated by Excel2LaTeX from sheet 'Sheet1'
\hadis{The results indicate that when no noise is present, interpolating is a good choice for surrogate fitting. When the objective function is noisy, interpolating may lead to a poor approximation and misleads the optimization. Note that No-Replication outperforms the replication approaches in some cases when we use non-interpolating surrogates. That is due to the fact that non-interpolating surrogates minimize the average error across all data points, which collectively reduces the impact of the noise. 
However, we realized non-interpolating models can capture the noise without replication when the objective function is not highly complex (highly oscillated). In such situations, where ``predictability'' is high, the model can use the value of other points to cancel out the noise for a specific point.
To further clarify this, let us consider the toy example in Figure \ref{fig:toy}. In Figure~\ref{fig:toy} (a), the underlying objective function is linear. We have highlighted a red star point, for which the function value is highly affected by the noise (dashed line). One can observe that minimizing the error for all the points, the non-interpolating model (linear regression) could cancel out the effect of the red star point.}

\hadis{This however, is not the case for ``unpredictable'' and highly-oscillating functions. That simply is because, due to the high oscillations, other points (slightly far from the star noisy point) cannot be used for estimating the correct value of the noisy point. We highlighted this in the toy example shown in \new{Figure~\ref{fig:toy}} (b) (Ackley function estimated by (non-interpolating) quadratic model). As one can confirm, assuming that the function had a peak (valley) at the noisy point, the existing data points cannot predict the depth of the peak. In such situations, replications are needed to cancel out the impact of noise. In \new{Figure~\ref{fig:toy}} (c), we added the replications for the star noisy point, the points shown as triangles. We can observe the non-interpolating model has a better estimation of the true function at the red star point, and was more effective in mitigating the noise effect of that point.}

\hadis{In summary, while non-interpolating methods may perform well for non-complex noisy environments, they fail to work well for complex noisy functions without replications.
We observed a similar behavior in our experiments where no-replication failed to work well for complex functions, such as Zakharov, Ackley, and Rastrigin.}

\begin{figure}[!tb]
\centering
    \begin{minipage}{\linewidth}
        \subfloat[Noise=0]{\includegraphics[width=0.55\textwidth]{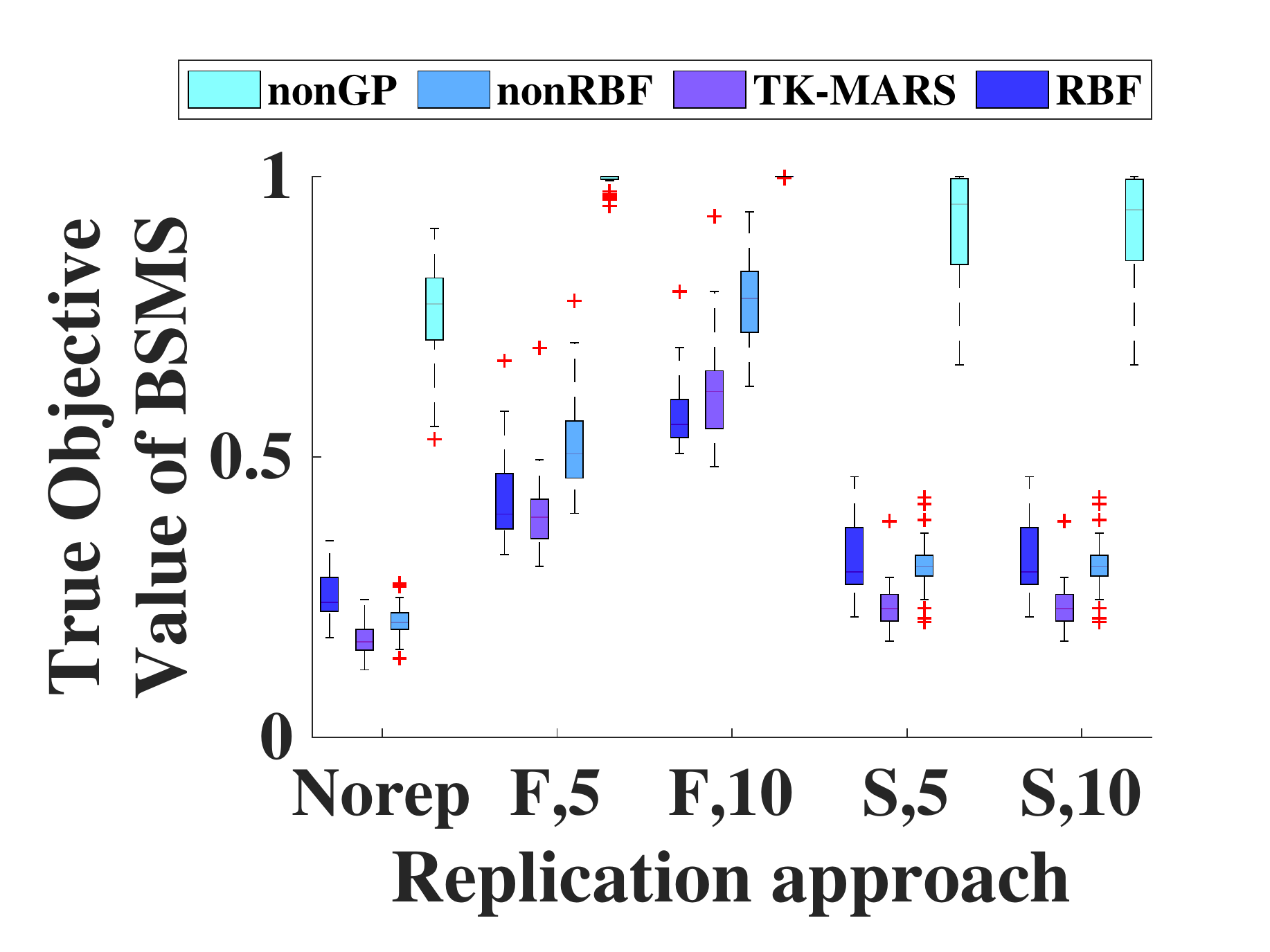}}
        \subfloat[Noise=0.05]{\includegraphics[width=0.55\textwidth]{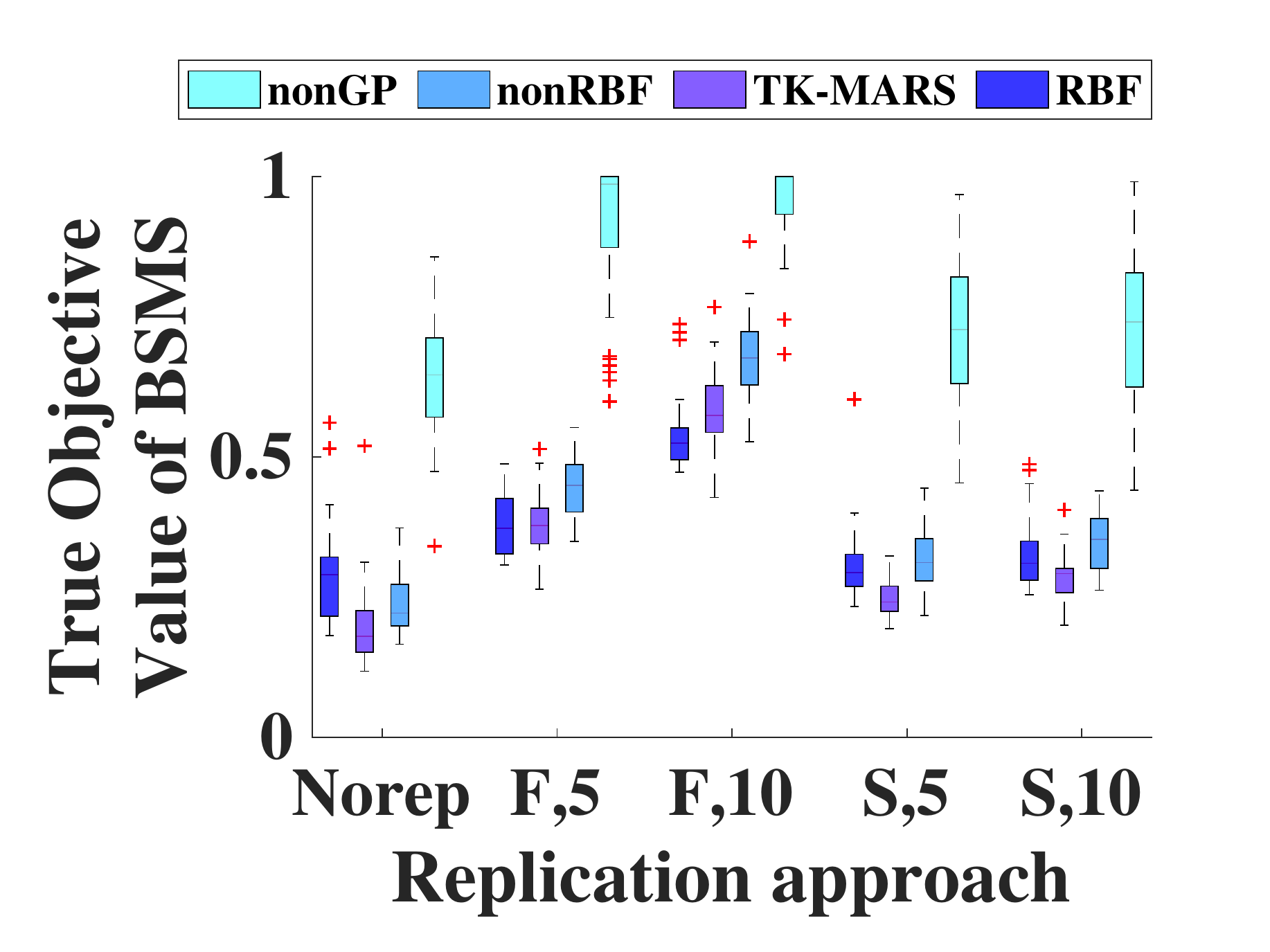}}
    \end{minipage}
    \begin{minipage}{\linewidth}
        \subfloat[Noise=0.1]{\includegraphics[width=0.55\textwidth]{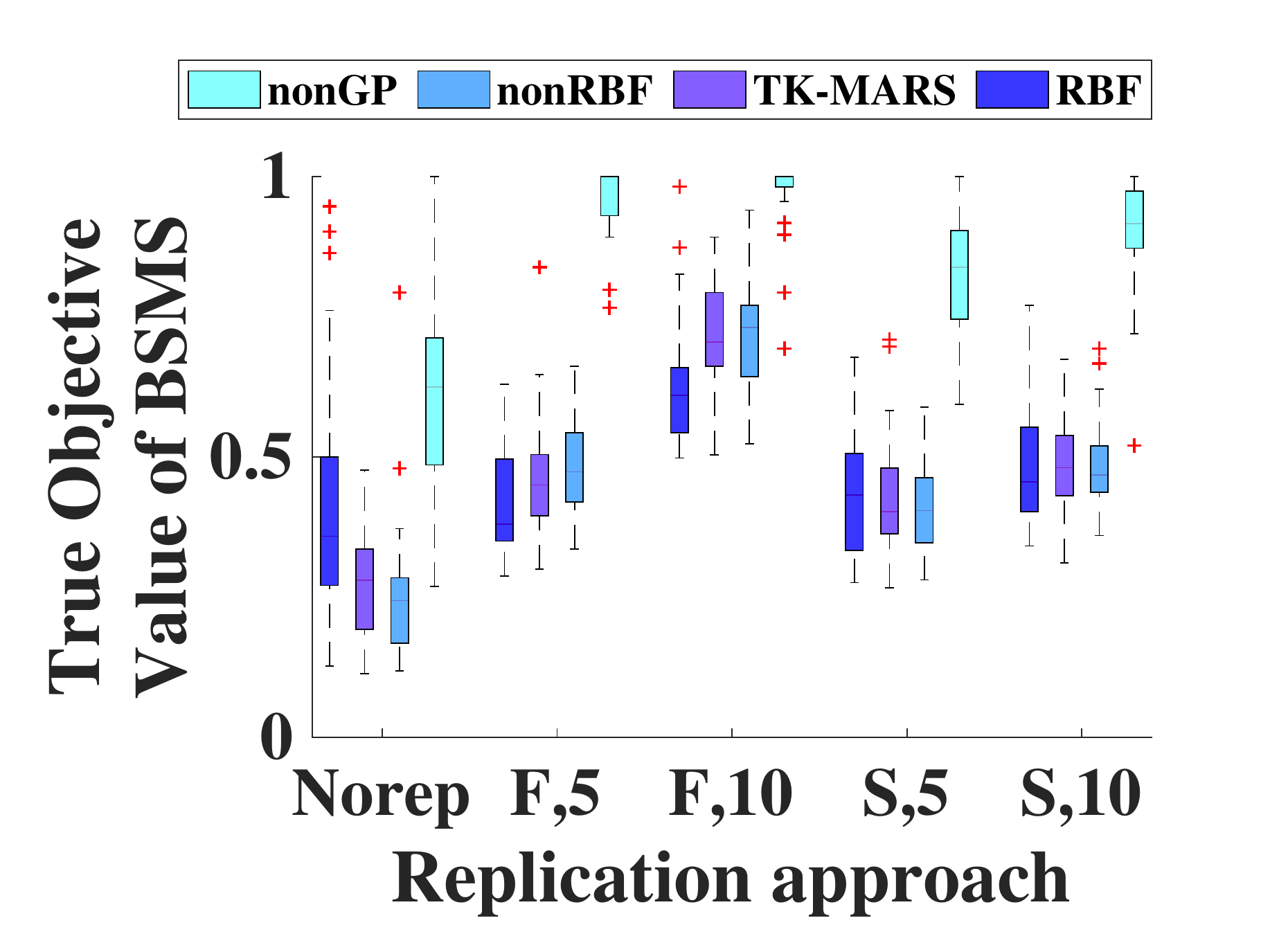}}
        \subfloat[Noise=0.25]{\includegraphics[width=0.55\textwidth]{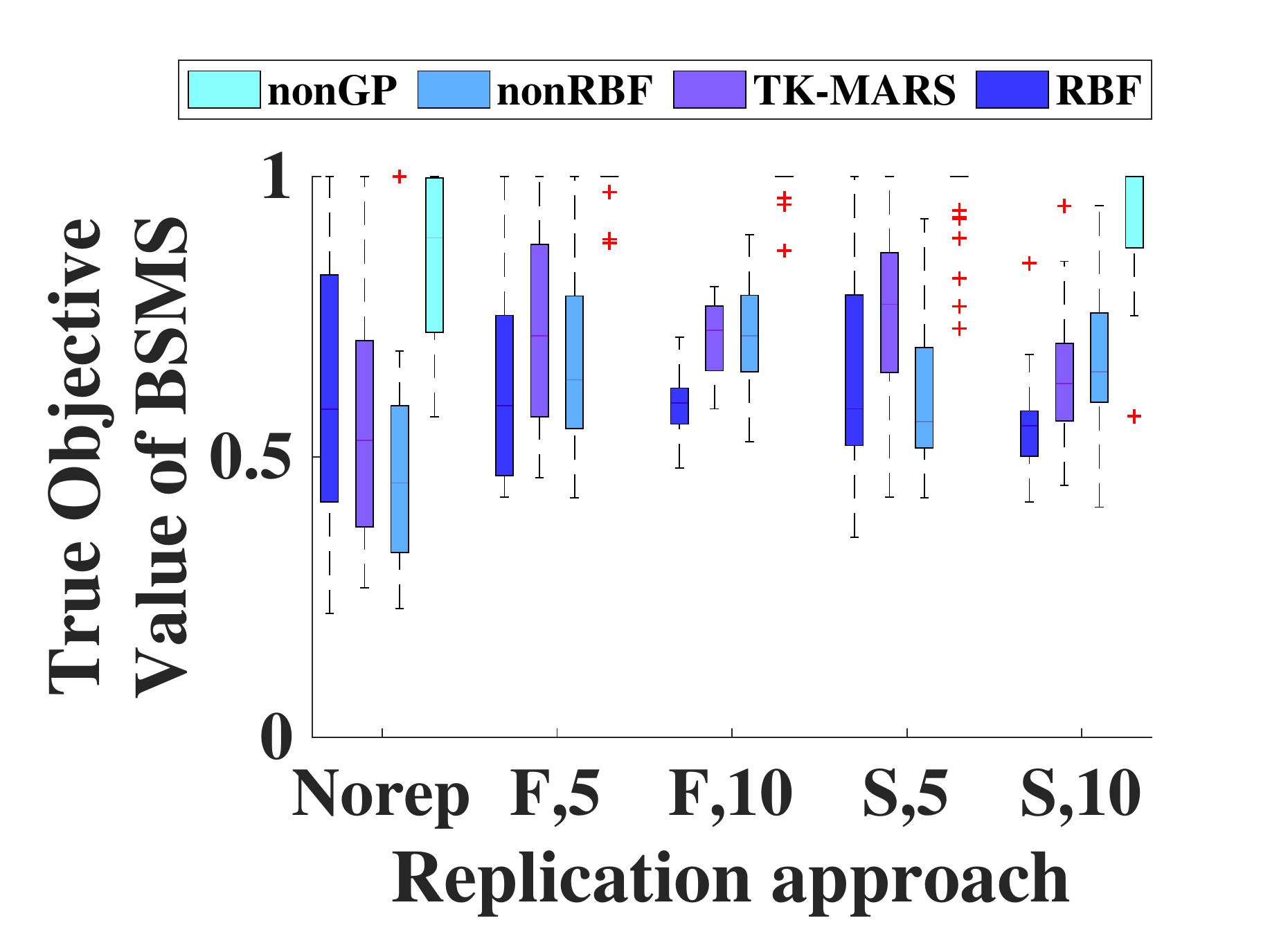}}
    \end{minipage}
    \caption{Box-plots of MTFAUC of surrogate optimization on the Rosenbrock function}
    \label{fig:auc_box_rosen}
    \vspace{-3mm}
\end{figure}

% \begin{figure}[!tb]
% \centering
%     \begin{minipage}{\linewidth}
%         \subfloat[]{\includegraphics[width=0.35\textwidth]{figures/box_auc_rosen_new_1.pdf}}
%         \subfloat[]{\includegraphics[width=0.35\textwidth]{figures/box_auc_rosen_new_2.pdf}}
%         \subfloat[]{\includegraphics[width=0.35\textwidth]{figures/box_auc_rosen_new_3.pdf}}
%     \end{minipage}
%     \begin{minipage}{\linewidth}
%         \subfloat[]{\includegraphics[width=0.35\textwidth]{figures/box_auc_rosen_new_4.pdf}}
%         \subfloat[]{\includegraphics[width=0.35\textwidth]{figures/box_auc_rosen_new_5.pdf}}
%     \end{minipage}
%     \caption{Box-plots of MTFAUC of surrogate optimization on the Rosenbrock function}
%     \label{fig:auc_box_rosen}
%     \vspace{-3mm}
% \end{figure}

\begin{figure}[!tb]
\centering
    \begin{minipage}{\linewidth}
        \subfloat[Noise=0]{\includegraphics[width=0.55\textwidth]{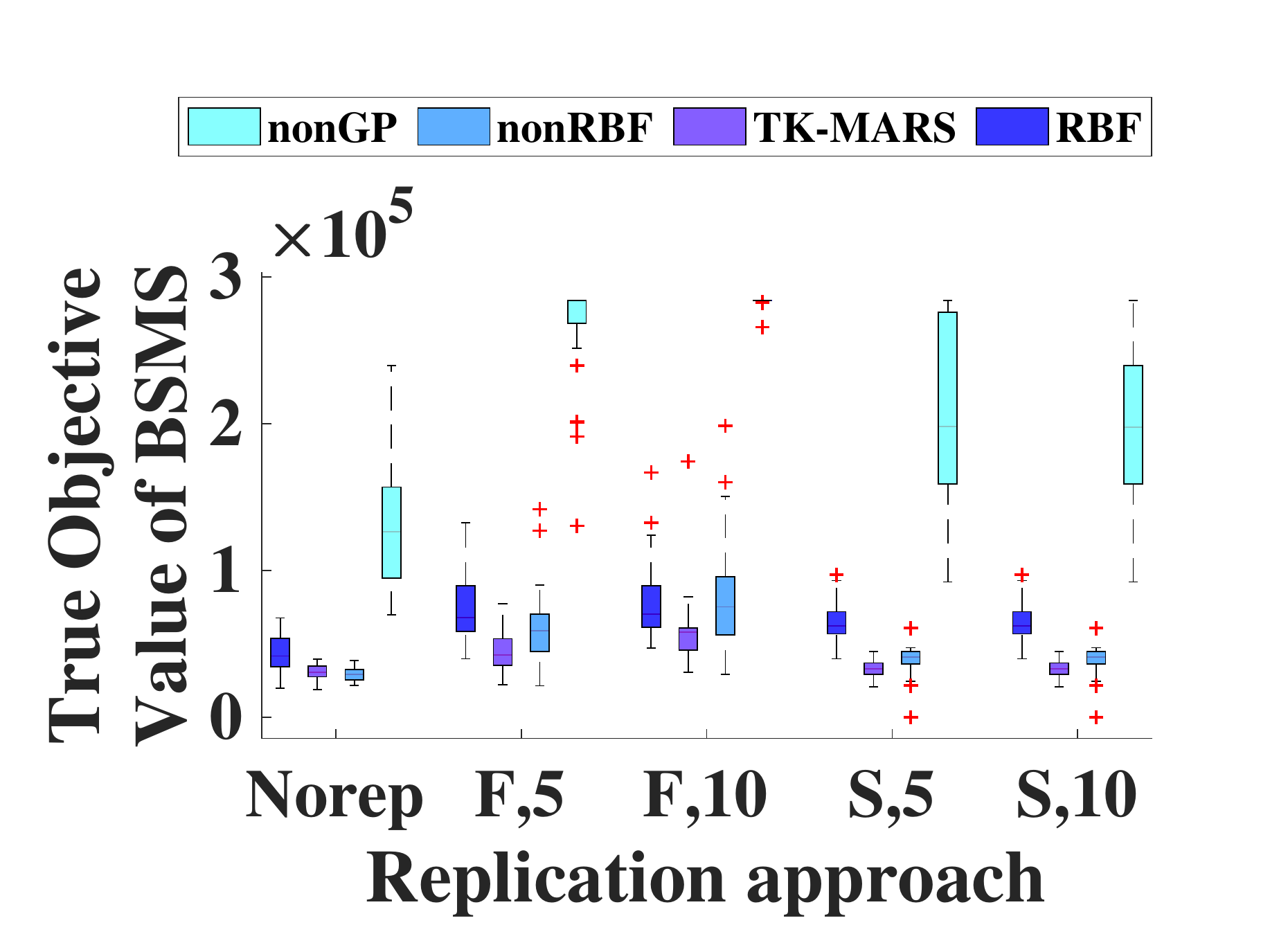}}
        \subfloat[Noise=0.05]{\includegraphics[width=0.55\textwidth]{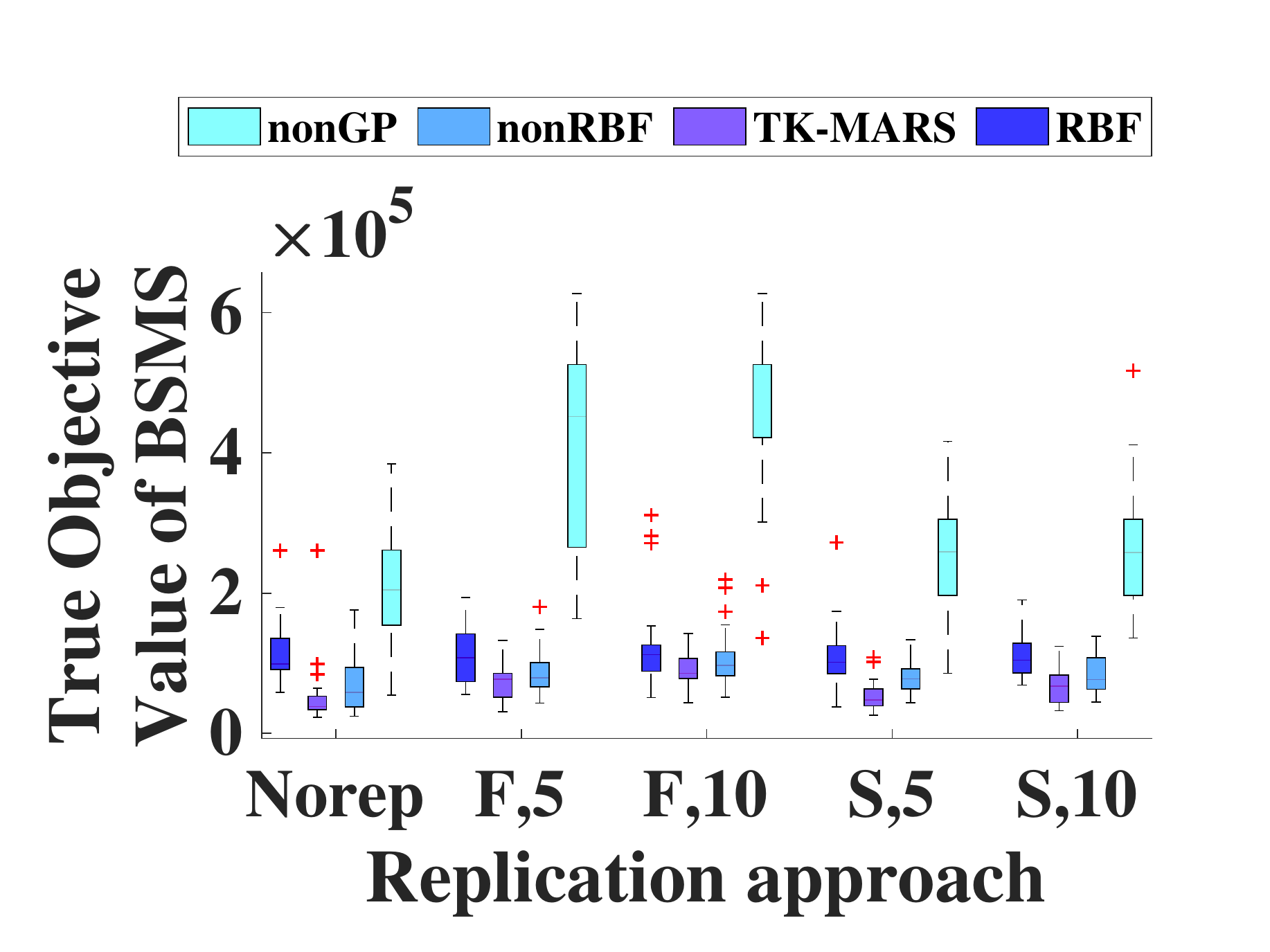}}
    \end{minipage}
    \begin{minipage}{\linewidth}
        \subfloat[Noise=0.1]{\includegraphics[width=0.55\textwidth]{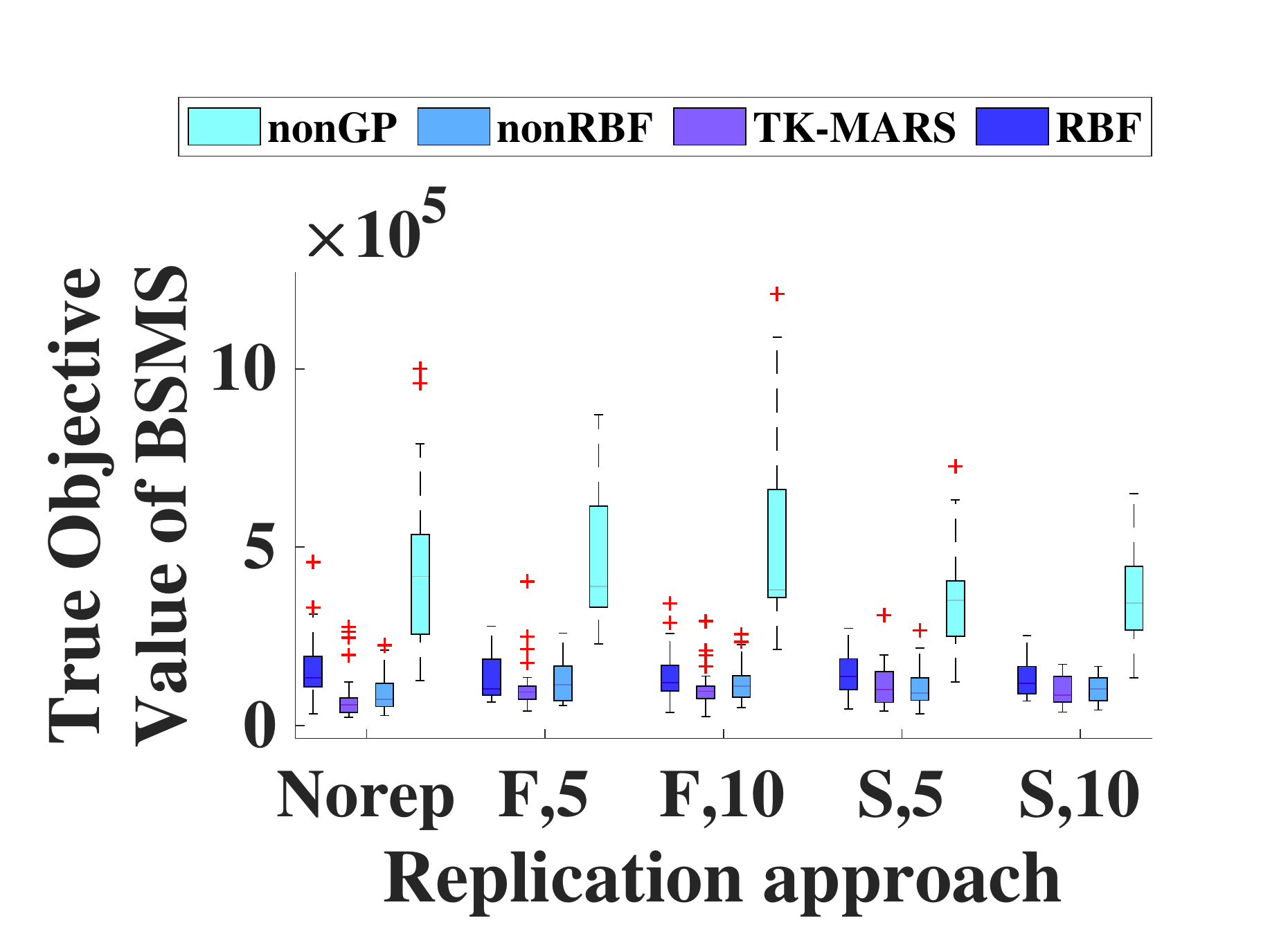}}
        \subfloat[Noise=0.25]{\includegraphics[width=0.55\textwidth]{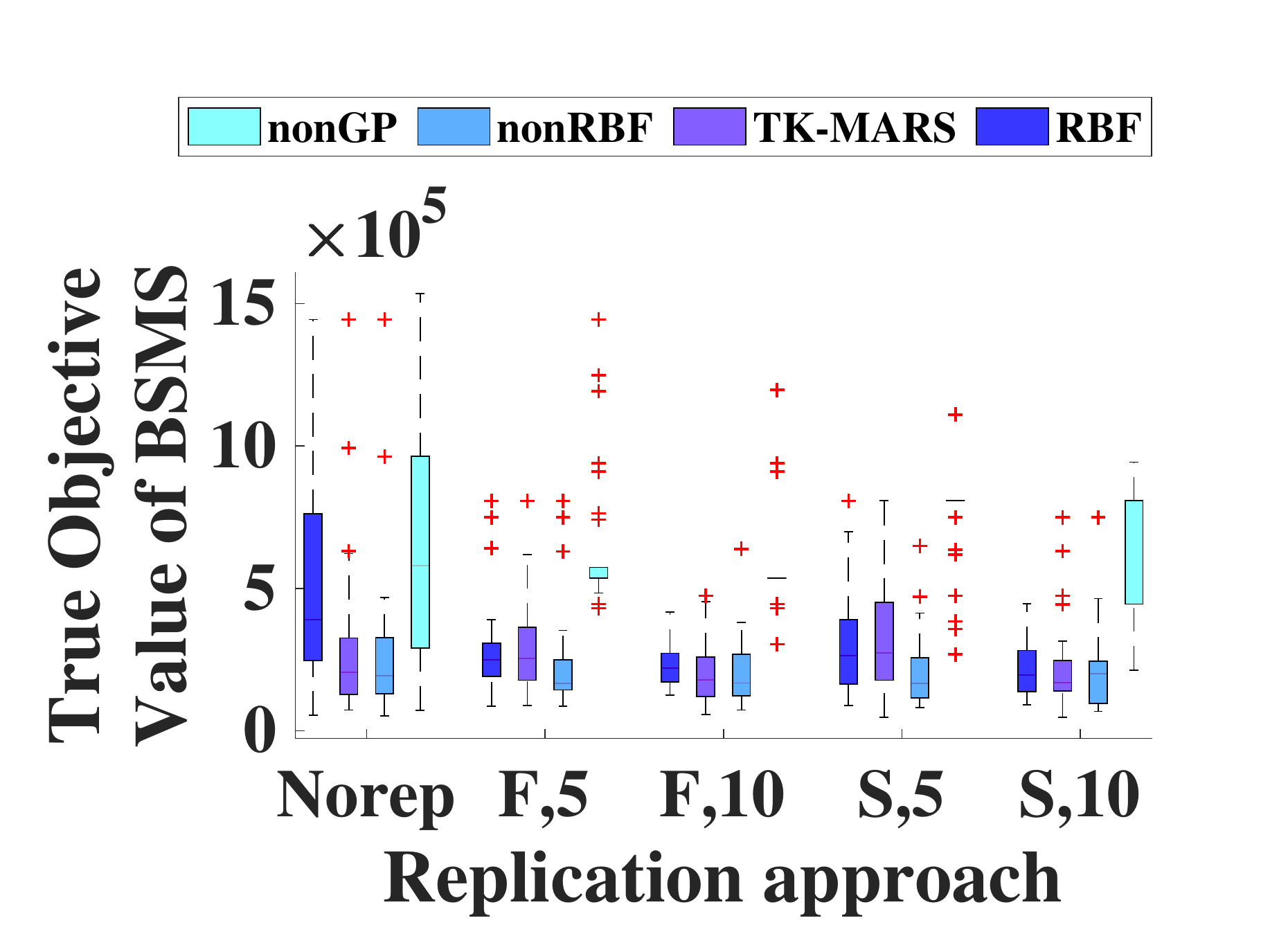}}
    \end{minipage}
    \caption{Box-plots of the true objective value of the BSMS after 1000 black-box function evaluations of surrogate optimization on the Rosenbrock function}
    \label{fig:bks_box_rosen}
    \vspace{-3mm}
\end{figure}

% \begin{figure}[!tb]
% \centering
%     \begin{minipage}{\linewidth}
%         \subfloat[]{\includegraphics[width=0.35\textwidth]{figures/box_bks_rosen_new_1.pdf}}
%         \subfloat[]{\includegraphics[width=0.35\textwidth]{figures/box_bks_rosen_new_2.pdf}}
%         \subfloat[]{\includegraphics[width=0.35\textwidth]{figures/box_bks_rosen_new_3.pdf}}
%     \end{minipage}
%     \begin{minipage}{\linewidth}
%         \subfloat[]{\includegraphics[width=0.35\textwidth]{figures/box_bks_rosen_new_4.pdf}}
%         \subfloat[]{\includegraphics[width=0.35\textwidth]{figures/box_bks_rosen_new_5.pdf}}
%     \end{minipage}
%     \caption{Box-plots of the true objective value of the BSMS after 1000 black-box function evaluations of surrogate optimization on the Rosenbrock function}
%     \label{fig:bks_box_rosen}
%     \vspace{-3mm}
% \end{figure}
\hadis{For a deeper evaluation of the performance of No-Replication and replication strategies using different surrogates, the robustness and the quality of BSMS must be assessed.} \hadis{
Figure~\ref{fig:auc_box_rosen}, show box-plots of the MTFAUC for 30 executions at different noise levels using different approaches on the Rosenbrock test function.}
%Although from Figure~\ref{fig:tkmars_mtfauc}(a) we observed that No-Replication outperforms Replication approaches, 
%Looking at the plots,
%Figure~\ref{fig:auc_box_rosen}(c) ($\mathit{Fixedrep},10$) and Figure~\ref{fig:auc_box_rosen}(e) ($\mathit{Smartrep},10$), 
We observe that the \emph{interquartile range} (the difference between the $25^{\mathit{th}}$ and $75^{\mathit{th}}$ percentiles, i.e., the height of the boxes) for Smart-Replication is smaller especially when the noise level is high. This suggests that Smart-Replication is more \emph{robust}. 
%, Figure~\ref{fig:auc_box_rosen}(a). As in Figures~\ref{fig:auc_box_rosen}(d) and (e), in lower levels of noise, Smart-Replication is competitive with No-Replication approach. This indicates the robustness of the Smart-Replication approach to randomness. 
\hadis{Although $\mathit{Fixedrep},10$ has also short interquartile range in higher levels of noise in some cases
%is robust to different noise levels, and there is no significant difference in the means across the noise levels
, the overall average MTFAUC is larger in lower noise levels than that of Smart-Replication. This is because the Fixed-Replication strategy makes unnecessary functional evaluations. $\mathit{Fixedrep},10$ outperforms $\mathit{Smartrep}$ when we use nonGP as a surrogate. This suggests that nonGP needs more replications to handle noise, and we perhaps need to increase $r_{max}$ of $\mathit{Smartrep}$ to observe its benefit. %compared with $\mathit{Smartrep}$. 
Note that as the noise level increases, $\mathit{Smartrep},10$ is more competitive with $\mathit{Norep}$ especially when we use the interpolating RBF model.}
%Besides, from the box-plots, note that as the level of noise increases, the variance of the No-Replication's $MTFAUC$ based on different pools increases.

Figure~\ref{fig:bks_box_rosen} presents the box-plots of the true objective value of BSMS after 1000 black-box function evaluations for the same 30 executions on the Rosenbrock test function as those in Figure \ref{fig:auc_box_rosen}. We observe that $\mathit{Smartrep}$ is competitive with $\mathit{Norep}$ in the lower levels of noise after 1000 black-box function evaluations. 
%$\mathit{Fixedrep},10$ and $\mathit{Smartrep},10$ are robust to different noise levels since the means are close compared to $\mathit{Norep}$. 
$\mathit{Smartrep},10$ has the smallest interquartile range with the highest level of noise, indicating its robustness. The box-plots for the other test functions are provided in Appendix~\ref{app:results}.

\hadis{It should be noted that while $Fixedrep$ and $Norep$ may have a competitive performance or outperformance compared to $Smartrep$ in some situations, the level of uncertainty of the black-box system is unknown apriori, $Smartrep$ may adapt to the unknown noise level to have the same performance of $Norep$ in the lower noise level and $Fixedrep$ in the higher noise level.}

\section{Conclusions}
%and Future Work}
\label{sec:conclusion}

In this paper, we developed a partitioning-based MARS model called TK-MARS, which is specifically designed for surrogate optimization of black-box functions. Moreover, we designed a smart replication approach to mitigate the impact of uncertainty associated with the black-box output. 
We demonstrated the performance of the proposed approaches using complex global optimization test functions. The performance of a surrogate optimization algorithm when using TK-MARS was compared with that when using original MARS. We observed that TK-MARS outperforms original MARS in the context of surrogate optimization. We also showed that TK-MARS is capable of correctly detecting the important variables.

No-Replication outperforms other approaches overall in a deterministic environment, but in a noisy environment, Smart-Replication yields more robust results. 
%However, Smart-Replication behaves similar to No-Replication when the uncertainty level is low and similar to Fixed-Replication when the uncertainty level is high.
The robustness and the quality of the final optimum solution found through Smart-Replication is competitive with that using no replications in environments with low levels of noise and using a fixed number of replications in highly noisy environments.
Moreover, a surrogate optimization algorithm with an interpolating surrogate, such as RBF, needs more replication than the one with a non-interpolating surrogate, such as TK-MARS.

\hadis{The results indicate that surrogate optimization using TK-MARS outperforms that of interpolating and non-interpolating RBF and non-interpolating GP in most of the cases, which are the prevalent surrogate models used in the literature. It is worth mentioning that non-interpolating RBF model proposed by Jakobsson in ~\cite{jakobsson2010method} is competitive with TK-MARS in a few instances. Non-interpolating surrogate models still need replications in the presence of noise when the underlying function is complex.}

% Overall, the experiments suggest a surrogate optimization algorithm using Smart-Replication and TK-MARS is the most robust approach for high-dimensional black-box optimization in both deterministic and noisy environments.

\hadis{
We consider using real datasets for evaluating our approaches as part of our future work. Our future work also includes investigating other feature selection techniques in the surrogate optimization framework to identify unimportant variables of a black-box function.
% In our future work, we aim to apply our proposed approaches using real datasets. We also would like to investigate other feature selection techniques in the surrogate optimization framework to identify unimportant variables of a black-box function.
}

\section{Acknowledgement}
The work was supported by the National Science Foundation Award CMMI–1434401.
% Non-BibTeX users please use
%\begin{thebibliography}{}
%

\bibliographystyle{ieeetr}
%\section*{References}
%\bibliographystyle{elsarticle-num}

\bibliography{mybibfile}
% and use \bibitem to create references. Consult the Instructions
% for authors for reference list style.

% \bibitem{RefJ}
% Format for Journal Reference
% Author, Article title, Journal, Volume, page numbers (year)
% Format for books
% \bibitem{RefB}
% Author, Book title, page numbers. Publisher, place (year)
% etc
%\end{thebibliography}
% \vspace{0.5cm}\noindent\LARGE{Appendix}
% \appendix

%\section{Rastrigin Box-plots}
%\section{Levy Box-plots}
% \appendix
% %\pagebreak
% {\bf \LARGE APPENDIX}

\appendix
%\appendixpage
\addappheadtotoc

{\bf \LARGE APPENDIX}
\section{Exploration-Exploitation Pareto Approach (EEPA)}\label{app:EEPA}
Dickson~\cite{dickson2014exploration} proposed an exploration-exploitation Pareto approach for sampling in surrogate optimization. Algorithm~\ref{alg:eepa} presents the EEPA pseudocode. EEPA starts with a set of input points $R$ representing the solution space. For each input point in $R$, it calculates the minimum distance, $\delta(x)$, from the already evaluated set of input points $I$. 
% For two input points $x$ and $\tilde{x}$, the distance metric $ \varDelta (x,\tilde{x})$ is defined as a euclidean distance given by 
% $$    \varDelta (x,\tilde{x}) ={\|x - \tilde{x}\|}$$
% or cosine distance given by
% $$    \varDelta (x,\tilde{x}) = 1 - \frac{x}{\|x\|}.\frac{\tilde{x}}{\|\tilde{x}\|}.$$
Next, it determines a non-dominated set of input points $F$ according to $\delta(x)$ and the predicted output values $\hat{f}(x)$ of the input points in $R$. To eliminate the close input points in $F$, EEPA applies a \emph{maximin} exploration measure to $F$ in order to select the most diverse candidates. The sampled input points are stored in $P$. EEPA limits the size of $P$ to a maximum $K'$ in each iteration due to the expensive black-box evaluation.

\begin{algorithm}[h]
\hadis{
\caption{{\bf EEPA,~\cite{dickson2014exploration}}\\
{\bf input}: $I$, $\hat{f}$, $R$ }
\begin{algorithmic}[1]
\label{alg:eepa}
%\STATE $I$ is the initial set, $R$ is the random data set
%\STATE $R=$ a fixed number of random data points from $D$
% \STATE For each $x\in I$, determine $f(x)$
%\WHILE {Termination criteria is not satisfied}
% * <jaymrosenberger@gmail.com> 2018-07-23T22:57:26.592Z:
% 
% > K
% Didn't you use N_Max of something for the iteration limit earlier?
% hadis: I removed that notation an replaced it with stopping criterion
    %\FOR{$i=1$ to $|R|$}
	\STATE $\delta(x)=\min_{\tilde{x}\in I} \lVert x-\tilde{x} \rVert, \forall x \in R$
    \STATE $F = \{x \in R~|~\nexists~\tilde{x} \in R, \hat{f}(\tilde{x}) \leq \hat{f}(x), \delta(\tilde{x})\geq \delta(x)\}$ 
	\STATE $P=\{\}; k^\prime= 1;$
    \STATE $x^\prime \in \arg\min \{\hat{f}(x)|x\in F\}$
 \STATE  $P = P \cup {x^\prime}$
    \WHILE {$k^\prime\leq K^\prime~\mathbf{and}~P$}
    \STATE $\delta^\prime (x)= \min_{\tilde{x}\in I \cup P}\lVert x-\tilde{x} \rVert, \forall x \in F$
    \STATE $x \in \arg\max \{\delta^\prime(x)|x\in P\}$
    \STATE $P =P\cup \{x\}$
    \STATE $k^\prime=k^\prime+1$
%    \ENDWHILE
\ENDWHILE
\STATE {\bf return} $P$
\end{algorithmic}}
\end{algorithm}

\section{Additional Experimental Results}\label{app:results}
This section provides experimental results in addition to those in Section \ref{sec:results}. 

\new{Table~\ref{tab:results-RBF} shows the preliminary results of using RBF with different kernels in surrogate optimization for a different fraction of important variables. The values in the table denote the MTFAUC measure and the total number of function evaluations before the optimum solution is found, respectively.}

\new{Table~\ref{tab:wrep-anova-oa} shows the ANOVA table of our preliminary analysis based on an orthogonal array design considering different parameters, as listed in Table~\ref{tab:params-OA}. We fixed the unimportant variables to their most significant level to study the impact of noise on different approaches, in a full factorial setting. }

\new{Tables~\ref{tab:mtfauc1}, \ref{tab:mtfauc2}, and \ref{tab:mtfauc3} show the average and variance of the performance of No-Replication, Fixed- Replication, and Smart-Replication for 30 different executions on the test functions listed in Table~\ref{tab:testfuncchar}. The highlighted cells show the best method for each test function under different uncertainty levels.}

\begin{figure}[!tb]
\centering
    \begin{minipage}{\linewidth}
        \subfloat[Noise=0]{\includegraphics[width=0.55\textwidth]{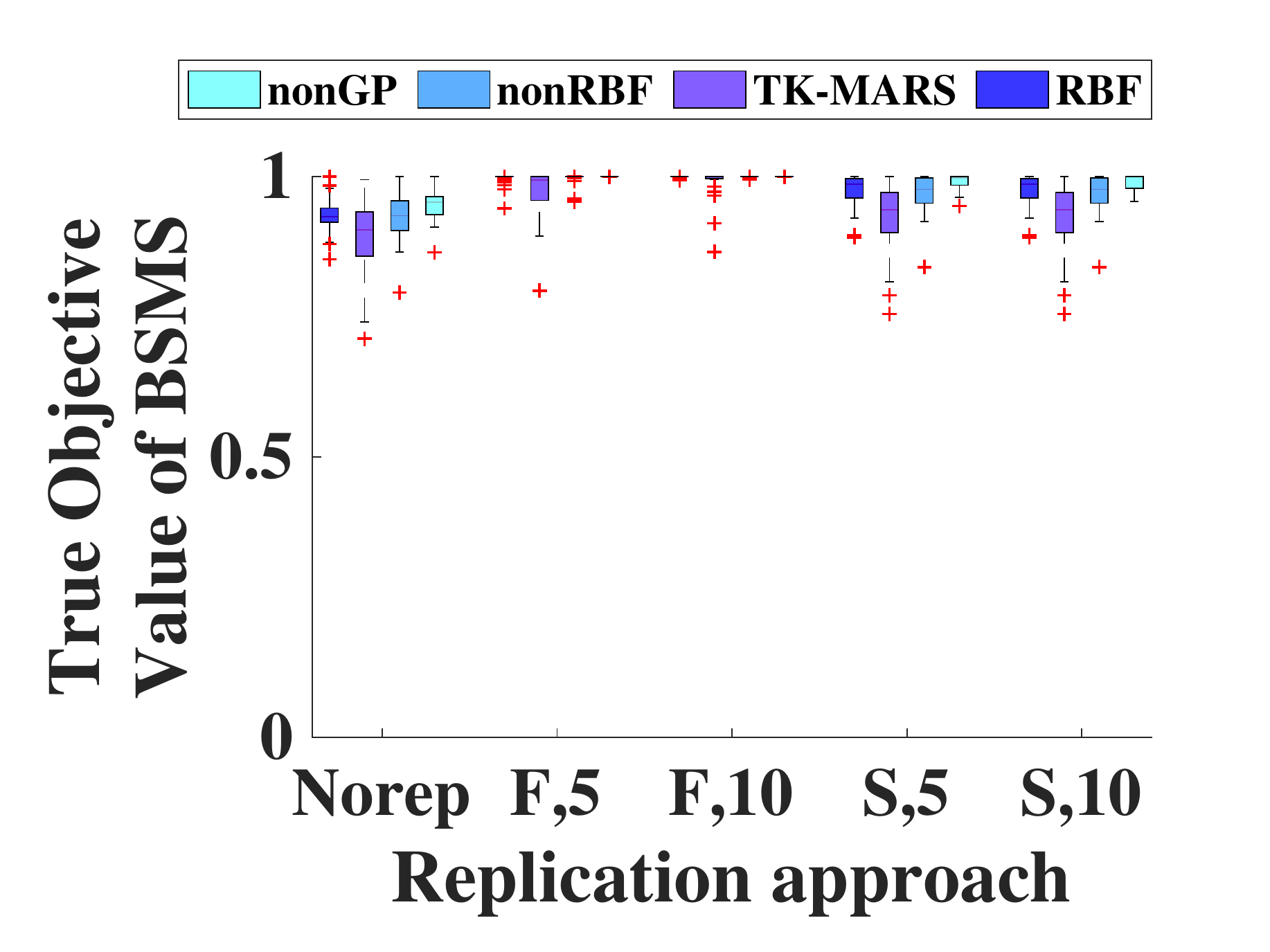}}
        \subfloat[Noise=0.05]{\includegraphics[width=0.55\textwidth]{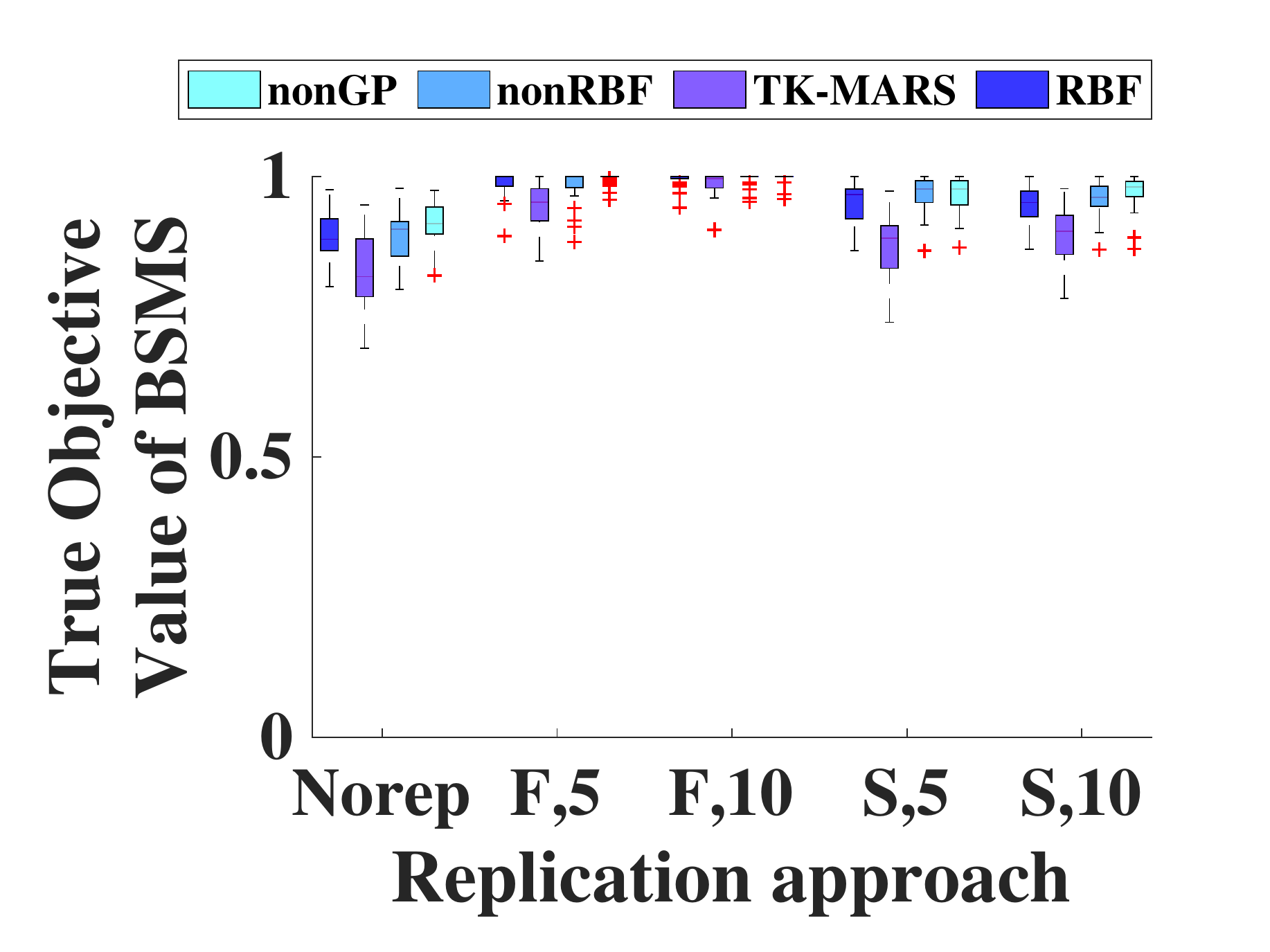}}
    \end{minipage}
    \begin{minipage}{\linewidth}
        \subfloat[Noise=0.1]{\includegraphics[width=0.55\textwidth]{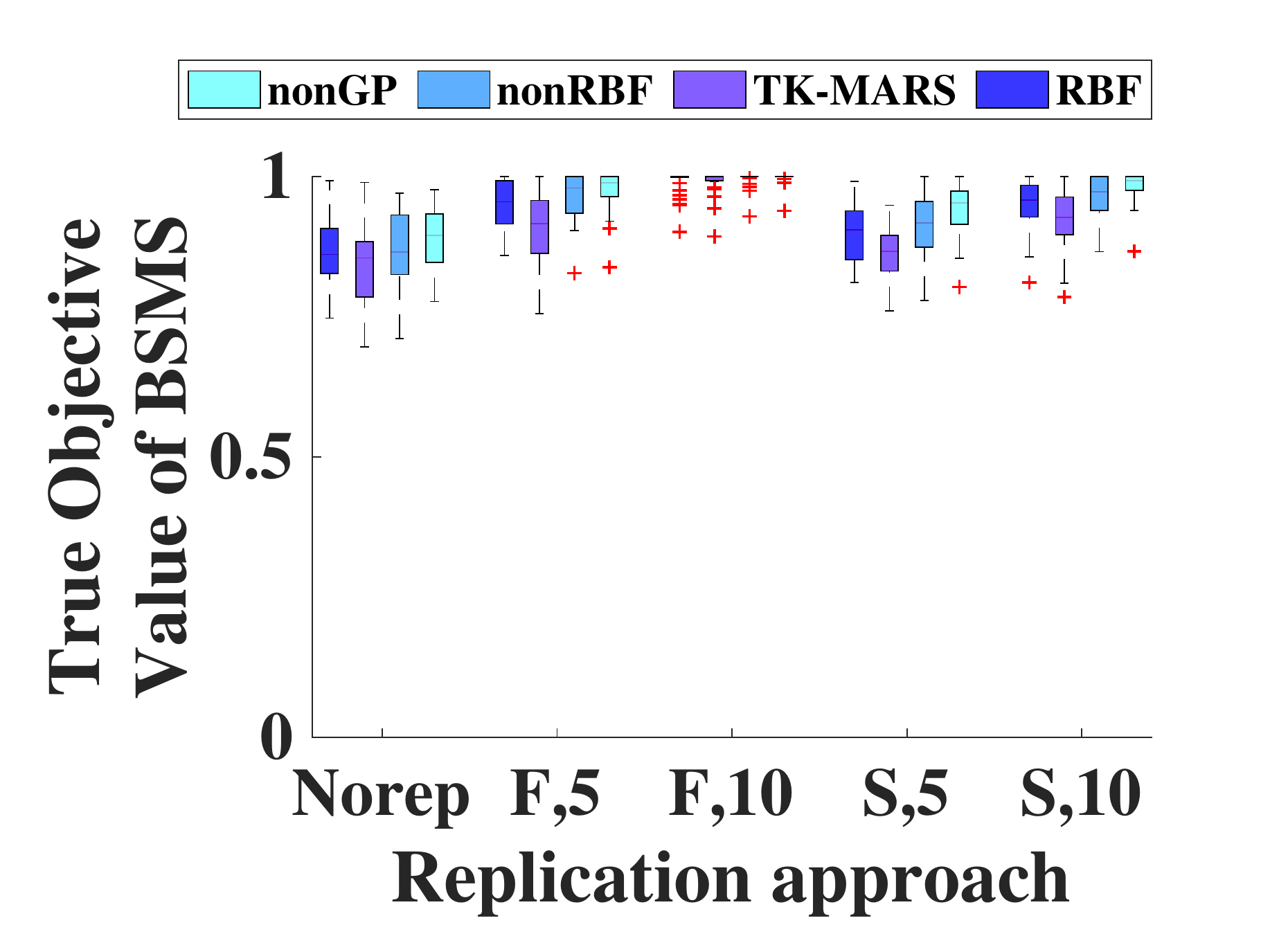}}
        \subfloat[Noise=0.25]{\includegraphics[width=0.55\textwidth]{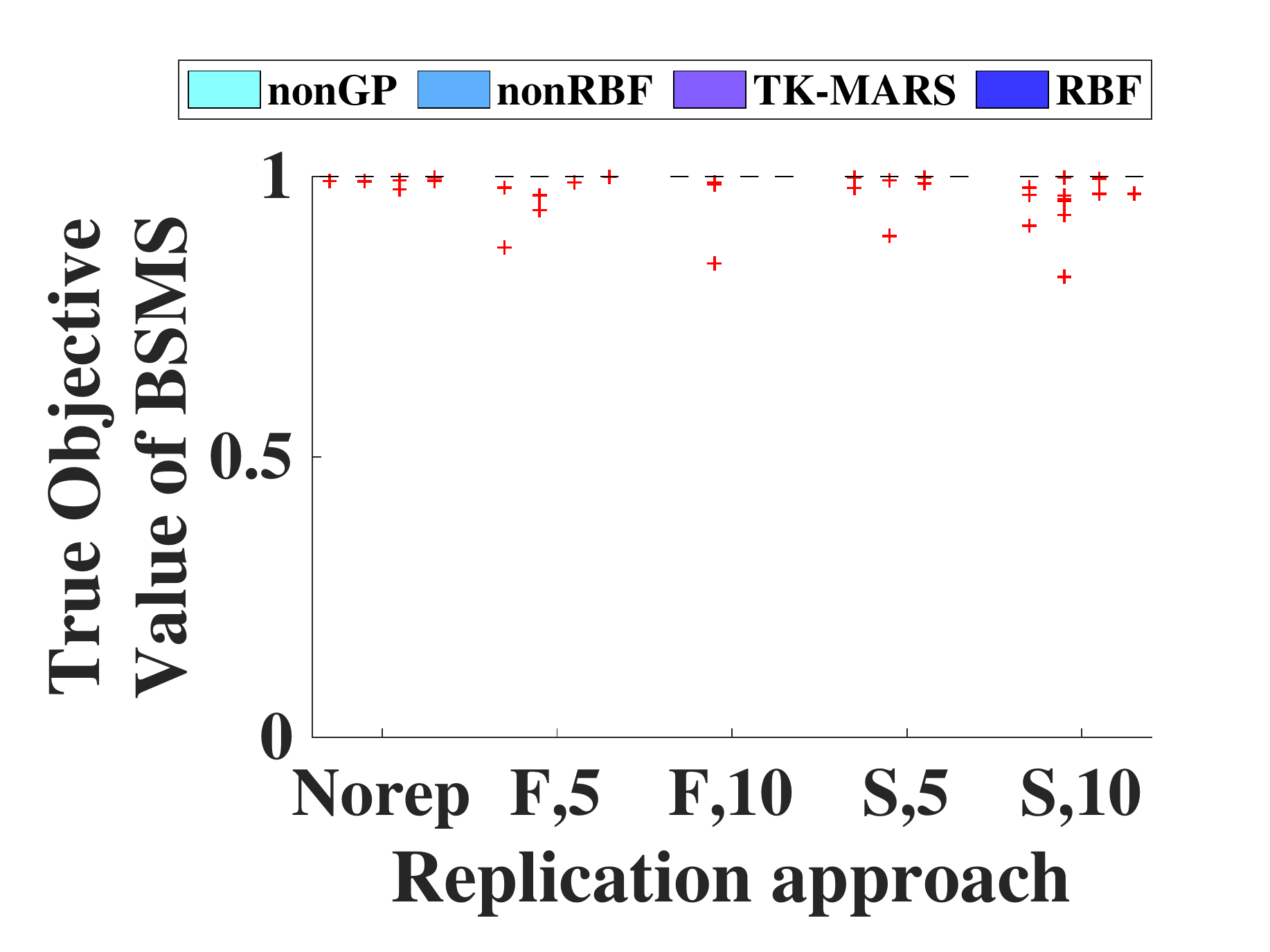}}
    \end{minipage}
    \caption{Box-plots of MTFAUC of surrogate optimization on the Rastrigin function}
    \label{fig:auc_box_rast}
    \vspace{-3mm}
\end{figure}

\hadis{Figure~\ref{fig:auc_box_rast} presents the box-plots of the MTFAUC for 30 executions using different approaches on the Rastrigin test function.
We can observe that $\mathit{Fixedrep},10$ has the smallest interquartile range and is the most robust approach to uncertainty.
The reason is likely due to the highly fluctuating behavior of the Rastrigin function with several local optima combined with the noise added to the function, which makes an optimum hard to obtain.}

\begin{figure}[!tb]
\centering
    \begin{minipage}{\linewidth}
        \subfloat[Noise=0]{\includegraphics[width=0.55\textwidth]{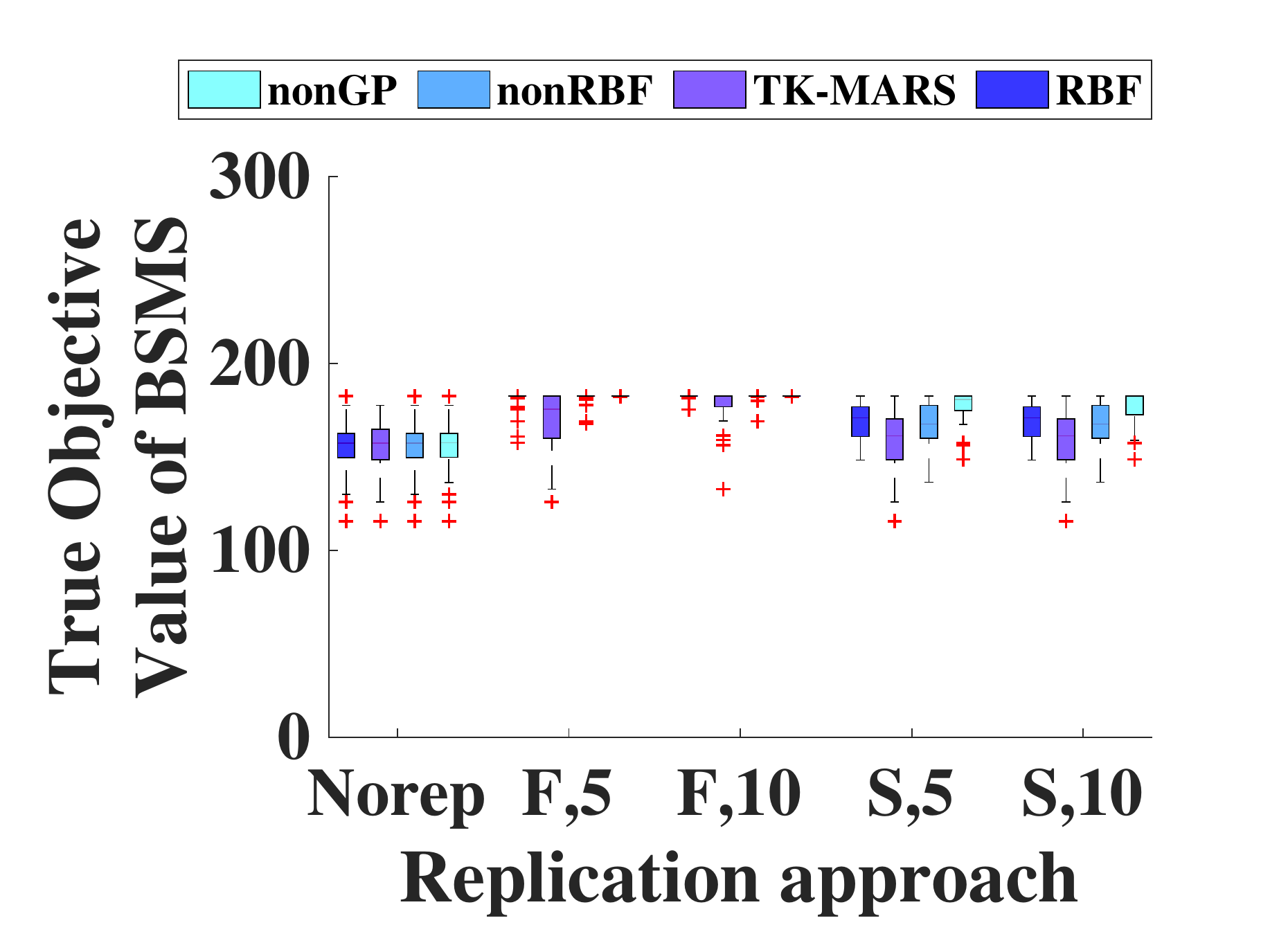}}
        \subfloat[Noise=0.05]{\includegraphics[width=0.55\textwidth]{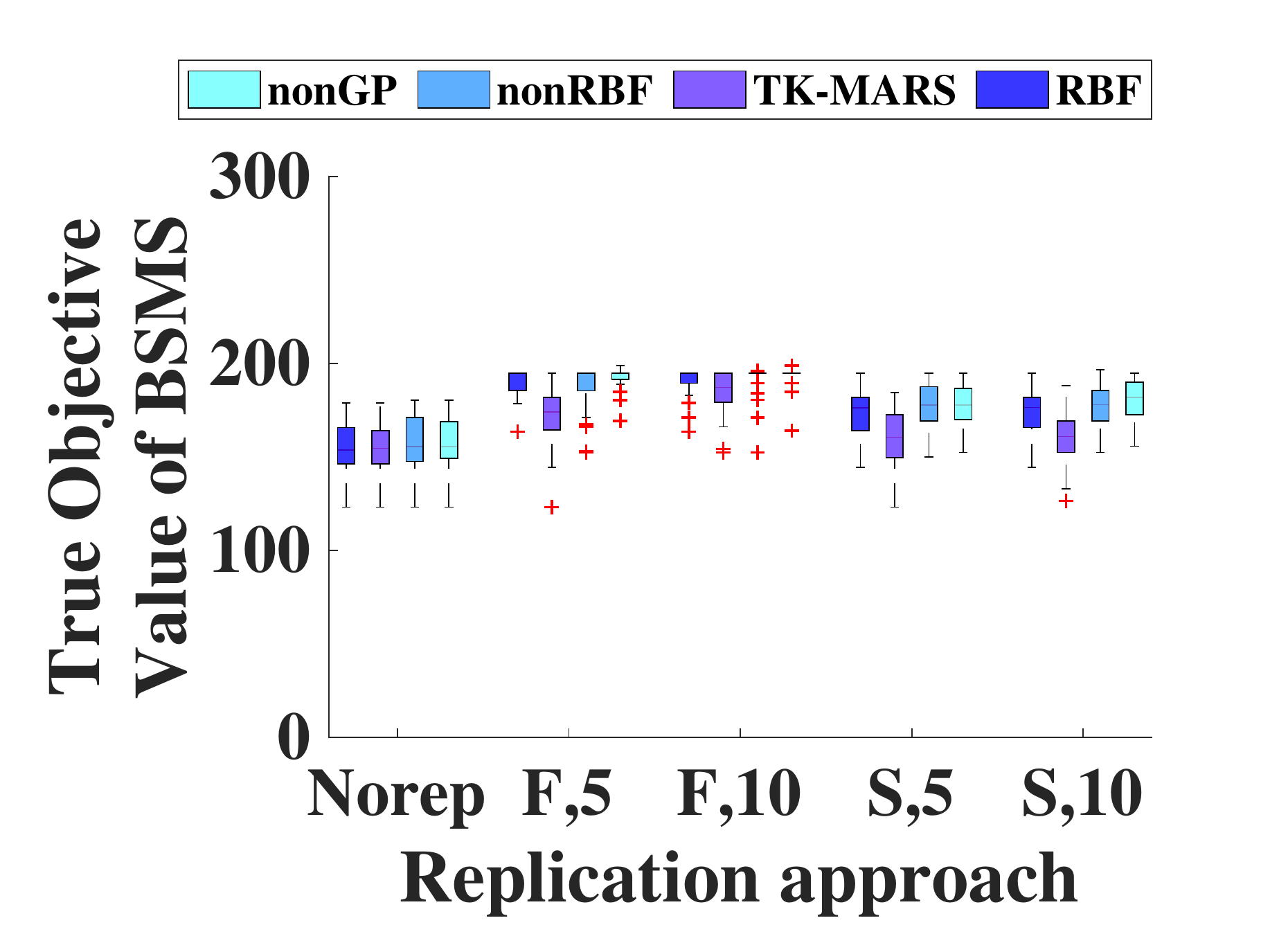}}
    \end{minipage}
    \begin{minipage}{\linewidth}
        \subfloat[Noise=0.1]{\includegraphics[width=0.55\textwidth]{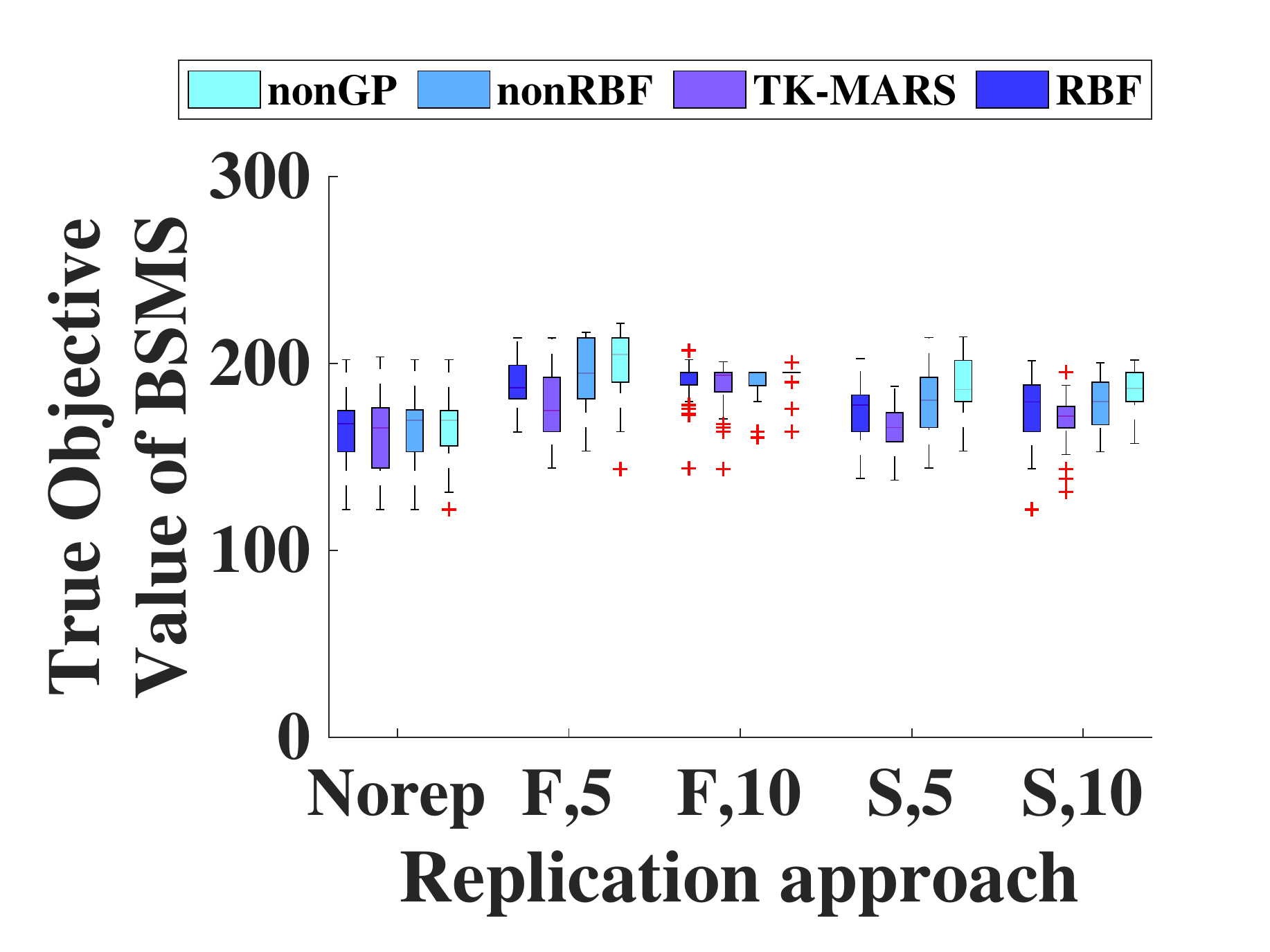}}
        \subfloat[Noise=0.25]{\includegraphics[width=0.55\textwidth]{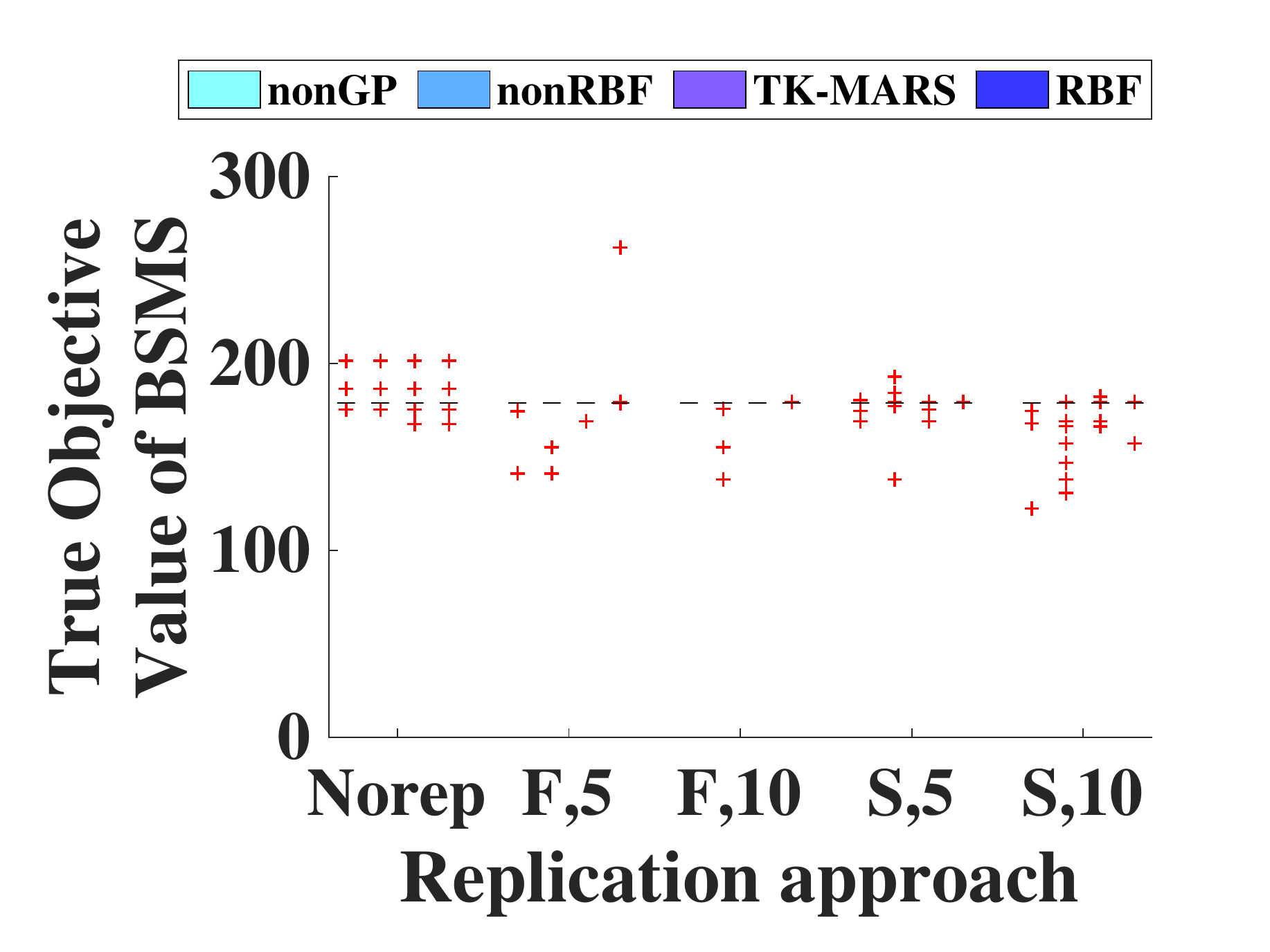}}
    \end{minipage}
    \caption{Box-plots of the true objective value of the BSMS after 1000 black-box function evaluations of surrogate optimization on the Rastrigin function}
    \label{fig:bks_box_rast}
    \vspace{-3mm}
\end{figure}

\hadis{Figure~\ref{fig:bks_box_rast} shows the box-plots of the true objective value of BSMS after 1000 black-box function evaluations for the same 30 executions on the Rastrigin test function. We observe the same pattern for the quality of the final BSMS after 1000 evaluations. $\mathit{Smartrep},10$ finds a slightly better BSMS across the noise levels with relative robustness. 
Note that $\mathit{Smartrep},10$ performance is independent of the surrogate model used. For example, we can observe that $\mathit{Fixedrep},10$ has higher BSMS value for nonRBF and nonGP at higher levels of noise, which is not the case for $\mathit{Smartrep},10$.}

\begin{figure}[!tb]
\centering
    \begin{minipage}{\linewidth}
        \subfloat[Noise=0]{\includegraphics[width=0.55\textwidth]{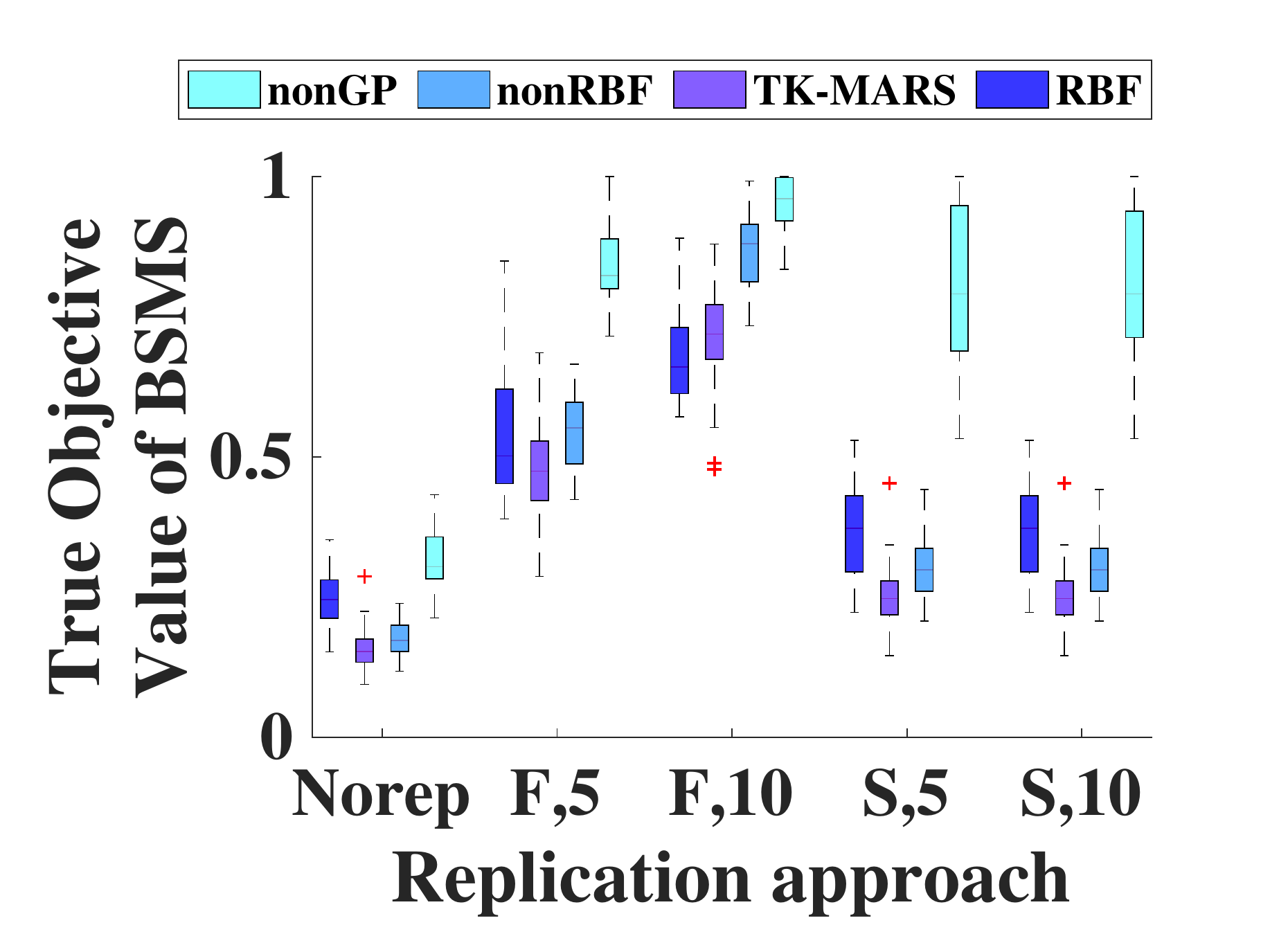}}
        \subfloat[Noise=0.05]{\includegraphics[width=0.55\textwidth]{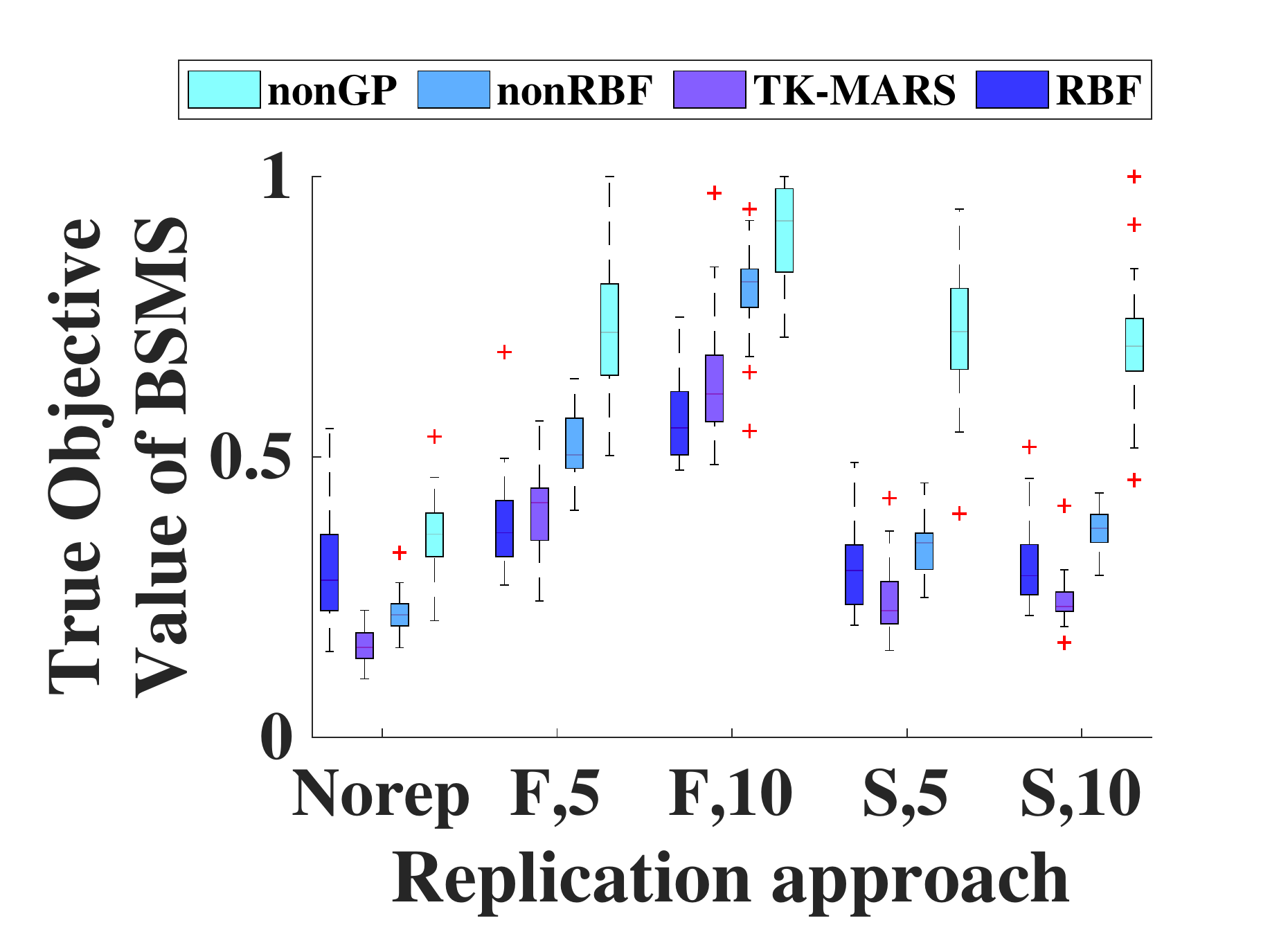}}
    \end{minipage}
    \begin{minipage}{\linewidth}
        \subfloat[Noise=0.1]{\includegraphics[width=0.55\textwidth]{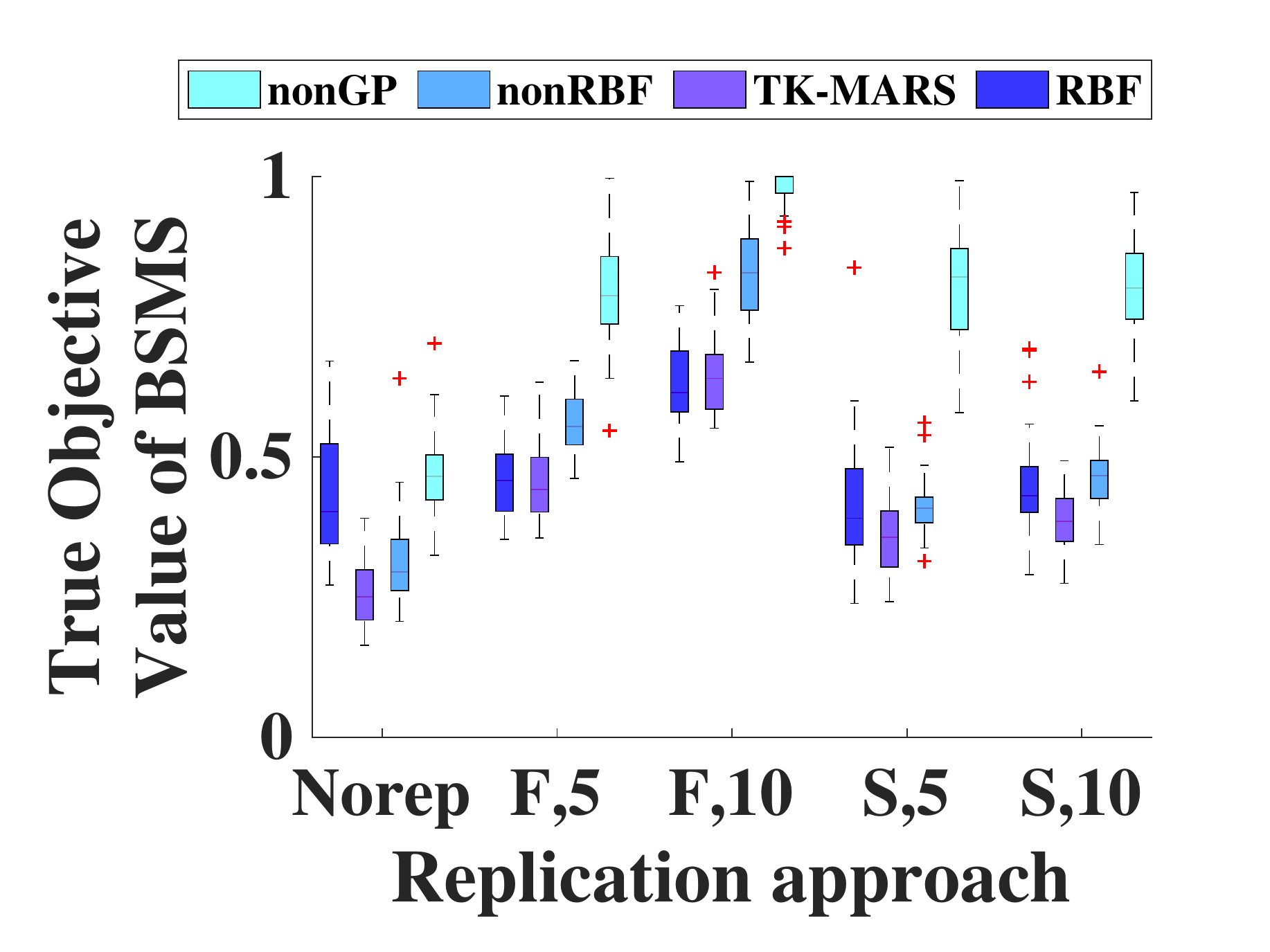}}
        \subfloat[Noise=0.25]{\includegraphics[width=0.55\textwidth]{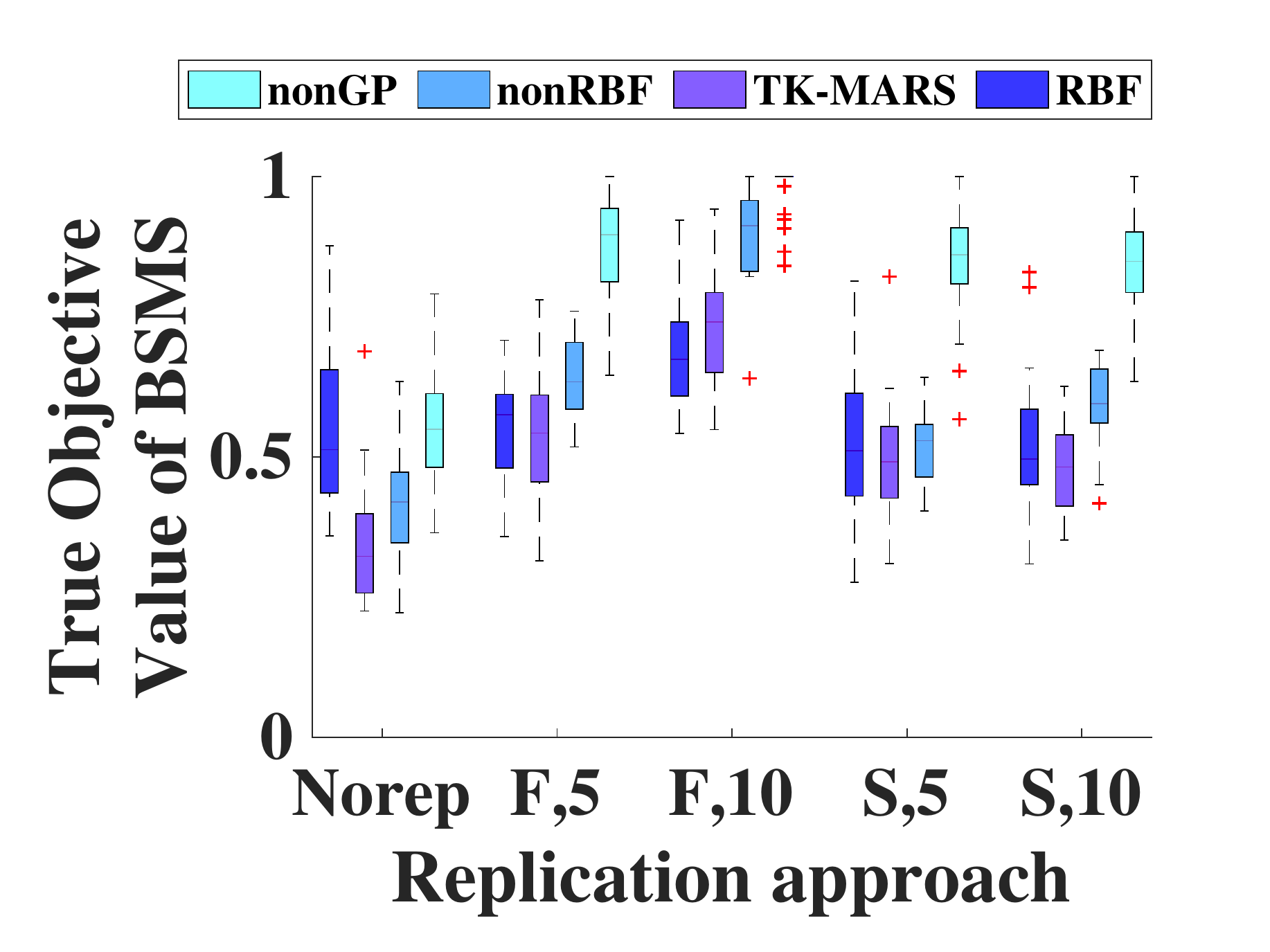}}
    \end{minipage}
    \caption{Box-plots of MTFAUC of surrogate optimization on the Levy function}
    \label{fig:auc_box_levy}
    \vspace{-3mm}
\end{figure}

\hadis{Figure~\ref{fig:auc_box_levy} displays the box-plots of the MTFAUC for 30 executions using different approaches on the Levy test function. We note that the interquartile range for $\mathit{Smartrep},10$ is comparatively small across different levels of noise. Note that in this instance, nonGP needs full replication, and the performance of $\mathit{Fixedrep},10$ outperforms $\mathit{Smartrep},10$ when we use nonGP as a surrogate.}

\begin{figure}[!tb]
\centering
    \begin{minipage}{\linewidth}
        \subfloat[Noise=0]{\includegraphics[width=0.55\textwidth]{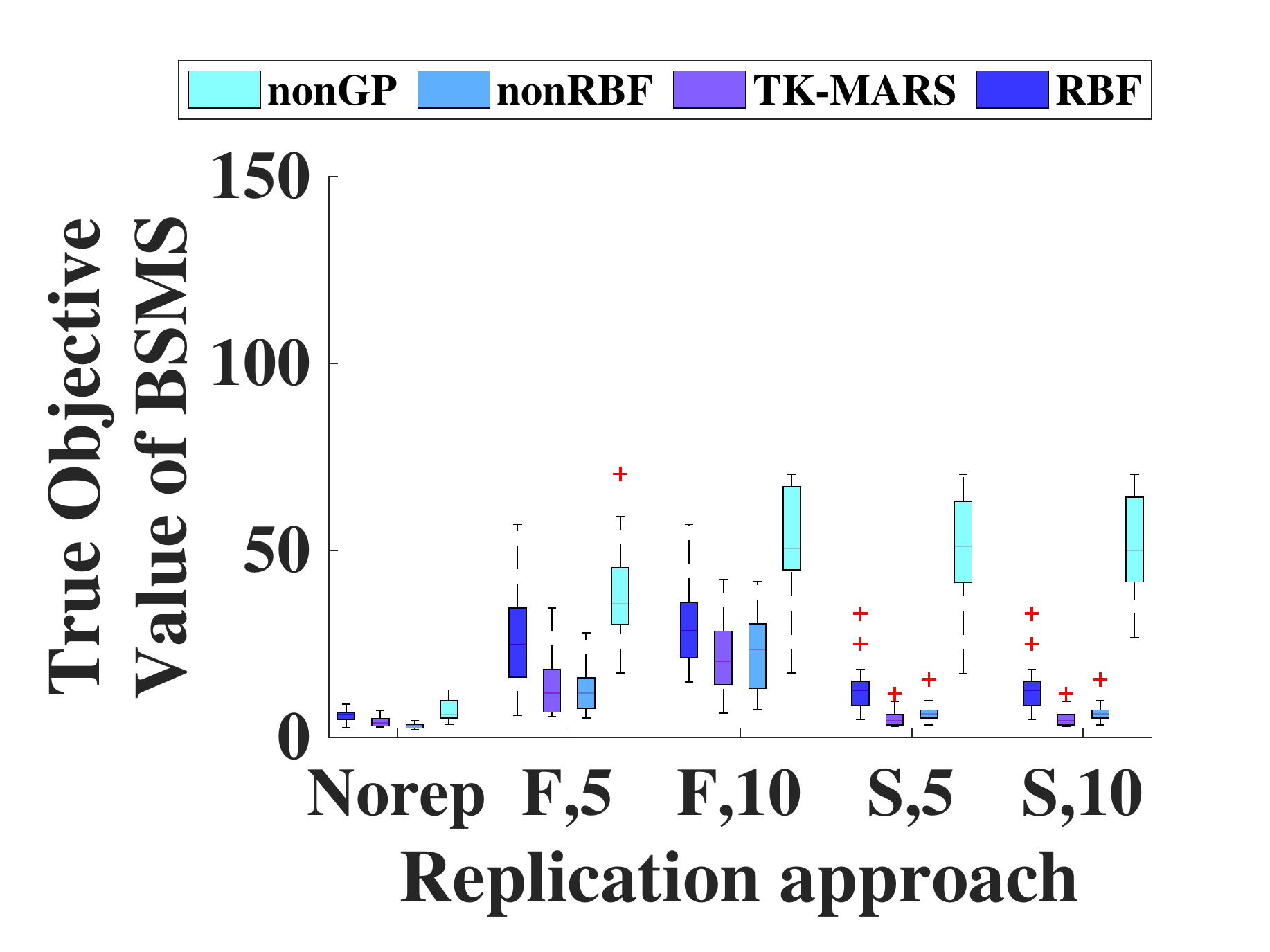}}
        \subfloat[Noise=0.05]{\includegraphics[width=0.55\textwidth]{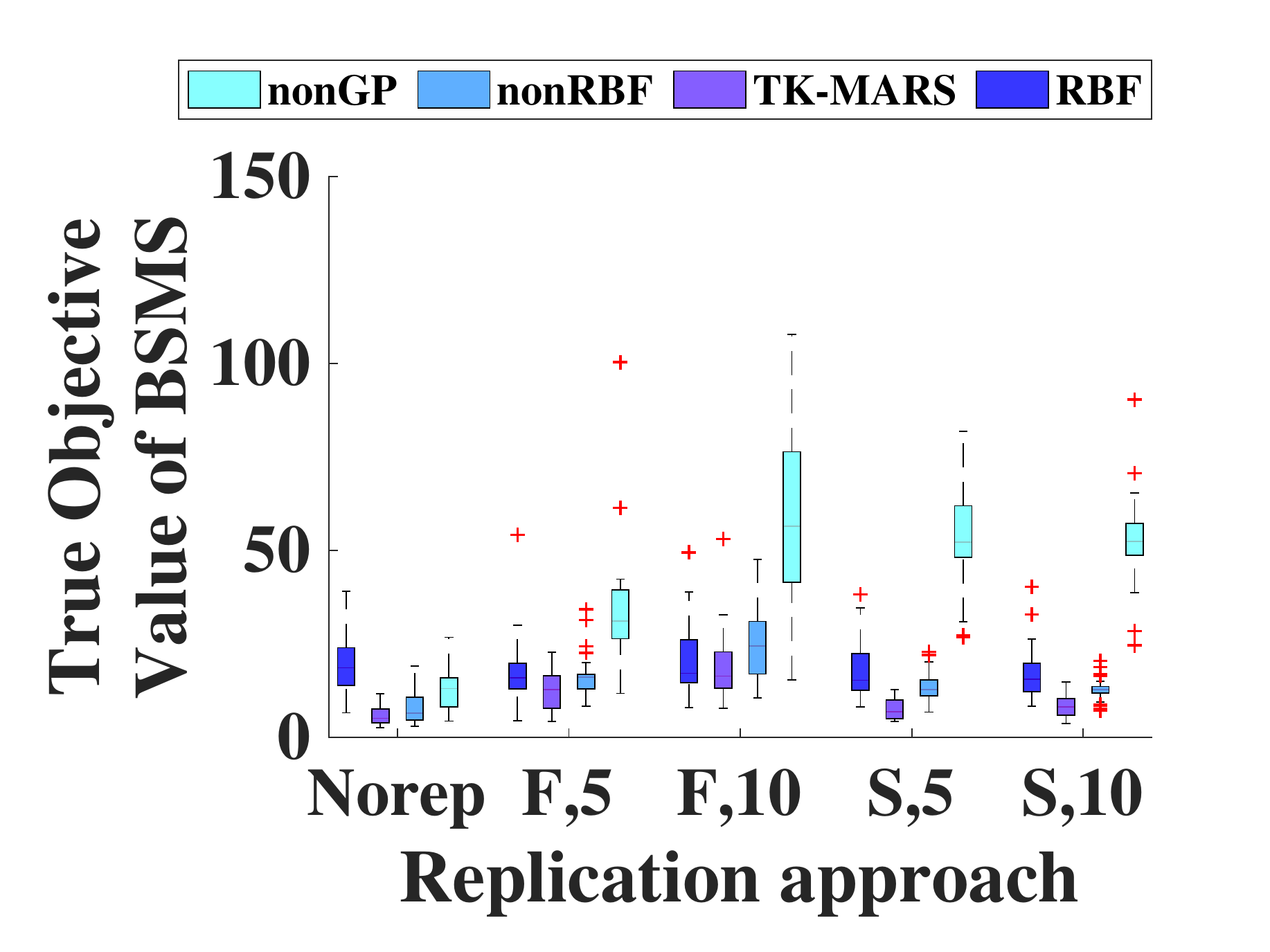}}
    \end{minipage}
    \begin{minipage}{\linewidth}
        \subfloat[Noise=0.1]{\includegraphics[width=0.55\textwidth]{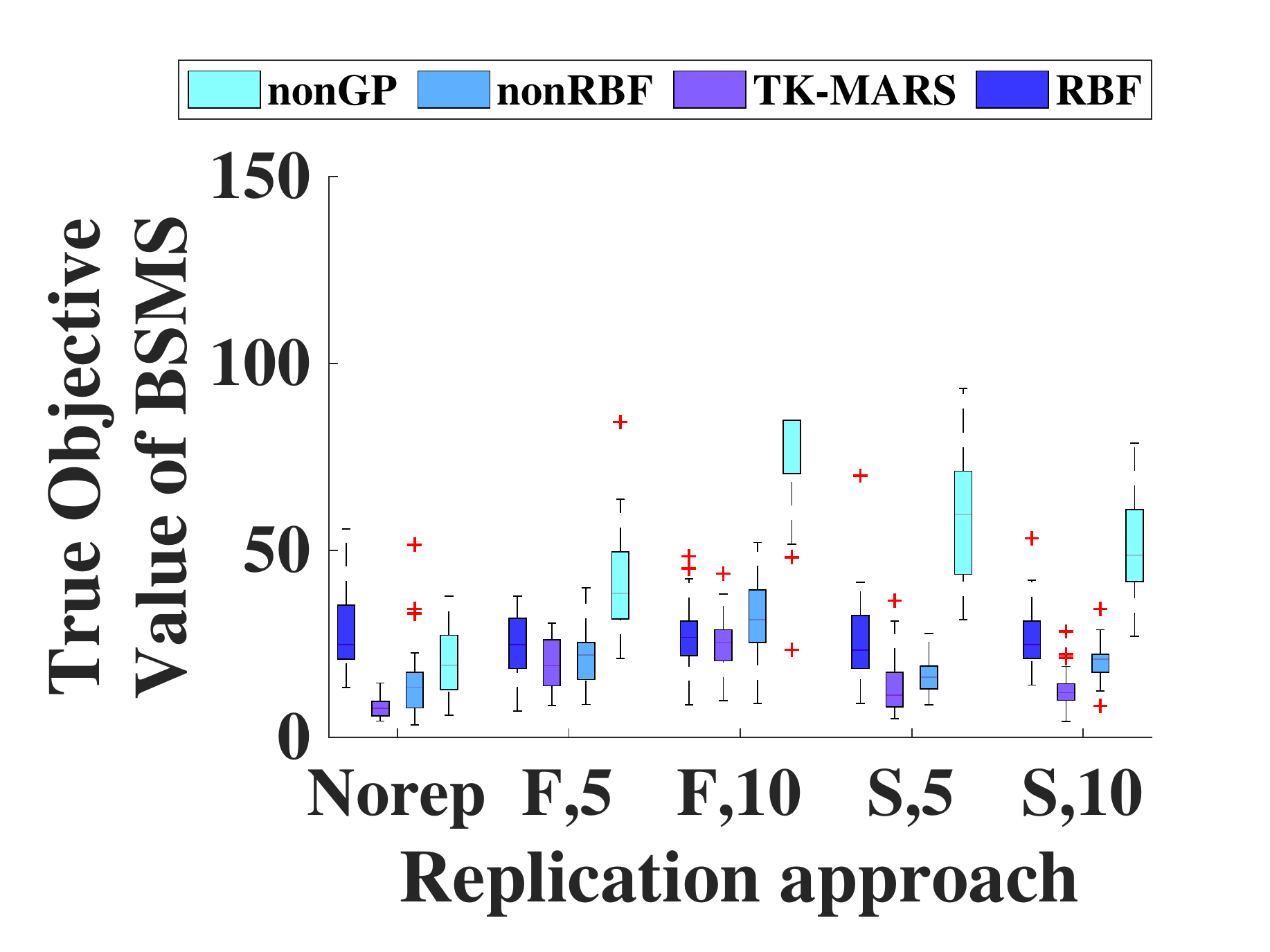}}
        \subfloat[Noise=0.25]{\includegraphics[width=0.55\textwidth]{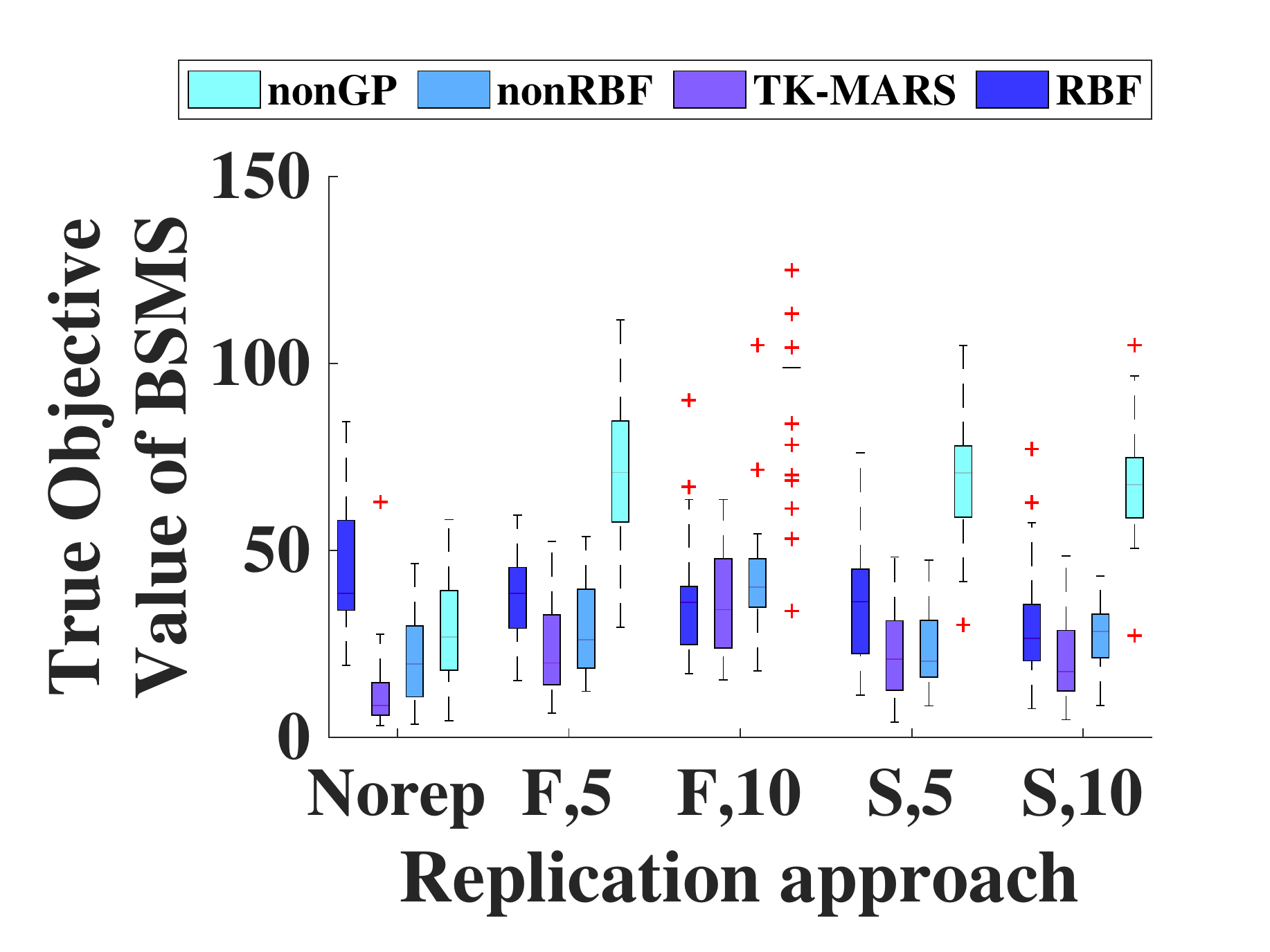}}
    \end{minipage}
    \caption{BBox-plots of the true objective value of the BSMS after 1000 black-box function evaluations of surrogate optimization on the Levy function}
    \label{fig:bks_box_levy}
    \vspace{-3mm}
\end{figure}

\hadis{Figure~\ref{fig:bks_box_levy} shows the box-plots of the true objective value of BSMS after 1000 black-box function evaluations for the same 30 executions on the Levy test function. $\mathit{Smartrep},10$ is the most robust approach to the different noise levels.}

\begin{figure}[!tb]
\centering
    \begin{minipage}{\linewidth}
        \subfloat[Noise=0]{\includegraphics[width=0.55\textwidth]{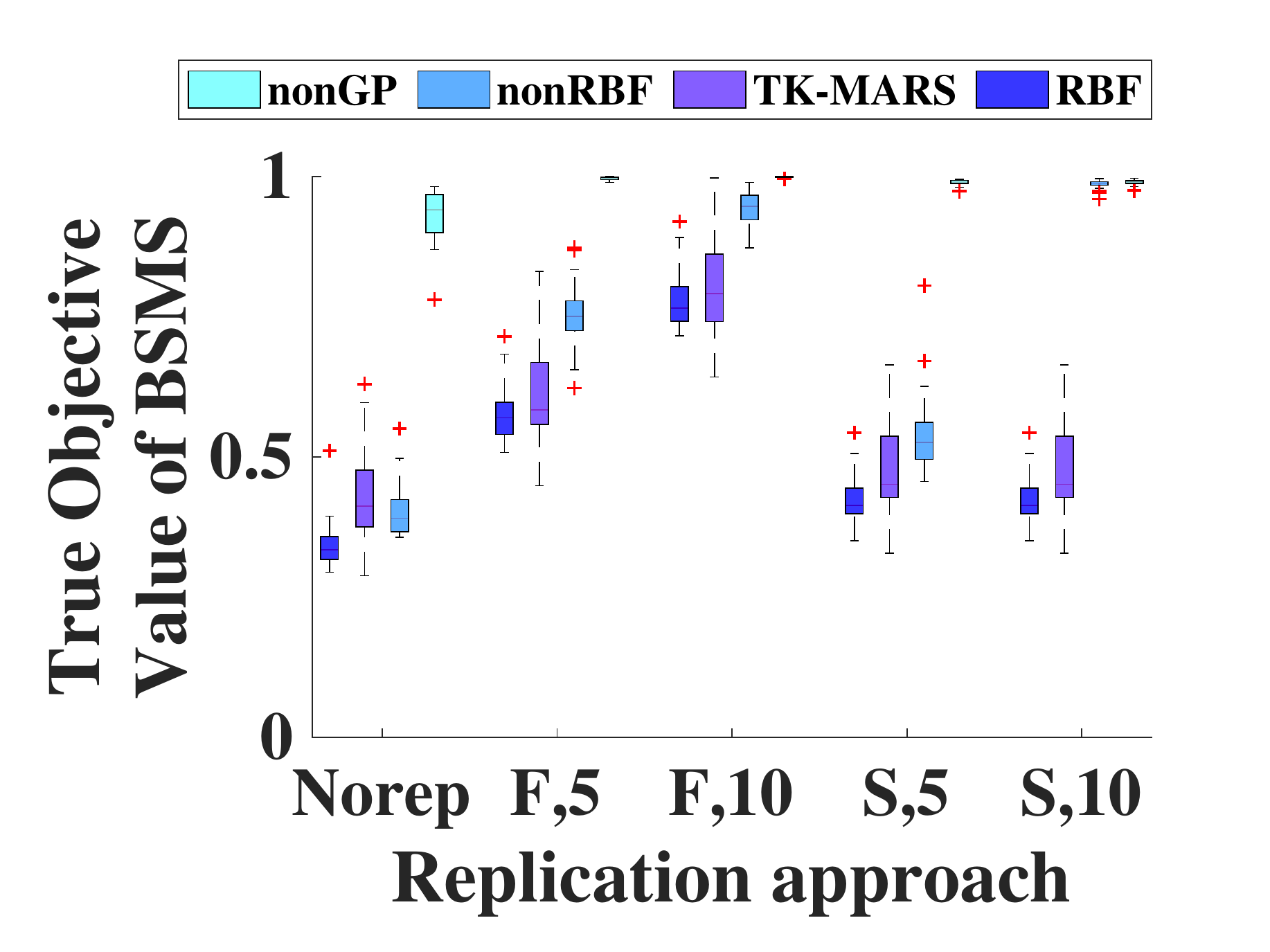}}
        \subfloat[Noise=0.05]{\includegraphics[width=0.55\textwidth]{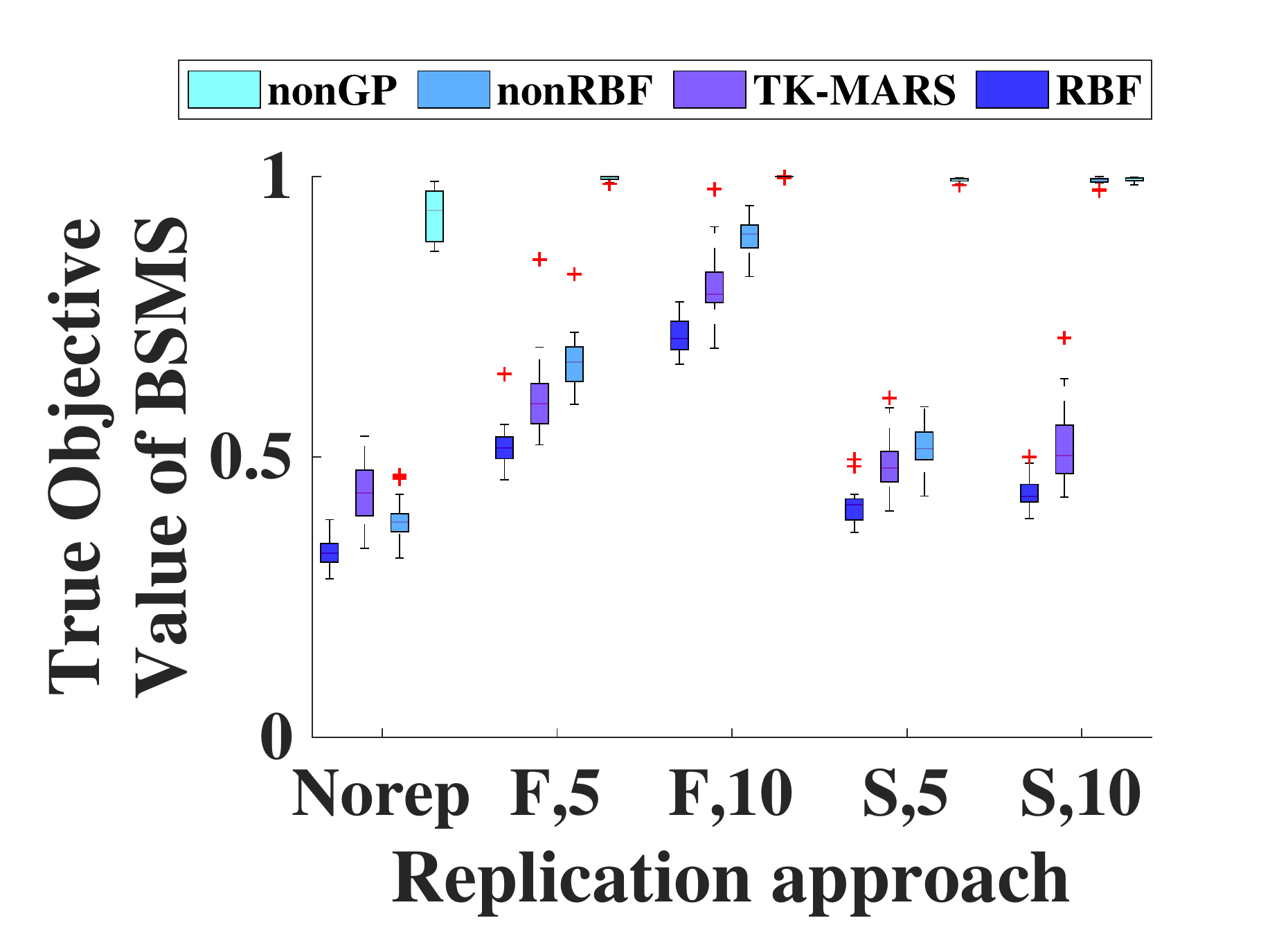}}
    \end{minipage}
    \begin{minipage}{\linewidth}
        \subfloat[Noise=0.1]{\includegraphics[width=0.55\textwidth]{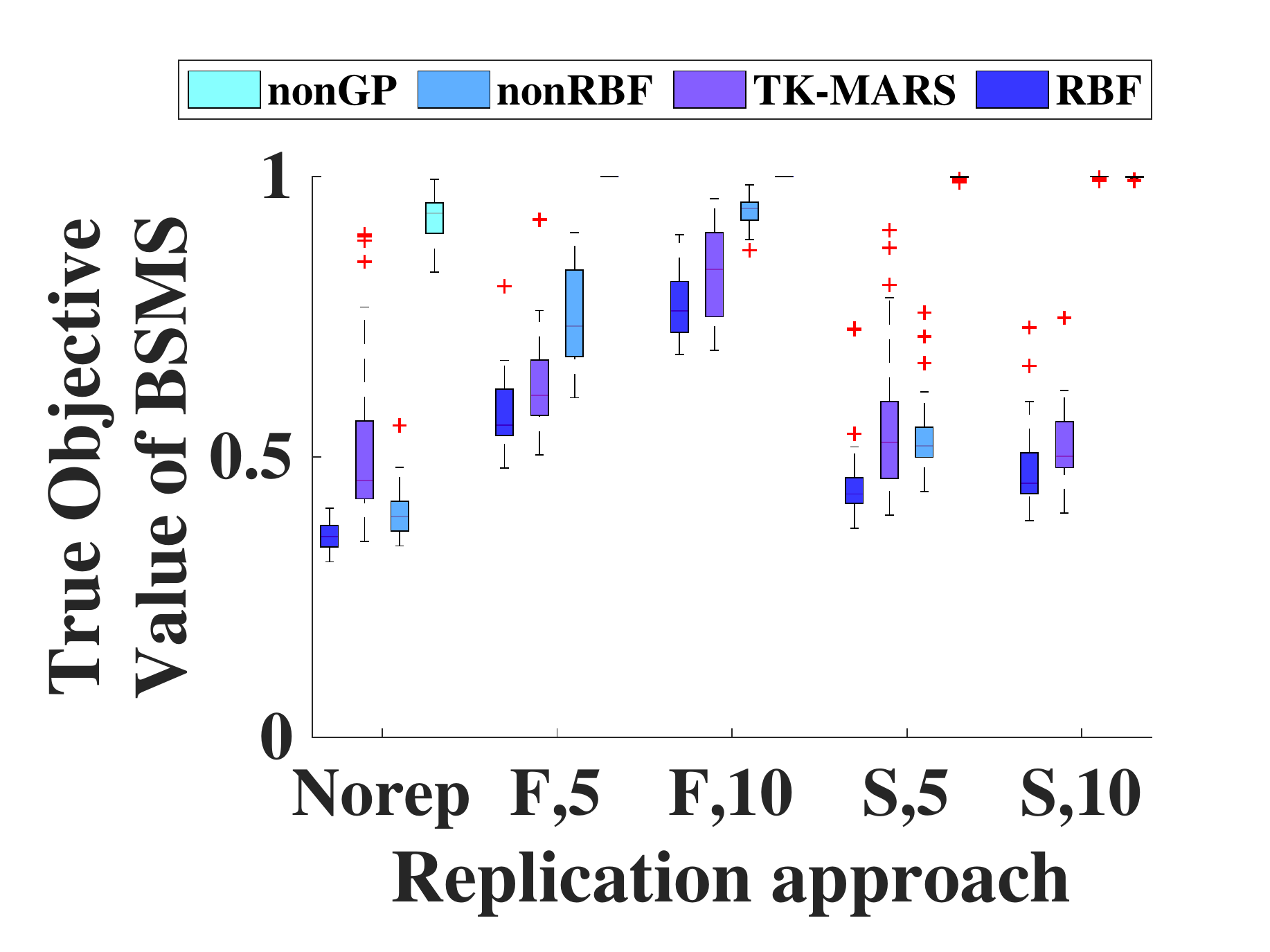}}
        \subfloat[Noise=0.25]{\includegraphics[width=0.55\textwidth]{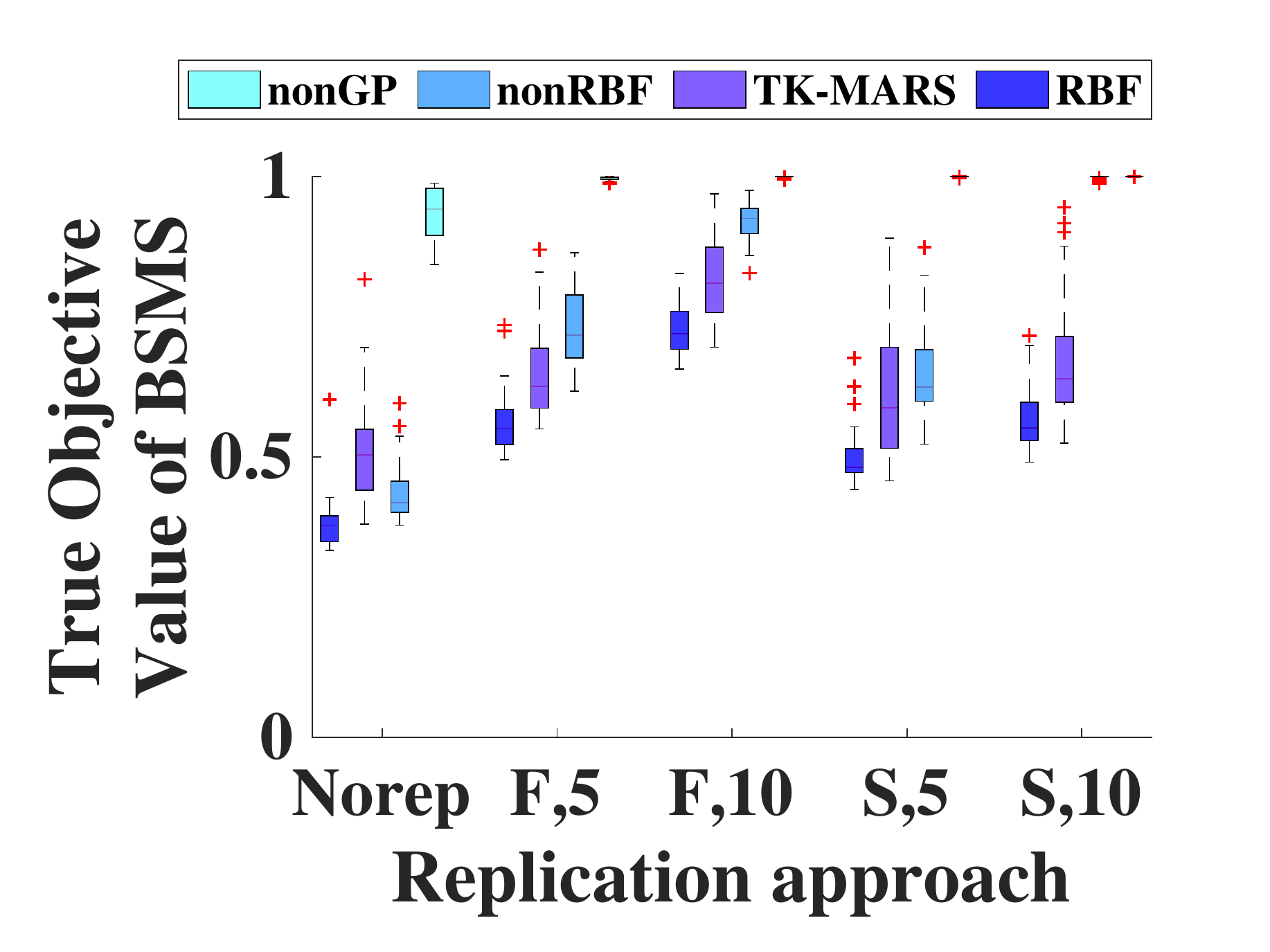}}
    \end{minipage}
    \caption{Box-plots of MTFAUC of surrogate optimization on the Ackley function}
    \label{fig:auc_box_ack}
    \vspace{-3mm}
\end{figure}

\hadis{Figure~\ref{fig:auc_box_ack} presents the box-plots of the MTFAUC for 30 executions using different approaches on the Ackley test function. We observe that $\mathit{Smartrep},10$ has the highest robustness to uncertainty. %Note that although $\mathit{Fixedrep},10$ performance is competitive, it has a larger true objective value of the BSMS, comparatively, as shown in Figure~\ref{fig:bks_box_ack}
.} 

\begin{figure}[!tb]
\centering
    \begin{minipage}{\linewidth}
        \subfloat[Noise=0]{\includegraphics[width=0.55\textwidth]{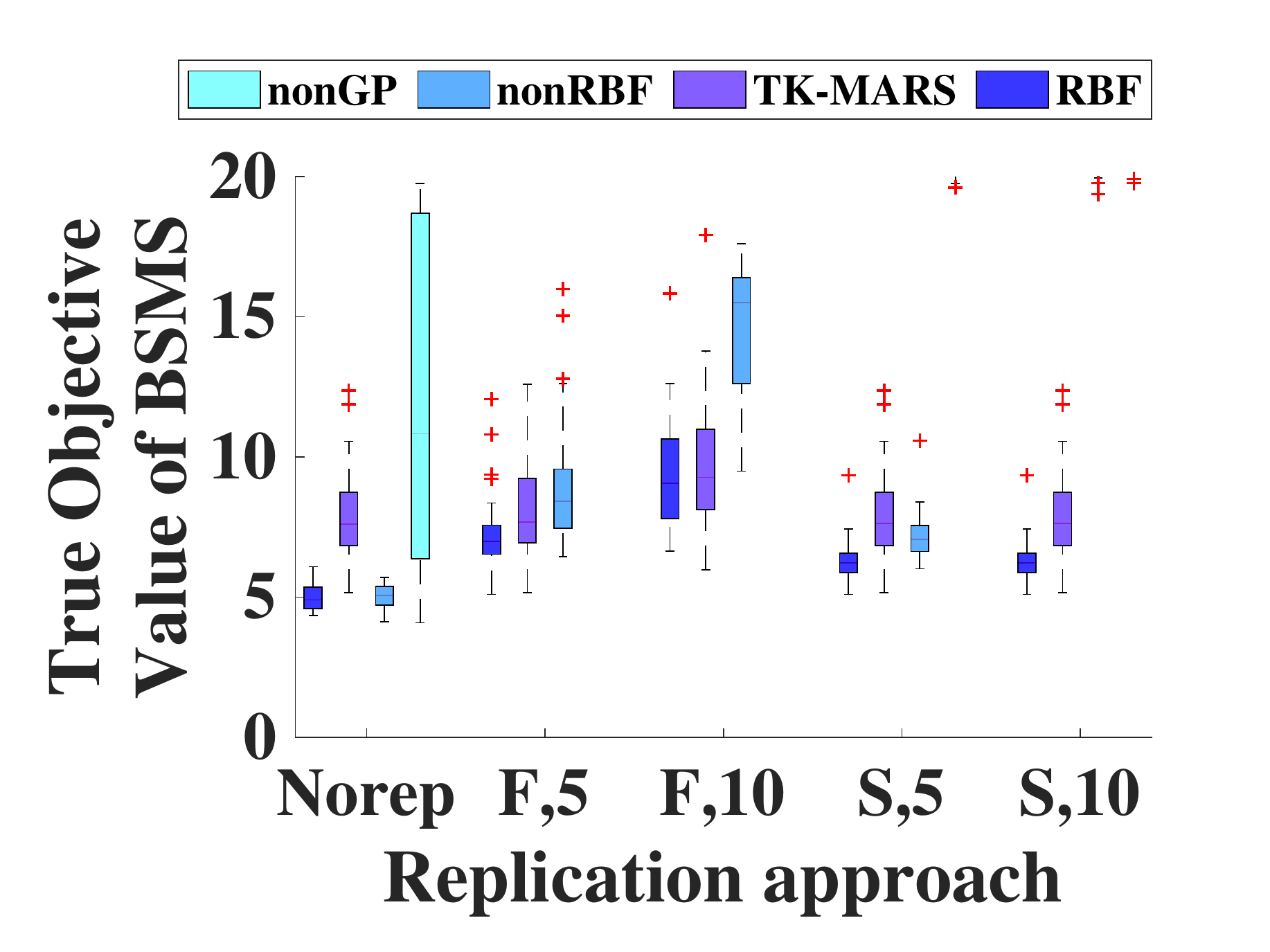}}
        \subfloat[Noise=0.05]{\includegraphics[width=0.55\textwidth]{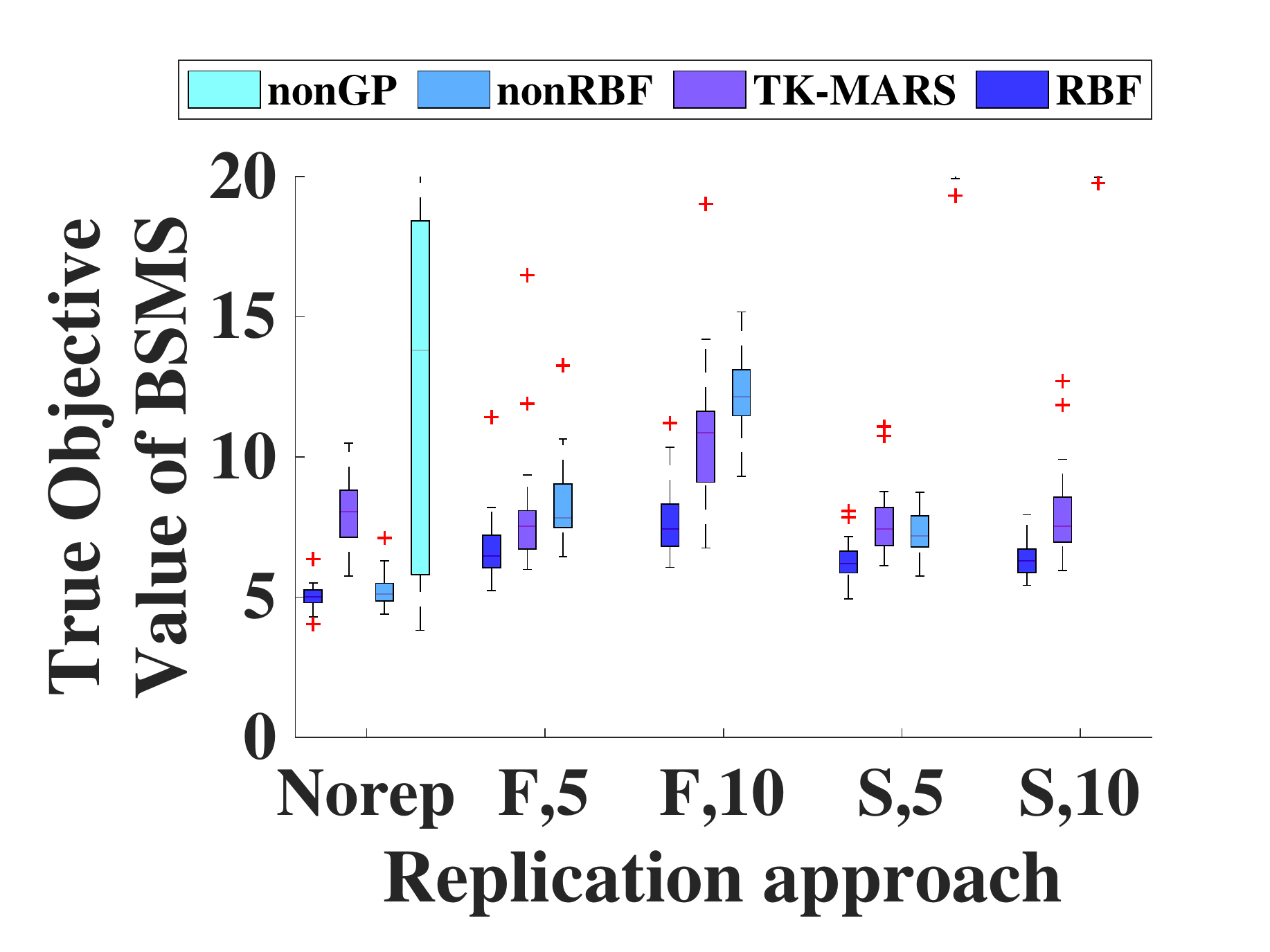}}
    \end{minipage}
    \begin{minipage}{\linewidth}
        \subfloat[Noise=0.1]{\includegraphics[width=0.55\textwidth]{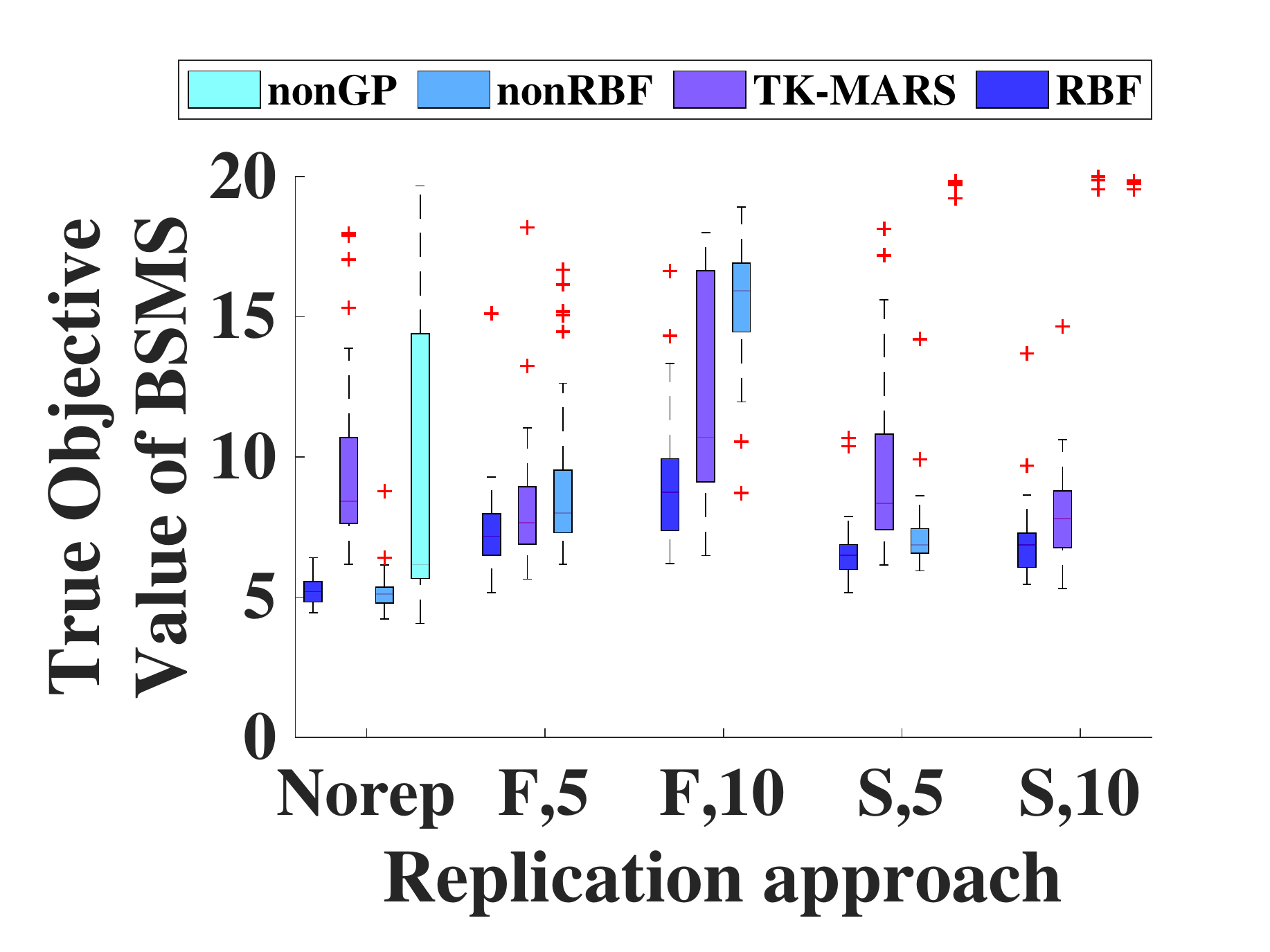}}
        \subfloat[Noise=0.25]{\includegraphics[width=0.55\textwidth]{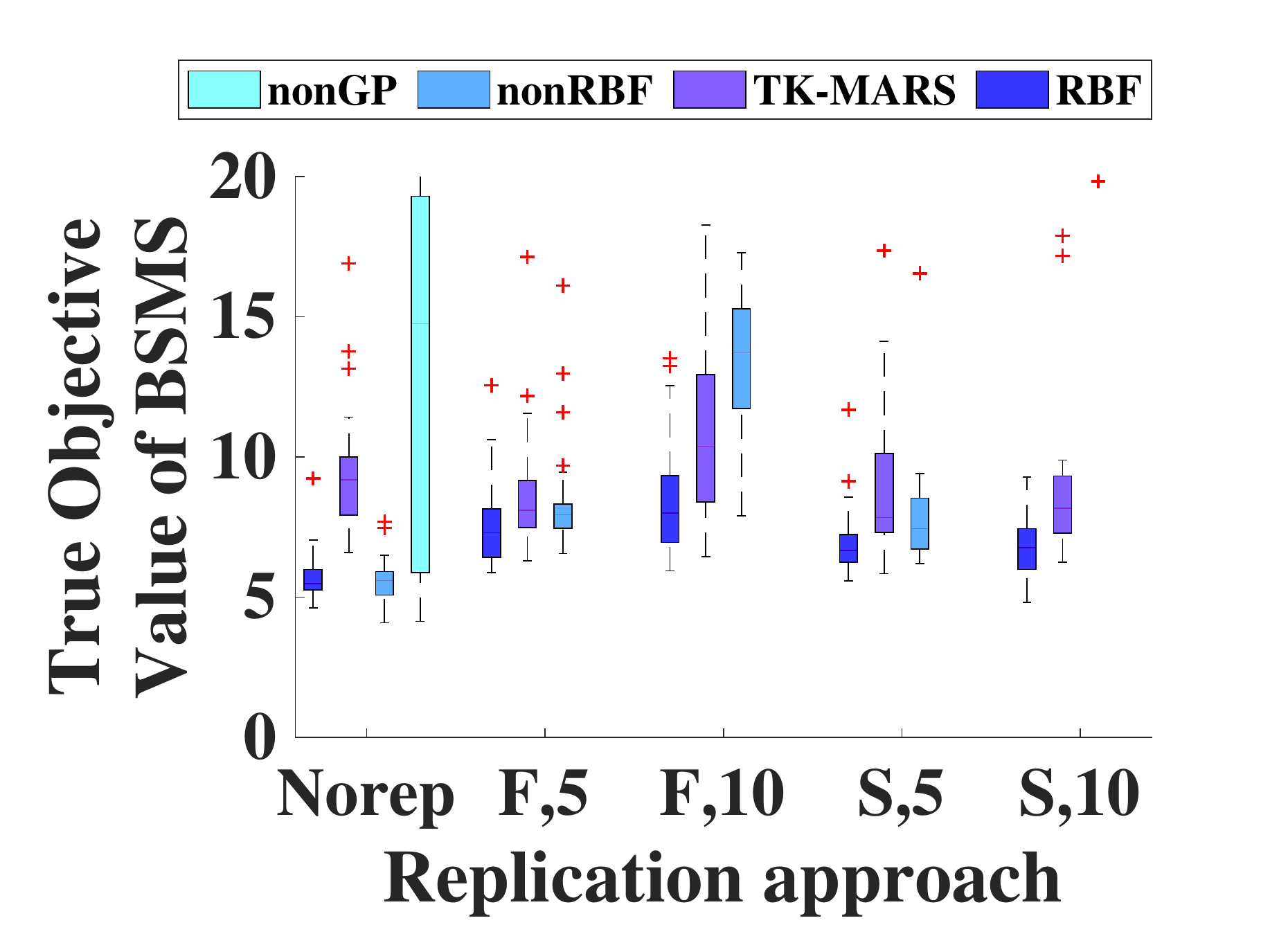}}
    \end{minipage}
    \caption{Box-plots of the true objective value of the BSMS after 1000 black-box function evaluations of surrogate optimization on the Ackley function}
    \label{fig:bks_box_ack}
    \vspace{-3mm}
\end{figure}

\hadis{Figure~\ref{fig:bks_box_ack} shows the box-plots of the true objective value of BSMS after 1000 black-box function evaluations for the same 30 executions.
On average, $\mathit{Smartrep},10$ finds better BSMS across different noise levels.}

\begin{figure}[!tb]
\centering
    \begin{minipage}{\linewidth}
        \subfloat[Noise=0]{\includegraphics[width=0.55\textwidth]{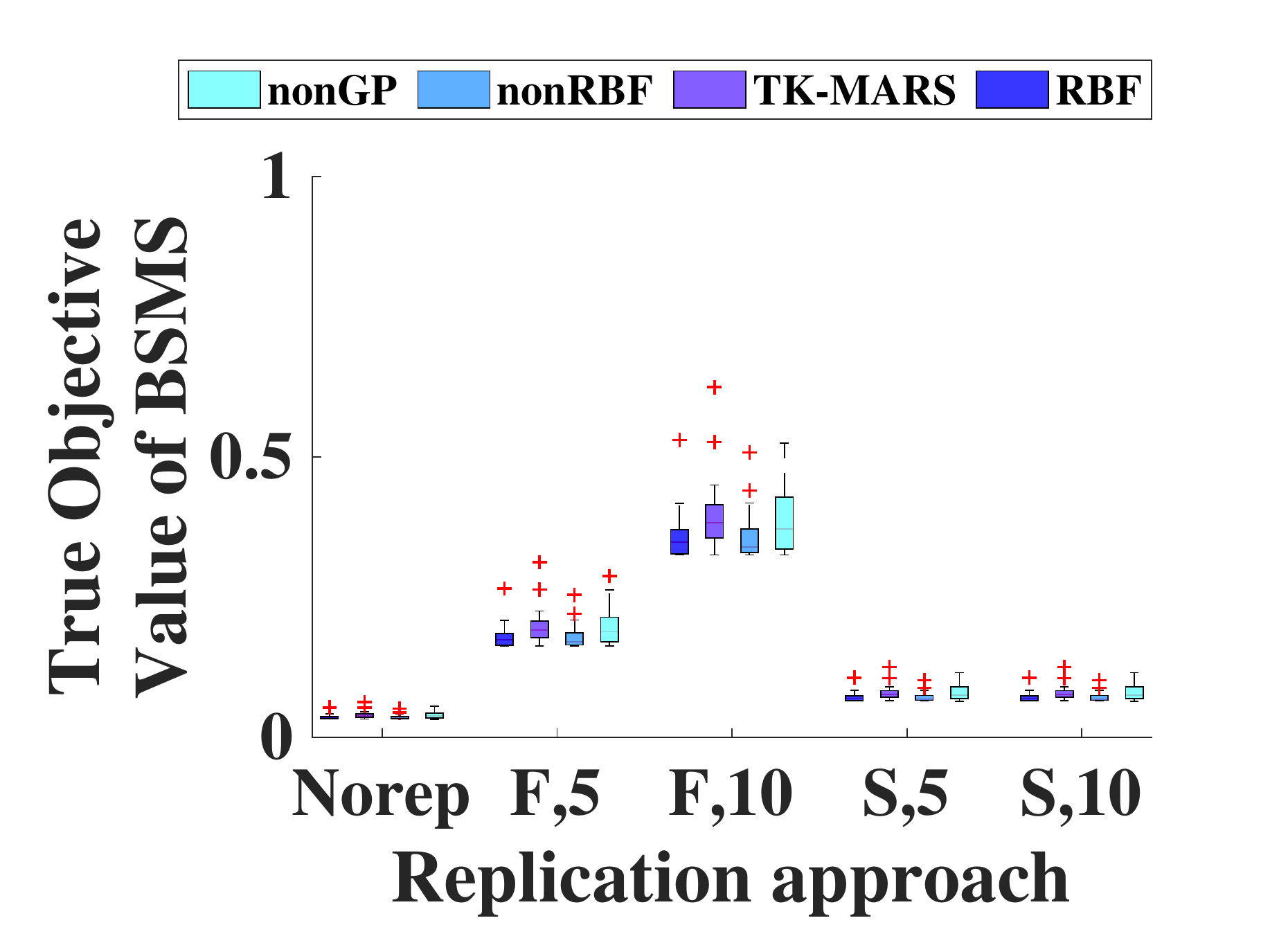}}
        \subfloat[Noise=0.05]{\includegraphics[width=0.55\textwidth]{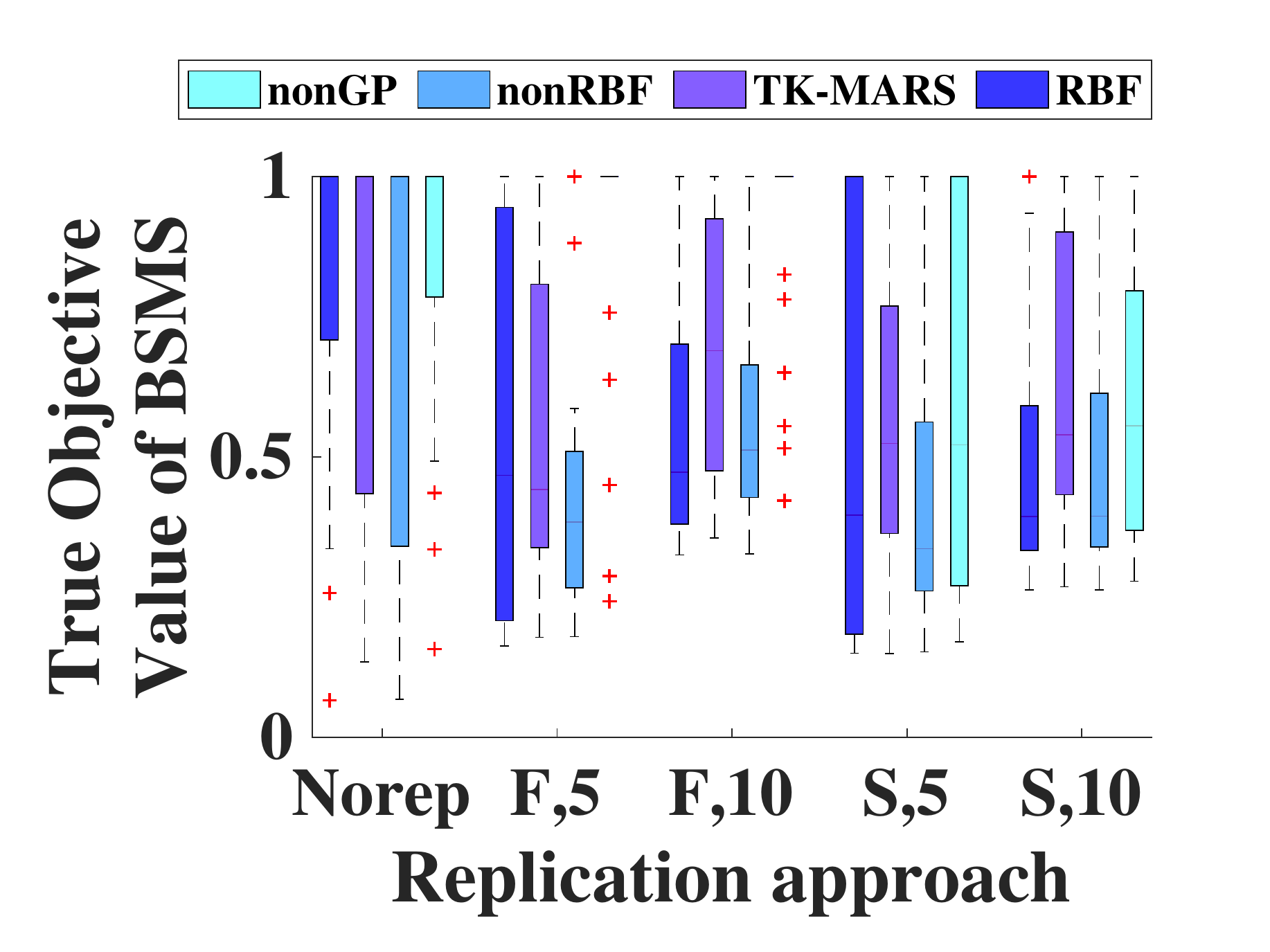}}
    \end{minipage}
    \begin{minipage}{\linewidth}
        \subfloat[Noise=0.1]{\includegraphics[width=0.55\textwidth]{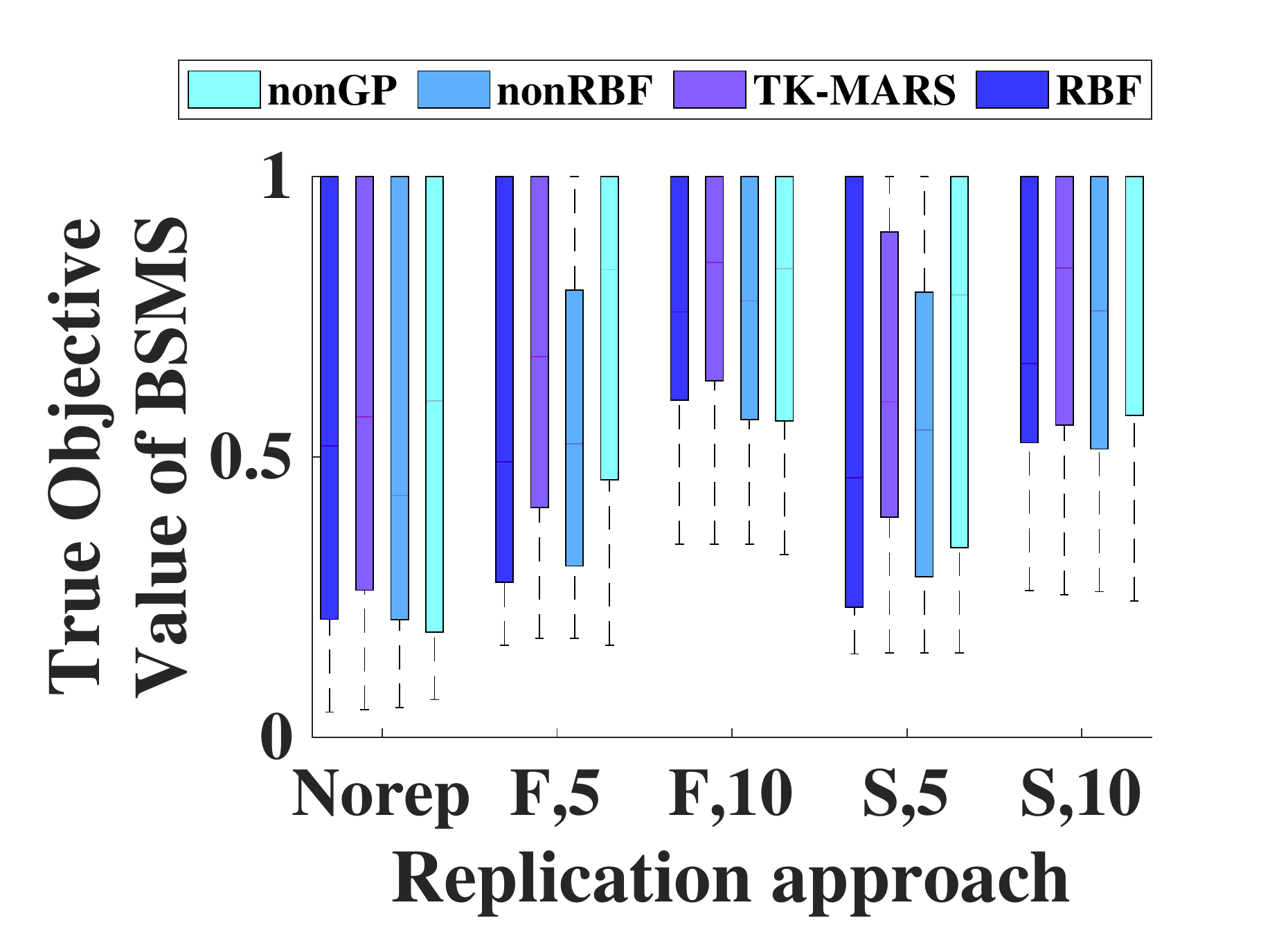}}
        \subfloat[Noise=0.25]{\includegraphics[width=0.55\textwidth]{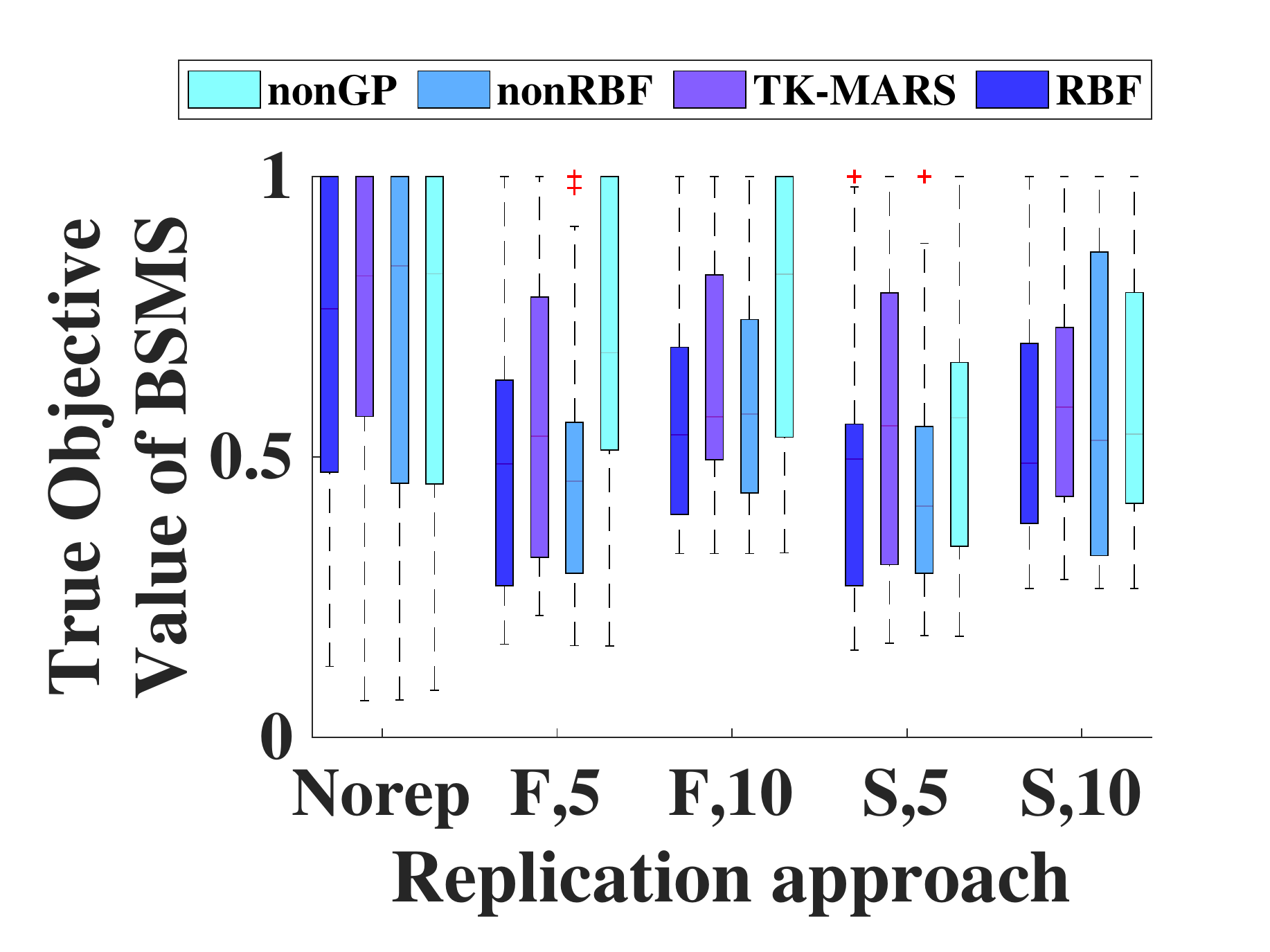}}
    \end{minipage}
    \caption{Box-plots of MTFAUC of surrogate optimization on the Zakharov function}
    \label{fig:auc_box_zak}
    \vspace{-3mm}
\end{figure}

\hadis{Figure~\ref{fig:auc_box_zak} presents the box-plots of the MTFAUC for 30 executions using different approaches on the Zakharov test function. The performance of $\mathit{Smartrep},10$ is noticeable at the higher noise levels. Although the performance of $\mathit{Fixedrep},10$ is competitive to some extent, it conducts unnecessary black-box function evaluations in a deterministic environment.}

\begin{figure}[!tb]
\centering
    \begin{minipage}{\linewidth}
        \subfloat[Noise=0]{\includegraphics[width=0.55\textwidth]{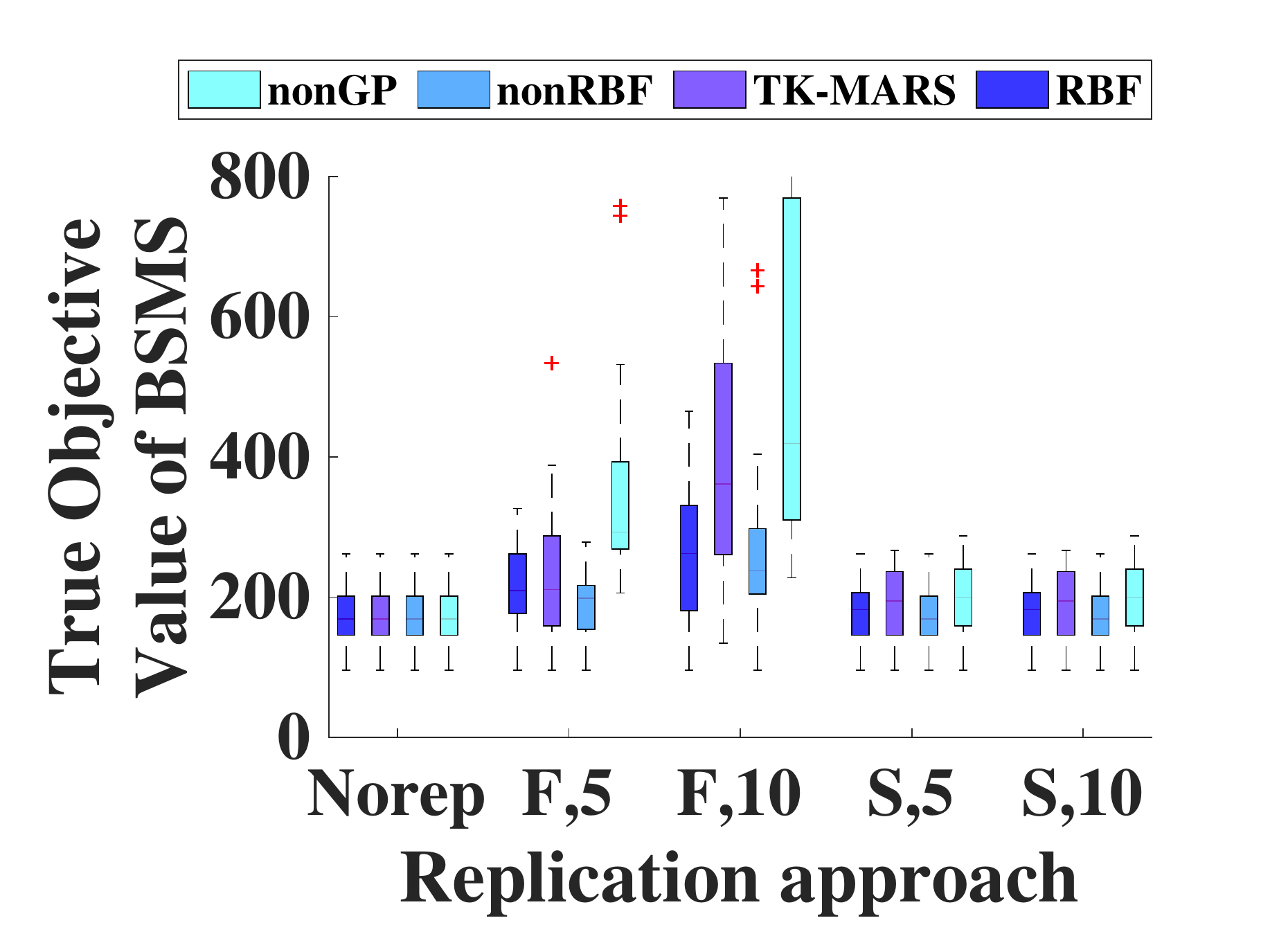}}
        \subfloat[Noise=0.05]{\includegraphics[width=0.55\textwidth]{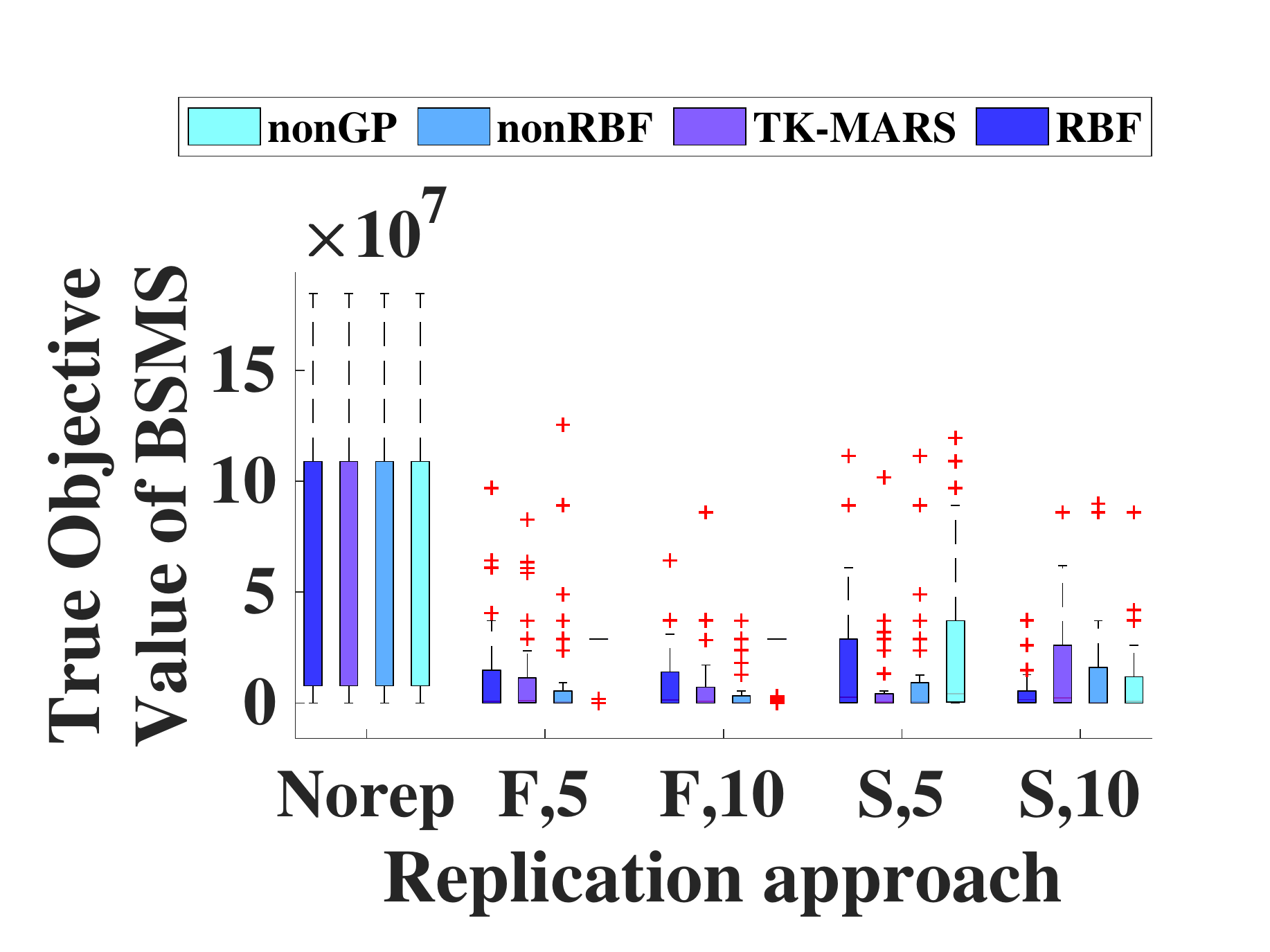}}
    \end{minipage}
    \begin{minipage}{\linewidth}
        \subfloat[Noise=0.1]{\includegraphics[width=0.55\textwidth]{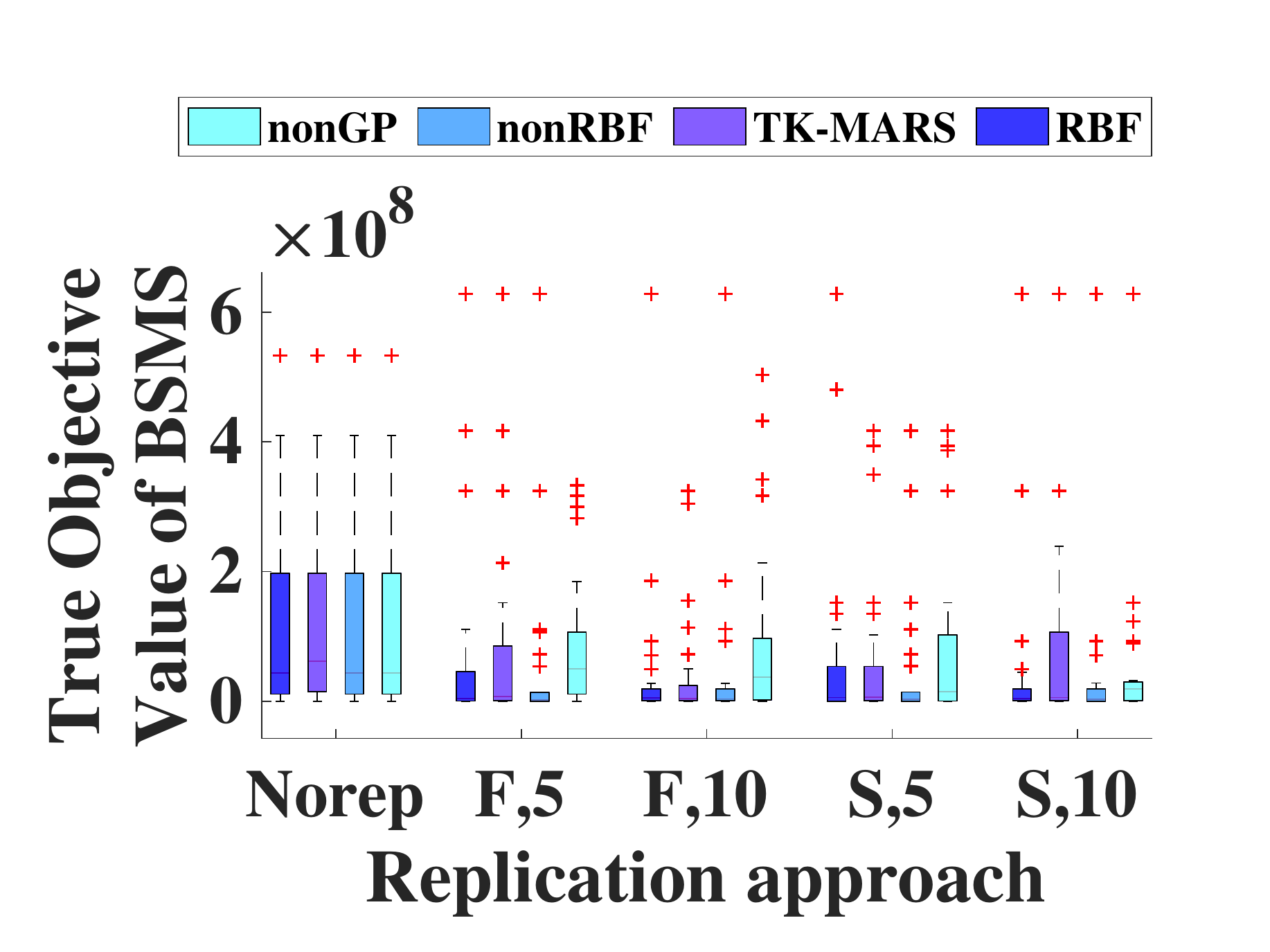}}
        \subfloat[Noise=0.25]{\includegraphics[width=0.55\textwidth]{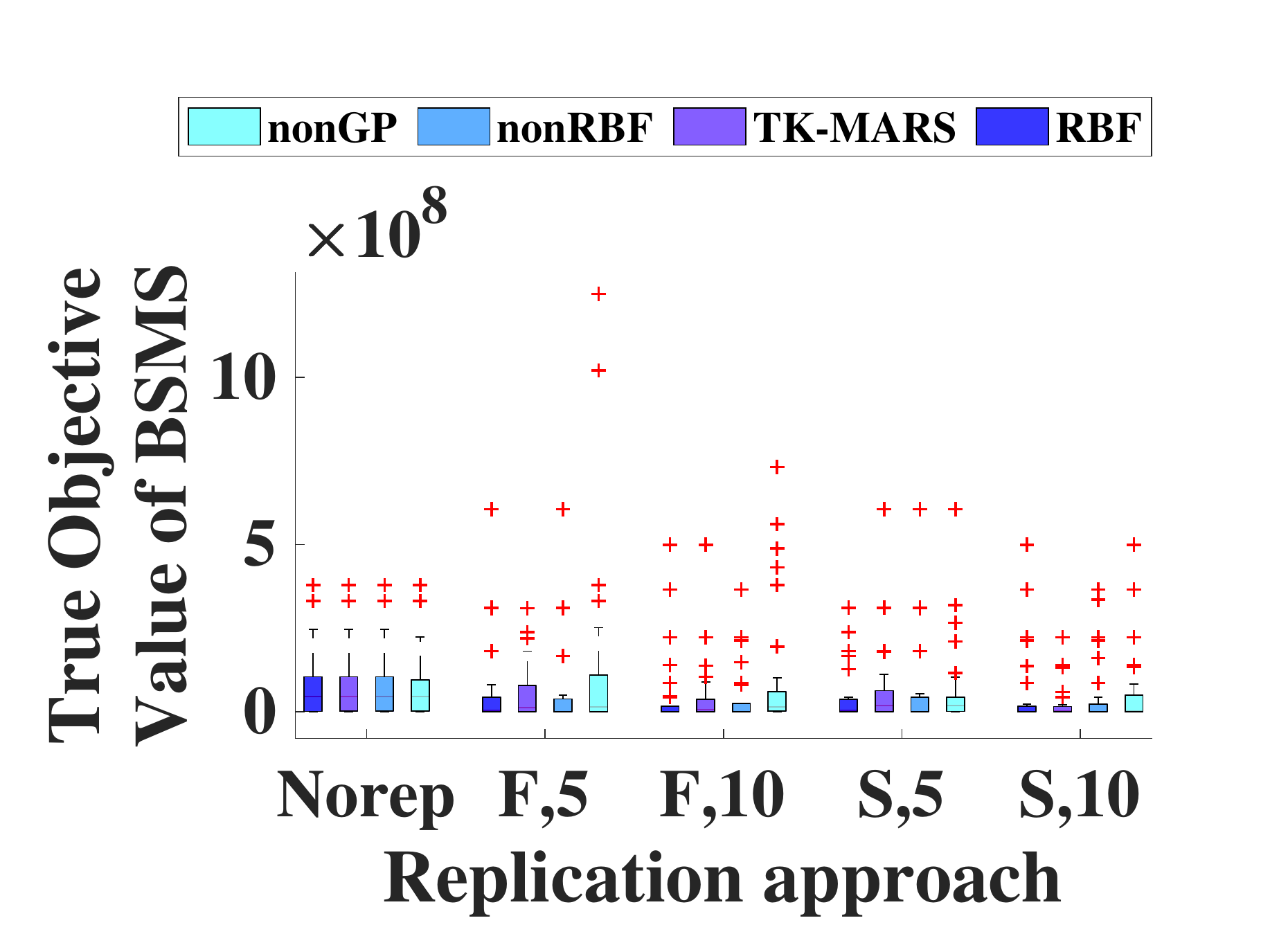}}
    \end{minipage}
    \caption{Box-plots of the true objective value of the BSMS after 1000 black-box function evaluations of surrogate optimization on the Zakharov function}
    \label{fig:bks_box_zak}
    \vspace{-3mm}
\end{figure}

\hadis{Figure~\ref{fig:bks_box_zak} shows the box-plots of the true objective value of the BSMS after 1000 black-box function evaluations for the same 30 executions on the Zakharov test function.
We verified that $\mathit{Smartrep},10$ has the smallest interquartile range overall in higher levels of noise compared with $\mathit{Fixedrep},10$, but at the lower levels of noise, it has larger range compared with $\mathit{Fixedrep},10$.}

% \begin{figure}[htb]
% \centering
%     \begin{minipage}{\linewidth}
%         \subfloat[]{\includegraphics[width=0.35\textwidth]{figures/box_auc_rast_new_1.pdf}}
%         \subfloat[]{\includegraphics[width=0.35\textwidth]{figures/box_auc_rast_new_2.pdf}}
%         \subfloat[]{\includegraphics[width=0.35\textwidth]{figures/box_auc_rast_new_3.pdf}}
%     \end{minipage}
%     \begin{minipage}{\linewidth}
%         \subfloat[]{\includegraphics[width=0.35\textwidth]{figures/box_auc_rast_new_4.pdf}}
%         \subfloat[]{\includegraphics[width=0.35\textwidth]{figures/box_auc_rast_new_5.pdf}}
%     \end{minipage}
%     \caption{Box-plots of MTFAUC of surrogate optimization on the Rastrigin function}
%     \label{fig:auc_box_rast}
%     \vspace{-3mm}

% \end{figure}

% \begin{figure}[htb]
% \centering
%     \begin{minipage}{\linewidth}
%         \subfloat[]{\includegraphics[width=0.35\textwidth]{figures/box_bks_rast_new_1.pdf}}
%         \subfloat[]{\includegraphics[width=0.35\textwidth]{figures/box_bks_rast_new_2.pdf}}
%         \subfloat[]{\includegraphics[width=0.35\textwidth]{figures/box_bks_rast_new_3.pdf}}
%     \end{minipage}
%     \begin{minipage}{\linewidth}
%         \subfloat[]{\includegraphics[width=0.35\textwidth]{figures/box_bks_rast_new_4.pdf}}
%         \subfloat[]{\includegraphics[width=0.35\textwidth]{figures/box_bks_rast_new_5.pdf}}
%     \end{minipage}
%     \caption{Box-plots of the true objective value of the BSMS after 1000 black-box function evaluations of surrogate optimization on the Rastrigin function}
%     \label{fig:bks_box_rast}
%     \vspace{-3mm}
% \end{figure}

% \begin{figure}[htb]
% \centering
%     \begin{minipage}{\linewidth}
%         \subfloat[]{\includegraphics[width=0.35\textwidth]{figures/box_auc_levy_new_1.pdf}}
%         \subfloat[]{\includegraphics[width=0.35\textwidth]{figures/box_auc_levy_new_2.pdf}}
%         \subfloat[]{\includegraphics[width=0.35\textwidth]{figures/box_auc_levy_new_3.pdf}}
%     \end{minipage}
%     \begin{minipage}{\linewidth}
%         \subfloat[]{\includegraphics[width=0.35\textwidth]{figures/box_auc_levy_new_4.pdf}}
%         \subfloat[]{\includegraphics[width=0.35\textwidth]{figures/box_auc_levy_new_5.pdf}}
%     \end{minipage}
%     \caption{Box-plots of MTFAUC of surrogate optimization on the Levy function}
%     \label{fig:auc_box_levy}
%     \vspace{-3mm}

% \end{figure}

% \begin{figure}[htb]
% \centering
%     \begin{minipage}{\linewidth}
%         \subfloat[]{\includegraphics[width=0.35\textwidth]{figures/box_bks_levy_new_1.pdf}}
%         \subfloat[]{\includegraphics[width=0.35\textwidth]{figures/box_bks_levy_new_2.pdf}}
%         \subfloat[]{\includegraphics[width=0.35\textwidth]{figures/box_bks_levy_new_3.pdf}}
%     \end{minipage}
%     \begin{minipage}{\linewidth}
%         \subfloat[]{\includegraphics[width=0.35\textwidth]{figures/box_bks_levy_new_4.pdf}}
%         \subfloat[]{\includegraphics[width=0.35\textwidth]{figures/box_bks_levy_new_5.pdf}}
%     \end{minipage}
%     \caption{Box-plots of the true objective value of the BSMS after 1000 black-box function evaluations of surrogate optimization on the Levy function}
%     \label{fig:bks_box_levy}
%     \vspace{-3mm}
% \end{figure}

% \begin{figure}[htb]
% \centering
%     \begin{minipage}{\linewidth}
%         \subfloat[]{\includegraphics[width=0.35\textwidth]{figures/box_auc_ack_new_1.pdf}}
%         \subfloat[]{\includegraphics[width=0.35\textwidth]{figures/box_auc_ack_new_2.pdf}}
%         \subfloat[]{\includegraphics[width=0.35\textwidth]{figures/box_auc_ack_new_3.pdf}}
%     \end{minipage}
%     \begin{minipage}{\linewidth}
%         \subfloat[]{\includegraphics[width=0.35\textwidth]{figures/box_auc_ack_new_4.pdf}}
%         \subfloat[]{\includegraphics[width=0.35\textwidth]{figures/box_auc_ack_new_5.pdf}}
%     \end{minipage}
%     \caption{Box-plots of MTFAUC of surrogate optimization on the Ackley function}
%     \label{fig:auc_box_ack}
%     \vspace{-3mm}

% \end{figure}

% \begin{figure}[htb]
% \centering
%     \begin{minipage}{\linewidth}
%         \subfloat[]{\includegraphics[width=0.35\textwidth]{figures/box_bks_ack_new_1.pdf}}
%         \subfloat[]{\includegraphics[width=0.35\textwidth]{figures/box_bks_ack_new_2.pdf}}
%         \subfloat[]{\includegraphics[width=0.35\textwidth]{figures/box_bks_ack_new_3.pdf}}
%     \end{minipage}
%     \begin{minipage}{\linewidth}
%         \subfloat[]{\includegraphics[width=0.35\textwidth]{figures/box_bks_ack_new_4.pdf}}
%         \subfloat[]{\includegraphics[width=0.35\textwidth]{figures/box_bks_ack_new_5.pdf}}
%     \end{minipage}
%     \caption{Box-plots of the true objective value of the BSMS after 1000 black-box function evaluations of surrogate optimization on the Ackley function}
%     \label{fig:bks_box_ack}
%     \vspace{-3mm}
% \end{figure}

% \begin{figure}[htb]
% \centering

%     \begin{minipage}{\linewidth}
%         \subfloat[]{\includegraphics[width=0.35\textwidth]{figures/box_auc_zak_new_1.pdf}}
%         \subfloat[]{\includegraphics[width=0.35\textwidth]{figures/box_auc_zak_new_2.pdf}}
%         \subfloat[]{\includegraphics[width=0.35\textwidth]{figures/box_auc_zak_new_3.pdf}}
%     \end{minipage}
%     \begin{minipage}{\linewidth}
%         \subfloat[]{\includegraphics[width=0.35\textwidth]{figures/box_auc_zak_new_4.pdf}}
%         \subfloat[]{\includegraphics[width=0.35\textwidth]{figures/box_auc_zak_new_5.pdf}}
%     \end{minipage}
%     \caption{Box-plots of MTFAUC of surrogate optimization on the Zakharov function}
%     \label{fig:auc_box_zak}
%     \vspace{-3mm}

% \end{figure}

% \begin{figure}[htb]
% \centering
%     \begin{minipage}{\linewidth}
%         \subfloat[]{\includegraphics[width=0.35\textwidth]{figures/box_bks_zak_new_1.pdf}}
%         \subfloat[]{\includegraphics[width=0.35\textwidth]{figures/box_bks_zak_new_2.pdf}}
%         \subfloat[]{\includegraphics[width=0.35\textwidth]{figures/box_bks_zak_new_3.pdf}}
%     \end{minipage}
%     \begin{minipage}{\linewidth}
%         \subfloat[]{\includegraphics[width=0.35\textwidth]{figures/box_bks_zak_new_4.pdf}}
%         \subfloat[]{\includegraphics[width=0.35\textwidth]{figures/box_bks_zak_new_5.pdf}}
%     \end{minipage}
%     \caption{Box-plots of the true objective value of the BSMS after 1000 black-box function evaluations of surrogate optimization on the Zakharov function}
%     \label{fig:bks_box_zak}
%     \vspace{-3mm}
% \end{figure}

% \newpage
% Table generated by Excel2LaTeX from sheet 'Sheet1'
\newpage
% \hspace{-30mm}

\begin{table}[h]
  \centering
  \caption{RBF Results}
  \begin{tiny}
    \begin{tabular}{|l|c|c|c|c|c|c|c|c|}
    \toprule
          & \multicolumn{8}{c|}{\new{\textbf{Multi-Quadric}}} \\
    \midrule
          &           \multicolumn{2}{c|}{\textbf{\new{$fiv=1 $}   }} & \multicolumn{2}{c|}{\textbf{\new{$fiv= 0.75$}}} & \multicolumn{2}{c|}{\textbf{\new{$fiv= 0.5$}}} & \multicolumn{2}{c|}{\textbf{\new{$fiv= 0.25$}}} \\
    \midrule
    \textbf{Rosenbrock} & 0.0087 & 47    & 0.0077 & 46    & 0.0126 & 52    & 0.0203 & 131 \\
    \midrule
    \textbf{Rastrigin} & 0.0087 & 47    & 0.0087 & 47    & 0.0145 & 53    & 0.0666 & 107 \\
    \midrule
    \textbf{Sphere} & 0.0087 & 47    & 0.0087 & 47    & 0.0079 & 47    & 0.0203 & 59 \\
    \midrule
    \textbf{Levy} & 0.0145 & 53    & 0.0145 & 53    & 0.0203 & 59    & 0.0366 & 83 \\ 
    \midrule 
          & \multicolumn{8}{c|}{\new{\textbf{Gaussian}}} \\
    \midrule
          &            \multicolumn{2}{c|}{\textbf{\new{$fiv=1 $}   }} & \multicolumn{2}{c|}{\textbf{\new{$fiv= 0.75$}}} & \multicolumn{2}{c|}{\textbf{\new{$fiv= 0.5$}}} & \multicolumn{2}{c|}{\textbf{\new{$fiv= 0.25$}}} \\
    \midrule
    \textbf{Rosenbrock} & 0.7376 & 895   & 0.7156 & 895   & 0.4293 & 895   & 0.2935 & 896 \\
    \midrule
    \textbf{Rastrigin} & 0.0087 & 47    & 0.0087 & 47    & 0.0087 & 47    & 0.0087 & 47 \\
    \midrule
    \textbf{Sphere} & 0.0083 & 47    & 0.0074 & 47    & 0.0087 & 47    & 0.0145 & 53 \\
    \midrule
    \textbf{Levy} & 0.8338 & 1035  & 0.8350 & 1035  & 0.6128 & 1035  & 0.4628 & 1035 \\
    \midrule
          &            \multicolumn{2}{c|}{\textbf{\new{$fiv=1 $}   }} & \multicolumn{2}{c|}{\textbf{\new{$fiv= 0.75$}}} & \multicolumn{2}{c|}{\textbf{\new{$fiv= 0.5$}}} & \multicolumn{2}{c|}{\textbf{\new{$fiv= 0.25$}}} \\
    \midrule
    \textbf{Rosenbrock} & 0.0183 & 57    & 0.0219 & 65    & 0.0287 & 69    & 0.0434 & 309 \\
    \midrule
    \textbf{Rastrigin} & 0.0241 & 63    & 0.0319 & 71    & 0.0550 & 95    & 0.3304 & 397 \\
    \midrule
    \textbf{Sphere} & 0.0193 & 58    & 0.0183 & 57    & 0.0203 & 59    & 0.0376 & 77 \\
    \midrule
    \textbf{Levy} & 0.0241 & 63    & 0.0261 & 65    & 0.0425 & 82    & 0.0531 & 93 \\
    \midrule
          & \multicolumn{8}{c|}{\new{\textbf{Thin Plate Spline}}} \\
    \midrule
          &            \multicolumn{2}{c|}{\textbf{\new{$fiv=1 $}   }} & \multicolumn{2}{c|}{\textbf{\new{$fiv= 0.75$}}} & \multicolumn{2}{c|}{\textbf{\new{$fiv= 0.5$}}} & \multicolumn{2}{c|}{\textbf{\new{$fiv= 0.25$}}} \\
    \midrule
    \textbf{Rosenbrock} & 0.0225 & 62    & 0.0219 & 65    & 0.0361 & 78    & 0.0245 & 171 \\
    \midrule
    \textbf{Rastrigin} & 0.0245 & 61    & 0.0261 & 65    & 0.0290 & 68    & 0.0705 & 111 \\
    \midrule
    \textbf{Sphere} & 0.0241 & 63    & 0.0241 & 63    & 0.0193 & 58    & 0.0319 & 71 \\
    \midrule
    \textbf{Levy} & 0.0357 & 75    & 0.0424 & 83    & 0.0338 & 73    & 0.0666 & 107 \\
    \bottomrule
    \end{tabular}%
    \end{tiny}
  \label{tab:results-RBF}%
\end{table}%

\begin{table}[h]
  \centering
  \caption{Parameters and levels for OA design in Preliminary Analysis}
  \begin{scriptsize}
   \begin{tabular}{|ll|}
    \toprule
	\textbf{Problem Parameters} & \textbf{levels} \\ \midrule
	Test function & Rosenbrock, Rastrigin, Levy \\
          Dimension & 10, 20, 30 \\
          Fraction of important variables & 0.25, 0.50, 0.75, 1 \\
          Noise level (\%) & 5, 10, 25 \\
	\midrule
    \textbf{Algorithm Parameters} & \textbf{levels} \\ \midrule
	Initial pool size & $d+1$, $2(d+1)$ \\
          DOE method & LHD, sobol \\
          EEPA distance & Euclidean, Cosine \\
          EEPA \# candidates & 3, 6 \\
          Replication type & fixed\_rep, smart\_rep \\
          Replication \# & 5, 10 \\
          Model  & RBF\_avrg, TK-MARS\_avrg, TK-MARS\_keepall \\ \bottomrule
    \end{tabular}%
    % \vspace{-10mm}
    \end{scriptsize}
  \label{tab:params-OA}%
\end{table}%
    % \vspace{-15mm}

\begin{table}[h]
  \centering
  \caption{Preliminary ANOVA table}
    \begin{scriptsize}
    \begin{tabular}{|lcccccc|}
    \toprule
    \multicolumn{1}{|r|}{} & \multicolumn{1}{l}{Estimate} & \multicolumn{1}{l}{Std. Error} & \multicolumn{1}{l}{t value} & \multicolumn{1}{l}{Pr(>|t|)} &       &  \\
    \midrule
    \multicolumn{1}{|l|}{(Intercept)} & \multicolumn{1}{r}{0.250806} & \multicolumn{1}{r}{0.086891} & \multicolumn{1}{r}{2.886} & \multicolumn{1}{r}{0.005626} & \multicolumn{1}{l}{**} &  \\
    \multicolumn{1}{|l|}{fiv=0.75} & \multicolumn{1}{r}{0.07965} & \multicolumn{1}{r}{0.056382} & \multicolumn{1}{r}{1.413} & \multicolumn{1}{r}{0.163596} &       &  \\
    \multicolumn{1}{|l|}{fiv=0.50} & \multicolumn{1}{r}{0.207822} & \multicolumn{1}{r}{0.056382} & \multicolumn{1}{r}{3.686} & \multicolumn{1}{r}{0.000537} & \multicolumn{1}{l}{***} &  \\
    \multicolumn{1}{|l|}{fiv=0.25} & \multicolumn{1}{r}{0.132775} & \multicolumn{1}{r}{0.056382} & \multicolumn{1}{r}{2.355} & \multicolumn{1}{r}{0.022263} & \multicolumn{1}{l}{*} &  \\
    \multicolumn{1}{|l|}{model=TK-MARS\_avrg} & \multicolumn{1}{r}{-0.04463} & \multicolumn{1}{r}{0.048828} & \multicolumn{1}{r}{-0.914} & \multicolumn{1}{r}{0.364833} &       &  \\
    \multicolumn{1}{|l|}{model=TK-MARS\_keepall} & \multicolumn{1}{r}{0.113446} & \multicolumn{1}{r}{0.048828} & \multicolumn{1}{r}{2.323} & \multicolumn{1}{r}{0.024028} & \multicolumn{1}{l}{*} &  \\
    \multicolumn{1}{|l|}{noise level=10\%} & \multicolumn{1}{r}{0.004939} & \multicolumn{1}{r}{0.048828} & \multicolumn{1}{r}{0.101} & \multicolumn{1}{r}{0.91982} &       &  \\
    \multicolumn{1}{|l|}{noise level=25\%} & \multicolumn{1}{r}{0.009702} & \multicolumn{1}{r}{0.048828} & \multicolumn{1}{r}{0.199} & \multicolumn{1}{r}{0.843259} &       &  \\
    \multicolumn{1}{|l|}{smart\_rep} & \multicolumn{1}{r}{-0.09262} & \multicolumn{1}{r}{0.039868} & \multicolumn{1}{r}{-2.323} & \multicolumn{1}{r}{0.02404} & \multicolumn{1}{l}{*} &  \\
    \multicolumn{1}{|l|}{Rastrigin} & \multicolumn{1}{r}{0.344011} & \multicolumn{1}{r}{0.048828} & \multicolumn{1}{r}{7.045} & \multicolumn{1}{r}{3.81E-09} & \multicolumn{1}{l}{***} &  \\
    \multicolumn{1}{|l|}{Levy} & \multicolumn{1}{r}{0.045359} & \multicolumn{1}{r}{0.048828} & \multicolumn{1}{r}{0.929} & \multicolumn{1}{r}{0.357128} &       &  \\
    \multicolumn{1}{|l|}{dimension=20} & \multicolumn{1}{r}{0.073677} & \multicolumn{1}{r}{0.048828} & \multicolumn{1}{r}{1.509} & \multicolumn{1}{r}{0.137262} &       &  \\
    \multicolumn{1}{|l|}{dimension=30} & \multicolumn{1}{r}{0.051091} & \multicolumn{1}{r}{0.048828} & \multicolumn{1}{r}{1.046} & \multicolumn{1}{r}{0.300151} &       &  \\
    \multicolumn{1}{|l|}{Replication =10} & \multicolumn{1}{r}{0.146338} & \multicolumn{1}{r}{0.048828} & \multicolumn{1}{r}{2.997} & \multicolumn{1}{r}{0.004142} & \multicolumn{1}{l}{**} &  \\
    \multicolumn{1}{|l|}{poolSize=2(d+1)} & \multicolumn{1}{r}{0.016954} & \multicolumn{1}{r}{0.039868} & \multicolumn{1}{r}{0.425} & \multicolumn{1}{r}{0.67238} &       &  \\
    \multicolumn{1}{|l|}{DOE=Sobol} & \multicolumn{1}{r}{-0.00669} & \multicolumn{1}{r}{0.039868} & \multicolumn{1}{r}{-0.168} & \multicolumn{1}{r}{0.867332} &       &  \\
    \multicolumn{1}{|l|}{EEPA distance=Cosine} & \multicolumn{1}{r}{-0.01711} & \multicolumn{1}{r}{0.039868} & \multicolumn{1}{r}{-0.429} & \multicolumn{1}{r}{0.669465} &       &  \\
    \multicolumn{1}{|l|}{EEPA number of candidates=6} & \multicolumn{1}{r}{-0.02138} & \multicolumn{1}{r}{0.039868} & \multicolumn{1}{r}{-0.536} & \multicolumn{1}{r}{0.594059} &       &  \\
    \multicolumn{1}{|r|}{} &       &       &       &       &       &  \\
    \midrule
    Signif. Codes: & 0 ‘***’ & 0.001 ‘**’ & 0.01 ‘*’ & 0.05 ‘.’ & 0.1 ‘ ’ & 1 \\
    \bottomrule
    \end{tabular}%
    \end{scriptsize}
  \label{tab:wrep-anova-oa}%
%   \vspace{100mm}
\end{table}%

\newpage

\begin{table}[htbp]
  \centering
  \caption{\new{Average and variance} of MTFAUC of different test functions and surrogates at different noise levels}
  \small
    % \begin{tabular}{|l|llll|llll|}
     \begin{tabular}{|l|@{\hspace{2pt}} *{9}{c|}}
    \toprule
          & \multicolumn{4}{p{15.745em}|}{Rosenbrock} & \multicolumn{4}{p{15.745em}|}{Rastrigin} \\
    \midrule
          & \multicolumn{1}{c}{RBF} & \multicolumn{1}{c}{TK-MARS} & \multicolumn{1}{c}{nonRBF} & \multicolumn{1}{c|}{nonGP} & \multicolumn{1}{c}{RBF} & \multicolumn{1}{c}{TK-MARS} & \multicolumn{1}{c}{nonRBF} & \multicolumn{1}{c|}{nonGP} \\
    \midrule
          & \multicolumn{4}{c|}{noise=0}  & \multicolumn{4}{c|}{noise=0} \\
    \midrule
    Norep & 0.25,0.04 & \textcolor{black}{\textbf{0.17,0.03}} & 0.21,0.03 & 0.76,0.09 & 0.93,0.03 & \textcolor{black}{\textbf{0.89,0.07}} & 0.93,0.04 & 0.95,0.03 \\
    Fixedrep,5 & 0.42,0.08 & 0.40,0.07 & 0.52,0.09 & 0.99,0.02 & 1.00,0.01 & 0.97,0.04 & 1.00,0.01 & 1.00,0.00 \\
    Fixedrep,10 & 0.58,0.06 & 0.62,0.09 & 0.78,0.07 & 1.00,0.00 & 1.00,0.00 & 0.99,0.03 & 1.00,0.00 & 1.00,0.00 \\
    Smartrep,5 & 0.32,0.06 & 0.23,0.04 & 0.31,0.05 & 0.92,0.09 & 0.97,0.03 & 0.92,0.06 & 0.97,0.03 & 0.99,0.01 \\
    Smartrep,10 & 0.32,0.06 & 0.23,0.04 & 0.31,0.05 & 0.91,0.09 & 0.97,0.03 & 0.92,0.06 & 0.97,0.03 & 0.99,0.02 \\
    \midrule
          & \multicolumn{4}{c|}{noise=0.05} & \multicolumn{4}{c|}{noise=0.05} \\
    \midrule
    Norep & 0.29,0.09 & \textcolor{black}{\textbf{0.20,0.08}} & 0.24,0.06 & 0.64,0.11 & 0.89,0.04 & \textcolor{black}{\textbf{0.83,0.07}} & 0.89,0.05 & 0.91,0.04 \\
    Fixedrep,5 & 0.38,0.05 & 0.38,0.06 & 0.45,0.06 & 0.90,0.14 & 0.99,0.02 & 0.95,0.04 & 0.98,0.03 & 1.00,0.01 \\
    Fixedrep,10 & 0.54,0.07 & 0.58,0.07 & 0.68,0.08 & 0.96,0.08 & 0.99,0.01 & 0.98,0.02 & 1.00,0.01 & 1.00,0.01 \\
    Smartrep,5 & 0.31,0.07 & 0.25,0.03 & 0.32,0.05 & 0.72,0.14 & 0.95,0.04 & 0.88,0.06 & 0.97,0.03 & 0.97,0.03 \\
    Smartrep,10 & 0.33,0.06 & 0.29,0.05 & 0.35,0.05 & 0.72,0.14 & 0.95,0.04 & 0.89,0.05 & 0.96,0.03 & 0.97,0.03 \\
    \midrule
          & \multicolumn{4}{c|}{noise=0.1} & \multicolumn{4}{c|}{noise=0.1} \\
    \midrule
    Norep & 0.43,0.23 & 0.27,0.10 & \textcolor{black}{\textbf{0.26,0.13}} & 0.62,0.19 & 0.87,0.06 & \textcolor{black}{\textbf{0.84,0.07}} & 0.86,0.07 & 0.89,0.06 \\
    Fixedrep,5 & 0.42,0.10 & 0.47,0.11 & 0.48,0.09 & 0.96,0.06 & 0.95,0.04 & 0.91,0.06 & 0.96,0.04 & 0.97,0.04 \\
    Fixedrep,10 & 0.63,0.11 & 0.72,0.10 & 0.72,0.10 & 0.97,0.07 & 0.99,0.02 & 0.99,0.02 & 0.99,0.01 & 1.00,0.01 \\
    Smartrep,5 & 0.43,0.11 & 0.43,0.10 & 0.41,0.08 & 0.82,0.11 & 0.90,0.05 & 0.87,0.05 & 0.91,0.06 & 0.94,0.05 \\
    Smartrep,10 & 0.48,0.11 & 0.49,0.09 & 0.49,0.08 & 0.89,0.11 & 0.95,0.05 & 0.93,0.05 & 0.96,0.04 & 0.98,0.03 \\
    \midrule
          & \multicolumn{4}{c|}{noise=0.25} & \multicolumn{4}{c|}{noise=0.25} \\
    \midrule
    Norep & 0.62,0.25 & 0.55,0.21 & \textcolor{black}{\textbf{0.48,0.20}} & 0.85,0.15 & 1.00,0.00 & 1.00,0.00 & 1.00,0.00 & 1.00,0.00 \\
    Fixedrep,5 & 0.63,0.18 & 0.74,0.18 & 0.67,0.16 & 0.99,0.03 & 0.99,0.02 & 1.00,0.01 & 1.00,0.00 & 1.00,0.00 \\
    Fixedrep,10 & 0.59,0.06 & 0.71,0.06 & 0.72,0.09 & 0.99,0.03 & 1.00,0.00 & \textcolor{black}{\textbf{0.99,0.03}} & 1.00,0.00 & 1.00,0.00 \\
    Smartrep,5 & 0.65,0.20 & 0.74,0.15 & 0.60,0.12 & 0.97,0.07 & 1.00,0.00 & 1.00,0.02 & 1.00,0.00 & 1.00,0.00 \\
    Smartrep,10 & 0.55,0.08 & 0.65,0.11 & 0.67,0.12 & 0.93,0.11 & 1.00,0.02 & 0.99,0.04 & 1.00,0.01 & 1.00,0.01 \\
    \bottomrule
    \end{tabular}%
  \label{tab:mtfauc1}%
\end{table}%

% Table generated by Excel2LaTeX from sheet 'Sheet4'
\begin{table}[htbp]
  \centering
    \caption{\new{Average and variance} of MTFAUC of different test functions and surrogates at different noise levels}
  \small
     \begin{tabular}{|l|@{\hspace{2pt}} *{9}{c|}}
    \toprule
          & \multicolumn{4}{p{15.745em}|}{Levy} & \multicolumn{4}{p{15.745em}|}{Ackley} \\
    \midrule
          & \multicolumn{1}{c}{RBF} & \multicolumn{1}{c}{TK-MARS} & \multicolumn{1}{c}{nonRBF} & \multicolumn{1}{c|}{nonGP} & \multicolumn{1}{c}{RBF} & \multicolumn{1}{c}{TK-MARS} & \multicolumn{1}{c}{nonRBF} & \multicolumn{1}{c|}{nonGP} \\
    \midrule
          & \multicolumn{4}{c|}{noise=0}  & \multicolumn{4}{c|}{noise=0} \\
    \midrule
    Norep & 0.25,0.05 & \textcolor{black}{\textbf{0.16,0.04}} & 0.17,0.03 & 0.32,0.06 & \textcolor{black}{\textbf{0.34,0.04}} & 0.43,0.08 & 0.40,0.05 & 0.93,0.05 \\
    Fixedrep,5 & 0.54,0.11 & 0.47,0.09 & 0.55,0.07 & 0.84,0.07 & 0.58,0.06 & 0.61,0.09 & 0.75,0.06 & 1.00,0.00 \\
    Fixedrep,10 & 0.68,0.08 & 0.71,0.10 & 0.86,0.07 & 0.95,0.05 & 0.78,0.05 & 0.80,0.08 & 0.94,0.03 & 1.00,0.00 \\
    Smartrep,5 & 0.37,0.09 & 0.25,0.06 & 0.30,0.05 & 0.80,0.14 & 0.43,0.04 & 0.48,0.08 & 0.54,0.07 & 0.99,0.00 \\
    Smartrep,10 & 0.37,0.09 & 0.25,0.06 & 0.30,0.05 & 0.80,0.14 & 0.43,0.04 & 0.48,0.08 & 0.99,0.01 & 0.99,0.01 \\
    \midrule
          & \multicolumn{4}{c|}{noise=0.05} & \multicolumn{4}{c|}{noise=0.05} \\
    \midrule
    Norep & 0.30,0.10 & \textcolor{black}{\textbf{0.17,0.03}} & 0.22,0.04 & 0.36,0.07 & \textcolor{black}{\textbf{0.33,0.02}} & 0.43,0.05 & 0.39,0.03 & 0.93,0.04 \\
    Fixedrep,5 & 0.38,0.08 & 0.40,0.07 & 0.52,0.06 & 0.72,0.13 & 0.52,0.04 & 0.60,0.07 & 0.67,0.05 & 1.00,0.00 \\
    Fixedrep,10 & 0.56,0.07 & 0.63,0.11 & 0.80,0.08 & 0.89,0.09 & 0.71,0.03 & 0.80,0.06 & 0.89,0.03 & 1.00,0.00 \\
    Smartrep,5 & 0.31,0.09 & 0.25,0.07 & 0.34,0.05 & 0.73,0.11 & 0.41,0.03 & 0.48,0.04 & 0.52,0.04 & 0.99,0.00 \\
    Smartrep,10 & 0.31,0.08 & 0.24,0.04 & 0.37,0.04 & 0.70,0.11 & 0.43,0.03 & 0.52,0.06 & 0.99,0.01 & 0.99,0.00 \\
    \midrule
          & \multicolumn{4}{c|}{noise=0.1} & \multicolumn{4}{c|}{noise=0.1} \\
    \midrule
    Norep & 0.44,0.11 & \textcolor{black}{\textbf{0.26,0.06}} & 0.32,0.09 & 0.47,0.07 & \textcolor{black}{\textbf{0.36,0.02}} & 0.54,0.16 & 0.40,0.05 & 0.93,0.04 \\
    Fixedrep,5 & 0.46,0.07 & 0.46,0.07 & 0.56,0.05 & 0.79,0.10 & 0.58,0.07 & 0.63,0.09 & 0.75,0.08 & 1.00,0.00 \\
    Fixedrep,10 & 0.63,0.07 & 0.65,0.08 & 0.83,0.08 & 0.98,0.03 & 0.77,0.06 & 0.83,0.09 & 0.94,0.03 & 1.00,0.00 \\
    Smartrep,5 & 0.41,0.11 & 0.35,0.06 & 0.42,0.05 & 0.80,0.10 & 0.46,0.08 & 0.56,0.14 & 0.54,0.07 & 1.00,0.00 \\
    Smartrep,10 & 0.45,0.10 & 0.39,0.05 & 0.47,0.06 & 0.80,0.09 & 0.48,0.08 & 0.52,0.07 & 1.00,0.00 & 1.00,0.00 \\
    \midrule
          & \multicolumn{4}{c|}{noise=0.25} & \multicolumn{4}{c|}{noise=0.25} \\
    \midrule
    Norep & 0.57,0.15 & \textcolor{black}{\textbf{0.34,0.10}} & 0.41,0.10 & 0.55,0.10 & \textcolor{black}{\textbf{0.38,0.05}} & 0.51,0.10 & 0.44,0.05 & 0.94,0.04 \\
    Fixedrep,5 & 0.56,0.10 & 0.54,0.11 & 0.64,0.07 & 0.88,0.09 & 0.56,0.06 & 0.65,0.09 & 0.73,0.07 & 1.00,0.00 \\
    Fixedrep,10 & 0.68,0.09 & 0.73,0.10 & 0.89,0.07 & 0.98,0.04 & 0.73,0.05 & 0.82,0.07 & 0.92,0.03 & 1.00,0.00 \\
    Smartrep,5 & 0.52,0.14 & 0.50,0.10 & 0.52,0.06 & 0.85,0.10 & 0.50,0.05 & 0.61,0.11 & 0.65,0.07 & 1.00,0.00 \\
    Smartrep,10 & 0.52,0.12 & 0.48,0.08 & 0.59,0.07 & 0.85,0.10 & 0.57,0.05 & 0.67,0.11 & 1.00,0.00 & 1.00,0.00 \\
    \bottomrule
    \end{tabular}%
  \label{tab:mtfauc2}%
\end{table}%

% Table generated by Excel2LaTeX from sheet 'Sheet4'
\begin{table}[htbp]
  \centering
 \caption{\new{Average and variance} of MTFAUC of different test functions and surrogates at different noise levels}
    % \begin{tabular}{|l|llll|}
     \begin{tabular}{|l|@{\hspace{2pt}} *{9}{c|}}
    \toprule
          & \multicolumn{4}{p{15.745em}|}{Zakharov} \\
    \midrule
          & \multicolumn{1}{c}{RBF} & \multicolumn{1}{c}{TK-MARS} & \multicolumn{1}{c}{nonRBF} & \multicolumn{1}{c|}{nonGP} \\
    \midrule
          & \multicolumn{4}{c|}{noise=0} \\
    \midrule
    Norep & \textcolor{black}{\textbf{0.04,0.00}} & 0.04,0.01 & \textcolor{black}{\textbf{0.04,0.00}} & 0.04,0.01 \\
    Fixedrep,5 & 0.18,0.02 & 0.19,0.03 & 0.18,0.02 & 0.20,0.03 \\
    Fixedrep,10 & 0.36,0.04 & 0.39,0.06 & 0.36,0.04 & 0.39,0.05 \\
    Smartrep,5 & 0.07,0.01 & 0.08,0.01 & 0.07,0.01 & 0.08,0.01 \\
    Smartrep,10 & 0.07,0.01 & 0.08,0.01 & 0.07,0.01 & 0.08,0.01 \\
    \midrule
          & \multicolumn{4}{c|}{noise=0.05} \\
    \midrule
    Norep & 0.84,0.26 & 0.79,0.30 & 0.75,0.36 & 0.85,0.24 \\
    Fixedrep,5 & 0.54,0.33 & 0.55,0.28 & \textcolor{black}{\textbf{0.47,0.28}} & 0.91,0.22 \\
    Fixedrep,10 & 0.57,0.23 & 0.68,0.22 & 0.56,0.17 & 0.92,0.17 \\
    Smartrep,5 & 0.52,0.34 & 0.55,0.25 & 0.47,0.30 & 0.60,0.34 \\
    Smartrep,10 & 0.50,0.23 & 0.64,0.25 & 0.50,0.24 & 0.58,0.24 \\
    \midrule
          & \multicolumn{4}{c|}{noise=0.1} \\
    \midrule
    Norep & 0.59,0.37 & 0.59,0.35 & \textcolor{black}{\textbf{0.54,0.37}} & 0.59,0.37 \\
    Fixedrep,5 & 0.58,0.32 & 0.68,0.31 & 0.56,0.28 & 0.73,0.28 \\
    Fixedrep,10 & 0.78,0.21 & 0.80,0.22 & 0.75,0.24 & 0.76,0.25 \\
    Smartrep,5 & 0.56,0.34 & 0.62,0.28 & 0.57,0.30 & 0.68,0.32 \\
    Smartrep,10 & 0.72,0.25 & 0.76,0.26 & 0.71,0.26 & 0.81,0.28 \\
    \midrule
          & \multicolumn{4}{c|}{noise=0.25} \\
    \midrule
    Norep & 0.71,0.30 & 0.74,0.31 & 0.71,0.32 & 0.69,0.34 \\
    Fixedrep,5 & 0.51,0.28 & 0.58,0.29 & 0.49,0.25 & 0.73,0.26 \\
    Fixedrep,10 & 0.60,0.23 & 0.64,0.21 & 0.63,0.24 & 0.79,0.22 \\
    Smartrep,5 & 0.49,0.27 & 0.57,0.28 & \textcolor{black}{\textbf{0.48,0.24}} & 0.56,0.25 \\
    Smartrep,10 & 0.58,0.26 & 0.60,0.21 & 0.59,0.27 & 0.61,0.25 \\
    \bottomrule
    \end{tabular}%
  \label{tab:mtfauc3}%
\end{table}%

\end{document}